\RequirePackage{fix-cm}
\documentclass[11pt,a4paper,twoside,openany,leqno]{book}

\usepackage{microtype}
\usepackage{amsmath,amssymb,amsthm}
\usepackage{mathrsfs}

\usepackage[
  inner=28mm,outer=25mm,top=26mm,bottom=28mm,
  headheight=15pt,headsep=8mm,footskip=12mm
]{geometry}
\usepackage{xcolor}
\definecolor{accent}{HTML}{243C54}
\definecolor{muted}{HTML}{616A73}
\usepackage{titlesec}
\usepackage{fancyhdr}
\usepackage{needspace}
\usepackage{emptypage}

\allowdisplaybreaks[2]

\titleformat{\chapter}[display]
  {\normalfont\raggedright}
  {\sffamily\small\bfseries\color{accent}\MakeUppercase{\chaptername}\enspace\thechapter}
  {10pt}
  {\Huge\bfseries\color{accent}}
  [\vspace{8pt}{\color{accent!30}\titlerule[0.4pt]}]
\titleformat{\section}
  {\Large\bfseries\color{accent}}{\thesection}{0.65em}{}
\titleformat{\subsection}
  {\large\bfseries}{\thesubsection}{0.65em}{}
\titlespacing*{\chapter}{0pt}{6pt}{25pt}
\titlespacing*{\section}{0pt}{20pt plus 3pt minus 2pt}{9pt plus 2pt}
\titlespacing*{\subsection}{0pt}{15pt plus 2pt minus 2pt}{7pt plus 1pt}
\fancypagestyle{plain}{%
  \fancyhf{}
  \fancyfoot[C]{\small\color{muted}\thepage}
  
}

\newtheoremstyle{lectureplain}
  {9pt plus 2pt minus 2pt}{8pt plus 2pt minus 2pt}
  {\itshape}{}{\bfseries}{.}{0.6em}{}
\newtheoremstyle{lecturedefinition}
  {9pt plus 2pt minus 2pt}{8pt plus 2pt minus 2pt}
  {\normalfont}{}{\bfseries}{.}{0.6em}{}
\newtheoremstyle{lectureremark}
  {8pt plus 2pt minus 2pt}{7pt plus 2pt minus 2pt}
  {\normalfont}{}{\itshape}{.}{0.6em}{}

\theoremstyle{lectureplain}
\newtheorem{theorem}{Theorem}[section]
\newtheorem{lemma}[theorem]{Lemma}
\newtheorem{corollary}[theorem]{Corollary}
\newtheorem{proposition}[theorem]{Proposition}
\newtheorem{conjecture}[theorem]{Conjecture}

\theoremstyle{lecturedefinition}
\newtheorem{definition}[theorem]{Definition}
\newtheorem{example}[theorem]{Example}
\newtheorem{problem}[theorem]{Problem}

\theoremstyle{lectureremark}
\newtheorem{remark}[theorem]{Remark}

\usepackage[
  unicode,colorlinks=true,linkcolor=accent,citecolor=accent,urlcolor=accent,
  bookmarksnumbered=true,bookmarksopen=true
]{hyperref}
\hypersetup{
  pdftitle={Lectures on Bergman  Theory},
  pdfauthor={Bo-Yong Chen}
}

\title{Lectures on Bergman  Theory}
\author{Bo-Yong Chen}
\date{October 16, 2026}

\begin{document}
\frontmatter

\begin{titlepage}
\raggedright
\vspace*{19mm}
{\sffamily\small\color{muted}THE FIRST HALF\par}
\vspace{10mm}
{\color{accent}\rule{0.18\textwidth}{1.2pt}\par}
\vspace{16mm}
{\fontsize{23}{29}\selectfont Lectures on\par}
\vspace{4mm}
{\fontsize{38}{44}\selectfont\bfseries\color{accent}
Bergman Theory\par}
\vspace{20mm}
{\Large Bo-Yong Chen\par}
\vfill
{\color{accent!30}\rule{\textwidth}{0.4pt}\par}
\vspace{5mm}
{\sffamily\small\color{muted}October 16, 2026\par}
\vspace*{12mm}
\end{titlepage}

\chapter*{Preface}
\addcontentsline{toc}{chapter}{Preface}
\markboth{Preface}{Preface}
The deepest concepts in mathematics connect one field with another. Descartes's coordinate system connects algebra with geometry, Newton's calculus connects geometry with dynamics, and Riemann surfaces connect analysis with geometry. Similarly, the kernel function introduced by Bergman connects complex analysis with functional analysis (Hilbert space theory). Bergman theory also serves as a testing ground for the $L^2$ theory of the Cauchy--Riemann equations, a fundamental tool in several complex variables.
Thus, this theory deserves continued exploration. As Hilbert said: We must know, we shall know!

In 2019 and 2025, the author taught a one-semester graduate course entitled \emph{Lectures on Bergman  Theory} at the School of Mathematical Sciences of Fudan University. These notes grew out of those courses and have been revised and expanded over the years, although progress was slow. A turning point came in September 2026: as AI became widely used, some problems in the original notes were resolved with its assistance. AI also supplied elegant new proofs of known theorems, such as a proof of the classical Riemann mapping theorem using only the Bergman kernel, which fully reveals the beauty of the Bergman kernel method. The author's tentative view of obtaining results through AI is that it resembles borrowing books from a library: the books are already there, and humans merely borrow them. Any dispute over priority is therefore meaningless. The meaning of mathematics lies in its supreme beauty. What humans should do, then, is to digest the proofs supplied by AI and reshape them through our uniquely human artistic sensibility.

Among monographs on the Bergman kernel, Bergman's own early work \cite{Bergman1970} is excellent, but the state of knowledge at the time prevented it from treating many deeper developments. The monograph by Jarnicki and Pflug  \cite{JarnickiPflug2013} and the survey of B\l ocki \cite{Blocki2014} also contain several chapters on the Bergman kernel and the Bergman metric. More recently, Krantz published a monograph on Bergman theory \cite{Krantz2013}. Rich in content, it centers on Kohn's theory of the $\bar{\partial}$-Neumann operator and Fefferman's asymptotic expansion theorem for the Bergman kernel on strongly pseudoconvex domains. However, it largely omits applications of weighted $L^2$ estimates for the $\bar{\partial}$-equation to the Bergman kernel and metric; in the author's view, this is precisely the most fascinating aspect of Bergman theory. These notes  aim to  lay a solid foundation for research on the Bergman kernel, upon which the $L^2$ method is widely used.

These notes are divided into two parts. The first part consists of the first three chapters and covers the basic theory of the Bergman kernel and metric on bounded domains, together with the parameter dependence of the Bergman kernel. The second part concerns Bergman theory on complex manifolds, Bergman kernels on hyperconvex manifolds, estimates for the Bergman distance, and mapping properties of the Bergman projection.
The choice of material largely reflects the author's preferences, so some important developments in Bergman theory are inevitably omitted. Interested readers are encouraged to consult the original literature. Results presented without attribution derive partly from the author's research over the past decade and partly from AI.

Finally,  the author strongly recommends Ohsawa's book \cite{Ohsawa100} for a popular introduction to the Bergman kernel.
	
\tableofcontents
\mainmatter
\chapter{The Bergman Kernel on Domains}
\section{Basic Properties}

Let $\Omega\subset \mathbb{C}^n$ be a domain, and let $\mathcal{O}(\Omega)$ denote the space of all holomorphic functions on $\Omega$. Consider the Hilbert space
$$
L^2(\Omega):=\left\{f\ \text{measurable}:\int_\Omega |f|^2<\infty\right\},
$$
with inner product $(f,g):=\int_\Omega f\cdot \bar{g}$ and norm $\|f\|:=\sqrt{(f,f)}$.

\begin{definition}
The space $ A^2(\Omega):=L^2(\Omega)\cap \mathcal{O}(\Omega) $ is called the Bergman space on $\Omega$.
\end{definition}

\begin{example}
$A^2(\mathbb{C}^n)=\{0\}$.
\end{example}

\begin{proof}
Let $P(z,r)$ be the polydisc centered at $z$ with radius $r$. Given $f\in A^2(\mathbb{C}^n)$ and $z\in \mathbb{C}^n$, the submean-value property yields
$$
|f(z)|^2 \le \frac1{|P(z,r)|}\int_{P(z,r)} |f|^2 \le \frac{\|f\|^2}{|P(z,r)|} \rightarrow 0\ \ \ (r\rightarrow \infty).
$$
\end{proof}

\begin{example}
For a bounded domain $\Omega\subset \mathbb{C}^n$, the space $A^2(\Omega)$ contains all polynomials.
\end{example}

Bergman theory begins with the following elementary property:

\begin{lemma}[Bergman's inequality {\cite{Bergman1970}}]\label{lm:BergmanIneq}
For every compact set $K\subset \Omega$, there exists a constant $C_K>0$ such that
  $$
  \sup_K |f|\le C_K \|f\|,\ \ \forall\,f\in A^2(\Omega).
  $$
 \end{lemma}

\begin{proof}
Choose a finite $r>0$ with $r<d(K,\partial\Omega)$ (if $\Omega=\mathbb{C}^n$, any $r$ will do). For every $z\in K$, the polydisc $P(z,r/n)$ is contained in $\Omega$. Since $|f|^2$ is plurisubharmonic,
the submean-value property gives
$$
|f(z)|^2 \le \frac{1}{|P(z,r/n)|}\int_{P(z,r/n)} |f|^2 \le \frac{1}{\pi^n (r/n)^{2n}} \|f\|^2.
$$
\end{proof}

\begin{proposition}
$A^2(\Omega)$ is a closed subspace of $L^2(\Omega)$; in particular, $A^2(\Omega)$ is a Hilbert space.
\end{proposition}

\begin{proof}
Let $\{f_j\}\subset A^2(\Omega)$ satisfy $f_j\rightarrow f\in L^2(\Omega)$.
By Lemma \ref{lm:BergmanIneq}, for every compact subset $K\subset \Omega$ we have
$$
\sup_K |f_j-f_k|\le C_K \|f_j-f_k\|,
$$
Hence $\{f_j\}$ converges uniformly on compact subsets to some $\tilde{f}\in \mathcal{O}(\Omega)$, and
$$
\int_K |f-\tilde{f}|^2\le 2\int_K |f_j-f|^2+2\int_K |f_j-\tilde{f}|^2\rightarrow 0,
$$
that is, $f=\tilde{f}$ almost everywhere on $K$. The conclusion follows since $K$ is arbitrary.
\end{proof}

\begin{proposition}
$A^2(\Omega)$ is separable.
\end{proposition}

\begin{proof}
Take the orthogonal projection $P:L^2(\Omega)\rightarrow A^2(\Omega)$. Since $L^2(\Omega)$ is separable, it has a countable dense subset
  $$
  \{e_1,e_2,\cdots\}.
  $$
For every $f\in A^2(\Omega)$, there is a subsequence $e_{j_k}\rightarrow f$. Thus
  $$
  \|f-P(e_{j_k})\|=\|P(f-e_{j_k})\|\le \|f-e_{j_k}\|\rightarrow 0.
  $$ 
Therefore $\{P(e_j)\}$ is dense in $A^2(\Omega)$.
  \end{proof}
  
 \begin{remark}
The Gram--Schmidt process produces an at most countable complete orthonormal basis of $A^2(\Omega)$.
\end{remark}
  
 \begin{proposition}
Let $\{h_j\}\subset A^2(\Omega)$ be an orthonormal basis. Then
  $$
  \sum  h_j(z)\overline{h_j(w)}
  $$
converges uniformly on compact subsets of $\Omega\times \Omega$.
  \end{proposition}

\begin{proof}
Let $K\subset \Omega$ be compact. By Lemma \ref{lm:BergmanIneq}, for every $w\in K$ and every $l>k$,
\begin{eqnarray}
\left(\sum_{j=k}^{l} |h_j(w)|^2\right)^2 & \le & C_K \left\| \sum_{j=k}^{l} h_j\cdot \overline{h_j(w)}\right\|^2  \nonumber\\
& = & C_K \sum_{j=k}^{l} \|h_j\|^2|h_j(w)|^2\nonumber\\
& = & C_K \sum_{j=k}^{l} |h_j(w)|^2,\nonumber
\end{eqnarray}
that is, $\sum_{j=k}^{l} |h_j(w)|^2\le C_K$. This shows that $\sum_j |h_j|^2$ converges pointwise on $\Omega$.

Moreover, for any $z,w\in K$,
$$
 \left| \sum_{j=k}^l h_j(z)\overline{h_j(w)} \right|^2\le \sum_{j=k}^l |h_j(z)|^2 \sum_{j=k}^l |h_j(w)|^2,  
$$ 
This shows that the sequence $\sum_{j=1}^k h_j(z)\overline{h_j(w)}$, holomorphic in $(z,\bar{w})$, converges pointwise. On the other hand, since
$$
 \left| \sum_{j=1}^k h_j(z)\overline{h_j(w)} \right|^2\le \sum_{j=1}^k |h_j(z)|^2 \sum_{j=1}^k |h_j(w)|^2\le C_K^2,\ \ \ \forall\,z,w\in K,
$$
Vitali's theorem implies uniform convergence on compact subsets.
\end{proof}

\begin{remark}
Surprisingly, even for a planar domain, the preceding proof uses properties of holomorphic functions of two variables!
\end{remark}

Let $w\in \Omega$. Define the linear functional
$$
F:A^2(\Omega)\rightarrow \mathbb{C},\ \ f\mapsto f(w).
$$
By Lemma \ref{lm:BergmanIneq}, $F$ is a bounded linear functional on $A^2(\Omega)$. The Riesz representation theorem therefore gives a \textbf{unique} element
$K_\Omega(\cdot,w)\in A^2(\Omega)$ with the following reproducing property:
$$
f(w)=(f,K_\Omega(\cdot,w))=\int_\Omega f\cdot \overline{K_\Omega(\cdot,w)},\ \ \ w\in \Omega.
$$ 

\begin{definition}
$K_\Omega(z,w)$ is called the Bergman kernel of $\Omega$, and $K_\Omega(z):=K_\Omega(z,z)$ is called its Bergman kernel function.
\end{definition}

\begin{proposition}\label{prop:CNB}
Let $\{h_j\}\subset A^2(\Omega)$ be a complete orthonormal basis. Then
	$$
	K_\Omega(z,w)=\sum_j h_j(z)\overline{h_j(w)}.
	$$
In particular, $K_\Omega(z)$ is a real-analytic, nonnegative plurisubharmonic function.
\end{proposition} 

\begin{proof}
The Fourier expansion gives
$$
K_\Omega(z,w)=\sum_j (K_\Omega(\cdot,w),h_j)\cdot h_j(z)=\sum_j \overline{h_j(w)}\cdot h_j(z).
$$
\end{proof}

\begin{corollary}
\begin{enumerate}
\item	$K_\Omega(z,w)=\overline{K_\Omega(w,z)}$.
\item $
|K_\Omega(z,w)|^2 \le K_\Omega(z)K_\Omega(w).
$
\end{enumerate}
\end{corollary}

\begin{proposition}
$P_\Omega(f)(z):= (f,K_\Omega(\cdot,z))$ is the orthogonal projection from $L^2(\Omega)$ onto $A^2(\Omega)$, called the Bergman projection of $\Omega$.
\end{proposition}

\begin{proof}
Let $\{h_j\}$ be a complete orthonormal basis of $A^2(\Omega)$. Then
$$
P_\Omega(f)(z)=(f,K_\Omega(\cdot,z))=\sum_j (f,h_j) h_j(z) \in A^2(\Omega).
$$
For any $g\in A^2(\Omega)$, convergence of the preceding series in $L^2(\Omega)$ and Parseval's identity give
\[
(P_\Omega f,g)=\sum_j(f,h_j)\overline{(g,h_j)}=(f,g),
\]
that is, $(f-P_\Omega(f))\, \bot\, g$.
\end{proof}

\begin{remark}
Let $u_0=f-P_\Omega(f)$. Then $\bar{\partial} u_0=\bar{\partial}f$ and $u_0\bot A^2(\Omega)$. If $u_1\in L^2(\Omega)$ also satisfies
$\bar{\partial} u_1=\bar{\partial}f$, then $u_1-u_0\in A^2(\Omega)$ and $\|u_1\|^2=\|u_0\|^2+\|u_1-u_0\|^2$. Hence $u_0$ is
the $L^2$-minimal solution of $\bar{\partial}u=\bar{\partial}f$.
\end{remark}

\begin{proposition}
Let $\Omega_1\subset \mathbb{C}^m$ and $\Omega_2\subset \mathbb{C}^n$ be domains. Then
$$
K_{\Omega_1\times \Omega_2}((z,\zeta),(w,\eta)) = K_{\Omega_1}(z,w)K_{\Omega_2}(\zeta,\eta).
$$ 
\end{proposition}

\begin{proof}
First, Fubini's theorem gives
$$
K_{\Omega_1}(\cdot,w)K_{\Omega_2}(\cdot,\eta)\in A^2(\Omega_1\times \Omega_2), \ \ \ \forall\,(w,\eta)\in \Omega_1\times \Omega_2.
$$
Let $f\in A^2(\Omega_1\times \Omega_2)$. Then
\begin{eqnarray*}
\iint_{\Omega_1\times \Omega_2} f(z,\zeta) \overline{K_{\Omega_1}(z,w)}\overline{K_{\Omega_2}(\zeta,\eta)}
&  = & \int_{\Omega_1} \overline{K_{\Omega_1}(z,w)}\int_{\zeta\in \Omega_2} f(z,\zeta) \overline{K_{\Omega_2}(\zeta,\eta)}\\
& = & \int_{z\in \Omega_1} f(z,\eta)\overline{K_{\Omega_1}(z,w)}=f(w,\eta).
\end{eqnarray*}
The conclusion follows from the uniqueness of the reproducing kernel.
\end{proof}

\begin{proposition}\label{prop:Transformation}
Let $F:\Omega_1\rightarrow \Omega_2$ be biholomorphic. Then
$$
K_{\Omega_1}(z,w)=K_{\Omega_2}(F(z),F(w))J_F(z)\overline{J_F(w)}
$$
where $J_F(z)=\mathrm{det}(\partial f_j/\partial z_k)$.
\end{proposition}

\begin{proof}
For $f\in A^2(\Omega_2)$ we have
$$
\int_{\Omega_2} |f|^2 = \int_{\Omega_1} |f\circ F|^2 |J_F|^2.
$$
Thus $f\in A^2(\Omega_2)\iff f\circ F\cdot J_F\in A^2(\Omega_1)$. Let $\{f_j\}\subset A^2(\Omega_2)$ be a complete orthonormal basis. We show that $\{g_j:=f_j\circ F\cdot J_F\}$ is a complete orthonormal basis of $A^2(\Omega_1)$. First,
$$
\int_{\Omega_1} g_j \bar{g}_k =\int_{\Omega_2} f_j \bar{f}_k =\delta_{jk},
$$
where $\delta_{jk}$ is the Kronecker delta. For any $g\in A^2(\Omega_1)$, choose $f\in A^2(\Omega_2)$ such that $g=f\circ F\cdot J_F$. Write $f=\sum c_j f_j$. Then
$$
g=\sum c_j f_j\circ F\cdot J_F=\sum c_j g_j,
$$
Hence $\{g_j\}$ is a complete orthonormal basis of $A^2(\Omega_1)$. Therefore
\begin{eqnarray*}
 K_{\Omega_1} (z,w)  =  \sum g_j(z)\overline{g_j(w)}
& = & \sum f_j( F(z) )\overline{f_j(F(w))} J_F(z) \overline{J_F(w)}\\
& = & K_{\Omega_2}(F(z),F(w))  J_F(z) \overline{J_F(w)}.
\end{eqnarray*}
\end{proof}

\begin{proposition}[Extremal properties]\label{prop:extremal}
Suppose that $A^2(\Omega)\neq \{0\}$. Then
\begin{enumerate}
\item $K_\Omega(z)=\max\{|f(z)|^2/\|f\|^2:f\in A^2(\Omega)\setminus \{0\}\}$; in other words, $K_\Omega(z)$ is the squared norm of the evaluation functional $f\mapsto f(z)$.
\item If $K_\Omega(z)>0$, then
$$
1/\sqrt{K_\Omega(z)} =\min\{\|f\|:f\in A^2(\Omega),f(z)=1\}=:m_\Omega(z);
$$
Moreover, the unique minimizer is $K_\Omega(\cdot,z)/{K_\Omega(z)}$.
\end{enumerate}
\end{proposition}

\begin{proof}
1. Take a complete orthonormal basis $\{h_j\}$ of $A^2(\Omega)$. For every $f\in A^2(\Omega)\setminus \{0\}$, we have $f=\sum_j(f,h_j)h_j$, so that
\begin{eqnarray*}
|f(z)|^2 & \le & \sum_j |(f,h_j)|^2\,\sum_j |h_j(z)|^2 \ \ \ \text{(Schwarz's inequality)}\\
& = & \|f\|^2\,K_\Omega(z)\ \ \ \text{(Parseval's identity)}.
\end{eqnarray*}
Note that
$$
\int_\Omega |K_\Omega(\cdot,z)|^2 = \int_\Omega K_\Omega(\cdot,z)\,\overline{K_\Omega(\cdot,z)}=K_\Omega(z)\ \ \ \text{(reproducing property)},
$$
If $K_\Omega(z)>0$, equality in the preceding inequality is attained at $f=K_\Omega(\cdot,z)$; if $K_\Omega(z)=0$, every nonzero $f\in A^2(\Omega)$ gives a ratio of zero.

2. Left-hand side $\ge$ right-hand side:
Take $f =K_\Omega(\cdot,z)/K_\Omega(z)$. Then $f\in A^2(\Omega)$, $f(z)=1$, and
$$
\|f\|^2=1/K_\Omega(z).
$$

Left-hand side $\le$ right-hand side: if $f\in A^2(\Omega)$ satisfies $f(z)=1$, the reproducing property gives
$$
1=f(z)=\left|\int_\Omega f\, \overline{K_\Omega(\cdot,z)}\right|\le \|f\| \sqrt{K_\Omega(z)}.
$$

Uniqueness: take a complete orthonormal basis $\{h_j\}$ with $h_1=K_\Omega(\cdot,z)/\sqrt{K_\Omega(z)}$. Then
$$
K_\Omega(z)=\sum |h_j(z)|^2 = K_\Omega(z)+\sum_{j=2}^\infty |h_j(z)|^2,
$$
Hence $h_j(z)=0$ for every $j\ge 2$. If $\tilde{f}=\sum_j c_j h_j$ is another minimizer, then $1=\tilde{f}(z)=c_1 h_1(z)$ implies $c_1=1/\sqrt{K_\Omega(z)}$. On the other hand, since
$$
\frac1{K_\Omega(z)}\ge \|\tilde{f}\|^2= \sum_j |c_j|^2 = \frac1{K_\Omega(z)} + \sum_{j\ge 2} |c_j|^2, 
$$ 
we have $c_j=0$ for every $j\ge 2$, and thus $\tilde{f}=h_1/\sqrt{K_\Omega(z)}=K_\Omega(\cdot,z)/K_\Omega(z)$.
\end{proof}

\begin{remark}
Conversely, fix a point $z$ with $K_\Omega(z)>0$. Classical variational arguments also derive the reproducing formula from the minimization problem defining $m_\Omega(z)$. First, a normal-family argument readily proves the existence of a minimizer. Second, the minimizer is unique: if $f_1,f_2$ are two minimizers, then $f_0:=\frac{f_1+f_2}2\in A^2(\Omega)$ satisfies $f_0(z)=1$ and
$$
m_\Omega(z)^2\le \int_\Omega |f_0|^2 \le  \int_\Omega \frac{|f_1|^2+|f_2|^2}2 = m_\Omega(z)^2.
$$
Equality holds in the second inequality if and only if $f_1=f_2$, proving uniqueness.

Now denote the corresponding minimizer by $m_\Omega(\cdot,z)$. Let $f\in A^2(\Omega)$ satisfy $f(z)=0$, and consider the holomorphic family
$$
f_t= m_\Omega(\cdot,z)+t f,  \ \ \ t\in \mathbb{C}.
$$
Since $f_t\in A^2(\Omega)$ and $f_t(z)=1$, we have
$$
m_\Omega(z)^2 \le \int_\Omega |f_t|^2 = m_\Omega(z)^2 + |t|^2\|f\|^2 + 2\mathrm{Re}\left\{t\int_\Omega f\cdot \overline{m_\Omega(\cdot,z)}\right\},\ \ \ \forall\, t\in \mathbb{C}.
$$
Thus
$$
\int_\Omega f\cdot \overline{m_\Omega(\cdot,z)}=0.
$$
For a general $f\in A^2(\Omega)$, substitute $f-f(z)m_\Omega(\cdot,z)$ into the preceding identity to obtain
$$
f(z)m_\Omega(z)^2=\int_\Omega f\cdot\overline{m_\Omega(\cdot,z)}.
$$
Therefore
$$
f(z)=\int_\Omega f\cdot\overline{m_\Omega(\cdot,z)/m_\Omega(z)^2}.
$$
By uniqueness of the reproducing kernel, $m_\Omega(\cdot,z)/m_\Omega(z)^2$ is precisely the Bergman kernel $K_\Omega(\cdot,z)$. This method extends to $p$-Bergman spaces:
$$
A^p(\Omega):=\left\{f\in \mathcal O(\Omega): \int_\Omega |f|^p <\infty\right\}, \ \ \ 1\le p<\infty
$$
(see \cite{ChenZhang,ChenXiong}).
\end{remark}

\begin{corollary}
If $\Omega_1\subset \Omega_2$, then $K_{\Omega_1}(z)\ge K_{\Omega_2}(z)$ for all $z\in \Omega_1$.
\end{corollary}

\begin{corollary}
Let $\Omega\subset \mathbb{C}^n$ be a domain, and write $\delta_\Omega(z)=\inf_{w\in \partial \Omega}|z-w|$. Then
$$
\frac1{|\Omega|} \le K_\Omega(z)\le \frac{n!}{\pi^n} \delta_\Omega(z)^{-2n}.
$$
\end{corollary}

\begin{proof}
If $|\Omega|<\infty$, the first inequality follows by applying the extremal property to the constant function $1$; if $|\Omega|=\infty$, we set $1/|\Omega|=0$, and the inequality is immediate. Moreover, for any $f\in A^2(\Omega)$, the submean-value property gives
$$
|f(z)|^2 \le \frac1{|B(z,\delta_\Omega(z))|}\int_{B(z,\delta_\Omega(z))}|f|^2\le \frac{n!}{\pi^n}\delta_\Omega(z)^{-2n} \|f\|^2,
$$
The second inequality follows from the extremal property.

\end{proof}

Introduce
\begin{eqnarray*}
R_\Omega^\pm(z,w) & := & K_\Omega(z) + K_\Omega(w) \pm  2 \mathrm{Re\,} K_\Omega(z,w)\\
I^{\pm}_\Omega(z,w) & := & K_\Omega(z) + K_\Omega(w) \pm  2 \mathrm{Im\,} K_\Omega(z,w).
\end{eqnarray*}
Similarly, we have the following extremal properties.

\begin{theorem}[Xu Wang]\label{th:Wang}
Suppose that $A^2(\Omega)\neq \{0\}$. Then
\begin{eqnarray*}
R^\pm_\Omega(z,w)  & = &  \max\left\{|f(z)\pm f(w)|^2: f\in A^2(\Omega),\,\|f\|=1\right\}\\
I^{\pm}_\Omega (z,w)& = & \max\left\{| f(z)\pm \mathrm{i}  f(w)|^2: f\in A^2(\Omega),\,\|f\|=1\right\}.
\end{eqnarray*}
\end{theorem}

\begin{proof}
First, note that
\begin{eqnarray*}
\|K_\Omega(\cdot,z)\pm K_\Omega(\cdot,w)\|^2 & = & \int_\Omega |K_\Omega(\cdot,z)|^2 + \int_\Omega |K_\Omega(\cdot,w)|^2\\
&& \pm 2\mathrm{Re} \int_\Omega K_\Omega(\cdot,w)\overline{K_\Omega(\cdot,z)}\\
& = & R^{\pm}_\Omega(z,w).
\end{eqnarray*}
The reproducing property gives
$$
f(z) \pm f(w) = \int_\Omega f\cdot \overline{K_\Omega(\cdot,z)\pm K_\Omega(\cdot,w)},\ \ \forall\,f\in A^2(\Omega).
$$
Combining this with Schwarz's inequality yields
\begin{eqnarray*}
|f(z)\pm f(w)|^2 & \le & \|f\|^2 \|K_\Omega(\cdot,z)\pm K_\Omega(\cdot,w)\|^2\\
& = & \|f\|^2 R_\Omega^\pm (z,w).  
\end{eqnarray*}
Observe that equality in the preceding inequality is attained when
$
f={K_\Omega(\cdot,z)\pm K_\Omega(\cdot,w)}
$
This proves the first identity.
 
Similarly, we have
 \begin{eqnarray*}
\|\mathrm{i} K_\Omega(\cdot,z)\pm K_\Omega(\cdot,w)\|^2 & = & \int_\Omega |K_\Omega(\cdot,z)|^2 + \int_\Omega |K_\Omega(\cdot,w)|^2\\
&& \pm 2\mathrm{Im} \int_\Omega K_\Omega(\cdot,w)\overline{K_\Omega(\cdot,z)}\\
& = & I^{\pm}_\Omega(z,w).
\end{eqnarray*}
The reproducing property gives
$$
 -\mathrm{i} f(z) \pm f(w) = \int_\Omega f\cdot \overline{ \mathrm{i} K_\Omega(\cdot,z)\pm K_\Omega(\cdot,w)},\ \ \forall\,f\in A^2(\Omega).
$$
Combining this with Schwarz's inequality yields
\begin{eqnarray*}
|f(z)\pm \mathrm{i} f(w)|^2 & \le & \|f\|^2 \| \mathrm{i} K_\Omega(\cdot,z)\pm K_\Omega(\cdot,w)\|^2\\
& = & \|f\|^2 I_\Omega^\pm (z,w),
\end{eqnarray*}
and equality is attained when
$
f={ \mathrm{i} K_\Omega(\cdot,z)\pm K_\Omega(\cdot,w)}
$
.
\end{proof}

\begin{corollary}[Bergman--Schiffer $(n=1)$; Xu Wang (general case)]\label{cor:Wang}
If $\Omega_1\subset \Omega_2$, then
$$
R_{\Omega_1}^\pm (z,w)\ge R_{\Omega_2}^\pm (z,w), \   I_{\Omega_1}^\pm (z,w)\ge I_{\Omega_2}^\pm (z,w), \   \   \forall\,z,w\in \Omega_1.
$$
\end{corollary}

\begin{proposition}
Let $\Omega\subset \mathbb{C}^n$ be a bounded Reinhardt domain, and put
 $\mathcal A=\{\alpha\in\mathbb{Z}^n:z^\alpha\in A^2(\Omega)\}$.
Then $\{z^\alpha/\|z^\alpha\|:\alpha\in\mathcal A\}$ is a complete orthonormal basis. In particular,
 $$
 K_\Omega(z,w)=\sum_{\alpha\in\mathcal A} \frac{z^\alpha\bar w^\alpha}{\|z^\alpha\|^2}.
 $$
\end{proposition}

\begin{proof}
Let $\alpha\neq \beta$. We may assume that $\alpha_1\neq \beta_1$. Make the holomorphic change of variables $z_1=\zeta_1 e^{i\theta}$, $z_j=\zeta_j$ for $j\ge 2$, where $\theta\in \mathbb{R}$. Then
$$
\int_\Omega z^\alpha \bar{z}^\beta =e^{i(\alpha_1-\beta_1)\theta} \int_\Omega \zeta^\alpha \bar{\zeta}^\beta. 
$$
It follows that the monomials $z^\alpha$ are mutually orthogonal. Every $f\in\mathcal O(\Omega)$ has a Laurent expansion. For $f\in A^2(\Omega)$, integration using orthogonality on circles shows that the coefficients of terms that are not square-integrable vanish, while the remaining series converges in $L^2$. Thus $\{z^\alpha/\|z^\alpha\|\}$ is a complete orthonormal basis.
\end{proof}

\begin{example}
Let $\mathbb{B}^n$ be the unit ball. Then
$$
K_{\mathbb{B}^n}(z,w)=\frac{n!}{\pi^n}(1-\langle z,w\rangle)^{-n-1}.
$$
Here $\langle z,w\rangle=\sum_j z_j \bar{w}_j$.
\end{example}

\begin{proof}
For $a\in \mathbb{D}$, put $\hat{a}=(a,0,\cdots,0)$ and
$$
F_a:(z_1,\cdots,z_n)\mapsto \left(\frac{z_1-a}{1-\bar{a}z_1},\frac{\sqrt{1-|a|^2}z_2}{1-\bar{a}z_1},\cdots,\frac{\sqrt{1-|a|^2}z_n}{1-\bar{a}z_1}\right).
$$
Then $F_a\in {\rm Aut}(\mathbb{B}^n)$ and $F_a(\hat{a})=0$. Since
$$
K_{\mathbb{B}^n}(z,w)=\sum_\alpha \frac{z^\alpha \bar{w}^\alpha}{\|z^\alpha\|^2},
$$
we obtain $K_{\mathbb{B}^n}(z,0)=|\mathbb{B}^n|^{-1}=n!/\pi^n$. The transformation formula for the Bergman kernel gives
$$
K_{\mathbb{B}^n}(z,\hat{a})  =  K_{\mathbb{B}^n}(z,F_a^{-1}(0)) = K_{\mathbb{B}^n}(F_a(z),0) J_{F_a}(z)\overline{J_{F_a}(\hat{a})}
=\frac{n!}{\pi^n} J_{F_a}(z)\overline{J_{F_a}(\hat{a})}.
$$
A direct calculation yields
$$
J_{F_a}(z) = \frac{(1-|a|^2)^{\frac{n+1}2}}{(1-\bar{a}z_1)^{n+1}},\ \ \ J_{F_a}(\hat{a}) = (1-|a|^2)^{-\frac{n+1}2}.
$$
Hence
$$
K_{\mathbb{B}^n}(z,\hat{a}) = \frac{n!}{\pi^n} (1-\bar{a}z_1)^{-n-1} =  \frac{n!}{\pi^n} (1-\langle z,\hat{a}\rangle)^{-n-1}.
$$
For a general $w\in \mathbb{B}^n$, choose a unitary transformation ${\bf U}$ taking $w$ to $(|w|,0,\cdots,0)$. Then
$$
K_{\mathbb{B}^n}(z,w) = K_{\mathbb{B}^n}({\bf U}z,{\bf U}w) =  \frac{n!}{\pi^n} (1-\langle {\bf U}z, {\bf U} w\rangle)^{-n-1} =  \frac{n!}{\pi^n} (1-\langle z,w\rangle)^{-n-1}. 
$$ 
\end{proof}

Combining this example with the product property of the Bergman kernel gives

\begin{example}
Let $\mathbb{D}^n$ be the unit polydisc. Then
$$
K_{\mathbb{D}^n}(z,w)=\frac{1}{\pi^n}\prod_{j=1}^n (1-z_j\bar{w}_j)^{-2}.
$$
\end{example}

\begin{example}\label{ex:Annuli}
Let $\mathcal A_r=\{z\in \mathbb{C}: r<|z|<1\}$. Then
$$
K_{\mathcal A_r}(z,w)=-\frac1{2\pi z\bar{w}\log r}+\frac1{\pi z\bar{w}}\sum_{n\in \mathbb{Z}\backslash \{0\}}\frac{n(z\bar{w})^n}{1-r^{2n}}.
$$
\end{example}

\begin{proof}
Clearly, $\{h_n(z)=z^n/\|z^n\|:n\in \mathbb{Z}\}$ is a complete orthonormal basis of $A^2(\mathcal A_r)$. Integration in polar coordinates gives
$$
\int_{\mathcal A_r} |z^n|^2 = \left\{
\begin{array}{ll}
2\pi \log 1/r, & n= -1\\
\frac{\pi}{n+1} (1-r^{2n+2}), & n\neq -1.
\end{array}
\right.
$$
Hence
$$
h_n(z) =z^n\cdot \left\{
\begin{array}{ll}
(2\pi \log 1/r)^{-1/2}, & n= -1\\
(\frac{\pi}{n+1} (1-r^{2n+2}))^{-1/2}, & n\neq -1,
\end{array}
\right.
$$
so that
\begin{eqnarray*}
K_{\mathcal A_r}(z,w) & = &  -\frac1{2\pi z\bar{w}\log r} + \sum_{n\in \mathbb{Z}\backslash \{-1\}} \frac{(n+1)(z\bar{w})^n}{\pi(1-r^{2n+2})}\\
& = & -\frac1{2\pi z\bar{w}\log r}+\frac1{\pi z\bar{w}}\sum_{n\in \mathbb{Z}\backslash \{0\}}\frac{n(z\bar{w})^n}{1-r^{2n}}.
\end{eqnarray*}
\end{proof}

\begin{definition}
A bounded domain $\Omega\subset \mathbb{C}^n$ is said to have $C^k$ boundary if, for every $p\in \partial \Omega$, there is a neighborhood $U$ such that, after an
affine change of coordinates,
$$
\Omega\cap U=\{z\in U: {\rm Im\,}z_n<g(z',{\rm Re\,}z_n)\},
$$
where $z'=(z_1,\cdots,z_{n-1})$ and $g$ is a $C^k$ function. Moreover, if $g$ is Lipschitz continuous, then $\Omega$ has Lipschitz boundary; if $0<k<1$, then $C^k$ means H\"older continuous of order $k$. As usual, $C^0$ means continuous.  
\end{definition}

\begin{proposition}
If a domain $\Omega\subset \mathbb{C}^n$ has $C^2$ boundary, then there exists a constant $C>0$ such that
$$
K_\Omega(z)\le C \delta_\Omega(z)^{-n-1}.
$$ 
\end{proposition}

\begin{proof}
First, there exists $r>0$ such that for every $a\in \partial \Omega$ there is a ball $B_a\subset \Omega$ of radius $r$ with $\partial B_a\cap \partial \Omega=\{a\}$. If $z$ is sufficiently close to the boundary, its normal projection $z^\ast\in \partial \Omega$ is defined. Thus
	$$
	K_\Omega(z)\le K_{B_{z^\ast}}(z)\le C\delta_\Omega(z)^{-n-1},
	$$
where $C$ depends only on $n,r$.
\end{proof}

\begin{proposition}[Density]\label{prop:density}
Let $E$ be a uniqueness set for $A^2(\Omega)$, meaning that $f|_E=0\Rightarrow f\equiv 0$ for every $f\in A^2(\Omega)$. Put
$
\Lambda=\{K_\Omega(\cdot,z):z\in E\}.
$ 
Then $\overline{\mathrm{Span\,}\Lambda}=
A^2(\Omega)$.
\end{proposition}

\begin{proof}
For any $f\in A^2(\Omega)$ with $f\bot\, \overline{\text{Span\,}\Lambda}$, we have
$$
0=\int_\Omega f\, \overline{K_\Omega(\cdot,z)}=f(z),\ \ \ \forall\,z\in E.
$$
Since $E$ is a uniqueness set, $f\equiv 0$.
\end{proof}

\section{Bergman's Construction of Complete Orthonormal Bases}
For multi-indices $\alpha=(\alpha_1,\cdots,\alpha_n)$ and $\beta=(\beta_1,\cdots,\beta_n)$, define the order $\beta \prec \alpha$ if and only if $|\beta|<|\alpha|$, or $|\beta|=|\alpha|$ and $\beta_j=\alpha_j$ for every $j<k$, but $\alpha_k<\beta_k$. Define $\alpha\preceq \beta$ if and only if $\alpha \prec \beta$ or $\alpha=\beta$. For example,
$$
(1,0,\cdots,0)\prec (0,1,0,\cdots,0) \prec \cdots \prec (0,\cdots,0,1).
$$

Given $w \in \Omega$, define
$$
H^{(\alpha)}_w(\Omega):=\left\{ f\in A^2(\Omega): f^{(\beta)}(w)=0,\forall\,\beta\prec \alpha,  f^{(\alpha)}(w)=1\right\}
$$
We always assume that $H^{(\alpha)}_w(\Omega)$ is nonempty (as is the case, for example, when $\Omega$ is bounded).
Consider the following higher-order minimization problem:
\begin{equation}\label{eq:mini}
m^{(\alpha)}(w)=\inf\left\{\|f\|: f\in H^{(\alpha)}_w(\Omega) \right\}. 
\end{equation}

\begin{lemma}
The minimization problem \eqref{eq:mini} has a unique minimizer $m^{(\alpha)}(\cdot,w)$.
\end{lemma}

\begin{proof}
We first prove existence. Choose $\{f_j\}\subset H^{(\alpha)}_w(\Omega)$ such that
$$
\|f_j\|\rightarrow m^{(\alpha)}(w)\ \ \  (j\rightarrow \infty).
$$
Bergman's inequality and Montel's theorem give a subsequence $\{f_{j_k}\}$ converging uniformly on compact subsets to some $f_0\in\mathcal O(\Omega)$, with $f_0^{(\beta)}(w)=0$ for every $\beta\prec\alpha$ and $f_0^{(\alpha)}(w)=1$. Moreover, Fatou's lemma gives
$$
\|f_0\| \le \liminf_{k\rightarrow \infty} \|f_{j_k}\| = m^{(\alpha)}(w),
$$
Hence $f_0\in H^{(\alpha)}_w(\Omega)$.
Since $\|f_0\|\ge m^{(\alpha)}(w)$ always holds, $\|f_0\|=m^{(\alpha)}(w)$, so $f_0$ is a minimizer.

We next prove uniqueness. Let $g_0$ be another minimizer. Then $\frac{f_0+g_0}2\in H^{(\alpha)}_w(\Omega)$ and
$$
m^{(\alpha)}(w) \le \left\| \frac{f_0+g_0}2 \right\|\le \frac{\|f_0\|+\|g_0\|}2=m^{(\alpha)}(w),
$$ 
Thus $f_0=cg_0$ for a constant $c$. Since $f_0^{(\alpha)}(w)=1=g_0^{(\alpha)}(w)$, we have $c=1$.
\end{proof}

\begin{lemma}\label{lm:orthogonal}
If $f\in A^2(\Omega)$ satisfies $f^{(\beta)}(w)=0$ for every $\beta\preceq\alpha$, then $f\bot m^{(\beta)}(\cdot,w)$ for every $\beta\preceq\alpha$. In particular, $m^{(\alpha)}(\cdot,w)\bot m^{(\beta)}(\cdot,w)$ whenever $\alpha\neq\beta$.
\end{lemma}

\begin{proof}
We may assume that $f\neq0$. For every constant $c$ and every $\beta\preceq\alpha$, we have $m^{(\beta)}(\cdot,w)+cf\in H^{(\beta)}_w(\Omega)$. If $(f,m^{(\beta)}(\cdot,w))\neq0$, choose $c=-(m^{(\beta)}(\cdot,w),f)/\|f\|^2$ to obtain
\begin{eqnarray*}
m^{(\beta)}(w)^2 & \le & \|m^{(\beta)}(\cdot,w)+ {c} f\|^2 \\
& = & \|m^{(\beta)}(\cdot,w)\|^2 + |c|^2\|f\|^2 + 2\mathrm{Re}\left\{{c}(f,m^{(\beta)}(\cdot,w)) \right\}\\
& = & m^{(\beta)}(w)^2-|c|^2\|f\|^2< m^{(\beta)}(w)^2,
\end{eqnarray*}
a contradiction.
\end{proof}

Put $h_\alpha:=m^{(\alpha)}(\cdot,w)/{m^{(\alpha)}(w)}$. We have

\begin{theorem}[Bergman {\cite{Bergman1970}}]
The family $\{h_\alpha\}$ is a complete orthonormal system in $A^2(\Omega)$.
\end{theorem}

\begin{proof}
Clearly, $h_\alpha\in A^2(\Omega)$ and $\|h_\alpha\|=1$. By Lemma \ref{lm:orthogonal}, it remains to prove completeness. Given $f\in A^2(\Omega)$, put $c_\alpha=(f,h_\alpha)$. Bessel's inequality gives
$$
\sum_\alpha |c_\alpha|^2\le \|f\|^2.
$$
Hence the series $s(z)=\sum_\alpha c_\alpha h_\alpha$ converges in norm. Bergman's inequality implies uniform convergence on compact subsets of $\Omega$, so $s\in\mathcal O(\Omega)$.

We claim that $f=s$. By the identity theorem, it suffices to check that $f^{(\alpha)}(w)=s^{(\alpha)}(w)$ for every $\alpha$. First, Lemma \ref{lm:orthogonal} gives $(f-f(w)m^{(0)}(\cdot,w),m^{(0)}(\cdot,w))=0$, and hence
$$
f(w)=\frac1{m^{(0)}(w)^2} (f,m^{(0)}(\cdot,w))=c_0 h_0(w)=s(w).
$$ 
Suppose that $f^{(\beta)}(w)=s^{(\beta)}(w)$ for all $\beta\prec\alpha$, that is,
$$
f^{(\beta)}(w)=\sum_{\gamma\preceq \beta} c_\gamma h_\gamma^{(\beta)}(w).
$$
Choose a complex number $c$ such that
$$
\sum_{\beta\prec \alpha} c_\beta h_\beta^{(\alpha)}(w) + c h^{(\alpha)}_\alpha(w) = f^{(\alpha)}(w).
$$
Put
$$
g=f-\sum_{\beta\prec \alpha} c_\beta h_\beta - c h_\alpha.
$$
Then Lemma \ref{lm:orthogonal} gives
$$
0=(g,h_\alpha)=(f,h_\alpha)-c,
$$
Thus $c=c_\alpha$, and consequently $f^{(\alpha)}(w)=s^{(\alpha)}(w)$.
\end{proof}

\section{Removable Singularities}

\begin{definition}
A subset $E$ of a domain $\Omega\subset\mathbb{R}^n$ is called polar if, for every $a\in E$, there exist a neighborhood $U\subset\Omega$ of $a$ and a subharmonic function $\psi$ on $U$, not identically $-\infty$, such that
$$
E\cap U\subset \{x\in U:\psi(x)=-\infty\}.
$$ 
If, in this definition, $\Omega$ is a domain in $\mathbb{C}^n$ and $\psi$ is required to be plurisubharmonic, then $E$ is called pluripolar.
\end{definition}

We always have
$$
\text{proper analytic subsets} \subsetneq \text{pluripolar sets} \subseteq \text{polar sets} \subsetneq \text{null sets}.
$$
When $n=1$, pluripolar and polar sets coincide; when $n\ge2$, the middle inclusion is also strict.

\begin{theorem}[Carleson $(n=1)$; Chen--Wu--Wang (general case) {\cite{Carleson1967,ChenWuWang2015}}]\label{th:polar}
Let $E$ be a closed polar subset of a domain $\Omega\subset\mathbb{C}^n$. Then $K_\Omega|_{\Omega\backslash E}=K_{\Omega\backslash E}$.
\end{theorem}

Let $W^{k,2}(\Omega)$ denote the usual Sobolev space, and let $W^{k,2}_0(\Omega)$ be the completion of $C^\infty_0(\Omega)$ in $W^{k,2}(\Omega)$. Clearly, $W^{k,2}_0(\Omega)\subset W^{k,2}(\Omega)$, but equality need not hold. Write
$$
W^{k,2}_{\rm loc}(\Omega)=\left\{f:\forall\,z\in\Omega,\ \exists\,\text{a neighborhood }U\subset\Omega\ \text{such that }f\in W^{k,2}(U)\right\}.
$$
We first prove the following lemma.

 \begin{lemma}\label{lm:PoincareIneq}
Let $\varphi<-c<0$ be a subharmonic function on a domain $\Omega\subset\mathbb{R}^n$, not identically $-\infty$. Then $\log(-\varphi)\in W^{1,2}_{\rm loc}(\Omega)$. Moreover,
\begin{equation}\label{eq:PoincareIneq-}
\int_{\Omega}\phi^2  |\nabla \log (-\varphi)|^2 \le 4 \int_\Omega  |\nabla\phi|^2,\ \ \ \phi\in C^1_0(\Omega).
\end{equation}
\end{lemma}

\begin{proof}
Since $\varphi$ is subharmonic, $\varphi\in L^1_{\rm loc}(\Omega)$, and hence $\log(-\varphi)\in L^p_{\rm loc}(\Omega)$ for every $p>0$. It therefore suffices to check that $\nabla\log(-\varphi)\in L^2_{\rm loc}(\Omega)$. Fix domains $\Omega'\subset\subset\Omega''\subset\subset\Omega$. Choose smooth subharmonic functions $\{\varphi_j\}$ defined on a neighborhood of $\overline{\Omega''}$, with $\varphi_j<-c$ and $\varphi_j\downarrow\varphi$. Let $\phi\in C^1_0(\Omega'')$. Since
\begin{eqnarray*}
\int_{\Omega''} \phi^2 \Delta [-\log(-\varphi_j)]  & = & 2 \int_{\Omega''} \phi \nabla \phi \cdot \nabla \log (-\varphi_j)\\
& \le & 2\left[\int_{\Omega''} |\nabla \phi|^2\right]^{1/2}\left[ \int_{\Omega''} \phi^2 |\nabla \log (-\varphi_j)|^2\right]^{1/2},
\end{eqnarray*}
and
$$
\Delta [-\log(-\varphi_j)]=\frac{\Delta\varphi_j}{-\varphi_j} + |\nabla \log(-\varphi_j)|^2\ge  |\nabla \log(-\varphi_j)|^2,
$$
we obtain
$$
\int_{\Omega''} \phi^2 |\nabla \log (-\varphi_j)|^2 \le 4 \int_{\Omega''} |\nabla \phi|^2.
$$
This shows that $\{\log(-\varphi_j)\}$ is uniformly bounded in $W^{1,2}(\Omega')$ (note that $\log(-\varphi_j)\le\log(-\varphi)$). The Banach--Alaoglu theorem gives a subsequence $\{\log(-\varphi_{j_k})\}$ converging weakly to some $\Phi\in W^{1,2}(\Omega')$. On the other hand, dominated convergence gives $\log(-\varphi_j)\rightarrow\log(-\varphi)$ in $L^2(\Omega')$. Thus $\Phi=\log(-\varphi)$ almost everywhere on $\Omega'$, so $\log(-\varphi)\in W^{1,2}(\Omega')$. Since $\Omega'$ is arbitrary, (\ref{eq:PoincareIneq-}) follows.
\end{proof}

\begin{proof}[Proof of Theorem \ref{th:polar}]
It suffices to check that $E$ is removable for $A^2(\Omega\setminus E)$. Given $z^0\in E$, choose a neighborhood $U\subset\subset\Omega$ and a subharmonic function $\psi$ on $U$ such that
$$
E\cap U\subset \{z\in U:\psi(z)=-\infty\}.
$$
Choose a further neighborhood $V\subset\subset U$ of $z^0$. We may assume that $\psi<-1$ on $\overline V$. By Lemma \ref{lm:PoincareIneq}, $\log(-\psi)$ belongs to $W^{1,2}$ on a neighborhood of $\overline V$.

Let $f\in A^2(\Omega\setminus E)$.  By Weyl's lemma, it suffices to verify
\begin{equation}\label{eq:distribution}
\int_V f\cdot \frac{\partial g}{\partial \bar{z}_j}=0,\ \ \ \forall\,g\in C_0^\infty(V),\,\forall\,1\le j\le n.
\end{equation}
 
Given $R>0$, define
$$
\chi_R(z)=\kappa(\log(-\psi)-\log R),
$$
where $\kappa$ is a cutoff function on $\mathbb{R}$ satisfying $\kappa|_{(-\infty,-\log2]}=1$ and $\kappa|_{[0,\infty)}=0$. The upper semicontinuity of $\psi$ implies that
$$
U_R=\{z\in U:\psi(z)<-R\}
$$ 
is open, $U\cap E\subset U_R$, and $\chi_R=0$ on $U_R$.
Thus, for $R>1$, the function $\chi_R\,g$ is compactly supported in $V\setminus E$, and
   \begin{equation}\label{eq:polar1}
 0  =  \int_V  \frac{\partial f}{\partial \bar{z}_j}\cdot \chi_R g
   =  \int_V \chi_R f \cdot \frac{\partial g}{\partial \bar{z}_j} + \int_V fg \cdot \frac{\partial  \chi_R}{\partial \bar{z}_j}.
 \end{equation}
Since
 \begin{eqnarray*}
 \left|\int_V f g \cdot \frac{\partial \chi_R}{\partial \bar{z}_j} \right|^2 & \le & 
 \sup|\kappa'|^2 \int_{V\setminus E} |fg|^2 \int_{U_{R/2}\setminus U_R}  |\nabla \log (-\psi)|^2\\
 & \rightarrow & 0\ \ \ (R\rightarrow \infty)
 \end{eqnarray*}
(note that $U_R\downarrow\psi^{-1}(-\infty)$ and $|\psi^{-1}(-\infty)|=0$), letting $R\rightarrow\infty$ in \eqref{eq:polar1} yields \eqref{eq:distribution}.
\end{proof}

\begin{remark}
The same method shows that if $u\in L^2_{\rm loc}(\Omega)$, $v\in L^1_{\rm loc}(\Omega)$, and $\bar\partial u=v$ holds distributionally on $\Omega\setminus E$, then the equation also holds on $\Omega$.
\end{remark}

On the other hand, we have the following result.
\begin{proposition}
Let $\Omega$ be a domain in $\mathbb{C}^n$ with $K_\Omega(z)>0$ for every $z\in\Omega$, and let $E$ be a closed subset of $\Omega$ of positive Lebesgue measure. Then, for every $z\in\Omega\setminus E$,
\begin{equation}\label{eq:<}
K_{\Omega}(z) < K_{\Omega\setminus E}(z).
\end{equation}
Conversely, still assuming $K_\Omega>0$, let $E$ be a compact subset of a domain $\Omega\subset\mathbb{C}^n\ (n\ge2)$ such that $\Omega\setminus E$ is connected. If \eqref{eq:<} holds at \textbf{some point} $z\in\Omega\setminus E$, then $E$ has positive measure. In particular, \eqref{eq:<} holds at \textbf{every} $z\in\Omega\setminus E$.
\end{proposition}

\begin{proof}
1. Note that, for every $z\in\Omega\setminus E$,
\begin{eqnarray*}
K_{\Omega\setminus E}(z)\ge \frac{K_\Omega(z)^2}{\int_{\Omega\setminus E} |K_\Omega(\cdot,z)|^2} & = & \frac{K_\Omega(z)^2}{\int_{\Omega} |K_\Omega(\cdot,z)|^2 - \int_E |K_\Omega(\cdot,z)|^2}\\
& = & \frac{K_\Omega(z)^2}{K_\Omega(z) - \int_E |K_\Omega(\cdot,z)|^2},
\end{eqnarray*}
It therefore suffices to check that $\int_E|K_\Omega(\cdot,z)|^2>0$. Otherwise, $K_\Omega(\cdot,z)=0$ almost everywhere on $E$. Since $|E|>0$, the properties of zero sets of holomorphic functions imply that $K_\Omega(\cdot,z)\equiv0$ on $\Omega$, which is impossible because
$
K_\Omega(z,z)= K_\Omega(z)>0.
$

2. By Hartogs's extension theorem, $K_{\Omega\setminus E}(\cdot,z)$ extends to an element of $A^2(\Omega)$, still denoted by $K_{\Omega\setminus E}(\cdot,z)$. Suppose that $|E|=0$. Then
\begin{eqnarray*}
K_{\Omega}(z)\ge \frac{K_{\Omega\setminus E}(z)^2}{\int_{\Omega} |K_{\Omega\setminus E}(\cdot,z)|^2} & = & \frac{K_{\Omega\setminus E}(z)^2}{\int_{\Omega\setminus E} |K_{\Omega\setminus E}(\cdot,z)|^2 + \int_E |K_{\Omega\setminus E} (\cdot,z)|^2}\\
& = & K_{\Omega\setminus E}(z),
\end{eqnarray*}
a contradiction.
\end{proof}

\begin{conjecture}
Let $E$ be a closed nonpolar subset of a domain $\Omega\subset\mathbb{C}$. Then
$$
K_{\Omega}(z) < K_{\Omega\setminus E}(z),\ \ \ \forall\,z\in \Omega\setminus E.
$$
\end{conjecture}

We now prove the conjecture for simply connected proper domains. The potential-theoretic fact used is that every nonpolar compact set supports a probability measure of finite logarithmic energy. We first establish a square-integrability estimate for its Cauchy transform.

\begin{lemma}\label{lem:finite-energy-cauchy}
Let $F\subset\mathbb{D}$ be compact, and let $\mu$ be a probability measure supported on $F$ such that
\[
I(\mu):=\iint_{F\times F}\log\frac4{|\zeta-\eta|}\,
 d\mu(\zeta)\,d\mu(\eta)<\infty.
\]
Then the Cauchy transform
\[
C_\mu(z):=\int_F\frac{d\mu(\zeta)}{z-\zeta},\qquad z\in\mathbb{D}\setminus F,
\]
belongs to $A^2(\mathbb{D}\setminus F)$ and satisfies
\[
\int_{\mathbb{D}\setminus F}|C_\mu(z)|^2\,dA(z)
\le 4\pi+2\pi I(\mu).
\]
Here $dA$ denotes ordinary planar area measure.
\end{lemma}

\begin{proof}
For distinct $\zeta,\eta\in\mathbb{D}$, put $d=|\zeta-\eta|$. On $\Delta(\zeta,d/2)$ we have $|z-\eta|\ge d/2$, so
\[
\int_{\mathbb{D}\cap \Delta(\zeta,d/2)}
\frac{dA(z)}{|z-\zeta|\,|z-\eta|}
\le\frac2d\int_{\Delta(\zeta,d/2)}\frac{dA(z)}{|z-\zeta|}=2\pi.
\]
The same estimate holds on the disc centered at $\eta$. On the remaining region $\mathcal R$, use
\[
\frac1{uv}\le\frac12\left(\frac1{u^2}+\frac1{v^2}\right)
\]
together with $d/2\le|z-\zeta|,|z-\eta|<2$ to obtain
\[
\int_{\mathcal R} \frac{dA(z)}{|z-\zeta|\,|z-\eta|}
\le 2\pi\int_{d/2}^2\frac{dr}{r}
=2\pi\log\frac4d.
\]
Combining these estimates gives
\[
\int_{\mathbb{D}}\frac{dA(z)}{|z-\zeta|\,|z-\eta|}
\le 4\pi+2\pi\log\frac4{|\zeta-\eta|}.
\]
Finite energy implies that $\mu$ has no atoms, so the diagonal has zero $\mu\times\mu$ measure. The preceding estimate and Tonelli's theorem give
\begin{align*}
\int_{\mathbb{D}\setminus F}|C_\mu(z)|^2\,dA(z)
&\le \int_{\mathbb{D}}
 \left(\int_F\frac{d\mu(\zeta)}{|z-\zeta|}\right)^2\,dA(z)\\
 & =  \int_{\mathbb{D}}
\int_F\frac{d\mu(\zeta)}{|z-\zeta|} \int_F\frac{d\mu(\eta)}{|z-\eta|}\,dA(z)\\
& = \int_F \int_F \int_{\mathbb{D}}
\frac{1}{|z-\zeta| |z-\eta|}\,dA(z)\,d\mu(\zeta)\,d\mu(\eta)\\
&\le 4\pi+2\pi I(\mu)<\infty.
\end{align*}
Since $C_\mu$ is holomorphic on $\mathbb{D}\setminus F$, the proof is complete.
\end{proof}

\begin{proposition}\label{prop:simply-connected-strict}
Let $\Omega\subsetneq\mathbb{C}$ be a simply connected domain, and let $E\subset\Omega$ be relatively closed and nonpolar. Then
\[
K_{\Omega\setminus E}(a)>K_\Omega(a),\qquad a\in\Omega\setminus E.
\]
\end{proposition}

\begin{proof}
Fix $a\in\Omega\setminus E$. By the Riemann mapping theorem, there is a biholomorphic map $\Phi:\Omega\to\mathbb{D}$ with $\Phi(a)=0$. Polarity is invariant under biholomorphic maps, and the transformation formula for the Bergman kernel gives
\[
\begin{aligned}
K_\Omega(a)&=|\Phi'(a)|^2K_{\mathbb{D}}(0),\\
K_{\Omega\setminus E}(a)&=|\Phi'(a)|^2K_{\mathbb{D}\setminus\Phi(E)}(0).
\end{aligned}
\]
Thus it suffices to consider $\Omega=\mathbb{D}$ and $a=0$. Continue to denote the removed set by $E$, and put $D=\mathbb{D}\setminus E$.

If $|E|>0$, using the constant function as a test function in the extremal problem for the Bergman kernel gives
\[
K_D(0)\ge\frac1{|D|}>\frac1\pi=K_{\mathbb{D}}(0).
\]
Now suppose that $|E|=0$. Since a countable union of polar sets is polar, a compact exhaustion of $E$ yields a nonpolar compact subset $F\subset E$. Classical potential theory gives a probability measure $\mu$ on $F$ with $I(\mu)<\infty$. Put
\[
g(z):=zC_\mu(z)=\int_F\frac{z}{z-\zeta}\,d\mu(\zeta).
\]
By Lemma~\ref{lem:finite-energy-cauchy}, $g\in A^2(D)$. Since $0\notin F$, we have $g(0)=0$.

For fixed $\zeta\in F$, averaging over circles gives
\[
\frac1{2\pi}\int_0^{2\pi}
\frac{re^{i\theta}}{re^{i\theta}-\zeta}\,d\theta
=\begin{cases}
0,&0<r<|\zeta|,\\
1,&|\zeta|<r<1.
\end{cases}
\]
Integration with respect to the radius then yields the key identity
\[
\int_{\mathbb{D}}\frac{z}{z-\zeta}\,dA(z)=\pi(1-|\zeta|^2).
\]
Note that $\int_{\mathbb{D}}|z-\zeta|^{-1}\,dA(z)\le4\pi$, so Fubini's theorem applies. Since $|E|=0$, we obtain
\[
A:=\int_D g(z)\,dA(z)
=\pi\int_F(1-|\zeta|^2)\,d\mu(\zeta)>0.
\]
For real $t>0$, put $f_t=1-tg$. Then $f_t(0)=1$ and
\[
\|f_t\|_{L^2(D)}^2
=\pi-2tA+t^2\|g\|_{L^2(D)}^2<\pi
\]
for sufficiently small $t$. The extremal property of the Bergman kernel therefore gives
\[
K_D(0)\ge\frac{|f_t(0)|^2}{\|f_t\|_{L^2(D)}^2}
>\frac1\pi=K_{\mathbb{D}}(0).
\]
This proves the required strict inequality.
\end{proof}

\section{Stability}

Does the Bergman kernel vary continuously when the domain varies? Approximation from inside and from outside leads to different answers.

\begin{theorem}[Ramadanov {\cite{Ramadanov1967}}]\label{th:Ramadanov}
If $\Omega_j\uparrow\Omega$, then $K_{\Omega_j}(z,w)$ converges to $K_\Omega(z,w)$ uniformly on compact subsets of $\Omega\times\Omega$.
\end{theorem}

\begin{proof}
We first prove that $K_{\Omega_j}(w)\downarrow K_\Omega(w)$ for every $w\in\Omega$. Since $\Omega_j\subset\Omega$, we have $K_{\Omega_j}(w)\ge K_\Omega(w)$. On the other hand,
the extremal property gives $f_j\in A^2(\Omega_j)$ with $\|f_j\|\le1$ and $K_{\Omega_j}(w)=|f_j(w)|^2$. By Bergman's inequality, $\{f_j\}$ is a normal family, so a subsequence $\{f_{j_k}\}$ converges uniformly on compact subsets to some $f\in\mathcal O(\Omega)$. Fatou's lemma gives, for every compact set $E\subset\Omega$,
$$
\int_E |f|^2 \le \liminf_{k\rightarrow \infty} \int_E |f_{j_k}|^2 \le 1,
$$
and hence $\|f\|\le1$.
Thus
$$
K_\Omega(w) \ge |f(w)|^2 = \lim_{k\rightarrow \infty} |f_{j_k}(w)|^2 =  \lim_{k\rightarrow \infty} K_{\Omega_{j_k}}(w) = \lim_{j\rightarrow \infty} K_{\Omega_{j}}(w), 
$$
where the last equality follows from the monotone decrease of $K_{\Omega_j}(w)$.

We next prove off-diagonal convergence. Note that
\begin{eqnarray*}
\int_{\Omega_j} |K_{\Omega_j}(\cdot,w)-K_\Omega(\cdot,w)|^2 & = & \int_{\Omega_j} |K_{\Omega_j}(\cdot,w)|^2 + \int_{\Omega_j} |K_\Omega(\cdot,w)|^2\\
&& -2{\rm Re} \int_{\Omega_j} K_\Omega(\cdot,w)\overline{K_{\Omega_j}(\cdot,w)}\\
& \le & K_{\Omega_j}(w)+ \int_{\Omega} |K_\Omega(\cdot,w)|^2-2K_\Omega(w)\\
& = & K_{\Omega_j}(w)-K_\Omega(w).
\end{eqnarray*}
Let $E$ be a compact subset of $\Omega$. For $j\ge j_0\gg1$, we have $E\subset\Omega_j$. The preceding inequality implies
$$
\iint_{E\times E}  |K_{\Omega_j}(\cdot,w)-K_\Omega(\cdot,w)|^2\le \int_E [K_{\Omega_j}(w)-K_\Omega(w)] \rightarrow 0,
$$
The last convergence follows by dominated convergence, since $K_{\Omega_j}(w)\downarrow K_\Omega(w)$ and
$$
K_{\Omega_j}(w)-K_\Omega(w)\le K_{\Omega_{j_0}}(w)-K_\Omega(w),\ \ \ \forall\,j\ge j_0.
$$
Viewing $K_{\Omega_j}(z,w)$ and $K_\Omega(z,w)$ as holomorphic functions of $(z,\bar w)\in\Omega\times\bar\Omega$, Bergman's inequality gives uniform convergence on compact subsets off the diagonal as well.
\end{proof}

\begin{remark}
Alternatively, we already know that $K_{\Omega_j}(z)\downarrow K_\Omega(z)$ for every $z\in\Omega$. Dini's theorem gives uniform convergence on compact subsets. Corollary \ref{cor:Wang} then gives
$$
|K_{\Omega_j}(z,w)-K_\Omega(z,w)|\le K_{\Omega_j}(z)-K_\Omega(z)+K_{\Omega_j}(w)-K_\Omega(w),
$$
and the conclusion follows.
\end{remark}

Let $\Omega\subset\mathbb{C}^n$ be a bounded domain, and put
$$
\Omega_\varepsilon:=\{z\in \Omega:\delta_\Omega(z)>\varepsilon\}, \ \ \ \varepsilon>0.
$$
Given a compact subset $E\subset\Omega$, one naturally hopes for quantitative estimates of
$$
\max_{z,w\in E}  | K_{\Omega_\varepsilon}(z,w)-K_\Omega(z,w)|
$$
as $\varepsilon\rightarrow0$. However, we have the following result.

\begin{proposition}\label{prop:slow-growth}
For every function $\eta:(0,1)\to(0,\infty)$ satisfying
$\lim_{\varepsilon\downarrow0}\eta(\varepsilon)=0$, there exist a fixed bounded  domain
$\Omega\subset\mathbb{D}$ and a fixed point $z_0=\frac12\in\Omega$ such that
\begin{equation}\label{eq:slow-main}
\lim_{\varepsilon\downarrow0}
\frac{K_{\Omega_\varepsilon}(z_0)-K_\Omega(z_0)}{\eta(\varepsilon)}
=+\infty.
\end{equation}
\end{proposition}

\begin{proof}
We first record a direct estimate on an annulus.
Let $\mathcal A_r=\{z\in\mathbb{C}:r<|z|<1\}$, where $0<r<1$.
By Example \ref{ex:Annuli},
\begin{equation}\label{eq:slow-annulus}
K_{\mathcal A_r}(z,z)-K_{\mathbb{D}}(z,z)
\ge\frac{1}{2\pi|z|^2\log(1/r)}, \quad z\in \mathcal A_r.
\end{equation}

Put $r_j=e^{-j-3}$ for $j\ge1$. Since $\eta(\varepsilon)\to0$, we can recursively choose a strictly decreasing
sequence of positive numbers $\varepsilon_j\downarrow0$ such that
\begin{equation}\label{eq:slow-eps}
\varepsilon_1<\frac18,\qquad
\varepsilon_{j+1}<\frac{\varepsilon_j}{2},\qquad
\sup_{0<\varepsilon\le\varepsilon_j}\eta(\varepsilon)
\le\frac{1}{(j+1)^3}.
\end{equation}
For each $j$, choose a finite set $F_j\subset\mathbb{D}_{r_j}$ such that
\begin{equation}\label{eq:slow-net}
\operatorname{dist}(z,F_j)\le\frac{\varepsilon_j}{2}
\quad z\in\overline{ \mathbb{D}}_{r_j}.
\end{equation}
Define
\[
F=\{0\}\cup\bigcup_{j=1}^\infty F_j,\qquad
\Omega=\mathbb{D}\setminus F.
\]
Since $r_j\to0$, the set $F$ is countable and compact, and all its points except $0$ are isolated.
Thus $\Omega$ is a bounded domain. Since $F\subset\mathbb{D}_{e^{-4}}$, we also have
$z_0=\frac12\in\Omega$ and
\[
\operatorname{dist}(z_0,\partial\Omega)\ge\frac12-e^{-4}>\frac18.
\]
In particular, $z_0\in\Omega_\varepsilon$ for every $0<\varepsilon\le\varepsilon_1$.

By Theorem \ref{th:polar}, the restriction map $A^2(\mathbb{D})\to A^2(\Omega)$
is an isometric isomorphism, and hence
\begin{equation}\label{eq:slow-removable}
K_\Omega(z,w)=K_{\mathbb{D}}(z,w),\qquad z,w\in\Omega.
\end{equation}
Fix $j$. Since $F\subset\partial\Omega$, \eqref{eq:slow-net} shows that, whenever
$\varepsilon\ge\varepsilon_j$, every point of $\Omega$ in the closed disc $\overline{\mathbb{D}}_{r_j}$
has distance at most $\varepsilon_j/2<\varepsilon$ from the boundary. Therefore
\begin{equation}\label{eq:slow-inclusion}
\Omega_\varepsilon\subset \mathcal A_{r_j}\qquad(\varepsilon\ge\varepsilon_j).
\end{equation}
Thus, for $\varepsilon_j\le\varepsilon\le\varepsilon_1$,
\begin{equation}\label{eq:slow-lower}
\begin{aligned}
K_{\Omega_\varepsilon}(z_0,z_0)-K_\Omega(z_0,z_0)
&\ge K_{\mathcal A_{r_j}}(z_0,z_0)-K_{\mathbb{D}}(z_0,z_0)\\
&\ge\frac{1}{2\pi|z_0|^2\log(1/r_j)}
=\frac{2}{\pi(j+3)}.
\end{aligned}
\end{equation}
For any $0<\varepsilon<\varepsilon_1$, choose the unique $j\ge2$ such that
$\varepsilon_j\le\varepsilon<\varepsilon_{j-1}$.
By \eqref{eq:slow-eps}, $\eta(\varepsilon)\le j^{-3}$. Combining this with
\eqref{eq:slow-lower} gives
\[
\frac{K_{\Omega_\varepsilon}(z_0,z_0)-K_\Omega(z_0,z_0)}{\eta(\varepsilon)}
\ge\frac{2j^3}{\pi(j+3)}.
\]
As $\varepsilon\downarrow0$, we have $j\to\infty$, and the right-hand side tends to $+\infty$, proving
\eqref{eq:slow-main}.
\end{proof}
\medskip

\begin{definition}
Let $E\subset\mathbb{C}^n$ be bounded. A decreasing sequence of domains $\{\Omega_j\}$ is said to approximate $E$ from outside if, for every open set $U\supset E$, we have $E\subset\Omega_j\subset U$ for all sufficiently large $j$. We write $\Omega_j\downarrow E$.
\end{definition}

\begin{problem}\label{prob:decreasing}
Let $\Omega,\Omega_j$ be bounded domains in $\mathbb{C}^n$ such that $\Omega_j\downarrow\overline\Omega$. Under what conditions does $K_{\Omega_j}(z,w)$ converge to $K_\Omega(z,w)$ uniformly on compact subsets?
\end{problem}

Under approximation from outside, the limit of the Bergman kernels need not equal the Bergman kernel of the original domain.

\begin{example}
Let $\Omega=\mathbb{D}\backslash[0,1)$ and $\Omega_j=\mathbb{D}_{1+1/j}$. Then $\Omega_j\downarrow\overline\Omega$, but
$$
K_{\Omega_j}(z)\uparrow K_\mathbb{D}(z) < K_\Omega(z),\ \ \ \forall\,z\in \Omega. 
$$
\end{example}

\begin{example}
Consider the Hartogs triangle
$$
\Omega=\{(z_1,z_2):|z_1|<|z_2|<1\}.
$$
If $\Omega_j\downarrow\overline\Omega$, then $K_{\Omega_j}\nrightarrow K_\Omega$. Note that $\Omega$ is fat, that is, $\overline\Omega^\circ=\Omega$.
\end{example}

\begin{proof}
By Theorem \ref{th:Skwarcz} given below, it suffices to verify that $\mathcal O(\overline\Omega)$ is not dense in $A^2(\Omega)$. Let $U\supset\overline\Omega$ be open and $f\in\mathcal O(U)$. Consider the Hartogs figure
\begin{eqnarray*}
\mathcal H_\varepsilon & := & \left\{(z_1,z_2): |z_1|<\varepsilon,\,|z_2|<1+\varepsilon\right\}\\
&&  \cup \left\{(z_1,z_2): |z_1|<1+\varepsilon,\,1-\varepsilon<|z_2|<1+\varepsilon\right\}\\
& \subset & U\ \ \ \text{(if $\varepsilon\ll1$)}.
\end{eqnarray*}
Then $f$ extends holomorphically to a neighborhood of $\overline{\mathbb{D}}\times\overline{\mathbb{D}}$. Put $f_0(z)=1/z_2\in\mathcal O(\Omega)$. Then
$$
\int_\Omega |f_0|^2 = \int_{|z_1|<|z_2|<1} |z_2|^{-2} =(2\pi)^2 \int_0^1 r_1(-\log r_1)dr_1<\infty.
$$
Suppose there exists $\{f_j\}\subset\mathcal O(\overline\Omega)$ with $\int_\Omega|f_j-f_0|^2\rightarrow0$ as $j\rightarrow\infty$. Since $f_j\in\mathcal O(\overline{\mathbb{D}}\times\overline{\mathbb{D}})$, we have
$$
\int_\Omega f_j \cdot \bar{f}_0 = \int_{|z_2|<1} 1/\bar{z}_2 \cdot \int_{|z_1|<|z_2|} f_j(z_1,z_2)=0,
$$
so that
$$
\int_\Omega |f_j-f_0|^2 = \int_\Omega |f_j|^2 + \int_\Omega |f_0|^2\ge  \int_\Omega |f_0|^2,
$$
a contradiction.
\end{proof}

\begin{theorem}[Skwarczy{\'n}ski {\cite{Skwarczynski1980}}]\label{th:Skwarcz}
Let $\Omega\subset\mathbb{C}^n$ be a bounded domain with $|\partial\Omega|=0$.
\begin{enumerate}
\item If $\mathcal O(\overline\Omega)$ is dense in $A^2(\Omega)$, then $\Omega_j\downarrow\overline\Omega$ implies that $K_{\Omega_j}(z,w)$ converges to $K_\Omega(z,w)$ uniformly on compact subsets.
\item If $\Omega_j\downarrow\overline\Omega$ and $K_{\Omega_j}(z)\rightarrow K_\Omega(z)$ for every $z\in\Omega$, then $\mathcal O(\overline\Omega)$ is dense in $A^2(\Omega)$.
\end{enumerate}
\end{theorem}

\begin{proof}
1. If $f\in\mathcal O(\overline\Omega)$, then $f\in A^2(\Omega_j)$ for all sufficiently large $j$. Thus, given $z\in\Omega$,
$$
|f(z)|^2 \le K_{\Omega_j}(z)\|f\|_{L^2(\Omega_j)}^2,\ \ \ j\gg 1.
$$ 
Since $\{K_{\Omega_j}(z)\}$ is increasing, its limit exists; denote it by $\widetilde K(z)$. Since $|\partial\Omega|=0$, letting $j\rightarrow\infty$ in the preceding inequality gives
$$
|f(z)|^2 \le \widetilde{K}(z)\|f\|_{L^2(\Omega)}^2.
$$ 
By density of $\mathcal O(\overline\Omega)$ in $A^2(\Omega)$, this inequality holds for all $f\in A^2(\Omega)$, yielding $K_\Omega(z)\le\widetilde K(z)$. Conversely, $\Omega\subset\Omega_j$ implies $K_\Omega(z)\ge\widetilde K(z)$. Thus $K_\Omega(z)=\widetilde K(z)$. As in the proof of Ramadanov's theorem, we obtain uniform convergence of $K_{\Omega_j}(z,w)$ to $K_\Omega(z,w)$ on compact subsets.

2. Let $f\in A^2(\Omega)$. By Proposition \ref{prop:density}, there exist $w_1^{(m)},\cdots,w_{l_m}^{(m)}\in\Omega$ and complex numbers $c_1^{(m)},\cdots,c_{l_m}^{(m)}$ such that
$$
f_m:=\sum_{k=1}^{l_m}c_k^{(m)}K_\Omega(\cdot,w_k^{(m)})\rightarrow f\ \text{in }L^2(\Omega),\ \ \ (m\rightarrow\infty).
$$
On the other hand, for fixed $m$,
$$
F_{m,j}:=\sum_{k=1}^{l_m} c_k^{(m)} K_{\Omega_j}(\cdot,w_k^{(m)}) \in \mathcal O(\overline{\Omega}). 
$$
The reproducing property and $\|K_{\Omega_j}(\cdot,w)\|_{L^2(\Omega)}^2\le K_{\Omega_j}(w)$ give
\[
\|K_{\Omega_j}(\cdot,w)-K_\Omega(\cdot,w)\|_{L^2(\Omega)}^2
\le K_\Omega(w)-K_{\Omega_j}(w)\longrightarrow0.
\]
Thus, for each $m$, as $j\to\infty$,
$$
K_{\Omega_j}(\cdot,w_k^{(m)})\rightarrow K_\Omega(\cdot,w_k^{(m)})\ \text{in }L^2(\Omega),
$$
and consequently $F_{m,j}\rightarrow f_m$ in $L^2(\Omega)$.
\end{proof}

\begin{theorem}[cf. \cite{ChenZhang2000}]\label{th:approximation}
Let $\Omega\subset\mathbb{C}^n$ be a bounded pseudoconvex domain with $C^0$ boundary. If there is a sequence of bounded pseudoconvex domains $\Omega_k\downarrow\overline\Omega$, then $\mathcal O(\overline\Omega)$ is dense in $A^2(\Omega)$. In particular, $K_{\Omega_k}(z,w)$ converges to $K_\Omega(z,w)$ uniformly on compact subsets.
\end{theorem}

\begin{proof}
Choose $p_1,\cdots,p_m\in\partial\Omega$, unit vectors $\mathcal N_1,\cdots,\mathcal N_m$, positive numbers $r_0,\varepsilon_0$, and $2n$-dimensional cubes centered at $p_j$,
$\mathcal C_j=\mathcal C_j'\times[(-r_0,r_0)\cdot\mathcal N_j]$, where $\mathcal C_j'$ is a $(2n-1)$-dimensional cube, such that
$$
\partial \Omega \cap \mathcal C_j -\varepsilon\cdot \mathcal N_j\subset \Omega,\ \ \  \partial \Omega \cap \mathcal C_j+\varepsilon\cdot \mathcal N_j \subset  \mathbb{C}^n\setminus \overline{\Omega},\ \ \ \forall\,0<\varepsilon\le \varepsilon_0,
$$
and $\{\mathcal C_j\}$ covers $\partial\Omega$. Choose a domain $U_0$ with $\Omega\supset\supset U_0\supset\Omega\setminus(\mathcal C_1\cup\cdots\cup\mathcal C_m)$.
Take a partition of unity $\{\chi_0,\chi_1,\cdots,\chi_m\}$ subordinate to $\{U_0,\mathcal C_1,\cdots,\mathcal C_m\}$.
Given $f\in A^2(\Omega)$ and $0<\varepsilon\ll1$, put
$$
f_\varepsilon:=\chi_0 f+\sum_{j=1}^m \chi_j f (\cdot-\varepsilon\cdot \mathcal N_j).
$$
Clearly, $f_\varepsilon\in C^\infty(\overline\Omega)$, and on $\Omega$ we have
$$
f_\varepsilon=f+\sum_{j=1}^m \chi_j \left(f(\cdot-\varepsilon\cdot \mathcal N_j)-f\right).
$$
For any $\eta>0$, choose a compact subset $E_\eta\subset\Omega$ with $|\Omega\backslash E_\eta|<\eta$. Note that, for each $j$,
\begin{eqnarray*}
&& \int_{\Omega\cap {\rm supp\,}\chi_j} |f(\cdot-\varepsilon\cdot \mathcal N_j)-f|^2\\
  & =&\left(\int_{\Omega\cap {\rm supp\,}\chi_j\cap E_\eta} + \int_{\Omega\cap {\rm supp\,}\chi_j\setminus E_\eta}\right)  |f(\cdot-\varepsilon\cdot \mathcal N_j)-f|^2\\
& =: & I_1(\varepsilon)+I_2(\varepsilon).
\end{eqnarray*}
Clearly, $I_1(\varepsilon)\rightarrow0$ as $\varepsilon\rightarrow0$, whereas
$$
I_2(\varepsilon) \le 2 \int_{\Omega\backslash E_\eta} |f|^2 + 2\int_{\Omega\cap {\rm supp\,}\chi_j\setminus E_\eta-\varepsilon\cdot \mathcal N_j} |f|^2
$$
can be made arbitrarily small when $\varepsilon\ll\eta\ll1$, by absolute continuity of the integral. Hence
$$
\|f_\varepsilon-f\|_{L^2(\Omega)}\rightarrow0\ \text{and }\ \|\bar\partial f_\varepsilon\|_{L^2(\Omega)}\rightarrow0\ \ \ (\varepsilon\rightarrow0).
$$
Since $\Omega$ has $C^0$ boundary (in particular, $|\partial\Omega|=0$) and $\Omega_k\downarrow\overline\Omega$, there exists $k_\varepsilon\in\mathbb{Z}^+$ such that $\|\bar\partial f_\varepsilon\|_{L^2(\Omega_{k_\varepsilon})}\rightarrow0$ as $\varepsilon\rightarrow0$. By pseudoconvexity of $\Omega_{k_\varepsilon}$ and H\"ormander's theorem,\footnote{For all $L^2$ estimates and extension tools used in these notes, including H\"ormander's estimates, the Donnelly--Fefferman estimates, and the Ohsawa--Takegoshi extension theorem, see the monograph of the author \cite{Chen2022}.} the equation $\bar\partial u=\bar\partial f_\varepsilon$ has a solution $u_\varepsilon$ satisfying
$$
\int_{\Omega_{k_\varepsilon}} |u_\varepsilon|^2 \le C \int_{\Omega_{k_\varepsilon}} |\bar{\partial} f_\varepsilon|^2,
$$
where $C$ is independent of $\varepsilon$.
Thus $F_\varepsilon:=f_\varepsilon-u_\varepsilon\in\mathcal O(\Omega_{k_\varepsilon})\subset\mathcal O(\overline\Omega)$ and
$$
\|F_\varepsilon-f\|_{L^2(\Omega)}\le \|f_\varepsilon-f\|_{L^2(\Omega)} + \|u_\varepsilon\|_{L^2(\Omega)}\rightarrow 0\ \ \ (\varepsilon\rightarrow 0).
$$
\end{proof}

\begin{theorem}[cf.  \cite{Chen1999}]\label{ch3:th:S}\footnote{The hypotheses can be substantially weakened; see \cite{ChenFornaessWu}.}
Let $\Omega\subset\mathbb{C}^n$ be a bounded domain, and let $\{\Omega_j\}$ be pseudoconvex domains satisfying $\Omega_j\downarrow\overline\Omega$. Suppose that there exist constants $0<\alpha<1$ and $C>0$ such that
\begin{enumerate}
\item[$(1)$] $\Lambda_j:=\max_{\partial \Omega} \delta_{\Omega_j}\rightarrow 0\ \ \ (j\rightarrow \infty)$;
\item[$(2)$] $\Lambda_j\le C\lambda_j^\alpha$, where $\lambda_j:=\min_{\partial\Omega}\delta_{\Omega_j}$.
\end{enumerate}
Then $K_{\Omega_j}\rightarrow K_\Omega$.
\end{theorem}

We first prove the following proposition.

\begin{proposition}\label{ch3:prop:S}
Let $\Omega\subset\mathbb{C}^n$ be a bounded domain, and let $\Omega_j\supset\overline\Omega$ be pseudoconvex domains. Suppose there exist negative continuous plurisubharmonic functions $\rho_j$ on $\Omega_j$ and a constant $C>1$ such that
\begin{enumerate}
\item[$(1)$] There exist $\varepsilon_j>0$ with $\rho_j>-\varepsilon_j$ on $\Omega_j\setminus \Omega$; 
\item[$(2)$]  $|\Omega\backslash \Omega_{j,-C\varepsilon_j}|\rightarrow 0\ \ \ (j\rightarrow \infty)$, where $\Omega_{j,-c}:=\{z\in\Omega:\rho_j(z)<-c\}$;.
\end{enumerate}
Then $K_{\Omega_j}\rightarrow K_\Omega$.
\end{proposition}

\begin{proof}
Choose a cutoff function $\chi:\mathbb{R}\rightarrow[0,1]$ with $\chi|_{(-\infty,-\log C]}=1$ and $\chi|_{[0,\infty)}=0$. Put
$$
\eta_j(z)=\chi(-\log(-\rho_j(z))+\log\varepsilon_j)\quad(z\in\Omega),\qquad\eta_j=0\quad\text{on }\Omega_j\setminus\Omega.
$$
By condition (1), this cutoff vanishes on a neighborhood of $\partial\Omega$, so the extension by zero is valid. We have
$$
{\rm supp\,}\eta_j\subset\{z\in\Omega:\rho_j(z)\le-\varepsilon_j\}\Subset\Omega,\ \ \ \{\eta_j=1\}\supset \Omega_{j,-C\varepsilon_j}
$$
and
$$
|\bar{\partial}\eta_j|^2_{i\partial\bar{\partial}(-\log(-\rho_j))}\le \sup |\chi'|^2, 
$$
i.e., 
$$
\sup |\chi'|^2 i \partial\bar{\partial}(-\log(-\rho_j))\ge i \partial \eta_j\wedge \bar{\partial}\eta_j
$$
 holds in the sense of distributions.  
Given $w\in\Omega$, the Donnelly--Fefferman estimate gives a solution $u_j$ of $\bar\partial u=K_\Omega(\cdot,w)\,\bar\partial\eta_j$ such that
$$
\int_{\Omega_j} |u_j|^2 \le C_0 \int_{\Omega_j} |K_\Omega(\cdot,w)|^2\,|\bar{\partial}\eta_j|^2_{i\partial\bar{\partial}(-\log(-\rho_j))}
\le C_0 \int_{\Omega\backslash \Omega_{j,-C\varepsilon_j}} |K_\Omega(\cdot,w)|^2.
$$
Put $E_j=\Omega\setminus\Omega_{j,-C\varepsilon_j}$ and $f=K_\Omega(\cdot,w)$. Since $|E_j|\to0$, absolute continuity of the integral gives $\|f\|_{L^2(E_j)}\to0$. Therefore
\[
f_j:=\eta_j f-u_j\in A^2(\Omega_j),\qquad
\|f_j-f\|_{L^2(\Omega)}\le (1+\sqrt{C_0})\|f\|_{L^2(E_j)}\longrightarrow0,
\]
where $\eta_jf$ is set equal to zero on $\Omega_j\setminus\Omega$. Also,
\[
\|f_j\|_{L^2(\Omega_j)}\le\|f\|_{L^2(\Omega)}+\|u_j\|_{L^2(\Omega_j)}
=\sqrt{K_\Omega(w)}+o(1).
\]
By Bergman's inequality, $f_j(w)\to f(w)=K_\Omega(w)$. The extremal property then gives
\[
\liminf_j K_{\Omega_j}(w)\ge K_\Omega(w).
\]
Together with $K_{\Omega_j}(w)\le K_\Omega(w)$, this proves the claim. Local uniform convergence of the two-variable kernels follows from normal families and uniqueness of polarization.
\end{proof}

\begin{proof}[Proof of Theorem \ref{ch3:th:S}]
Discarding finitely many terms, we may assume that $\Omega_j$ is bounded. Put
\[
\widehat\Omega_j=\{z\in\Omega_j:\delta_{\Omega_j}(z)>\lambda_j/2\},\qquad
\rho_j=-\log\delta_{\Omega_j}+\log(\lambda_j/2).
\]
By Oka's lemma, $\widehat\Omega_j$ is pseudoconvex, $\rho_j$ is negative, continuous, and plurisubharmonic there, and $\overline\Omega\subset\widehat\Omega_j\Subset\Omega_j$. If these sets are disconnected, kernels and estimates are understood componentwise. Put
\[
\varepsilon_j=\log\frac{2\Lambda_j}{\lambda_j/2}.
\]
Then $\rho_j\ge-\varepsilon_j+\log2$ on $\widehat{\Omega_j}\setminus\Omega$, so condition (1) of Proposition \ref{ch3:prop:S} holds. The sublevel sets below are taken only inside $\Omega$.
Choose $1<\beta<\frac1{1-\alpha}$. Then
\begin{eqnarray*}
\widehat{\Omega}_{j,-\beta\varepsilon_j} & = & \left\{z\in\Omega: \delta_{\Omega_j}> (\lambda_j/2)^{1-\beta}(2\Lambda_j)^\beta=:\tau_j\right\}\\
& \supset & \left\{z\in\Omega:\delta_\Omega>\tau_j\right\}\ \ \ (\text{since }\delta_{\Omega_j}\ge\delta_\Omega).
\end{eqnarray*}
Since
$$
\tau_j\lesssim \Lambda_j^{\frac{1-\beta}\alpha+\beta} \asymp  \Lambda_j^{\frac{1-\beta(1-\alpha)}\alpha}\rightarrow 0\ \ \ (j\rightarrow \infty),
$$
we obtain $|\Omega\backslash\widehat\Omega_{j,-\beta\varepsilon_j}|\rightarrow0$. Proposition \ref{ch3:prop:S} gives $K_{\widehat\Omega_j}\rightarrow K_\Omega$. Finally, use cofinality of the outer approximations: for each fixed $j$, the set $\widehat\Omega_j$ is an open neighborhood of $\overline\Omega$, so $\Omega_k\subset\widehat\Omega_j$ for all sufficiently large $k$. Therefore
\[
K_\Omega\ge\lim_kK_{\Omega_k}\ge K_{\widehat\Omega_j}.
\]
Letting $j\to\infty$ proves the conclusion.
\end{proof}

In fact, the second condition in Theorem \ref{th:Skwarcz} can be substantially weakened:

\begin{theorem}\label{th:outer-nonpluripolar}
Let $\Omega\subset\mathbb{C}^n$ be a bounded domain with $|\partial\Omega|=0$, and suppose
$\Omega_j\downarrow\overline\Omega$. Assume there is a uniqueness set $E\subset\Omega$ for $A^2(\Omega)$ (that is, $f\in A^2(\Omega)$ and $f|_E=0$ imply $f=0$) such that
\[
K_{\Omega_j}(a)\longrightarrow K_\Omega(a),\qquad a\in E,
\]
Then $\mathcal O(\overline\Omega)$ is dense in $A^2(\Omega)$, and
\[
K_{\Omega_j}(z,w)\longrightarrow K_\Omega(z,w)
\quad\text{uniformly on compact subsets of }\Omega\times\Omega.
\]
\end{theorem}

We first prove the following result.

\begin{lemma}\label{lem:outer-limit-space}
Let $\Omega\subset\mathbb{C}^n$ be a bounded domain with $|\partial\Omega|=0$, and suppose
$\Omega_j\downarrow\overline\Omega$. Let
$
M
$
be the closure of $\mathcal O(\overline\Omega)$ in $A^2(\Omega)$, put
 \[
 N:=A^2(\Omega)\ominus M,
\]
and denote by $K_M$ the reproducing kernel of $M$ for the $L^2(\Omega)$ inner product. Then
\begin{equation}\label{eq:outer-limit-space}
K_{\Omega_j}(z,w)\longrightarrow K_M(z,w)
\quad\text{uniformly on compact subsets of }\Omega\times\Omega.
\end{equation}
In particular, the limit is independent of the chosen outer approximation. If $\{g_\nu\}$ is a complete orthonormal basis of $N$, then
\begin{equation}\label{eq:outer-defect-kernel}
K_\Omega(z)-K_M(z)=\sum_\nu |g_\nu(z)|^2.
\end{equation}
Consequently, $K_{\Omega_j}(a)\to K_\Omega(a)$ at a point $a\in\Omega$
if and only if every $g\in N$ vanishes at $a$.
\end{lemma}

\begin{proof}
Point evaluation is continuous on $M$, so $K_M$ exists and has the usual extremal characterization.
For sufficiently large $j$, the restriction of every $f\in A^2(\Omega_j)$ to $\Omega$ belongs to $M$, and
$\|f\|_{L^2(\Omega)}\le\|f\|_{L^2(\Omega_j)}$. Thus
\[
K_{\Omega_j}(z)\le K_M(z)\le K_\Omega(z).
\]
Conversely, fix $f\in\mathcal O(\overline\Omega)$. For sufficiently large $j$, the function $f$
is holomorphic on a neighborhood of $\overline{\Omega}_j$, and hence
\[
|f(z)|^2\le K_{\Omega_j}(z)\|f\|_{L^2(\Omega_j)}^2.
\]
By outer approximation and $|\partial\Omega|=0$, we have
$\|f\|_{L^2(\Omega_j)}\to\|f\|_{L^2(\Omega)}$.
Letting $j\to\infty$ and using the density of $\mathcal O(\overline\Omega)$ in $M$
gives $K_M(z)\le\lim_jK_{\Omega_j}(z)$. This proves convergence on the diagonal.
The estimate
\[
|K_{\Omega_j}(z,w)|^2\le K_\Omega(z)K_\Omega(w)
\]
shows that these kernels form a normal family. Every subsequential limit equals $K_M$ on the diagonal, so uniqueness of polarization yields
\eqref{eq:outer-limit-space}.

Finally, combining complete orthonormal bases of $M$ and $N$ gives one for $A^2(\Omega)$.
The series expansion of the kernel then gives \eqref{eq:outer-defect-kernel}. All terms on the right are nonnegative,
which proves the final equivalence. If $N=\{0\}$, the right-hand side is understood to be zero.
\end{proof}

\begin{proof}[Proof of Theorem \ref{th:outer-nonpluripolar}]
Use the notation of Lemma \ref{lem:outer-limit-space}. The hypothesis and \eqref{eq:outer-defect-kernel} imply that every $g\in N$ vanishes on $E$. Uniqueness gives $N=\{0\}$, hence $M=A^2(\Omega)$, and the result follows.
\end{proof}

Note that if $E$ is not contained in any analytic subset of $\Omega$ then it is a uniqueness set for $\mathcal O(\Omega)$.  Thus the following example shows that, when $n\ge2$, the uniqueness-set condition in Theorem \ref{th:outer-nonpluripolar} is somewhat the best possible.

\begin{example}\label{ex:outer-one-point-counterexample}
Put
\begin{equation}\label{eq:outer-log-hartogs}
\Omega=\left\{(z,w)\in\mathbb{C}^2:
0<|w|<1,\quad |z|<\frac{1}{\sqrt{\log(e/|w|)}}\right\}.
\end{equation}
Then $\Omega$ is a bounded domain, $|\partial\Omega|=0$, and
$\Omega=\operatorname{int}\overline\Omega$. For every outer approximation
$\Omega_j\downarrow\overline\Omega$, we have
\begin{equation}\label{eq:outer-counterexample-gap}
K_\Omega(\zeta)-\lim_{j\to\infty}K_{\Omega_j}(\zeta)
=\frac{|z|^2}{\pi^2|w|^2(1-|z|^2)^3},
\qquad \zeta=(z,w)\in\Omega.
\end{equation}
Thus the kernel functions converge to $K_\Omega$ on $\{(0,w):0<|w|<1\}$, but fail to converge to
$K_\Omega$ at every other point; moreover, $\mathcal O(\overline\Omega)$ is not dense in $A^2(\Omega)$.
\end{example}

\begin{proof}
Put $R(r)=[\log(e/r)]^{-1/2}$. Connect a point first to $(0,w)$ in the $z$ direction, and then move along a path in the punctured disc.
This shows that $\Omega$ is connected. Its boundary consists of
\[
\{|z|=R(|w|),\ 0<|w|<1\},\qquad
\{|w|=1,\ |z|\le1\},\qquad \{(0,0)\}
\]
and therefore has zero four-dimensional Lebesgue measure. Since $R(r)\to0$ as $r\to0$, it is easy to see that
$\operatorname{int}\overline\Omega=\Omega$. Outer approximations always exist; for example, take
\[
\Omega_j=\{\xi\in\mathbb{C}^2:
\operatorname{dist}(\xi,\overline\Omega)<1/j\}.
\]

Since $\Omega$ is a Reinhardt domain, the square-integrable Laurent monomials form
a complete orthogonal system in $A^2(\Omega)$. Since $\Omega$ contains points with $z=0$, the exponent of $z$ must be nonnegative.
For $m\ge0$ and $k\in\mathbb{Z}$, direct integration gives
\begin{equation}\label{eq:outer-monomial-norm}
\|z^m w^k\|_{L^2(\Omega)}^2
=\frac{2\pi^2}{m+1}
\int_0^1\frac{r^{2k+1}}{[\log(e/r)]^{m+1}}\,dr.
\end{equation}
The integral is finite if $k\ge0$. For $k=-1$, the substitution $t=\log(e/r)$ reduces it to
$\int_1^\infty t^{-m-1}\,dt$, which is finite exactly when $m\ge1$. For $k\le-2$, it diverges.
Thus, if $P$ denotes the closure of the polynomials in $A^2(\Omega)$, then
\begin{equation}\label{eq:outer-counterexample-splitting}
A^2(\Omega)=P\mathbin{\oplus}
\overline{\operatorname{span}\{z^m/w:m\ge1\}}^{\,A^2(\Omega)}.
\end{equation}

We claim that $M=P$. Since polynomials belong to $\mathcal O(\overline\Omega)$, we have $P\subset M$.
Conversely, if $f\in\mathcal O(\overline\Omega)$, then $(0,0)\in\overline\Omega$ implies that
$f$ is holomorphic in a polydisc neighborhood of the origin. Computing Laurent coefficients on a small torus in this neighborhood
shows that all coefficients of negative powers of $w$ vanish. By uniqueness of Laurent expansions, these coefficients also vanish in the expansion on
all of $\Omega$. Completeness and \eqref{eq:outer-counterexample-splitting}
then give $f\in P$, hence $M\subset P$.

In particular, $z/w$ is a nonzero function in $A^2(\Omega)$ orthogonal to $M$,
so $\mathcal O(\overline\Omega)$ is not dense.

By \eqref{eq:outer-monomial-norm}, for $m\ge1$ we have
\[
\|z^m/w\|_{L^2(\Omega)}^2=\frac{2\pi^2}{m(m+1)}.
\]
Combining Lemma \ref{lem:outer-limit-space} with the preceding orthogonal decomposition gives, for $\zeta=(z,w)\in\Omega$,
\begin{align*}
K_\Omega(\zeta)-\lim_jK_{\Omega_j}(\zeta)
&=\frac{1}{2\pi^2|w|^2}
\sum_{m=1}^\infty m(m+1)|z|^{2m}\\
&=\frac{|z|^2}{\pi^2|w|^2(1-|z|^2)^3}.
\end{align*}
Here we used $\sum_{m\ge1}m(m+1)t^m=2t/(1-t)^3$; the series converges since $|z|<1$.
The right-hand side vanishes exactly when $z=0$, proving all the assertions.
\end{proof}

\begin{problem}\label{prob:outer-one-point}
Let $\Omega\subset\mathbb{C}$ be a bounded domain with $|\partial\Omega|=0$, and suppose
$\Omega_j\downarrow\overline\Omega$. If, for {\bf some} $a\in\Omega$,
$K_{\Omega_j}(a)\to K_\Omega(a)$, must
$\mathcal O(\overline\Omega)$ be dense in $A^2(\Omega)$, and hence must the kernel functions converge throughout the domain?
\end{problem}

For bounded planar domains, the following theorem characterizes convergence of Bergman kernels along decreasing sequences of domains in terms of Sobolev spaces.

\begin{theorem}[Skwarczy{\'n}ski]\label{ch3:th:Skw}
Let $\Omega$ be a bounded domain in $\mathbb{C}$ satisfying $\Omega=\overline{\Omega}^\circ$ and $|\partial \Omega|=0$. If $\Omega_j\downarrow \overline{\Omega}$, then $K_{\Omega_j} \rightarrow K_\Omega \iff W^{1,2}_0(\Omega) = W^{1,2}_0(\overline{\Omega})$.
\end{theorem}

Here $W^{1,2}_0(G)$, for an open set $G\subset\mathbb{R}^N$, denotes the closure of $C^\infty_0(G)$ in $W^{1,2}(G)$, with its elements extended by zero to $\mathbb{R}^N$. For a compact set $E$, we use the convention
$$
W^{1,2}_0(E)=\{u\in W^{1,2}(\mathbb{R}^N):u=0\ \text{a.e. on }\mathbb{R}^N\backslash E\}.
$$
These spaces are compared as subspaces of $W^{1,2}(\mathbb{R}^N)$.

The proof of Theorem \ref{ch3:th:Skw} is based on the following result.

\begin{theorem}[Havin;Hedberg]\label{ch3:th:Havin}
Let $E$ be a compact set in $\mathbb{C}$, and write $L^2_a(E)=L^2(E)\cap \mathcal O(E^\circ)$. Let $R^2(E)$ be the closure in $L^2(E)$ of the restrictions of functions holomorphic on a neighborhood of $E$. Then
$$
R^2(E)=L^2_a(E) \iff W^{1,2}_0(E^\circ) = W^{1,2}_0(E).
$$
\end{theorem}

\begin{proof}[Theorem \ref{ch3:th:Havin} $\Rightarrow$ Theorem \ref{ch3:th:Skw}]
Take $E=\overline{\Omega}$ in Theorem \ref{ch3:th:Havin} and combine it with Skwarczy{\'n}ski's theorem (Theorem \ref{th:Skwarcz}).
\end{proof}

The proof of Theorem \ref{ch3:th:Havin} will be postponed to the next section.

\section{The Case of Planar Domains}

The Bergman kernel on planar domains has its own distinctive aspects. We begin with Bergman's representation of the Riemann map through kernel functions.

\begin{proposition}\label{prop:SC}
Let $\Omega\subsetneq\mathbb{C}$ be simply connected, let $a\in\Omega$, and let $F:\Omega\mapsto\mathbb{D}$ be a Riemann map normalized by $F(a)=0$ and $F'(a)>0$. Then
\begin{enumerate}
\item[$(1)$] $F'(z)=\sqrt{\frac{\pi}{K_\Omega(a)}} K_\Omega(z,a)$;
\item[$(2)$] $F(z)= \sqrt{\frac{\pi}{K_\Omega(a)}} \int_a^z K_\Omega(\zeta,a) d\zeta$.  
\end{enumerate}
\end{proposition}

\begin{proof}
1. The transformation formula for the Bergman kernel gives
$$
K_\Omega(z,a)=K_{\mathbb{D}}(F(z),0)F'(z)\overline{F'(a)}=\frac1\pi F'(z)\overline{F'(a)}.
$$
Setting $z=a$ gives $F'(a)=\sqrt{\pi K_\Omega(a)}$; substitution proves the assertion.

2. Integrate the identity in (1).
\end{proof}

Proposition \ref{prop:SC} provides an effective numerical method for approximating Riemann maps. If $\Omega$ is a bounded simply connected domain with $C^0$ boundary, Theorem \ref{th:approximation} shows that $\mathcal O(\overline\Omega)$ is dense in $A^2(\Omega)$. Moreover, $\mathbb{C}\setminus\overline\Omega$ is connected, so Runge's approximation theorem implies that every function in $\mathcal O(\overline\Omega)$ can be approximated uniformly on $\overline\Omega$ by polynomials. Hence polynomials are dense in $A^2(\Omega)$. The Gram--Schmidt
process produces a complete orthonormal basis $\{h_j\}$ of $A^2(\Omega)$ consisting of polynomials. Thus
 \begin{eqnarray*}
 K_\Omega(a) & = & \lim_{n\rightarrow \infty} \sum_{j=1}^n |h_j(a)|^2,\\
 F(z) & = & \lim_{n\rightarrow \infty}\sum_{j=1}^n \sqrt{\frac{\pi}{K_\Omega(a)}} \int_a^z h_j(\zeta) \overline{h_j(a)} d\zeta.
 \end{eqnarray*}
 
Conversely, Bergman asked whether the kernel function could provide
a new proof of the Riemann mapping theorem.
Following the Garabedian--Schiffer kernel method, we give such a proof, starting from an extremal problem for area.

\begin{theorem}\label{th:kernel-extremal-riemann}
Let $\Omega\subsetneq\mathbb{C}$ be simply connected and $a\in\Omega$.
Then $K_\Omega(a)>0$, and
\begin{equation}\label{eq:kernel-extremal-riemann-map}
F(z)=\sqrt{\frac{\pi}{K_\Omega(a)}}
\int_a^z K_\Omega(\zeta,a)\,d\zeta
\end{equation}
is a biholomorphic map from $\Omega$ onto $\mathbb{D}$, with $F(a)=0$ and $F'(a)>0$.
\end{theorem}

We first establish the following Cauchy boundary pairing:

\begin{lemma}\label{lm:kernel-boundary-pairing}
Let $D\subset\subset\mathbb{C}$ be a Jordan domain with  real-analytic boundary, and let $a\in D$.
Suppose $f\in\mathcal O(D)$, $f'\in A^2(D)$, and, for some $c\in\mathbb{C}$,
\begin{equation}\label{eq:kernel-pairing-hypothesis}
\int_D h(z)\overline{f'(z)}\,dA(z)=c\,h(a),
\qquad h\in\mathcal O(\overline D).
\end{equation}
Then $f$ extends holomorphically to a neighborhood of $\overline D$, and there is a meromorphic function $g$ on that neighborhood such that
\begin{equation}\label{eq:kernel-pairing-conclusion}
g(z)-\frac{c}{\pi(z-a)}\in\mathcal O(\overline D),\qquad
g|_{\partial D}=\overline f.
\end{equation}
\end{lemma}

\begin{proof}
Complexifying a positively oriented real-analytic parametrization of the boundary and shrinking the neighborhood gives a holomorphic univalent parametrization
\[
\gamma:\{r <|\zeta|<r^{-1}\}\longrightarrow\mathbb{C},
\qquad 0<r<1,
\]
taking the unit circle to $\partial D$ and its inner side to the inner side of $D$.
This uses only analytic continuation of the parametrization, the local inverse function theorem, and compactness of the boundary.
On the inner annulus, write
\[
f_1(\zeta):=f(\gamma(\zeta))=\sum_{n\in\mathbb{Z}}b_n\zeta^n.
\]
The change-of-variables formula for area gives $\int_{r<|\zeta|<1}|f_1'|^2\,dA<\infty$.
Fix $r<r_0<1$. By orthogonality on circles and passage to the limit as the outer radius tends to $1$,
\[
\int_{r_0<|\zeta|<1}|f_1'(\zeta)|^2\,dA(\zeta)
=\pi\sum_{n\ne0}n|b_n|^2(1-r_0^{2n}).
\]
Every term is nonnegative, so $\sum_{n\ge1}n|b_n|^2<\infty$.
On the other hand, Cauchy's estimates on any fixed circle $|\zeta|=s_0\in(r,1)$ give
$|b_{-n}|\le M_{s_0}s_0^n$ ($n\ge1$).
Thus $\sum_{n\in\mathbb{Z}}|b_n|^2<\infty$, and, as $r\uparrow1$,
$f_1(re^{it})$ converges in $L^2(0,2\pi)$ to a boundary function.
Changing variables by $\gamma$ gives the $L^2$ boundary values of $f$ on $\partial D$.

Let $D_{1-\varepsilon}$ be the interior of $\gamma(\{|\zeta|=1-\varepsilon\})$, for sufficiently small $\varepsilon>0$.
For $h\in\mathcal O(\overline D)$, apply the Cauchy--Green formula on $D_{1-\varepsilon}$ and let $\varepsilon\downarrow0$.
The boundary integrals converge by the preceding $L^2$ convergence, and the area integrals converge since $f'\in L^2(D)$. Thus
\begin{equation}\label{eq:kernel-boundary-integral}
\frac1{2\pi i}\int_{\partial D}\overline{f(\zeta)}h(\zeta)\,d\zeta
=\frac1\pi\int_D h\,\overline{f'}\,dA
=\frac c\pi h(a).
\end{equation}
Here the boundary integral is interpreted using the $L^2$ boundary values.

Next,
define a holomorphic function on an exterior annular neighborhood by
\[
f^*(\gamma(\zeta))=\overline{f_1(1/\overline\zeta)},\qquad
1<|\zeta|<r^{-1}.
\]
Its exterior $L^2$ boundary values are exactly $\overline f$.
For $1<s<r^{-1}$, let $\Gamma_s=\gamma(\{|\zeta|=s\})$ be positively oriented. Fix $1<s_0<r^{-1}$ and put
\[
G(z)=\frac1{2\pi i}\int_{\Gamma_{s_0}}\frac{f^*(\zeta)}{\zeta-z}\,d\zeta,
\qquad z\in\operatorname{Int}\Gamma_{s_0}.
\]
Then $G$ is holomorphic inside $\Gamma_{s_0}$, in particular on a neighborhood of $\overline D$.
For $z\in\operatorname{Int}\Gamma_{s_0}\setminus\overline D$ and $s$ sufficiently close to $1$, apply Cauchy's formula on the annulus between $\Gamma_s$ and $\Gamma_{s_0}$ to obtain
\begin{eqnarray*}
f^*(z) & = & \frac1{2\pi i} \left(\int_{\Gamma_{s_0}}-\int_{\Gamma_s}\right) \frac{f^\ast(\zeta)}{\zeta-z}\,d\zeta\\
& = & G(z) - \frac1{2\pi i} \int_{\Gamma_s} \frac{f^\ast(\zeta)}{\zeta-z}\,d\zeta\\
& \to & G(z)-\frac1{2\pi i} \int_{\partial D} \frac{\overline{f(\zeta)}}{\zeta-z}\,d\zeta\quad (s\to 1)\\
& = &  G(z) - \frac{c}{\pi(a-z)},
\end{eqnarray*}
where the last equality follows from \eqref{eq:kernel-boundary-integral}. 
Thus 
$$g(z):=G(z)+\frac{c}{\pi(z-a)}
$$
 equals $f^*(z)$ on $\operatorname{Int}\Gamma_{s_0}\setminus\overline D$.

The boundary reflection
\[
\sigma(\gamma(\zeta))=\gamma(1/\overline\zeta)
\]
is an antiholomorphic involution fixing every boundary point.
Shrinking the neighborhood to avoid $a$, the function $\overline{g(\sigma(z))}$ is holomorphic across the boundary
and agrees with the original $f(z)$ on the inner side.
This gives the holomorphic extension of $f$. Letting $z$ approach the boundary yields $g=\overline f$.
\end{proof}

\begin{proof}[Proof of Theorem \ref{th:kernel-extremal-riemann}]
The proof has three steps. First assume that $\Omega$ is a Jordan domain with  real-analytic boundary.
The extremal property of the Bergman kernel gives
\[
m_\Omega(a)^2:=\min\left\{\int_\Omega |h|^2\,dA:h\in A^2(\Omega),\ h(a)=1\right\}
=\frac1{K_\Omega(a)}
\]
and
$
h_0(z)=\frac{K_\Omega(z,a)}{K_\Omega(a)}
$
is the extremal function.
By simple connectivity, we may define
\[
f(z)=\int_a^z h_0(\zeta)\,d\zeta.
\]
Then $f(a)=0$, $f'(a)=1$, and $\int_\Omega|f'|^2\,dA=m_\Omega(a)^2$.
For $h\in A^2(\Omega)$ with $h(a)=0$, minimizing the squared norm of $h_0+th$ at $t=0$ gives
$\int_\Omega h\,\overline{f'}\,dA=0$.
Decomposing a general $h$ into $h-h(a)f'$ and $h(a)f'$ yields
\begin{equation}\label{eq:kernel-extremal-orthogonality}
\int_\Omega h\,\overline{f'}\,dA=m_\Omega(a)^2 h(a),\qquad h\in A^2(\Omega).
\end{equation}

Put $R^2=m_\Omega(a)^2/\pi$. By Lemma \ref{lm:kernel-boundary-pairing}, $f$ extends holomorphically across the boundary,
and there is a meromorphic function
\[
g(z)=\frac{R^2}{z-a}+G(z),\qquad G\in\mathcal O(\overline \Omega),
\qquad g|_{\partial \Omega}=\overline f.
\]
Since $f(a)=0$, the singularity of $fg$ at $a$ is removable, and
\[
(fg)(a)=R^2f'(a)=R^2.
\]
On the boundary, $fg=|f|^2$ is real. Thus $\operatorname{Im}(fg)$ has zero boundary values
and vanishes identically by the maximum principle. The holomorphic function $fg$ is real-valued and hence constant. Therefore
\begin{equation}\label{eq:kernel-extremal-product}
fg\equiv R^2\quad\hbox{on }\Omega.
\end{equation}
This also shows that $|f|=R$ on $\partial\Omega$ and that $a$ is the only zero of $f$ in $\Omega$,
with multiplicity one.
By the maximum modulus principle, $f(\Omega)\subset\{|w|<R\}$.
For any $|w|<R$, we have $|w|<|f|$ on the boundary, so Rouch\'e's theorem implies that $f-w$ and $f$ have the same number of zeros,
namely one, counted with multiplicity.
Thus $f$ maps $\Omega$ biholomorphically onto the disc of radius $R$.
Normalization gives
\[
F(z)=\frac{f(z)}R
=\sqrt{\frac\pi{K_\Omega (a)}}\int_a^zK_\Omega (\zeta,a)\,d\zeta,
\qquad F'(a)=\sqrt{\pi K_\Omega(a)}>0.
\]

In the second step, consider a general bounded
domain $\Omega$. Choose an exhaustion by Jordan domains with  real-analytic boundary:
\[
a\in \Omega_1\subset\subset  \Omega_2\subset\subset \cdots,\qquad \bigcup_j \Omega_j=\Omega.
\]
Such an exhaustion can be constructed geometrically in the plane: surround a prescribed compact set by a connected finite union of small squares, fill in the bounded holes,
smooth the boundary, and approximate its parametrization sufficiently closely in $C^1$ by trigonometric polynomials.
Simple connectivity ensures that the filled holes remain inside $\Omega$; including the closure of the preceding domain at each stage ensures nesting.

For each $\Omega_j$, there is a normalized biholomorphic map $F_j:\Omega_j\to\mathbb{D}$.
By Theorem \ref{th:Ramadanov}, the kernels converge locally uniformly. Integration along compact paths therefore gives
\[
F_j\longrightarrow F(z):=\sqrt{\frac\pi{K_\Omega(a)}}
\int_a^zK_\Omega(\zeta,a)\,d\zeta
\]
locally uniformly on $\Omega$.
Since $\Omega$ is bounded, $K_\Omega(a)>0$, so $F'(a)=\sqrt{\pi K_\Omega(a)}>0$.
By Hurwitz's theorem, a nonconstant limit of univalent functions is univalent. The maximum modulus principle also gives $F(\Omega)\subset\mathbb{D}$.

To verify surjectivity, put $Q_j=F_j^{-1}:\mathbb{D}\to\Omega_j\subset\Omega$.
Boundedness and Montel's theorem give a locally uniformly convergent subsequence $Q_j\to Q$, with $Q(0)=a$.
For every $c\notin\Omega$, the functions $Q_j-c$ have no zeros, and $Q(0)-c\ne0$.
Hurwitz's theorem implies that $Q-c$ has no zeros, so $Q(\mathbb{D})\subset\Omega$.
Passing to the limit in $F_j(Q_j(\zeta))=\zeta$ gives $F(Q(\zeta))=\zeta$.
Thus $F$ is surjective, proving the bounded case.

Finally, let $\Omega\subsetneq\mathbb{C}$ be an unbounded simply connected domain.
Choose $b\notin\Omega$. Simple connectivity gives a holomorphic branch $s(z)=\sqrt{z-b}$.
This function is univalent and satisfies $s(\Omega)\cap(-s(\Omega))=\varnothing$.
Choose $B(c,r)\subset s(\Omega)$. Then $B(-c,r)$ is disjoint from $s(\Omega)$, and hence
\[
T(z)=\frac1{s(z)+c}
\]
is bounded and univalent. Apply the second step to the bounded simply connected domain $T(\Omega)$, compose with $T$, and rotate the disc
to obtain a biholomorphic map $H:\Omega\to\mathbb{D}$ with $H(a)=0$ and $H'(a)>0$.

The existence of the map has now been established independently. Only at this point do we invoke the transformation formula for the Bergman kernel, obtaining
\[
K_\Omega(z,a)=\frac1\pi H'(z)\overline{H'(a)},\qquad
K_\Omega(a)=\frac{H'(a)^2}{\pi}>0.
\]
Integration gives $H=F$, where $F$ is defined by \eqref{eq:kernel-extremal-riemann-map}.
\end{proof}

\needspace{5\baselineskip}

  \begin{proposition}\label{prop:simply-connected}
If $\Omega\subsetneq\mathbb{C}$ is simply connected, then
        $$
 \frac1{16\pi}\cdot \delta_\Omega(z)^{-2} \le K_\Omega(z)\le \frac1{\pi}\cdot \delta_\Omega(z)^{-2}.
$$
        \end{proposition}
        
        \begin{proof}
It suffices to prove the first inequality. Given $z\in\Omega$, consider a Riemann map $F:\Omega\rightarrow\mathbb{D}$ with
$
 F(z)=0,\,F'(z)>0.
$
The classical Koebe one-quarter theorem states that if $f$ is univalent on $\mathbb{D}$, $f(0)=0$, and $f'(0)=1$, then $f(\mathbb{D})\supset\Delta(0,1/4)$. Put $G=F^{-1}$. If we take
$$
f(\zeta)=\frac{G(\zeta)-G(0)}{G'(0)},
$$
Koebe's theorem immediately gives
$$
\Delta\bigl(z,|G'(0)|/4\bigr)\subset\Omega.
$$
Therefore $\delta_\Omega(z)\ge|G'(0)|/4$, that is,
$$
\frac14\le\delta_\Omega(z)|F'(z)|.
$$
Thus
$$
K_\Omega(z)=K_{\mathbb{D}}(0) |F'(z)|^2\ge \frac1{16 \pi}\cdot \delta_\Omega(z)^{-2}.
$$  
        \end{proof}
 
\begin{definition}
A domain $\Omega\subset\mathbb{C}$ is called hyperbolic in the potential-theoretic sense (or simply hyperbolic in this section) if $E:=\mathbb{C}\setminus\Omega$ is nonpolar.
\end{definition} 

Every hyperbolic domain $\Omega$ has a Green function $g_\Omega<0$ such that
\begin{equation}\label{eq:GreenDef}
g_\Omega(z,w)=\log|z-w|+h_\Omega(z,w),
\end{equation}
where $h_\Omega\in C^\infty(\Omega\times\Omega)$ is harmonic separately in $z$ and $w$, and $h_\Omega(z,w)=h_\Omega(w,z)$.

\begin{theorem}\label{th:hyperbolicBergman}
Let $\Omega$ be hyperbolic. For every $w\in\Omega$ and every nonnegative integer $k$, there exists $f\in A^2(\Omega)$ such that
$$
f^{(j)}(w)=0,\ \forall\,j<k;\ \ f^{(k)}(w)=1.
$$
In particular, $K_\Omega(z)>0$ for all $z\in\Omega$, and $A^2(\Omega)$ is infinite-dimensional.
\end{theorem}

\begin{proof}
Put $\Omega^\ast:=\Omega\setminus\{w\}$. Given $\varepsilon>0$ and a nonnegative integer $k$, define
$$
\varphi_{k,\varepsilon}(z):=2(k+1)g_\Omega(z,w)-\log(-g_\Omega(z,w)+1) +\varepsilon |z|^2.  
$$
On $\Omega^\ast$ we have
\begin{equation}\label{eq:hyperbolic}
\frac{\partial^2 \varphi_{k,\varepsilon}}{\partial z\partial\bar{z}} = \frac{|{\partial g_\Omega(z,w)}/{\partial z}|^2}{(-g_\Omega(z,w)+1)^2} +\varepsilon,
\end{equation}
In particular, $\varphi_{k,\varepsilon}$ is strictly subharmonic on $\Omega^\ast$. Choose a small disc $\Delta\Subset\Omega$ about $w$, then choose $R\gg1$ so that $g_\Omega(\cdot,w)>-R$ on $\partial\Delta$. The component $U_R$ of $\{g_\Omega(\cdot,w)<-R\}$ containing $w$ is then relatively compact in $\Delta$.
Choose a smooth cutoff $\kappa:\mathbb{R}\to[0,1]$ with $\kappa|_{(-\infty,1]}=0$ and $\kappa|_{[2,\infty)}=1$.
Below, $\kappa(-g_\Omega/R)$ is used only in $U_R$ and is set to zero on $\Omega\setminus U_R$; the transition region in the integrals is also restricted to $U_R$. By
H\"ormander's theorem, the equation
$$
\bar{\partial} u = (z-w)^k\bar{\partial} \kappa(-g_\Omega(\cdot,w)/R)
$$ 
has a solution $u_{k,\varepsilon}$ on $\Omega^\ast$ satisfying
\begin{eqnarray*}
\int_{\Omega^\ast} |u_{k,\varepsilon}|^2 e^{-\varphi_{k,\varepsilon}} & \le & \int_{\Omega^\ast} \left(\frac{\partial^2 \varphi_{k,\varepsilon}}{\partial z\partial\bar{z}}\right)^{-1} |z-w|^{2k} \left|\frac{\partial \kappa(-g_\Omega(\cdot,w)/R)}{\partial z}\right|^2 e^{-\varphi_{k,\varepsilon}} \\
& \le & \frac{\sup|\kappa'|^2}{R^2}\int_{-2R<g_\Omega(\cdot,w)<-R} (-g_\Omega(z,w)+1)^2 |z-w|^{2k}  e^{-\varphi_{k,\varepsilon}} \\
& \le & C=C(k,w,R).
\end{eqnarray*}
Since $u_{k,\varepsilon}$ is holomorphic for $0<|z-w|\ll1$, the preceding integrability condition implies that it extends holomorphically across $w$ and satisfies
$$
u_{k,\varepsilon}^{(j)}(w)=0,\ \forall\,j\le k.
$$
Thus
$$
f_{k,\varepsilon}:= (z-w)^k \kappa(-g_\Omega(\cdot,w)/R)-u_{k,\varepsilon}\in \mathcal O(\Omega)
$$
and satisfies
$$
f^{(j)}_{k,\varepsilon}(w)=0,\ \forall\,j<k;\ \ f^{(k)}_{k,\varepsilon}(w)=k!,
$$
as well as
$$
\int_\Omega |f_{k,\varepsilon}|^2 e^{-\varepsilon |z|^2} \le C
$$
(note that $\varphi_{k,\varepsilon}-\varepsilon|z|^2<0$). A subsequential limit of $\{\frac1{k!}f_{k,\varepsilon}\}$ as $\varepsilon\downarrow0$ gives the desired $f$.
\end{proof}

Write
\[\frac{\partial^2g_\Omega(z,w)}{\partial{z}\partial\bar{w}}:=\frac{\partial^2h_\Omega(z,w)}{\partial{z}\partial\bar{w}},\]
where $h_\Omega$ is as in \eqref{eq:GreenDef}.
Note that $\partial^2g_\Omega(z,w)/\partial z\partial\bar w$ is holomorphic separately in $z,\bar w$.

\begin{theorem}[Schiffer {\cite{Schiffer1946}}]\label{th:Schiffer}
\begin{equation}\label{eq:Schiffer}
K_\Omega(z,w)=\frac{2}{\pi}\frac{\partial^2g_\Omega(z,w)}{\partial{z}\partial\bar{w}}.
\end{equation}
\end{theorem}

\begin{proof}
First assume that $\Omega$ is bounded with real-analytic boundary. Since $g_\Omega(z,w)$ is harmonic in $z\neq w$ and vanishes on $\partial\Omega$, Schwarz reflection gives
$$
g_\Omega(z,w)\in{C^\infty(\overline{\Omega}\times\Omega\setminus \{z=w\})}.
$$
By symmetry, $g_\Omega\in C^\infty(\Omega\times\overline\Omega\setminus\{z=w\})$. Given $z\in\Omega$, put
 $$
 \Omega_\varepsilon:=\Omega\setminus\Delta(z,\varepsilon),
 $$
where $\varepsilon$ is sufficiently small. Let $f\in\mathcal O(\Omega)\cap C^1(\overline\Omega)$. Stokes's formula gives
\begin{eqnarray*}
\frac{i}{2}\int_{\Omega_\varepsilon}f(w)\frac{\partial^2g_\Omega(z,w)}{\partial{z}\partial\bar{w}}dw\wedge{d\bar{w}}
&=&-\frac{i}{2}\int_{\Omega_\varepsilon}d_w\bigg(f(w)\frac{\partial{g_\Omega(z,w)}}{\partial{z}}dw\bigg)\\
&=&-\frac{i}{2}\int_{\partial\Omega}f(w)\frac{\partial{g_\Omega(z,w)}}{\partial{z}}dw\\
& &+\frac{i}{2}\int_{\partial\Delta(z,\varepsilon)}f(w)\frac{\partial{g_\Omega(z,w)}}{\partial{z}}dw.
\end{eqnarray*}
Since $g_\Omega(z,\cdot)|_{\partial\Omega}\equiv0$, we have $\frac{\partial g_\Omega(z,\cdot)}{\partial z}|_{\partial\Omega}\equiv0$, and hence
\begin{eqnarray*}
\frac{i}{2}\int_{\Omega_\varepsilon}f(w)\frac{\partial^2g_\Omega(z,w)}{\partial{z}\partial\bar{w}}dw\wedge{d\bar{w}}
&=&\frac{i}{2}\int_{\partial\Delta(z,\varepsilon)}f(w)\frac{\partial{g_\Omega(z,w)}}{\partial{z}}dw\\
&=&\frac{i}{2}\int_{\partial\Delta(z,\varepsilon)}f(w)\frac{\partial\log|z-w|}{\partial{z}}dw\\
& &+\frac{i}{2}\int_{\partial\Delta(z,\varepsilon)}f(w)\frac{\partial{h(z,w)}}{\partial{z}}dw\\
&=&\frac{i}{4}\int_{\partial\Delta(z,\varepsilon)}\frac{f(w)}{z-w}dw+O(\varepsilon)\\
&=&\frac{\pi}{2}f(z)+O(\varepsilon).
\end{eqnarray*}
Letting $\varepsilon\rightarrow0$ gives
$$
f(z) = \frac{i}{2}\int_{\Omega} f(w) \left[\frac{2}{\pi}\frac{\partial^2g_\Omega(z,w)}{\partial{z}\partial\bar{w}}\right]dw\wedge{d\bar{w}}.
$$
By Theorem \ref{th:approximation}, this identity holds for every $f\in A^2(\Omega)$. Uniqueness of the reproducing kernel gives \eqref{eq:Schiffer}.

For a general hyperbolic domain $\Omega$, take domains $\{\Omega_j\}$ with real-analytic boundary such that $\Omega_j\subset\subset\Omega$ and $\Omega_j\uparrow\Omega$ as $j\rightarrow\infty$. By Ramadanov's theorem (Theorem \ref{th:Ramadanov}), $K_{\Omega_j}(z,w)\rightarrow K_\Omega(z,w)$. Since $g_{\Omega_j}(z,w)$ converges uniformly on compact subsets to $g_\Omega(z,w)$, we have $h_{\Omega_j}(z,w)\downarrow h_\Omega(z,w)$. Both $h_{\Omega_j}$ and $h_\Omega$ are harmonic in each variable; successive applications of the Poisson integral formula give
\begin{eqnarray*}
&& 4\pi^2 h_{\Omega_j}(z,w)\\
& = & \int_0^{2\pi}\int_0^{2\pi} h_{\Omega_j}(z_0+\varepsilon e^{i\theta},w_0+\varepsilon e^{i\vartheta})\frac{\varepsilon^2-|z-z_0|^2}{|\varepsilon e^{i\theta}-(z-z_0)|^2} \frac{\varepsilon^2-|w-w_0|^2}{|\varepsilon e^{i\vartheta}-(w-w_0)|^2}d\theta d\vartheta.
\end{eqnarray*}
The same formula holds for $h_\Omega$. Dominated convergence then gives
\[
\frac{\partial^2h_{\Omega_j}(z,w)}{\partial{z}\partial\bar{w}}\rightarrow\frac{\partial^2h_\Omega(z,w)}{\partial{z}\partial\bar{w}},\ \ \ (j\rightarrow\infty).
\]
and \eqref{eq:Schiffer} follows.

\end{proof}

We now give some applications of Schiffer's theorem. First, we give a new proof of Theorem \ref{th:polar} when $n=1$. Let $E$ be a closed polar subset of a plane domain $\Omega$. For $w\in\Omega\setminus E$, we have
\[
g_{\Omega\setminus{E}}(\cdot,w)\in{SH^-(\Omega\setminus{E})=SH^-(\Omega)},
\]
where $SH(\Omega)$ denotes the subharmonic functions on $\Omega$ and $SH^-(\Omega)$ the negative subharmonic functions.
The Green function has the extremal characterization
$$
g_\Omega(z,w) =\sup\left\{u(z): u\in SH^-(\Omega),\,\limsup_{z\rightarrow w}(u(z)-\log |z-w|)<\infty\right\},
$$
so $g_\Omega(z,w)\ge g_{\Omega\setminus E}(z,w)$ for every $z\in\Omega\setminus E$. Conversely, $\Omega\setminus E\subset\Omega$ implies $g_\Omega(z,w)\le g_{\Omega\setminus E}(z,w)$. Hence equality holds. Schiffer's theorem gives $K_\Omega(z,w)=K_{\Omega\setminus E}(z,w)$ for $z,w\in\Omega\setminus E$.

Next we give an application to discrete groups. Let $\Omega\subset\mathbb{C}$ be hyperbolic. By uniformization, its universal covering space is $\mathbb{D}$. Thus there is a discrete subgroup $\Gamma\subset\mathrm{Aut\,}\mathbb{D}$, isomorphic to $\pi_1(\Omega)$, such that $\Omega=\mathbb{D}/\Gamma$ (such groups are historically called Fuchsian groups). Denote the covering map by $\varpi:\mathbb{D}\rightarrow\Omega$.

\begin{theorem}[Summation formula]\label{th:BergmanSum}
\begin{equation}\label{eq:BergmanSum}
K_{\Omega}(\varpi(z),\varpi(w))\, \varpi'(z)\overline{\varpi'(w)}=\sum_{T\in \Gamma} K_{\mathbb{D}}(z,T(w))\,\overline{T'(w)}.
\end{equation}
\end{theorem}

We need the following theorem of Myrberg.

\begin{theorem}[Myrberg {\cite{Myrberg1933}}]\label{th:Myrberg}
The following properties are equivalent:
\begin{enumerate}
\item[$(1)$] $\Omega=\mathbb{D}/\Gamma$ is hyperbolic;
\item[$(2)$] $\sum_{T\in\Gamma}(1-|T(w)|)<\infty$ for some $w\in\mathbb{D}$;
\item[$(3)$] $\sum_{T\in\Gamma}(1-|T(w)|)<\infty$ for every $w\in\mathbb{D}$.
\end{enumerate}
In this case, the Green function satisfies the summation formula
\begin{equation}\label{eq:GreenSum}
g_{\Omega}(\varpi(z),\varpi(w)) = \sum_{T\in \Gamma} g_{\mathbb{D}}(z,T(w)).  
\end{equation}
\end{theorem}

\begin{remark}
Poincar\'e proved that, for every Fuchsian group $\Gamma$ and every $\alpha>1$,
$$
\sum_{T\in \Gamma} \left(1-|T(w)|\right)^\alpha <\infty,\ \ \ \forall\,w\in \mathbb{D}.
$$
\end{remark}

The proof of Theorem \ref{th:Myrberg} uses the following classical result.

\begin{theorem}[Blaschke]\label{th:Blaschke_1}
Suppose $0<|a_n|<1$ and $\sum_{n=1}^\infty(1-|a_n|)<\infty$. Then the Blaschke product associated with $\{a_n\}$,
$$
B(z)=\prod_{n=1}^\infty \frac{|a_n|}{a_n}\cdot \frac{a_n-z}{1-\bar{a}_nz}
$$ 
is bounded and holomorphic on $\mathbb{D}$, and satisfies $|B(e^{i\theta})|=1$ for almost every $\theta$.

Conversely, if $\{a_1,a_2,\cdots\}$ lists all zeros, counted with multiplicity, of a nonzero bounded holomorphic function $f$ on $\mathbb{D}$, then
$$
(1-|a_1|)+(1-|a_2|)+\cdots <\infty.
$$
\end{theorem}

\begin{proof}[Proof of Theorem \ref{th:Myrberg}]
$(2)\Rightarrow(3)$. Conformal invariance of the Poincar\'e distance gives
$$
\log \frac{1+\left|\frac{T(w)-T(0)}{1-\overline{T(0)}T(w)}\right|}{1-\left|\frac{T(w)-T(0)}{1-\overline{T(0)}T(w)}\right|} = \log \frac{1+|w|}{1-|w|},
$$
Hence
\begin{eqnarray*}
1-|w|^2  & = & 1- \left|\frac{T(w)-T(0)}{1-\overline{T(0)}T(w)}\right|^2\ \ =\ \   \frac{(1-|T(w)|^2)(1-|T(0)|^2)}{\left|1-\overline{T(0)} T(w)\right|^2}\\
& \le & \min\left\{\frac{4(1-|T(0)|)(1-|T(w)|)}{(1-|T(0)|)^2},  \frac{4(1-|T(0)|)(1-|T(w)|)}{(1-|T(w)|)^2}\right\},
\end{eqnarray*}
so that
$$
\frac{(1-|w|^2)(1-|T(0)|)}4 \le 1-|T(w)| \le \frac{4(1-|T(0)|)}{1-|w|^2}.
$$
and the implication follows.

$(3)\Rightarrow(2)$ is immediate.

$(3)\Rightarrow(1)$. Put
$$
g^\ast(z):=\sum_{T\in \Gamma} g_\mathbb{D}(z,T(w)),\ \ \ z\in \mathbb{D}.  
$$
Note that
\begin{eqnarray}\label{eq:GreenConverge}
2g_\mathbb{D}(z,T(w)) & = & \log \frac{|z-T(w)|^2}{|1-\overline{T(w)}z|^2}\ \ =\ \ \log\left(1-\frac{(1-|z|^2)(1-|T(w)|^2)}{|1-\overline{T(w)}z|^2}\right)\nonumber\\
& \ge & \log \left(1-\frac{(1+|z|)(1-|T(w)|^2)}{1-|z|}\right)\nonumber\\
& \sim & -\frac{(1+|z|)(1-|T(w)|^2)}{1-|z|}\ \ \ (\text{as }T(w)\rightarrow\partial\mathbb{D}),
\end{eqnarray}
so $g^\ast\in SH^-(\mathbb{D})$ and is harmonic on $\mathbb{D}\backslash\{T(w):T\in\Gamma\}$. Moreover, for every $S\in\Gamma$,
$$
g^\ast(S(z))=\sum_{T\in \Gamma} g_\mathbb{D}(S(z),T(w)) = \sum_{T\in \Gamma} g_\mathbb{D}(z,S^{-1}\circ T(w))=g^\ast(z),
$$
Thus $g^\ast$ is $\Gamma$-invariant and induces on the quotient $\Omega$ a negative subharmonic function
$g_\ast(\varpi(z)):=g^\ast(z)$, harmonic on $\Omega\setminus\{\varpi(w)\}$ with a logarithmic singularity at $\varpi(w)$. The extremal property of the Green function gives
\begin{equation}\label{eq:Myrberg_1}
g_\Omega(\varpi(z),\varpi(w))\ge g_\ast(\varpi(z))=g^\ast(z).  
\end{equation}
In particular, $\Omega$ is hyperbolic.

We next verify \eqref{eq:GreenSum}. Choose a sufficiently small closed disc $E$ centered at $w$ so that the sets $T(E)$, $T\in\Gamma$, are pairwise disjoint. Put $D=\mathbb{D}\setminus\bigcup_{T\in\Gamma}T(E)$. The maximum principle and $\Gamma$-invariance of $g^\ast$ give
\begin{equation}\label{eq:Myrberg_2}
0<-g^\ast(z)\le \sup_{\bigcup\partial T(E)} (-g^\ast)=\sup_{\partial E}(-g^\ast),\ \ \ \forall\,z\in D.
\end{equation}
Thus, if we put
$$
H(z)=g_\Omega(\varpi(z),\varpi(w))-g^\ast(z),\ \ \ z\in \mathbb{D},
$$
then \eqref{eq:Myrberg_1} and \eqref{eq:Myrberg_2} imply that $H$ is bounded on $D$. Since $H$ is $\Gamma$-invariant and harmonic on $\bigcup_{T\in\Gamma}T(E)$, it is bounded and harmonic on $\mathbb{D}$. Fatou's theorem implies that $\lim_{r\rightarrow1-}H(re^{i\theta})=H(e^{i\theta})$ exists for almost every $\theta$. Since $\sum_{T\in\Gamma}(1-|T(w)|)<\infty$ and
$$
0\le H(z)<-g^\ast(z)=-\log \prod_{T\in \Gamma} \left|\frac{z-T(w)}{1-\overline{T(w)}z}\right|,
$$ 
Blaschke's theorem gives $H(e^{i\theta})=0$ for almost every $\theta$. The Poisson formula then gives $H\equiv0$, proving \eqref{eq:GreenSum}.

$(1)\Rightarrow(3)$. Since $g_\Omega$ exists, $g^\ast(z):=g_\Omega(\varpi(z),\varpi(w))$ is negative and subharmonic on $\mathbb{D}$, harmonic away from $\{T(w):T\in\Gamma\}$ and having logarithmic singularities there. The monodromy theorem allows us to construct $f\in\mathcal O(\mathbb{D})$ such that
$$
|f(z)|=e^{g^\ast(z)},\ \ \ z\in \mathbb{D}.
$$
Note that $|f|<1$ and its zeros are exactly $T(w)$, $T\in\Gamma$. Blaschke's theorem proves the implication.
\end{proof}

\begin{proof}[Proof of Theorem \ref{th:BergmanSum}]
By \eqref{eq:GreenConverge}, the series in \eqref{eq:GreenSum} converges uniformly on compact subsets of $\mathbb{D}\setminus\{T(w):T\in\Gamma\}$. Applying $\frac{\partial^2}{\partial z\partial\bar w}$ to both sides of \eqref{eq:GreenSum} and using \eqref{eq:Schiffer} gives \eqref{eq:BergmanSum}.
\end{proof}

\begin{definition}
Let $\Omega\subset\mathbb{C}$ be hyperbolic. Given $w\in\Omega$ and $m\in\mathbb{Z}^+$, the function
$$
K_{\Omega,w}^{(m)}(z):=\sup\left\{|f(z)|^2: f\in A^2(\Omega), \|f\|=1, \mathrm{ord}_w\, f\ge m\right\}
$$
is called the partial Bergman kernel of order $m$ on $\Omega$.
\end{definition}

\begin{theorem}[Nivoche \cite{Nivoche},  Xu Wang]
$$
g_\Omega(z,w)=\lim_{m\rightarrow \infty} \frac1{2m} \log K_{\Omega,w}^{(m)}(z). 
$$
\end{theorem}

\begin{proof}
Let $f\in A^2(\Omega)$ satisfy $\|f\|=1$ and $\mathrm{ord}_w\,f\ge m$. Choose domains
$$
\Omega_1\subset\subset \Omega_2\subset\subset \cdots, \ \Omega=\bigcup \Omega_j.
$$
Bergman's inequality gives
$$
\sup_{\Omega_j} |f|^2\le C_j \|f\|^2 =C_j,
$$
so that
$$
\sup_{\Omega_j}  K_{\Omega,w}^{(m)}\le C_j.
$$
Thus
$$
\frac1{2m}\log \frac{K_{\Omega,w}^{(m)}}{C_j}\in SH^-(\Omega_j)
$$
and
$$
\frac1{2m}\log \frac{K_{\Omega,w}^{(m)}(z)}{C_j} = \log |z-w| + O(1)\ \ \ (z\rightarrow w).
$$
The extremal property of the Green function gives
$$
g_{\Omega_j}(z,w)\ge \frac1{2m}\log \frac{K_{\Omega,w}^{(m)}(z)}{C_j}.
$$
Letting first $m\rightarrow\infty$ and then $j\rightarrow\infty$ yields
$$
g_\Omega(z,w)\ge \limsup_{m\rightarrow \infty} \frac1{2m}\log K_{\Omega,w}^{(m)}(z).
$$

On the other hand, the Ohsawa--Takegoshi extension theorem \cite{Chen2022} gives, for each $z\in\Omega\setminus\{w\}$, a function $f\in\mathcal O(\Omega)$ with $f(z)=1$ and
$$
\int_{\Omega} |f|^2 e^{-2m g_{\Omega}(\cdot,w)} \le C_z e^{-2m g_{\Omega}(z,w)},
$$
where $C_z$ depends on the fixed domain and the point $z$, but not on $m$.
Thus $f^{(k)}(w)=0$ for every $k<m$, and
$$
\int_{\Omega} |f|^2  \le C_z e^{-2m g_{\Omega}(z,w)},
$$ 
Hence
$$
 \frac1{2m}\log K_{{\Omega},w}^{(m)}(z)\ge  \frac1{2m}\log \frac{|f(z)|^2}{\|f\|^2}\ge g_{\Omega}(z,w)-\frac1{2m}\log C_z,
$$
so that
$$
\liminf_{m\rightarrow \infty}  \frac1{2m}\log K_{{\Omega},w}^{(m)}(z)\ge g_\Omega(z,w).
$$
\end{proof}

If $\Omega$ is hyperbolic, Theorem \ref{th:hyperbolicBergman} allows us to define two K\"ahler metrics on it:
$$
ds^2_{1}=K_\Omega(z)dz\otimes d\bar{z},\ \ ds_{2}^2:=\frac{\partial^2\log K_\Omega(z)}{\partial z\partial \bar z}dz\otimes d\bar{z}.
$$
By Proposition \ref{prop:Transformation}, both $ds_1^2$ and $ds_2^2$ are invariant under biholomorphic maps. They are called the first and second Bergman invariant metrics, respectively.

Note that the Gaussian curvature of $ds_1^2$ is
$$
\mathscr{K}= -\frac2{K_\Omega(z)} \,\frac{\partial^2\log K_\Omega(z)}{\partial z\partial \bar z}=-2\,\frac{ds^2_{2}}{ds^2_{1}}.
$$
Thus $\mathscr K<0$ everywhere.

\begin{problem}
Under what conditions is $ds_1^2$ complete? What geometric condition characterizes the property that $\mathscr K$ lies between two negative constants?
\end{problem}

Finally,  we give the proof of Theorem \ref{ch3:th:Havin}.  Orthogonal complements below are taken in $L^2(\Omega)$ or $L^2(E)$, as indicated; functions are extended by zero when considered on $\mathbb{C}$.

\begin{lemma}\label{ch3:lm:Havin1}
Let $\Omega\subset\mathbb{C}$ be bounded and open, and let $g\in L^2(\Omega)$. Then $g\in A^2(\Omega)^\bot\subset L^2(\Omega)$ if and only if there exists a unique $\phi\in W^{1,2}_0(\Omega)$ such that $g=\partial\phi/\partial z$ in the sense of distributions. The equality also holds on $\mathbb{C}$ after zero extension.
\end{lemma}

\begin{proof}
By Weyl's lemma, a function $f\in L^2(\Omega)$ is holomorphic on $\Omega$ if and only if $\bar{\partial} f=0$ holds on $\Omega$ in the sense of distributions, that is,
$$
\int_\Omega f \cdot \overline{\partial \phi/\partial z}=0,\ \ \ \phi\in C^\infty_0(\Omega).
$$
Taking orthogonal complements, $g\in A^2(\Omega)^\bot \subset L^2(\Omega)$ if and only if there exists $\{\phi_j\}\subset C^\infty_0(\Omega)$ such that
$$
\int_\Omega \left|g-\partial \phi_j/\partial z\right|^2 \rightarrow 0\ \ \ (j\rightarrow \infty).
$$
Since
\begin{equation}\label{ch3:eq:dVSdbar}
\begin{aligned}
\int_\Omega |\partial \phi/\partial {z}|^2
&=\int_\Omega \partial \phi/\partial z \cdot \overline{\partial \phi/\partial z}\\
&= - \int_\Omega \frac{\partial^2 \phi}{\partial z\partial \bar{z}}\cdot \bar{\phi}\\
&= \int_\Omega |\partial \phi/\partial \bar{z}|^2,\ \ \ \forall\,\phi\in C^\infty_0(\Omega),
\end{aligned}
\end{equation}
applying this identity to $\phi_j-\phi_k$ shows that both $\{\partial\phi_j/\partial z\}$ and $\{\partial \phi_j/\partial \bar{z}\}$ are Cauchy sequences in $L^2(\Omega)$. Hence the vector gradients $\{\nabla\phi_j\}$ are Cauchy in $L^2(\Omega;\mathbb{C}^2)$. Poincar\'e's inequality shows that $\{\phi_j\}\subset W^{1,2}(\Omega)$ is a Cauchy sequence. Thus $\phi_j\rightarrow \phi\in W^{1,2}_0(\Omega)$, with
$\partial \phi/\partial z=\lim_{j\rightarrow \infty} \partial \phi_j/\partial z=g$.

Next, we verify uniqueness. Suppose that there is another $\widetilde{\phi}\in W^{1,2}_0(\Omega)$ such that $g=\partial \widetilde{\phi}/\partial z$. Then $\psi:=\phi-\widetilde{\phi}\in W^{1,2}_0(\Omega)$ and $\partial\psi/\partial z=0$. By density, (\ref{ch3:eq:dVSdbar}) also holds for $\psi$, and gives
$$
\int_\Omega |\partial \psi/\partial z|^2 = \int_\Omega |\partial\psi/\partial\bar z|^2=0.
$$
Thus $\nabla\psi=0$, and Poincar\'e's inequality yields $\psi\equiv 0$.

Finally, suppose that $g=\partial \phi/\partial z$, where $\phi\in W^{1,2}_0(\Omega)$. Choose $\phi_n\in C^\infty_0(\Omega)$ such that $\partial \phi_n/\partial z\rightarrow \partial\phi/\partial z$ in $L^2(\Omega)$. The first identity in the proof then gives $g\in A^2(\Omega)^\bot\subset L^2(\Omega)$.
\end{proof}

\begin{lemma}\label{ch3:lm:Havin2}
Let $E$ be a compact set in $\mathbb{C}$ and $g\in L^2(E)$, extended by zero outside $E$. Then $g\in R^2(E)^\bot \subset L^2(E)$ if and only if $g=\partial\phi/\partial z$ on $\mathbb{C}$ for some $\phi\in W^{1,2}_0(E)$.
\end{lemma}

\begin{proof}
First, for any bounded open set $U\supset E$, the restrictions to $E$ of functions in $A^2(U)$ belong to $R^2(E)$. Thus, if $g\in R^2(E)^\bot$, its zero extension satisfies
$$
\int_U g\cdot\bar f=0,\qquad f\in A^2(U).
$$
Conversely, suppose these orthogonality relations hold for every such $U$. For every $f\in R^2(E)$, there exists a sequence $f_j\in \mathcal O(E)$ such that $\int_E |f-f_j|^2\rightarrow 0\ (j\rightarrow\infty)$.

Choose bounded neighborhoods $U_j\downarrow E$ on which $f_j\in A^2(U_j)$. Since $g=0$ outside $E$,
$$
\int_E g\cdot \bar{f} = \lim_{j\rightarrow \infty} \int_{U_j} g\cdot \overline{f_j} =0,
$$
that is, $g\in R^2(E)^\bot$.

For $g\in R^2(E)^\bot$, Lemma \ref{ch3:lm:Havin1} therefore gives $g=\partial\phi_U/\partial z$ with $\phi_U\in W^{1,2}_0(U)$ for each bounded open neighborhood $U$ of $E$. Any two such zero extensions agree: apply uniqueness in Lemma \ref{ch3:lm:Havin1} on a ball containing both neighborhoods. Denote this common function by $\phi$. Shrinking $U$ to $E$ gives $\mathrm{supp}\,\phi\subset E$, that is, $\phi\in W^{1,2}_0(E)$.

On the other hand, if $g=\partial \phi/\partial z$ with $\phi\in W^{1,2}_0(E)$, then $\phi\in W^{1,2}_0(U)$ for every open neighborhood $U$ of $E$, by mollification. Lemma \ref{ch3:lm:Havin1} gives the preceding orthogonality relations for every bounded such $U$, and hence $g\in R^2(E)^\bot$.
\end{proof}

\begin{proof}[Proof of Theorem \ref{ch3:th:Havin}]
Since $R^2(E),\,L^2_a(E)$ are both closed subspaces of $L^2(E)$,
$$
R^2(E) = L^2_a(E) \iff R^2(E)^\bot = L^2_a(E)^\bot.
$$
Every element of $L^2_a(E)^\bot$ vanishes a.e. on $E\backslash E^\circ$, since $L^2_a(E)$ contains all $L^2(E)$ functions supported there. Lemma \ref{ch3:lm:Havin1}, applied on $E^\circ$, therefore gives
$$
L^2_a(E)^\bot
=\{\partial\phi/\partial z:\phi\in W^{1,2}_0(E^\circ)\},
$$
where the derivatives are extended by zero. Lemma \ref{ch3:lm:Havin2} gives the same description of $R^2(E)^\bot$ with $W^{1,2}_0(E)$ in place of $W^{1,2}_0(E^\circ)$. The map $\phi\mapsto\partial\phi/\partial z$ is injective on these spaces by uniqueness in Lemma \ref{ch3:lm:Havin1} on a ball containing $E$. This proves the claimed equivalence.
\end{proof}

\begin{definition}
For a compact set $K\subset \mathbb{R}^N$, define its Sobolev $2$-capacity by
$$
C_2(K) = \inf\left\{ \int_{\mathbb{R}^N} (|u|^2+|\nabla u|^2):
u\in C^\infty_0(\mathbb{R}^N;\mathbb{R}),\ u\ge1\ \text{near }K\right\}.
$$
If $\Omega\subset \mathbb{R}^N$ is open, define
$$
C_2(\Omega)=\sup\left\{C_2(K): K\subset \Omega\ \text{compact}\right\}.
$$
For an arbitrary set $E\subset \mathbb{R}^N$, define
$$
C_2(E)=\inf\left\{C_2(\Omega):  \Omega\supset E\ \text{open}\right\}.
$$
\end{definition}

The $|u|^2$ term is needed here: in dimension two, the infimum of the Dirichlet energy alone over the whole plane is zero for every compact set.

We list some basic properties of capacity:

(1) $E_1\subset E_2$ $\Rightarrow$ $C_2(E_1)\le C_2(E_2)$.

(2) If $K_j,\,K$ are compact sets with $K_j\downarrow K$, then $C_2(K)=\lim_{j\rightarrow \infty} C_2(K_j)$.

(3) If $E_j\uparrow E$, then $C_2(E)=\lim_{j\rightarrow \infty} C_2(E_j)$.

(4) Subadditivity: $C_2(\bigcup E_j)\le \sum C_2(E_j)$.

In pointwise statements about Sobolev functions, we use their quasi-continuous representatives; q.e. means outside a set of zero capacity. We use the standard characterization that $u\in W^{1,2}(\mathbb{R}^N)$ belongs to $W^{1,2}_0(G)$, for an open set $G$, if and only if $u=0$ q.e. on $\mathbb{R}^N\backslash G$. We also use inner regularity of capacity on Borel sets.

\begin{proposition}\label{ch3:prop:capacity}
Let $E\subset \mathbb{R}^N$ be compact and satisfy $W^{1,2}_0(E^\circ) = W^{1,2}_0(E)$. Then, for every open set $U\subset \mathbb{R}^N$,
$$
C_2(U\backslash E) = C_2(U\backslash E^\circ).
$$ 
\end{proposition}

\begin{proof}
Suppose that there exists an open set $U$ such that
\begin{equation}\label{ch3:eq:capacity}
C_2(U\backslash E)<C_2(U\backslash E^\circ).
\end{equation}
The strict inequality implies $C_2(U\backslash E)<\infty$. Choose its capacitary minimizer $u\in W^{1,2}(\mathbb{R}^N)$, so that $u\ge1$ q.e. on $U\backslash E$ and
$$
\int_{\mathbb{R}^N} (|u|^2+|\nabla u|^2)=C_2(U\backslash E).
$$
There must then exist a compact set $K\subset U\cap\partial E$ such that $C_2(K)>0$ and $u<1$ on $K$. Indeed, if $u\ge1$ q.e. on $U\cap\partial E$, the identity $U\backslash E^\circ=(U\backslash E)\cup(U\cap\partial E)$ would imply
$$
C_2(U\backslash E)=\int_{\mathbb{R}^N} (|u|^2+|\nabla u|^2)
\ge C_2(U\backslash E^\circ),
$$
contradicting (\ref{ch3:eq:capacity}). Thus the Borel set $\{u<1\}\cap U\cap\partial E$ has positive capacity, and inner regularity gives such a $K$.

Set $v=\min\{u,1\}\in W^{1,2}(\mathbb{R}^N)$. Then $v=1$ q.e. on $U\backslash E$. Choose $\chi\in C^\infty_0(U)$ such that $\chi=1$ near $K$. The function $\chi(1-v)$ vanishes a.e. outside $E$, so it belongs to $W^{1,2}_0(E)$. However, it is positive on $K\subset\partial E$, a set of positive capacity. The preceding characterization therefore gives $\chi(1-v)\notin W^{1,2}_0(E^\circ)=W^{1,2}_0(E)$, a contradiction!
\end{proof}

Next, we give an example in which stability of the Bergman kernel fails.

\begin{example}[Skwarczy{\'n}ski]
Let $\{\overline{\Delta(z_k,r_k)}\}$ be a sequence of pairwise disjoint closed discs in $\mathbb{D}\backslash [0,1]$ such that the set of accumulation points of $\{z_k\}$ is exactly $[0,1]$. Choose the radii so small that
$$
r_k<\tfrac14 d(z_k,[0,1]),\qquad
\sum_{k=1}^{\infty}C_2(\overline{\Delta(z_k,r_k)})<C_2([0,1]).
$$
This is possible because points have zero Sobolev $2$-capacity in the plane, whereas the nondegenerate segment $[0,1]$ has positive capacity. Set
$$
\Omega=(\mathbb{D}\backslash [0,1])\left\backslash \bigcup \overline{\Delta(z_k,r_k)}\right..
$$
The discs are locally finite off $[0,1]$. Paths in $\mathbb{D}\backslash[0,1]$ can therefore be detoured around the finitely many discs they meet, so $\Omega$ is a domain. Every point of $[0,1]$ is approached by disc boundaries, which belong to $\overline{\Omega}$, and also by disc interiors. Hence $\overline{\Omega}^\circ=\Omega$, while $\partial\Omega$ is contained in $[0,1]$, the unit circle, and the countably many disc boundaries; thus $|\partial\Omega|=0$. Since
$$
\mathbb{D}\backslash\overline{\Omega}=\bigcup_k\Delta(z_k,r_k),
\qquad
\mathbb{D}\backslash\Omega=[0,1)\cup\bigcup_k\overline{\Delta(z_k,r_k)},
$$
and $C_2([0,1))=C_2([0,1])$ because $C_2(\{1\})=0$, subadditivity gives
\begin{align*}
C_2(\mathbb{D}\backslash\overline{\Omega})
&\le\sum_k C_2(\overline{\Delta(z_k,r_k)})\\
&<C_2([0,1))\le C_2(\mathbb{D}\backslash\Omega).
\end{align*}
Set $\Omega_j=\{z\in \mathbb{C}: d(z,\overline{\Omega})<1/j\}$. Clearly, $\Omega_j\downarrow \overline{\Omega}$. On the other hand, Theorem \ref{ch3:th:Skw}
and Proposition \ref{ch3:prop:capacity} imply $K_{\Omega_j}\nrightarrow K_\Omega$.
\end{example}

\section{Exhaustiveness}

\begin{definition}
A bounded domain $\Omega\subset\mathbb{C}^n$ is called an $L^2$-domain of holomorphy if there exists $f\in A^2(\Omega)$ such that
\begin{equation}\label{eq:L2-domain-holomorphy}
\limsup_{z\rightarrow a} |f(z)|=\infty,\ \ \ \forall\,a\in \partial \Omega.
\end{equation}
\end{definition}

\begin{proposition}[Pflug-Zwonek {\cite{PflugZwonek2002}}]\label{prop:L^2-domain-holomorphy}
A bounded domain $\Omega$ is an $L^2$-domain of holomorphy if and only if
\begin{equation}\label{eq:L^2-domain-holomorphy2}
\limsup_{z\rightarrow a} K_\Omega(z)=\infty,  \ \ \ \forall\,a\in \partial \Omega.
\end{equation}
\end{proposition}

\begin{proof}
First suppose there is $f\in A^2(\Omega)$ satisfying \eqref{eq:L2-domain-holomorphy}. Since
$$
K_\Omega(z)\ge \frac{|f(z)|^2}{\|f\|^2},
$$
\eqref{eq:L^2-domain-holomorphy2} follows.

Conversely, assume \eqref{eq:L^2-domain-holomorphy2}. Given $a\in\partial\Omega$, choose $\Omega\ni a_j\rightarrow a$ such that
$$
\sup_j K_\Omega(a_j)=\infty.
$$
We first show that there is $f_a\in A^2(\Omega)$ such that
$$
\sup_j  |f_a(a_j)|=\infty.
$$
Otherwise,
$$
\sup_j |f(a_j)|<\infty,  \ \ \ \forall\,f\in A^2(\Omega).  
$$ 
Applying the Banach--Steinhaus theorem to the functionals $A^2(\Omega)\ni f\mapsto f(a_j)$ gives $\sup_jK_\Omega(a_j)<\infty$, a contradiction.

Next, take all open balls $V_\nu$ with rational centers and radii that meet $\partial\Omega$, and put
\[
\mathcal U_{\nu,m}=\{f\in A^2(\Omega):\sup_{\Omega\cap V_\nu}|f|>m\},\qquad m\in\mathbb{N}.
\]
Point evaluation is continuous, so these sets are open. By what we have proved, there is $h_\nu\in A^2(\Omega)$ unbounded on $\Omega\cap V_\nu$.
If $f\notin\mathcal U_{\nu,m}$, then $f$ is bounded there, whereas $f+\varepsilon h_\nu$ is unbounded there for every nonzero $\varepsilon$. Letting $\varepsilon\to0$ shows that $\mathcal U_{\nu,m}$ is dense.
Baire's theorem gives $f\in\bigcap_{\nu,m}\mathcal U_{\nu,m}$. Every neighborhood of every boundary point contains one of the above balls meeting the boundary, so this $f$ satisfies \eqref{eq:L2-domain-holomorphy}.
\end{proof}

\begin{corollary}\label{cor:fat_L2HoloConvex}
If $\Omega\subset\mathbb{C}^n$ is bounded and pseudoconvex with $\overline\Omega^\circ=\Omega$, then it is an $L^2$-domain of holomorphy.
\end{corollary}

\begin{proof}
Given $a\in\partial\Omega$, choose $\{a_j\}\subset\mathbb{C}^n\setminus\overline\Omega$ with $a_j\rightarrow a$. Put $r_j=d(a_j,\overline\Omega)$, and choose $a_j^\ast\in\partial\Omega$ with $|a_j-a_j^\ast|=r_j$. By Example \ref{ex:Annuli}, there exists $b_j\in\Omega\cap\Delta(a_j^\ast,r_j)$ such that
$$
K_{\Omega\cap \mathcal L_j}(b_j) > j,
$$
where $\mathcal L_j$ is the complex line through $a_j,b_j$. The Ohsawa--Takegoshi theorem gives
$$
K_\Omega(b_j)\ge C K_{\Omega\cap \mathcal L_j}(b_j) > C j,
$$
where $C$ depends only on the diameter of $\Omega$. Apply Proposition \ref{prop:L^2-domain-holomorphy}.

\end{proof}

\begin{definition}
For a domain $\Omega\subset\mathbb{C}^n$, let $\mathrm{Aut}(\Omega)$ denote its holomorphic automorphism group. A subgroup $\Gamma\subset\mathrm{Aut}(\Omega)$ is called properly discontinuous if, for every compact set $E\subset\Omega$, only finitely many $T\in\Gamma$ satisfy $E\cap T(E)\neq\emptyset$.
\end{definition}

By a theorem of H.~Cartan, the quotient $\Omega/\Gamma$ is a complex space. If $\Gamma$ acts freely, meaning every $T\in\Gamma\setminus\{\mathrm{id}\}$ has no fixed points, then $\Omega/\Gamma$ is a complex manifold.

Let $\Gamma_a:=\{T(a):T\in\Gamma\}$ be the orbit of $a\in\Omega$. Let $L(\Gamma_a)$ denote its set of limit points, and define
$$
L(\Gamma):=\bigcup_{a\in \Omega} L(\Gamma_a).
$$
Proper discontinuity implies $L(\Gamma)\subset\partial\Omega$.

\begin{proposition}\label{prop:proper-discontinuous}
Let $\Omega\subset\mathbb{C}^n$ be bounded and let $\Gamma\subset\mathrm{Aut}(\Omega)$ be properly discontinuous. If $L(\Gamma)$ is dense in $\partial\Omega$, then $\Omega$ is an $L^2$-domain of holomorphy.
\end{proposition}

We first prove the following classical result.

\begin{lemma}[Poincar\'e-Siegel]
$$
\sum_{T\in \Gamma} |J_T(a)|^2<\infty,\ \ \ \forall\,a\in \Omega.
$$
Here $J_T$ is the complex Jacobian determinant of $T$.
\end{lemma}

\begin{proof}
Given $a\in\Omega$, choose $r_a<\delta_\Omega(a)$ and let $\overline B_a$ be the closed ball centered at $a$ with radius $r_a$. Proper discontinuity gives $N\in\mathbb{Z}^+$ such that each $T(\overline B_a)$ intersects at most $N$ other sets $T'(\overline B_a)$. Thus
$$
\infty>(N+1)|\Omega| \ge \sum_{T\in \Gamma} |T(\overline{B}_a)|=\sum_{T\in \Gamma} \int_{\overline{B}_a} |J_T|^2\ge |B_a| \sum_{T\in \Gamma} |J_T(a)|^2,
$$
where the last inequality is the submean-value inequality.
\end{proof}

\begin{proof}[Proof of Proposition \ref{prop:proper-discontinuous}]
Given $b\in L(\Gamma)$, choose $a\in\Omega$ and $\{T_j\}\subset\Gamma$ with $T_j(a)\rightarrow b$. Since
$$
K_\Omega(a)=K_\Omega(T_j(a)) |J_{T_j}(a)|^2,
$$
we have
$$
\sum_{j=1}^\infty \frac{K_\Omega(a)}{K_\Omega(T_j(a))} =\sum_{j=1}^\infty |J_{T_j}(a)|^2<\infty,
$$
so that
$$
\lim_{j\rightarrow \infty} K_\Omega(T_j(a))=\infty.
$$

Since $L(\Gamma)$ is dense in $\partial\Omega$, Proposition \ref{prop:L^2-domain-holomorphy} proves the claim.

\end{proof}

\begin{corollary}[Siegel]
Let $\Omega\subset\mathbb{C}^n$ be bounded and let $\Gamma\subset\mathrm{Aut}(\Omega)$ be properly discontinuous. If $\Omega/\Gamma$ is compact, then $\Omega$ is an $L^2$-domain of holomorphy.
\end{corollary}

\begin{definition}
A bounded domain $\Omega\subset\mathbb{C}^n$ is called Bergman exhaustive if $K_\Omega(z)$ is an exhaustion function on $\Omega$, that is,
$\lim_{z\rightarrow \partial \Omega} K_\Omega(z)=+\infty$.
\end{definition}

\begin{definition}
A domain $\Omega\subset\mathbb{C}^n$ is called homogeneous if, for every pair of points $a,b\in\Omega$, there exists $F\in\mathrm{Aut}(\Omega)$ such that
$F(a)=b$.
\end{definition}

\begin{theorem}
Every bounded homogeneous domain is Bergman exhaustive; in particular, it is pseudoconvex.
\end{theorem}

\begin{proof}
Let $\Omega$ be bounded and homogeneous. Suppose there is a sequence $\{a_j\}\subset\Omega$ such that $a_j\rightarrow b\in\partial\Omega$ and
$$
\sup_j K_\Omega(a_j)\le M<+\infty.
$$
Fix $a\in\Omega$ and choose $F_j\in\mathrm{Aut}(\Omega)$ with $F_j(a)=a_j$. Since
$$
K_\Omega(a)=K_\Omega(a_j) |J_{F_j}(a)|^2,
$$
we have $|J_{F_j}(a)|\ge c>0$ for all $j$. By Montel's theorem, a subsequence, still denoted by $\{F_j\}$,
converges uniformly on compact subsets to some $F\in\mathcal O(\Omega,\mathbb{C}^n)$. Weierstrass's theorem gives
$$
|J_F(a)|=\lim_{j\rightarrow \infty} |J_{F_j}(a)|\ge c>0.
$$ 
Thus there is a neighborhood $U\ni a$ such that $F(U)$ is a neighborhood of $F(a)$. On the other hand, $F(\Omega)\subset\overline\Omega$ and
$$
F(a)=\lim_{j\rightarrow \infty} F_j(a)=\lim_{j\rightarrow \infty} a_j = b\in \partial  \Omega,
$$
Take a sufficiently small neighborhood $U\Subset\Omega$ of $a$ on which $F$ is univalent. Since $J_F(a)\ne0$ and $F_j\to F$ locally uniformly, the local inverse mapping theorem (or the multidimensional Rouch\'e theorem) gives $F(a)\in F_j(U)$ for sufficiently large $j$. Hence $b\in\Omega$, contradicting $b\in\partial\Omega$.
\end{proof}

\begin{proposition}\label{prop:cone}
Let $\Omega\subset\mathbb{C}^n$ be bounded and pseudoconvex. Suppose that, for every $a\in\partial\Omega$, there is a sequence $\{a_j\}\subset\mathbb{C}^n\backslash\overline\Omega$ with $a_j\rightarrow a$ and
\begin{equation}\label{eq:cone}
|a_j-a|^2 \log d(a_j,\overline{\Omega})\rightarrow 0\ \ \ (j\rightarrow +\infty).
\end{equation}
Then $\Omega$ is Bergman exhaustive.
\end{proposition}

\begin{corollary}[Pflug \cite{Pflug1975}]
Let $\Omega\subset\mathbb{C}^n$ be bounded and pseudoconvex. Suppose that, for every $a\in\partial\Omega$, there exist $r,\alpha>0$ and a sequence $\{a_j\}\subset\mathbb{C}^n\backslash\overline\Omega$
such that $a_j\rightarrow a$ and
$$
d(a_j,\overline{\Omega}) \ge r |a_j-a|^\alpha.
$$
Then $\Omega$ is Bergman exhaustive.
\end{corollary}

\begin{proof}[Proof of Proposition \ref{prop:cone}]
Fix $a\in\partial\Omega$. For each $j$ and $b\in\Omega\cap B(a,|a_j-a|)$, let $\mathcal L_{a_j,b}$ be the complex line through $a_j,b$. Put $r_j=d(a_j,\overline\Omega)$. After an affine change of coordinates, we may assume
$$
 \mathcal L_{a_j,b} =\{z\in \mathbb{C}^n: z_1=\cdots = z_{n-1}=0\}.
$$
Let $a_j=(0,\cdots,0,t_j)$ and put $f(z_n)=(z_n-t_j)^{-1}$. Then $f\in\mathcal O(\Omega\cap\mathcal L_{a_j,b})$ and
$$
\int_{\Omega\cap \mathcal L_{a_j,b}} |f|^2 \le \int_{r_j<|z_n-t_j|<2R}  |z_n-t_j|^{-2}=2\pi \log \frac{2R}{r_j},
$$
where $R$ is the diameter of $\Omega$. Thus
$$
K_{\Omega\cap \mathcal L_{a_j,b}} (b) \ge \frac{|b-a_j|^{-2}}{2\pi \log \frac{2R}{r_j}}.
$$
The Ohsawa--Takegoshi extension theorem gives
$$
K_\Omega(b)\ge CK_{\Omega\cap \mathcal L_{a_j,b}} (b) \ge \frac{C'}{|a_j-a|^2 \log (2R/r_j)},
$$
so $K_\Omega(z)$ is an exhaustion function.

\end{proof}

The following example shows that Proposition \ref{prop:cone} is nearly optimal.

\begin{definition}
A domain $D\subset\mathbb{D}$ is said to be of Zalcman type $(L)$ if
$$
D=\mathbb{D}\backslash\left(\{0\}\cup\bigcup_{k=1}^\infty\overline{\Delta(x_k,r_k)}\right),
$$ 
where $\{x_k\}\subset(0,1)$ is strictly decreasing to $0$ and the closed discs $\overline{\Delta(x_k,r_k)}\subset\mathbb{D}$ are pairwise disjoint.
\end{definition}

\begin{proposition}\label{prop:Zalcman}
Let $D$ be a domain of type $(L)$ satisfying
\begin{enumerate}
\item  $x_{k-1}-2\sqrt{r_{k-1}}> x_k+2\sqrt{r_k}$,  $\forall\,k\ge 2$,
\item $\sum_k (x_k^2 \log 1/r_k)^{-1}<+\infty$.
\end{enumerate}
Then
$
\sup_{x\in (-1/2,0)} K_D(x)<+\infty.
$
\end{proposition}

\begin{remark}
By Corollary \ref{cor:fat_L2HoloConvex}, every domain of type $(L)$ is an $L^2$-domain of holomorphy.
\end{remark}

\begin{proof}
The idea is essentially due to Zwonek \cite{Zwonek}.  
For $n\in\mathbb{Z}^+$, put
$$
D_n=\mathbb{D}\setminus \bigcup_{k=1}^n \overline{\Delta(x_k,r_k)}.
$$ 
Choose a cutoff $\chi:\mathbb{R}\rightarrow[0,1]$ with $\chi|_{(-\infty,-\log2]}=1$ and $\chi|_{[0,+\infty)}=0$. Put
$$
\eta_n=\chi_0\cdot \prod_{k=1}^n \chi(\log(-\log|z-x_k|)-\log (-\log 2r_k)),
$$
where $\chi_0\in C^\infty_0(\mathbb{D})$ satisfies $\chi_0|_{\mathbb{D}_{3/4}}=1$. Each factor in the product is defined to be $1$ when $|z-x_k|\ge\sqrt{2r_k}$ and $0$ when $|z-x_k|\le2r_k$. Then $\eta_n\in C^\infty_0(D_n)$ and $\eta_n(x)=1$ for all $x\in(-1/2,0)$. By localization of the metric and kernel (Lemma \ref{lm:LocalizationPrinciple}), we may localize near $0$, rescale, and discard finitely many discs away from $0$ to assume $x_1+r_1<1/2$. Values on intervals away from $0$ are controlled by compact-set estimates.
Let $f\in A^2(D_n)$ with $\|f\|_{L^2(D_n)}=1$. The Cauchy--Pompeiu formula gives, for every $x\in(-1/2,0)$,
$$
f(x)=\eta_n(x) f(x) = \frac1{2\pi i} \int_{\mathbb{C}}  \frac{\partial(\eta_n f)/\partial \bar{z}}{z-x} dz\wedge d\bar{z} =  \frac1{2\pi i} \int_{\mathbb{C}}  \frac{f\partial \eta_n /\partial \bar{z}}{z-x} dz\wedge d\bar{z},
$$
so that
\begin{eqnarray*}
|f(x)|  & \le & \frac1{\pi} \int_{\chi_0\neq 1} \frac{|f|\, |\partial \chi_0/\partial \bar{z}|}{|z-x|}\\
 && +\frac{\sup|\chi'|}{\pi} \sum_{k=1}^n \int_{2r_k<|z-x_k|<\sqrt{2r_k}}\frac{|f|}{|z-x||z-x_k||\log |z-x_k||}.
 \end{eqnarray*}
The Cauchy--Schwarz inequality gives
\begin{eqnarray*}
&& \sum_{k=1}^n \int_{2r_k<|z-x_k|<\sqrt{2r_k}}\frac{|f|}{|z-x||z-x_k||\log |z-x_k||}\\
& \le &  \frac12 \sum_{k=1}^n \int_{2r_k<|z-x_k|<\sqrt{2r_k}} |f|^2\\
&& +\frac12 \sum_{k=1}^n \int_{2r_k<|z-x_k|<\sqrt{2r_k}} 
\frac1{|z-x|^2|z-x_k|^2|\log |z-x_k||^2} 
\\
& \lesssim & 1+  \sum_{k=1}^n \frac1{(x_k-\sqrt{2r_k})^2 \log 1/r_k}\\
& \lesssim &   1+ \sum_{k=1}^n \frac1{x_k^2 \log 1/r_k}.
\end{eqnarray*}
Thus
$$
K_{D_n}(x)\le C_1 \left[ 1+\sum_{k=1}^\infty \frac1{x_k^2 \log 1/r_k}\right]^2=:C_2,\ \ \ \forall\,x\in (-1/2,0).
$$
Since $K_{D_n}(x)\uparrow K_D(x)$ (exercise), we obtain $K_D(x)\le C_2$.
\end{proof}

The next example shows that condition (\ref{eq:cone}) cannot be weakened to
$$
|a_k-a|^2 |\log |a_k-a||^{-\varepsilon} \log d(a_k,\overline{\Omega}) \rightarrow 0\ \ \ (k\rightarrow +\infty).
$$

\begin{example}
In Proposition \ref{prop:Zalcman}, take $x_k=e^{-e^k}$ and $r_k=e^{-k^2/x_k^2}$. Then, for every $\varepsilon>0$,
$$
x_k^2 |\log x_k|^{-\varepsilon
} |\log r_k|=k^2 e^{-\varepsilon k}\rightarrow 0\ \ \ (k\rightarrow +\infty),
$$
but
$$
\sum_k \frac1{x_k^2 \log 1/r_k}=\sum_k \frac1{k^2}<\infty,
$$
so $D$ is not Bergman exhaustive.
\end{example}

B\l ocki \cite{Blocki2014}
posed a question on whether Bergman exhaustiveness is conformally invariant for bounded planar domains.  The following result gives an affirmative answer to this question.  

\begin{theorem}\label{th:exhaustiveness-conformal}
Let $\Omega_1,\Omega_2\subset\mathbb{C}$ be bounded domains, and let
$F:\Omega_1\to\Omega_2$ be biholomorphic. Then
\[
\Omega_1\text{ is Bergman exhaustive}
\quad\Longleftrightarrow\quad
\Omega_2\text{ is Bergman exhaustive}.
\]
\end{theorem}

We first prove the following density property.

\begin{lemma}[cf. \cite{Chen2000}]\label{lm:point_dense}
For every bounded plane domain $D$ and every $a\in\partial D$,
the space
\[
\mathcal H_a(D):=\{g\in A^2(D):
g\text{ is bounded on }D\cap U\text{ for some neighborhood }U\text{ of }a\}
\]
is dense in $A^2(D)$.
\end{lemma}

\begin{proof}
Let $f\in A^2(D)$. For simplicity, assume $\sup_{z\in D}|z-a|<1$. Put
$$
\psi(z)=-\log (-\log |z-a|).
$$
Clearly, it satisfies the Donnelly--Fefferman condition
$$
i\partial\bar{\partial}\psi \ge i\partial \psi\wedge \bar{\partial} \psi.
$$
Choose a cutoff $\chi:\mathbb{R}\rightarrow[0,1]$ with $\chi|_{(-\infty,-\log2]}=1$ and $\chi|_{[0,\infty)}=0$. For $0<\varepsilon<1$, put $D_\varepsilon=D\cup\Delta(a,\varepsilon)$.
Then $\chi(-\psi-\log(-\log\varepsilon))f\in C^\infty(D_\varepsilon)$. By the Donnelly--Fefferman estimate \cite{Chen2022}, the equation
$$
\bar{\partial} u = \bar{\partial}(\chi(-\psi-\log(-\log \varepsilon))f)=:v_\varepsilon
$$ 
has a solution $u_\varepsilon$ on $D_\varepsilon$ such that
$$
\int_{D_\varepsilon} |u_\varepsilon|^2 \le C_0 \int_{D_\varepsilon} |v_\varepsilon|^2_{i\partial\bar{\partial}\psi}\le C_0\sup|\chi'|^2 \int_{D\cap \Delta(a,\sqrt{\varepsilon})}|f|^2,
$$
Here $C_0$ is a universal constant. Thus
$$
f_\varepsilon:=\chi(-\psi-\log(-\log \varepsilon))f-u_\varepsilon\in \mathcal O(D_\varepsilon),\qquad f_\varepsilon|_D\in\mathcal H_a(D)
$$
and satisfies
\begin{eqnarray*}
\int_D |f_\varepsilon -f|^2 & \le & 2\int_{D\cap \Delta(a,\sqrt{\varepsilon})}|f|^2 + 2 \int_D |u_\varepsilon|^2\\
& \le  & 2(1+C_0\sup|\chi'|^2) \int_{D\cap \Delta(a,\sqrt{\varepsilon})}|f|^2\\
& \rightarrow & 0 \ \ \ (\varepsilon\rightarrow 0).
\end{eqnarray*}
\end{proof}

\begin{proof}[Proof of Theorem \ref{th:exhaustiveness-conformal}]
Let $D$ be a bounded Bergman exhaustive planar domain. We first prove
\begin{equation}\label{eq:exhaustiveness-normalized-evaluation}
\lim_{z\to\partial D}\frac{|f(z)|^2}{K_D(z)}=0,
\qquad f\in A^2(D).
\end{equation}
Fix $a\in\partial D$ and $\varepsilon>0$. By the preceding density result, choose
$g\in\mathcal H_a(D)$ with $\|f-g\|_{L^2(D)}<\varepsilon$.
The extremal property of the Bergman kernel gives
\[
\frac{|f(z)|}{\sqrt{K_D(z)}}
\le \frac{|g(z)|}{\sqrt{K_D(z)}}
   +\frac{|f(z)-g(z)|}{\sqrt{K_D(z)}}
\le \frac{|g(z)|}{\sqrt{K_D(z)}}+\varepsilon.
\]
As $D\ni z\to a$, the function $g$ remains bounded while $K_D(z)\to\infty$. Hence
\[
\limsup_{D\ni z\to a}\frac{|f(z)|}{\sqrt{K_D(z)}}
\le\varepsilon.
\]
Letting $\varepsilon\downarrow0$ gives a zero limit at every boundary point.
Since $D$ is bounded, every sequence approaching its boundary has a subsequence converging to a boundary point;
a contradiction argument therefore proves \eqref{eq:exhaustiveness-normalized-evaluation}.

Suppose that $\Omega_1$ is Bergman exhaustive. Since $F$ is univalent and $\Omega_2$ is bounded,
the area change-of-variables formula gives
\[
\int_{\Omega_1}|F'(z)|^2\,dA(z)=|\Omega_2|<\infty.
\]
Thus $F'\in A^2(\Omega_1)$. The kernel transformation formula gives
\[
K_{\Omega_1}(z)=K_{\Omega_2}(F(z))\,|F'(z)|^2,
\]
and \eqref{eq:exhaustiveness-normalized-evaluation} yields
\begin{equation}\label{eq:exhaustiveness-conformal-limit}
\frac{1}{K_{\Omega_2}(F(z))}
=\frac{|F'(z)|^2}{K_{\Omega_1}(z)}
\longrightarrow0\qquad(z\to\partial\Omega_1).
\end{equation}

Finally, consider boundary-approaching sequences. If $w_j\in\Omega_2$ approaches $\partial\Omega_2$,
then $z_j:=F^{-1}(w_j)$ leaves every compact subset of $\Omega_1$. Otherwise a subsequence would satisfy
$z_{j_k}\to z_0\in\Omega_1$, giving
$w_{j_k}\to F(z_0)\in\Omega_2$, a contradiction.
By \eqref{eq:exhaustiveness-conformal-limit},
$K_{\Omega_2}(w_j)\to\infty$, so $\Omega_2$ is Bergman exhaustive.
Apply the same argument to $F^{-1}$ for the converse.
\end{proof}

\begin{remark}
Theorem \ref{th:exhaustiveness-conformal} does not extend to higher dimensions: the Hartogs triangle and $\mathbb{D}\times\mathbb{D}^\ast$
give a simple counterexample.
\end{remark}

\begin{definition}
A domain $\Omega\subset\mathbb{C}^n$ is called hyperconvex if it admits a negative continuous plurisubharmonic function $\rho$ such that
$$
\{\rho< -c\}\subset\subset \Omega,\ \ \ \forall\,c>0.
$$ 
\end{definition}

\begin{theorem}[Ohsawa {\cite{Ohsawa1993}}]\label{th:Ohsawa_2}
Every bounded hyperconvex domain $\Omega\subset\mathbb{C}^n$ is Bergman exhaustive.
\end{theorem}

\begin{proof}
For simplicity, assume the diameter of $\Omega$ is less than $1/2$. Given $\varepsilon\ll1$, there exists $C_\varepsilon\gg1$ such that, whenever $\rho(w)\le-2\varepsilon$,
$$
u_w(z):= 
\left\{
\begin{array}{ll}
 \log |z-w|, & \rho(z)\le -\varepsilon\\
 \max\{C_\varepsilon \rho(z),\log |z-w| \}, & \rho(z)>-\varepsilon
 \end{array}
 \right.
$$
defines a negative plurisubharmonic function on $\Omega$ with a logarithmic singularity at $w$. Fix $w'\in\Omega$ and let $\mathcal L$ be any complex line through $w'$. For $w\in\mathcal L$, the restriction of $u_w$ to $\Omega\cap\mathcal L$, still denoted by $u_w$, is negative and subharmonic with a logarithmic singularity at $w$. Since $g_{\Omega\cap\mathcal L}(z,w)\ge u_w(z)$, there exists $\delta=\delta(\varepsilon)>0$ such that
$$
g_{ \Omega\cap \mathcal L }(z,w) >-1, \ \ \ z,w\in \Omega\cap \mathcal L,\, \rho(z)>-\delta,\ \rho(w)\le -2\varepsilon.
$$
Symmetry of $g_{\Omega\cap\mathcal L}$ gives
$$
g_{\Omega\cap \mathcal L} (z,w) >-1, \ \ \ z,w\in \Omega\cap \mathcal L,\, \rho(w)>-\delta,\ \rho(z)\le -2\varepsilon.
$$
Thus, whenever $\rho(w')>-\delta$,
$$
A_{\Omega\cap \mathcal L}(w',-1):=\{z\in \Omega\cap \mathcal L:  g_{\Omega\cap \mathcal L} (z,w') \le -1\} \subset \{\rho>-2\varepsilon\}\cap \mathcal L,
$$
so that\footnote{Put $E_\varepsilon=\{\rho>-2\varepsilon\}$, and let $d\mu$ be normalized Fubini--Study volume on the space of complex lines through $w'$. Polar integration gives
$\int |E_\varepsilon\cap\mathcal L|\,d\mu(\mathcal L)=c_n\int_{E_\varepsilon}|z-w'|^{2-2n}\,dV(z)\le C_n|E_\varepsilon|^{1/n}$.
The last estimate follows by splitting the integral into $|z-w'|<R$ and its complement and choosing $R=|E_\varepsilon|^{1/(2n)}$ (for $n=1$ it is immediate). Since $|E_\varepsilon|\to0$, the infimum of the section areas also tends to zero.}
$$
\inf_{\mathcal L} |A_{\Omega\cap \mathcal L}(w',-1)| \rightarrow 0\ \ \ (w'\rightarrow \partial \Omega),
$$
Here the infimum is over all complex lines $\mathcal L$ through $w'$.

On the other hand, the proof of Theorem \ref{th:hyperbolicBergman} gives an absolute constant $C_0>0$ such that
$$
K_{\Omega\cap \mathcal L}(w')\ge \frac{C_0}{|A_{\Omega\cap \mathcal L}(w',-1)|},
$$
The Ohsawa--Takegoshi extension theorem gives
a constant $C$, depending only on the diameter of $\Omega$, such that
$$
K_\Omega(w')\ge C \sup_{\mathcal L} K_{\Omega\cap \mathcal L}(w')\ge \frac{C_0C}{\inf_\mathcal L |A_{\Omega\cap \mathcal L}(w',-1)|}\rightarrow \infty\ \ \ (w'\rightarrow \partial \Omega).
$$
\end{proof}

A longstanding open problem is the following.

\begin{problem}
Is every bounded pseudoconvex domain with $C^0$ boundary hyperconvex?
\end{problem}

\begin{remark}
This problem has  been solved only very recently (see \cite{WangYang}).
\end{remark}

Independently, Bergman exhaustiveness can be established for bounded pseudoconvex domains with $C^0$ boundary. More precisely, we have the following sharp lower bound for the Bergman kernel function:

\begin{theorem}[cf.  {\cite{ChenCapacities}}]\label{th:ChenC^0}\footnote{The case of $C^2$ boundary is known (see \cite{Fu}).}
If $\Omega\subset\mathbb{C}^n$ is bounded and pseudoconvex with $C^0$ boundary, then there exists $C>0$ such that
\begin{equation}\label{eq:delta^-2}
K_\Omega(z)\ge C\delta_\Omega(z)^{-2}.
\end{equation}
\end{theorem}

We first introduce a large class of domains for which \eqref{eq:delta^-2} holds. Note that a domain $\Omega\subset\mathbb{R}^n$ is convex if and only if its intersection with every real line is empty or an interval, while a domain $\Omega\subset\mathbb{C}^n$ is called $\mathbb{C}$-convex if its intersection with every complex line is empty or a simply connected domain. On the other hand, \textbf{every} domain
$\Omega\subset\mathbb{R}^n$ has the elementary property that each component of its nonempty intersection with a real line is an interval. This suggests the following definition.
    
    \begin{definition}
A domain $\Omega\subset\mathbb{C}^n$ is said to be \textbf{Runge along complex lines} if, for every complex line $\mathcal L$, the intersection $\Omega\cap\mathcal L$ is either empty or a Runge open set, that is, every connected component is simply connected.
\end{definition}

\begin{proposition}\label{prop:Runge}
If $\Omega\subset\mathbb{C}^n$ is pseudoconvex and Runge, meaning that every holomorphic function on $\Omega$ can be approximated uniformly on compact subsets by polynomials, then $\Omega$ is Runge along complex lines.
\end{proposition}

\begin{proof}
Let $\mathcal L$ be a complex line with $\Omega\cap\mathcal L\neq\emptyset$, and let $D$ be a component of $\Omega\cap\mathcal L$. For $f\in\mathcal O(D)$, define $\tilde f$ holomorphically on $\Omega\cap\mathcal L$ by $\tilde f=f$ on $D$ and $\tilde f=0$ on the other components. Since $\Omega$ is pseudoconvex, $\tilde f$ extends to $F\in\mathcal O(\Omega)$ with $F|_{\Omega\cap\mathcal L}=\tilde f$.

Since $\Omega$ is Runge, there are polynomials $P_j$ converging to $F$ uniformly on compact subsets of $\Omega$. In particular,
$
P_j|_{\mathcal L} 
$
converges to $f$ uniformly on compact subsets of $D$. Runge's theorem therefore implies that $D$ is simply connected.

\end{proof}

\begin{remark}
Conversely,  Bremermann asks whether a pseudoconvex domain is Runge if it is Runge along complex lines.  Joi\c{t}a \cite{Joita} gave a negative answer: there exists a bounded strongly pseudoconvex domain in $\mathbb{C}^2$ with real-analytic boundary,  which is Runge along complex lines but  is not Runge.
\end{remark}

Every $\mathbb{C}$-convex domain is pseudoconvex, whereas there are nonpseudoconvex domains which are Runge along complex lines.

\begin{example}\label{L:nonpsc}
Write $\mathbb{D}_r=\{\zeta\in\mathbb{C}:|\zeta|<r\}$. Then
\[
U=(\mathbb{D}_2\times\mathbb{D}_1)\cup(\mathbb{D}_1\times\mathbb{D}_2)
\]
is Runge along complex lines
 but is not pseudoconvex.
\end{example}

\begin{proof}
Each complex-line section is the union of two convex planar open sets. If they are disjoint, each component is convex. If they intersect, their union is star-shaped with respect to any point in the intersection, hence simply connected. Thus $U$ is Runge along complex lines.

We verify nonpseudoconvexity directly using power series. Let $f\in\mathcal O(U)$ and expand at the origin as $\sum_{j,k\ge0}a_{jk}z^jw^k$. Fix $r=3/2$ and $s=3/4$. Cauchy's estimates on the tori with radii $(r,s)$ and $(s,r)$ give $M_f<\infty$ such that
\[
|a_{jk}|\le M_f\min\{r^{-j}s^{-k},s^{-j}r^{-k}\}
\le M_f(rs)^{-(j+k)/2}.
\]
Thus the series converges on $\mathbb{D}_{\sqrt{rs}}^2$. It agrees with $f$ on the intersection, so every $f\in\mathcal O(U)$ extends to the strictly larger connected open set $U\cup\mathbb{D}_{\sqrt{rs}}^2$, where $\sqrt{rs}=\sqrt{9/8}>1$. Hence $U$ is not a domain of holomorphy and, by the solution of the Levi problem, is not pseudoconvex.
\end{proof}

\begin{example}\label{L:many}
For an integer $m\ge2$ and $0<\varepsilon<1$, put
\[
\Omega_{m,\varepsilon}=\{(z,w)\in\mathbb{C}^2:|z|<2,\ |w-z^m|<\varepsilon\}.
\]
This is a bounded pseudoconvex Runge domain and therefore is Runge along complex lines. However, its section at $w=1$ has exactly $m$ components. In particular, $\Omega_{m,\varepsilon}$ is not $\mathbb{C}$-convex.
\end{example}

\begin{proof}
The polynomial automorphism $(z,w)\mapsto(z,w-z^m)$ maps the domain onto $\mathbb{D}_2\times\mathbb{D}_\varepsilon$. Polynomial automorphisms preserve the Runge property, proving the first assertion.

The section at $w=1$ is given by $|z^m-1|<\varepsilon$. The disc $\mathbb{D}(1,\varepsilon)$ avoids zero and has $m$ distinct holomorphic branches of the $m$th root. Their images are pairwise disjoint and lie in $\mathbb{D}_2$. Thus the section has exactly $m$ simply connected components, and $\Omega_{m,\varepsilon}$ is not $\mathbb{C}$-convex.
\end{proof}

\begin{proposition}\label{prop:propertyL}
If  a bounded pseudoconvex domain $\Omega$ is   Runge along complex lines, then \eqref{eq:delta^-2} holds.
\end{proposition}

\begin{proof}
Given $w\in\Omega$, choose $w^\ast\in\partial\Omega$ with $|w-w^\ast|=\delta_\Omega(w)$, and let $\mathcal L$ be the complex line through $w,w^\ast$. Let $D$ be the component of $\Omega\cap\mathcal L$ containing $w$. By the Ohsawa--Takegoshi theorem, for every $f\in A^2(D)$ there exists $F\in A^2(\Omega)$
with $F=f$ on $D$, $F=0$ on the other components of $\Omega\cap\mathcal L$, and
$$
\int_\Omega |F|^2 \le C\int_{D} |f|^2,
$$
where $C$ depends only on the diameter of $\Omega$. Thus
$$
K_\Omega(w)\ge C K_D(w) \ge C\delta_{D}(w)^{-2}=C\delta_\Omega(w)^{-2},
$$
where the second inequality follows from Proposition \ref{prop:simply-connected}.
\end{proof}

\begin{proposition}\label{prop:C^0Runge}
Let $B'\subset\mathbb{C}^{n-1}\times\mathbb{R}$ be a ball, and let $g(z',\mathrm{Re}\,z_n)$ be a positive continuous function on a neighborhood of $\overline{B'}$. If
$$
\Omega=\left\{z\in B'\times \mathbb{R}: 0<\mathrm{Im}\,z_n < g(z',\mathrm{Re}\,z_n) \right\}
$$
is pseudoconvex, then it is Runge.
  \end{proposition}
  
The proof requires the following fundamental result.

   \begin{theorem}[Docquier-Grauert {\cite{DocquierGrauert1960}}]\label{th:DG}  
Let $\{\Omega_t\}_{t\in[0,1]}$ be an increasing family of bounded pseudoconvex domains in $\mathbb{C}^n$, continuous in $t$ in the sense that $d(z,\partial\Omega_t)$ is continuous in $(z,t)\in\mathbb{C}^n\times[0,1]$. Then, for every $0\le t_1<t_2\le1$, the pair $(\Omega_{t_1},\Omega_{t_2})$ is Runge: holomorphic functions on $\Omega_{t_1}$ can be approximated uniformly on compact subsets by holomorphic functions on $\Omega_{t_2}$.
   \end{theorem}  
   
   \begin{proof}
Let $E\subset\Omega_{t_1}$ be compact and put $r_1:=\frac12d(E,\partial\Omega_{t_1})$. Continuity in $t$ gives $0<\varepsilon\ll1$ such that, for $t_1\le s<t\le t_2$ with $t-s<\varepsilon$,
   $$
   E \subset \{z\in \Omega_t: \delta_{\Omega_{t}}(z) > r_1\} =: \widetilde{\Omega}_t\subset\subset \Omega_s.  
   $$
By the classical Oka--Weil approximation theorem, $(\widetilde\Omega_t,\Omega_t)$ is a Runge pair. Hence every $f\in\mathcal O(\Omega_s)$ can be uniformly approximated on $E$ by holomorphic functions on $\Omega_t$. Finally, choose a partition
   $$
t_1=s_1<s_2<\cdots<s_m=t_2,\ \ \text{with }s_{j+1}-s_j<\varepsilon,
   $$
and apply induction.
   
   \end{proof}  
   
\begin{proof}[Proof of Proposition \ref{prop:C^0Runge}]
For $t\ge0$, put
\[
\Omega_t=\{z\in B'\times\mathbb{R}:-t<\operatorname{Im}z_n<g(z',\operatorname{Re}z_n)+t\}.
\]
The upper boundary is a translate of the original upper boundary, the lower boundary is a real hyperplane, and the side boundary comes from the convex cylinder $B'\times\mathbb{R}$. By the local nature of pseudoconvexity, every $\Omega_t$ is pseudoconvex. These bounded domains form an increasing continuous family, with $\Omega_0=\Omega$ and union $B'\times\mathbb{R}$.
By Theorem \ref{th:DG}, $(\Omega,B'\times\mathbb{R})$ is a Runge pair. The convex cylinder itself is Runge, so transitivity of approximation shows that $\Omega$ is Runge.
\end{proof}
        
Propositions \ref{prop:Runge}, \ref{prop:propertyL}, and \ref{prop:C^0Runge}, together with the next lemma, immediately imply Theorem \ref{th:ChenC^0}.
 
 \begin{lemma}[Localization lemma]\label{lm:LocalizationPrinciple}
Let $\Omega\subset\mathbb{C}^n$ be bounded and pseudoconvex, let $a\in\partial\Omega$, and let $V\subset\subset U$ be neighborhoods of $a$. For every $k\in\mathbb{Z}^+$, there exists $C=C(k,U,V,\Omega)$ such that, for every $f\in A^2(\Omega\cap U)$ and $w\in\Omega\cap V$, there is $F\in A^2(\Omega)$ satisfying
$$
F^{(j)}(w)=f^{(j)}(w),\ \forall\,j\le k;\ \ \ \int_\Omega |F|^2 \le C\int_{\Omega\cap U} |f|^2.
$$
In particular,
$$
K_\Omega(w)\ge C' K_{\Omega\cap U}(w),\ \ \ \forall\,w\in \Omega\cap V,
$$
where $C'$ depends only on $U,V,\Omega$.
   \end{lemma}  
   
   \begin{proof}
Choose $\chi\in C_0^\infty(U)$ equal to $1$ on a neighborhood of $\overline V$. By H\"ormander's theorem, the equation
   $$
   \bar{\partial} u = \bar{\partial}(\chi f)
   $$
has a solution satisfying
   \begin{eqnarray*}
   \int_\Omega |u|^2 e^{-2(n+k)\log |z-w|} & \le &  C_1 \int_\Omega |f|^2 |\bar{\partial}\chi|^2 e^{-2(n+k)\log |z-w|}\\
   & \le & C \int_{\Omega\cap U} |f|^2,
   \end{eqnarray*}
Here $C_1$ depends only on the diameter of $\Omega$, while $C$ denotes a \textbf{generic} constant depending on $k,\Omega,U,V$. Thus
   $$
   F:=\chi f-u\in \mathcal O(\Omega)
   $$
satisfies $F^{(j)}(w)=f^{(j)}(w)$ for all $j\le k$, and
   $$
   \int_\Omega |F|^2 \le 2\int_\Omega |\chi f|^2 + 2\int_\Omega |u|^2 \le C \int_{\Omega\cap U} |f|^2
   $$
(note that $2(n+k)\log|z-w|\le C$ on $\Omega$).
\end{proof}   

Another class of domains satisfying \eqref{eq:delta^-2} arises in complex dynamics.
For $F\in\mathrm{Aut}(\mathbb{C}^n)$, denote its $k$th iterate by $F^k$. A fixed point $a$ is attracting if all eigenvalues of $F'(a)$ lie in $\mathbb{D}$. Rosay and Rudin \cite{RosayRudin} proved that if
$n>1$ and $a$ is an attracting fixed point of $F\in\mathrm{Aut}(\mathbb{C}^n)$, then
    $$
    \Omega_{F,a}:=\left\{z\in \mathbb{C}^n: \lim_{k\rightarrow \infty} F^k(z) = a\right\}
    $$
is biholomorphic to $\mathbb{C}^n$ and is Runge. If $\Omega_{F,a}\ne\mathbb{C}^n$, it is  a Fatou--Bieberbach domain; we consider this case below.

Clearly, $K_{\Omega_{F,a}}\equiv0$. Since $\Omega_{F,a}$ is pseudoconvex, choose a smooth plurisubharmonic exhaustion $\rho$ on it. Given $b\in\partial\Omega$ and $r>0$, for every $R>0$,
$$
\left\{B(b,tr)\cap \{\rho<R+t\}: t\ge 1\right\}
$$
is an increasing continuous family of bounded pseudoconvex open sets. Theorem \ref{th:DG} shows that $B(b,r)\cap\{\rho<R+1\}$ and $\Omega_{F,a}$ form a Runge pair. Since $\Omega_{F,a}$ is Runge, so is $\Omega_{F,a}\cap B(b,r)$. Propositions \ref{prop:Runge} and \ref{prop:propertyL} therefore give
$$
K_{ \Omega_{F,a}\cap B(b,r) }(z) \gtrsim \delta_{\Omega_{F,a}}(z)^{-2}.
$$   

\begin{problem}   
Does the preceding estimate also hold for Fatou--Bieberbach domains that are not Runge?
\end{problem}

\begin{remark}\label{rem:WoldLocalBergman}
Wold \cite{Wold} constructed a Fatou--Bieberbach domain $\Omega_W\subset\mathbb{C}^*\times\mathbb{C}$ that is Runge in $\mathbb{C}^*\times\mathbb{C}$ but not in $\mathbb{C}^2$.
Nevertheless, this example is Runge along complex lines. Thus the preceding question has an affirmative answer for Wold's non-Runge example.
\end{remark}

\begin{proof}
First consider its complex-line sections. Let $\mathcal L$ be a complex line, let $\gamma\subset\Omega_W\cap \mathcal L$ be a piecewise smooth Jordan curve, and let $G\subset \mathcal L$ be its bounded interior. Since $\Omega_W$ is biholomorphic to $\mathbb{C}^2$, $\gamma$ is contractible in $\Omega_W$, so $\int_\gamma dz_1/z_1=0$. If the intersection of $\mathcal L$ with $\{z_1=0\}$ lay in $G$, the argument principle would give $\pm2\pi i$, a contradiction. Hence $\overline G\subset\mathbb{C}^*\times\mathbb{C}$ (the claim is immediate if $\mathcal L$ is parallel to $\{z_1=0\}$). By the maximum modulus principle, $\overline G$ lies in the holomorphic hull $\widehat{\gamma}_{\mathcal O(\mathbb C^\ast\times \mathbb C)}$  of $\gamma$ in $\mathbb{C}^*\times\mathbb{C}$. The Runge property and pseudoconvexity of $\Omega_W$ then imply 
$$
\widehat{\gamma}_{\mathcal O(\mathbb C^\ast\times \mathbb C)}=\widehat{\gamma}_{\mathcal O(\Omega_W)}\subset \Omega_W,
$$
so $G\subset\Omega_W$. The Jordan-curve criterion proves simple connectivity of each section component.

Put $D=\Omega_W\cap B(b,r)$. For every Jordan curve in $D\cap \mathcal L$, the preceding argument places its interior in $\Omega_W$. Since $B(b,r)\cap \mathcal L$ is a disc, the interior also lies in the ball. Thus every component of $D$ is Runge along complex lines  and is bounded and pseudoconvex.
\end{proof}

The following proposition gives bounded pseudoconvex domains $\Omega\subset \mathbb C^2$ that are \textbf{not hyperconvex} but satisfy $K_\Omega(z)\ge C\delta_\Omega(z)^{-2}$.\footnote{For bounded planar domains,  $K_\Omega(z)\ge C\delta_\Omega(z)^{-2}$ would imply that $\partial \Omega$ is uniformly perfect; in particular,  $\Omega$ is hyperconvex (cf. Zheng \cite{Zheng}). }

\begin{proposition}[Kerzman-Rosay {\cite{KerzmanRosay1981}}]\label{prop:KR}
Let $\varphi:\mathbb{D}\rightarrow(0,1)$ be subharmonic and discontinuous at $0$. Then the Hartogs domain
$$
\Omega:=\left\{(z_1,z_2)\in \mathbb{D}\times \mathbb{C}: |z_2|<e^{-\varphi(z_1)}\right\}
$$
is bounded, pseudoconvex, and Runge, but not hyperconvex.
\end{proposition}

For example, take
$$
\psi(z_1)=\sum_{j=1}^\infty 2^{-j}\log\frac{|z_1-1/j|}{2},
\qquad \varphi=\frac{1+e^\psi}{2}.
$$
Since $\sum_j2^{-j}\log(2j)<\infty$, we have $\psi(0)>-\infty$, whereas $\psi(1/j)=-\infty$ for $j\ge2$. Thus $\varphi$ is subharmonic with values in $(0,1)$, and $\varphi(1/j)=1/2<\varphi(0)$, so it is discontinuous at $0$.

To prove the proposition, we introduce the following concept.

\begin{definition}
A domain $\Omega\subset\mathbb{C}^n$ is called taut if every sequence $\{f_j\}\subset\mathcal O(\mathbb{D},\Omega)$ either

(a) has a subsequence $\{f_{j_k}\}$ converging in $\mathcal O(\mathbb{D},\Omega)$,
or

(b) for every compact $K\subset\mathbb{D}$ and compact $L\subset\Omega$, there is $j_0\in\mathbb{Z}^+$ such that
$$
f_j(K)\cap L=\emptyset,\ \ \ \forall\,j>j_0.
$$
\end{definition}

\begin{proposition}\label{prop:hyperconvexTaut}
Every bounded hyperconvex domain $\Omega\subset\mathbb{C}^n$ is taut.
\end{proposition}

\begin{proof}
Choose a negative continuous plurisubharmonic exhaustion $\rho$ and extend it by zero to $\partial\Omega$, obtaining a continuous function on $\overline\Omega$. Given $\{f_j\}\subset\mathcal O(\mathbb{D},\Omega)$, if (b) fails, take a subsequence and $\zeta_j\in K\Subset\mathbb{D}$ with $f_j(\zeta_j)\in L\Subset\Omega$. By Montel's theorem, pass to a further subsequence with $f_j\to f$ uniformly on compact subsets and $\zeta_j\to\zeta\in K$. Then $f(\zeta)\in L$.

Since $\rho\circ f_j$ converges uniformly on compact subsets to $\rho\circ f$, the latter is nonpositive and subharmonic, and strictly negative at $\zeta$. The maximum principle gives $\rho\circ f<0$ on $\mathbb{D}$, hence $f(\mathbb{D})\subset\Omega$. Thus (a) holds.
\end{proof}

\begin{proof}[Proof of Proposition \ref{prop:KR}]
Clearly, $\Omega$ is bounded and pseudoconvex. By Proposition \ref{prop:hyperconvexTaut}, it suffices to show that $\Omega$ is not taut. Since $\varphi$ is not lower semicontinuous at $0$, there are $c_j\in\mathbb{D}$ tending to $0$ such that
$$
0< \varphi(c_j)<c<\varphi(0),\ \ \ \forall\,j.
$$
Define
$$
f_j(\zeta):=\left(c_j,e^{-c}\zeta\right)\in \mathcal O(\mathbb{D},\Omega),\ \ \ j\in \mathbb{Z}^+.
$$
Note that $\lim_{j\rightarrow\infty}f_j(0)=(0,0)\in\Omega$, whereas
$$
\lim_{j\rightarrow \infty} f_j\left(\exp(c-\varphi(0))\right)=\left(0,\exp(-\varphi(0))\right)\notin \Omega,
$$
This shows that $\Omega$ is not taut.

Finally, we verify the Runge property. Every $h\in\mathcal O(\Omega)$ has an expansion converging uniformly on compact subsets,
$$
h(z_1,z_2)=\sum_{k=0}^\infty h_k(z_1)z_2^k,\qquad h_k\in\mathcal O(\mathbb{D}).
$$
Truncate the series, then approximate the finitely many coefficients $h_k$ by polynomials. This approximates $h$ uniformly on compact subsets of $\Omega$ by polynomials in $\mathbb{C}^2$.
\end{proof}

Finally, we would like to point out that the Bers embedding $\Omega_T\subset\mathbb{C}^{3g-3}$ of the Teichm\"uller space of compact Riemann surfaces of genus $g\ge2$ is a bounded pseudoconvex Runge domain (see \cite{Krushkal1995} or \cite{Dujardin2017}). Thus $K_{\Omega_T}(z)\gtrsim\delta_{\Omega_T}(z)^{-2}$. Note that the boundary of $\Omega_T$ is extremely complicated.

\chapter{The Bergman Metric on Domains}
\section{Basic Properties}

If $\Omega\subset\mathbb{C}^n$ is bounded, the extremal property gives $K_\Omega(z)\ge1/|\Omega|>0$. We may therefore define the Hermitian form
$$
ds^2_\Omega=\sum_{j,k=1}^n \frac{\partial^2 \log K_\Omega(z)}{\partial z_j\partial \bar{z}_k} dz_j \otimes d\bar{z}_k. 
$$

\begin{proposition}
The form $ds_\Omega^2$ is positive definite. For $z\in\Omega$ and a complex tangent vector $X=\sum_jX_j\partial/\partial z_j$, we have
\begin{eqnarray*}
B_\Omega(z;X) & :=& \sqrt{\sum_{j,k=1}^n \frac{\partial^2 \log K_\Omega(z)}{\partial z_j\partial \bar{z}_k}  X_j\overline{X}_k}\\
& = & \frac1{\sqrt{K_\Omega(z)}} \sup\{|Xf(z)|:f\in A^2(\Omega), f(z)=0, \|f\|=1\}.
\end{eqnarray*}
\end{proposition}

\begin{proof}
Given $z\in\Omega$ and $X\neq0$, put
$$
A_0=\{f\in A^2(\Omega): f(z)=0\},\ \ \ A_1=\{f\in A_0: Xf(z)=0\}.
$$
Clearly, $A_0\subset A^2(\Omega)$ and $A_1\subset A_0$ are closed subspaces, so there are orthogonal decompositions
$$
A^2(\Omega) = A_0 \oplus A_0^\bot,\ \ \ A_0=A_1\oplus A_1^\bot.
$$
Since $1\in A^2(\Omega)\backslash A_0$, there exists $0\neq h_0\in A_0^\bot$ with $\|h_0\|=1$. As $h_0(z)\neq0$, for any $f\in A^2(\Omega)$,
$$
f-\frac{f(z)}{h_0(z)} h_0\in A_0,
$$
that is, $A^2(\Omega)=\mathrm{Span}\{h_0,A_0\}$. Similarly, choose $h_1\in A_1^\bot\subset A_0$ with $\|h_1\|=1$ and $Xh_1(z)\neq0$. We thus obtain a complete orthonormal basis $\{h_0,h_1,h_2,\cdots\}$ of $A^2(\Omega)$ with $h_j\in A_1$ for every $j\ge2$. Since $K_\Omega(\zeta)=\sum_{j=0}^\infty|h_j(\zeta)|^2$, we have
\begin{eqnarray*}
B_\Omega (z;X)^2 & = & \frac{X\overline{X}K_\Omega(z)}{K_\Omega(z)} -\frac{|XK_\Omega(z)|^2}{K_\Omega(z)^2}\\
& = & \frac{ |Xh_0(z)|^2+|Xh_1(z)|^2}{K_\Omega(z)}-\frac{|h_0(z)|^2|X h_0 (z)|^2}{K_\Omega(z)^2}\\
&  = & \frac1{K_\Omega(z)}|Xh_1(z)|^2>0.
 \end{eqnarray*} 
On the other hand, for $f\in A_0$, consider its Fourier series
$$
f=\sum_{j=1}^\infty c_jh_j,\ \text{where }\sum_{j=1}^\infty|c_j|^2=\|f\|^2.
$$
Thus
$$
|Xf(z)|^2 = |c_1|^2 |Xh_1(z)|^2 \le \|f\|^2 |Xh_1(z)|^2.
$$
\end{proof}

\begin{definition}
The form $ds_\Omega^2$ is called the Bergman metric of $\Omega$.
\end{definition}

Let $\omega_\Omega$ denote its K\"ahler form. When no confusion arises, we also call $\omega_\Omega$ the Bergman metric of $\Omega$.

\begin{proposition}\label{prop:metric-trans}
If $F:\Omega_1\rightarrow\Omega_2$ is biholomorphic, then $\omega_{\Omega_1}=F^*\omega_{\Omega_2}$.
\end{proposition}

\begin{proof}
Since $K_{\Omega_1}(z)=K_{\Omega_2}(F(z))|J_F(z)|^2$, we have
$$
\partial\bar{\partial} \log K_{\Omega_1}(z)=\partial\bar{\partial}\log K_{\Omega_2}(F(z)).
$$
Put $w=F(z)$. Direct calculation gives
\begin{eqnarray*}
\omega_{\Omega_1} &  =& \frac{i}2 \sum_{j,k}\sum_{\alpha,\beta} \frac{\partial^2 \log K_{\Omega_2}(w)}{\partial w_\alpha \partial \bar{w}_\beta} \frac{\partial w_\alpha}{\partial z_j}
\overline{\frac{\partial w_\beta}{\partial z_k}} dz_j\wedge d\bar z_k\\
& = & \frac{i}2 \sum_{\alpha,\beta} \frac{\partial^2 \log K_{\Omega_2}(w)}{\partial w_\alpha \partial \bar{w}_\beta} \left(\sum_j  \frac{\partial w_\alpha}{\partial z_j} dz_j\right)
\wedge \overline{\sum_k  \frac{\partial w_\beta}{\partial z_k} dz_k}\\
& = &  \frac{i}2 \sum_{\alpha,\beta} \frac{\partial^2 \log K_{\Omega_2}(w)}{\partial w_\alpha \partial \bar{w}_\beta}  dw_\alpha\wedge d\bar{w}_\beta=F^*\omega_{\Omega_2}.
\end{eqnarray*}
\end{proof}

The product property of the Bergman kernel immediately implies

\begin{proposition}\label{prop:prod}
For bounded domains $\Omega_1\subset\mathbb{C}^m$ and $\Omega_2\subset\mathbb{C}^n$,
$$
\omega_{\Omega_1\times \Omega_2}=\omega_{\Omega_1} \oplus \omega_{\Omega_2}.
$$ 
\end{proposition}

\begin{example}
The Bergman metrics of the unit ball $\mathbb{B}^n$ and the unit polydisc $\mathbb{D}^n$ are, respectively,
$$
\omega_{\mathbb{B}^n} = (n+1)\frac{i}2\partial\bar{\partial}[-\log(1-|z|^2)], \ \ \ \omega_{\mathbb{D}^n}=\sum_{j=1}^n  {i}\partial\bar{\partial}[-\log(1-|z_j|^2)].
$$
\end{example}

Let $\Gamma(a,b)$ be the collection of piecewise $C^1$ curves in $\Omega$ joining $a$ to $b$.

\begin{definition}
The quantity
$$
d_\Omega(a,b):=\inf_{\gamma\in \Gamma(a,b)} \int_0^1 B_{\Omega}({\gamma}(t);{\gamma}'(t))dt 
$$
is called the Bergman distance from $a$ to $b$.
\end{definition}

\begin{corollary}
If $F:\Omega_1\rightarrow\Omega_2$ is biholomorphic, then
 $$
 d_{\Omega_1}(a,b)=d_{\Omega_2}(F(a),F(b)),\ \ \ \forall\,a,b\in \Omega_1,
 $$
that is, biholomorphic maps are isometries for the Bergman distance.
\end{corollary}

\begin{proof}
By Proposition \ref{prop:metric-trans},
\begin{eqnarray*}
d_{\Omega_1}(a,b) &  = & \inf_{\gamma\in \Gamma(a,b)} \int_0^1 B_{\Omega_2}(F\circ \gamma(t);(F\circ \gamma)'(t)) dt\\
& \ge & \inf_{\widetilde{\gamma}\in \Gamma(F(a),F(b))} \int_0^1 B_{\Omega_2}(\widetilde{\gamma}(t);(\widetilde{\gamma}'(t))dt\\
&  \ge & d_{\Omega_2}(F(a),F(b)).
\end{eqnarray*}
Repeat the argument with $F^{-1}$ in place of $F$ to obtain the reverse inequality.
\end{proof}

\begin{definition}
For a bounded domain $\Omega\subset\mathbb{C}^n$, the quantity
$$
C_\Omega(z;X)  =  \sup\left\{|Xf(z)|: f\in\mathcal{O}(\Omega,\mathbb{D}),\ f(z)=0 \right\},
$$
is called the Carath\'eodory metric on $\Omega$.
\end{definition}

\begin{theorem}[Lu Qi-Keng]
$$
B_\Omega(z;X) \ge C_\Omega(z;X)
$$
Moreover, the inequality is strict for every $X\neq0$.
\end{theorem}

\begin{proof}
This simple proof is due to K.~T.~Hahn \cite{Hahn1978}. Put $b_z=K_\Omega(\cdot,z)/\sqrt{K_\Omega(z)}$. Since $\mathcal O(\Omega,\mathbb{D})$ is a normal family, there exists $c_{z,X}\in\mathcal O(\Omega)$ with $c_{z,X}(z)=0$, $|c_{z,X}|<1$, and $Xc_{z,X}(z)=C_\Omega(z;X)$. Using $b_zc_{z,X}$ as a test function in the extremal problem for the Bergman metric gives
	$$
B_\Omega(z;X) \ge \frac{1}{ \sqrt{K_\Omega(z)} } \frac{ |b_z(z)| |Xc_{z,X}(z)|}{\|b_z c_{z,X} \|_{L^2(\Omega)}}=\frac{C_\Omega(z;X)}{\|b_z c_{z,X}\|_{L^2(\Omega)}}.
$$
Note that
 $$
 \|b_z c_{z,X} \|_{L^2(\Omega)}< \|b_z\|_{L^2(\Omega)}=1,
 $$
and the conclusion follows.
\end{proof}

\begin{remark}
Chen, Xiong, and Zhang \cite{CXZ} further proved the following. Write
\[
S_{z,X}:=\{f\in{A^2(\Omega)}:  f(z)=0,\ Xf(z)=0\},
\]
and let $P_{z,X}$ be the orthogonal projection from $A^2(\Omega)$ onto $S_{z,X}^\perp$. Then
$$
B_\Omega (z;X)=\frac{ C_\Omega(z;X)}{\|P_{z,X}(b_z c_{z,X})\|_{L^2(\Omega)}}.
$$
\end{remark}

Lemma \ref{lm:LocalizationPrinciple} immediately gives

\begin{lemma}[Localization principle]\label{lm:local}
Let $\Omega$ be bounded and pseudoconvex, and let $V\subset\subset U$ be neighborhoods of a boundary point $a$. Then there exists $C>0$ such that
$$
C^{-1} B_\Omega(z;X) \le  B_{\Omega\cap U}(z;X) \le  C B_\Omega(z;X), \ \ \ \forall\,z\in \Omega\cap V,\ X\in T^{1,0}(\mathbb{C}^n).
$$
\end{lemma}

\begin{definition}
A domain $\Omega$ is called Bergman complete if $(\Omega,d_\Omega)$ is a complete metric space, that is, every Cauchy sequence converges.
\end{definition}

\begin{remark}
By the Hopf--Rinow theorem, $\Omega$ is Bergman complete if and only if there exists $a\in\Omega$ such that $d_\Omega(a,z)$ is an exhaustion function on $\Omega$.
\end{remark}

\begin{theorem}[Bremermann]\label{th:Brem}
Every Bergman complete domain $\Omega$ is pseudoconvex.
\end{theorem}

\begin{proof}
If $\Omega$ is not pseudoconvex, there exist $a\in\partial\Omega$ and $r>0$ such that every $f\in\mathcal O(\Omega)$ extends holomorphically to $\widehat f\in\mathcal O(B(a,r))$.
Consider the linear space
$$
H:=\{(f,\widehat{f}):f\in A^2(\Omega)\}\subset A^2(\Omega) \times A^2(B(a,r/2)),
$$
with norm $\|(f,\widehat f)\|:=\|f\|_{L^2(\Omega)}+\|\widehat f\|_{L^2(B(a,r/2))}$.
Clearly, $H$ is a closed subspace of $A^2(\Omega)\times A^2(B(a,r/2))$, hence a Banach space. Consider the natural projection
$$
\varpi:H\rightarrow A^2(\Omega),\ \ \ (f,\widehat{f})\mapsto f.
$$
Clearly, $\varpi$ is surjective. If $f_1=f_2$, the identity theorem gives $\widehat f_1=\widehat f_2$, so $\varpi$ is also injective. Since
$$
\|\varpi(f,\widehat{f})\|=\|f\|_{L^2(\Omega)}\le \|(f,\widehat{f})\|,
$$
we have $\|\varpi\|\le1$. The bounded inverse theorem implies that $\varpi^{-1}$ is bounded, so there exists $C>0$ such that
$$
\|\widehat{f}\|_{L^2(B(a,r/2))} \le \|(f,\widehat{f})\| \le C\|f\|_{L^2(\Omega)}.
$$
Let $\{h_j\}\subset A^2(\Omega)$ be a complete orthonormal basis. For $w\in B(a,r/3)$, put
$$
g=\sum_{j=k}^l \overline{\widehat{h}_j(w)}\, \widehat{h}_j,\ \ \  l>k.
$$
Bergman's inequality gives
$$
 |g(w)|^2  \le   C_r \|g\|^2_{L^2(B(a,r/2))}
 \le  C_r C^2 \|g\|_{L^2(\Omega)}^2 
 =  C_r C^2 \sum_{j=k}^l |\widehat{h}_j(w)|^2,
$$
that is, $\sum_{j=k}^l|\widehat h_j(w)|^2\le C_rC^2$. Thus $\widehat K(w):=\sum_{j=1}^\infty|\widehat h_j(w)|^2$ is real-analytic on $B(a,r/3)$. Since $\widehat K(w)=K_\Omega(w)$ on $\Omega\cap B(a,r/3)$, the metric $ds_\Omega^2$ is incomplete, a contradiction.
\end{proof}

\begin{corollary}[Siegel]
If a bounded domain $\Omega\subset\mathbb{C}^n$ is a regular unramified covering of a compact complex manifold $M$, then $\Omega$ is pseudoconvex.
\end{corollary}

\begin{proof}
Since $ds_\Omega^2$ is invariant under $\mathrm{Aut}(\Omega)$, it descends to a K\"ahler metric $ds_M^2$ on $M$. Compactness of $M$ implies completeness of $ds_M^2$, hence also of $ds_\Omega^2$. By Theorem \ref{th:Brem}, $\Omega$ is pseudoconvex.
\end{proof}

On the other hand, we have

\begin{theorem}[Kodaira]
If a bounded domain $\Omega\subset\mathbb{C}^n$ is a regular unramified covering of a compact complex manifold $M$, then $M$ is algebraic: it admits a holomorphic embedding into a complex projective space.
\end{theorem}

\begin{proof}
Write $dz=dz_1\wedge\cdots\wedge dz_n$ and $d\bar z=d\bar z_1\wedge\cdots\wedge d\bar z_n$. The form $K_\Omega(z)\,dz\wedge d\bar z$ is invariant under $\mathrm{Aut}(\Omega)$, giving a smooth Hermitian metric $|dz|^2=K_\Omega(z)^{-1}$ on the canonical line bundle, which descends to $M$. Its curvature is exactly the positive $(1,1)$-form on $M$ induced by the Bergman metric of $\Omega$, itself invariant under $\mathrm{Aut}(\Omega)$. Kodaira's embedding theorem proves the claim.
\end{proof}

By $L^2$ removability, $ds_{\mathbb{D}^\ast}^2$ is incomplete, although $\mathbb{D}^\ast$ is pseudoconvex. This leads naturally to the question:

\begin{problem}[Kobayashi {\cite{Kobayashi1959}}]
Which pseudoconvex domains are Bergman complete?
\end{problem}

\section{Kobayashi's Criterion}

Blow-up of the kernel at the boundary describes growth of point-evaluation functionals; completeness of the metric describes the length needed to approach the boundary. How are these two boundary behaviors related?

We first introduce infinite-dimensional complex projective space. Write
$$
l^2_{\mathbb{C}}:=\left\{(\zeta_0,\zeta_1,\cdots):\zeta_j\in \mathbb{C},\,\sum |\zeta_j|^2<\infty\right\}.
$$
Define an equivalence relation on $l^2_\mathbb{C}\setminus\{0\}$ by
$$
(\zeta_0',\zeta_1',\cdots)\sim(\zeta_0'',\zeta_1'',\cdots)\iff\exists\,c\in\mathbb{C}^\ast\ \text{such that }\zeta_j'=c\cdot\zeta_j'',\ j=0,1,\cdots.
$$
The quotient $\mathbb{P}^\infty:=(l^2_\mathbb{C}\setminus\{0\})/\sim$ is called infinite-dimensional complex projective space.

Let $\varpi:l^2_\mathbb{C}\setminus\{0\}\rightarrow\mathbb{P}^\infty$ be the natural projection. It induces the quotient topology: $U\subset\mathbb{P}^\infty$ is open if and only if $\varpi^{-1}(U)$ is open in $l^2_\mathbb{C}\setminus\{0\}$. For $\zeta=(\zeta_0,\zeta_1,\cdots)\in l^2_\mathbb{C}\setminus\{0\}$,
$$
\varpi(\zeta)=:[\zeta_0:\zeta_1:\cdots]
$$
are called homogeneous coordinates on $\mathbb{P}^\infty$. Define
$$
U_j=\{[\zeta_0:\zeta_1:\cdots]: \zeta_j\neq 0\},\ \ \ j=0,1,\cdots.
$$
Then $\{U_j\}_{j=0}^\infty$ is an open cover of $\mathbb{P}^\infty$. For each $j$, define the homeomorphism
$$
\Psi_j: U_j\rightarrow l^2_\mathbb{C},\ \ \ [\zeta_0:\zeta_1:\cdots]\mapsto \left(\frac{\zeta_0}{\zeta_j},\cdots,\frac{\zeta_{j-1}}{\zeta_j},\frac{\zeta_{j+1}}{\zeta_j},\cdots\right).
$$
This makes $\mathbb{P}^\infty$ an infinite-dimensional topological manifold. The coordinates
\[
\left(\frac{\zeta_0}{\zeta_j},\cdots,\frac{\zeta_{j-1}}{\zeta_j},\frac{\zeta_{j+1}}{\zeta_j},\cdots\right)
\]
are called inhomogeneous coordinates. Note that
$$
\Psi_k\circ \Psi_j^{-1}:\Psi_j(U_j\cap U_k)\rightarrow \Psi_k(U_j\cap U_k)
$$
$$
\left(\frac{\zeta_0}{\zeta_j},\cdots,\frac{\zeta_{j-1}}{\zeta_j},\frac{\zeta_{j+1}}{\zeta_j},\cdots\right)\mapsto 
\left(\frac{\zeta_0}{\zeta_k},\cdots,\frac{\zeta_{k-1}}{\zeta_k},\frac{\zeta_{k+1}}{\zeta_k},\cdots\right)
$$
is biholomorphic, so $\mathbb{P}^\infty$ is an infinite-dimensional complex manifold.

The classical Fubini--Study metric $ds_{FS}^2$ on $\mathbb{P}^\infty$ has K\"ahler form
$$
\omega_{FS}=\frac{i}2\partial \bar{\partial} \log \sum_{j=0}^\infty |\zeta_j|^2
$$
which, in the inhomogeneous coordinates on $U_0$, is
$$
\omega_{FS}=\frac{i}2\partial \bar{\partial} \log \left(1+\sum_{j=1}^\infty |\xi_j|^2\right), \ \ \ \xi_j=\zeta_j/\zeta_0.
$$

Let $\Omega\subset\mathbb{C}^n$ be bounded, and let $\{h_j\}$ be a complete orthonormal basis of $A^2(\Omega)$. Define the holomorphic map
$$
\iota:\Omega \rightarrow \mathbb{P}^\infty,\ \ \ z\mapsto [h_0(z):h_1(z):\cdots].
$$

\begin{lemma}\label{lm:imbed}
The map $\iota$ is a holomorphic embedding.
\end{lemma}

\begin{proof}
First we prove that $\iota$ is injective. Let $a,b\in\Omega$ be distinct. Since $\Omega$ is bounded, we can choose $h\in A^2(\Omega)$ with $h(a)=0$ and $h(b)\ne0$: take $h(z)=z_k-a_k$ for a coordinate with $a_k\ne b_k$. Write $h=\sum_{j=0}^{\infty}c_jh_j$. If $\iota(a)=\iota(b)$, there is a nonzero scalar $\lambda$ such that $h_j(a)=\lambda h_j(b)$ for every $j$. Hence
\[
0=h(a)=\sum_{j=0}^{\infty}c_jh_j(a)
=\lambda\sum_{j=0}^{\infty}c_jh_j(b)
=\lambda h(b)\ne0,
\]
a contradiction. Here we used the locally uniform convergence of the basis expansion, which follows from Bergman's inequality.

Next we prove that $\iota$ is an immersion. Fix $a\in\Omega$ and put
\[
A_0=\{f\in A^2(\Omega):f(a)=0\},\qquad
A_1=\{f\in A_0:\partial f(a)=0\}.
\]
Point evaluation and first derivatives at $a$ are continuous on $A^2(\Omega)$, so these are closed subspaces. Changing the orthonormal basis induces a unitary projective automorphism of $\mathbb{P}^{\infty}$. We may therefore choose a basis adapted to these subspaces without affecting the assertion.

Choose a unit vector $h_0$ spanning $A^2(\Omega)\ominus A_0$; then $h_0(a)\ne0$. To choose the next $n$ basis functions, set $S_0=A_0$ and
\[
S_k=\left\{f\in S_{k-1}:
\frac{\partial f}{\partial z_k}(a)=0\right\},
\qquad 1\le k\le n.
\]
Since $z_k-a_k\in S_{k-1}$ and its $z_k$-derivative is $1$, the space $S_k$ has codimension one in $S_{k-1}$. Choose a unit vector $h_k$ spanning $S_{k-1}\ominus S_k$. Finally, complete the basis with an orthonormal basis $\{h_j\}_{j\ge n+1}$ of $S_n=A_1$.

By construction, for $1\le j\le n$,
\[
h_j(a)=0,\qquad
\frac{\partial h_j}{\partial z_k}(a)=0\quad(k<j),\qquad
\frac{\partial h_j}{\partial z_j}(a)\ne0.
\]
In the projective chart where $h_0\ne0$, the map $\iota$ has coordinates
$h_1/h_0,h_2/h_0,\ldots$. Since
\[
\left.\frac{\partial(h_j/h_0)}{\partial z_k}\right|_a
=\frac{1}{h_0(a)}\frac{\partial h_j}{\partial z_k}(a),
\]
the $n\times n$ matrix formed from the first $n$ coordinate functions is upper triangular with nonzero diagonal entries. It is therefore invertible, and $\iota$ is an immersion at $a$.

Finally, we check that the inverse on the image is continuous. Since $\Omega$ is bounded, $z_kh_0\in A^2(\Omega)$, so we may expand
\[
z_kh_0=\sum_{j=0}^{\infty}d_{kj}h_j,\qquad
\sum_{j=0}^{\infty}|d_{kj}|^2<\infty.
\]
Wherever $h_0(z)\ne0$, this gives
\[
z_k=d_{k0}+\sum_{j=1}^{\infty}d_{kj}\frac{h_j(z)}{h_0(z)}.
\]
By Cauchy--Schwarz, the right-hand side is a continuous function of the $l^2$ chart coordinates. Thus the coordinates of $z$ depend continuously on $\iota(z)$ near $\iota(a)$. Since $a$ was arbitrary, $\iota^{-1}$ is continuous on its image, and $\iota$ is a holomorphic embedding.
\end{proof}

\begin{proposition}
$ds_\Omega^2=\iota^\ast(ds_{FS}^2)$ and $d_\Omega(z,w)\ge d_{FS}(\iota(z),\iota(w))$.
\end{proposition}

\begin{proof}
Note that
\begin{eqnarray*}
\omega_{\Omega} & = & \frac{i}2 \sum_{\alpha,\beta} \frac{\partial^2 \log \sum_{l=0}^\infty |h_l(z)|^2 }{\partial z_\alpha \partial \bar{z}_\beta} dz_\alpha \wedge d\bar{z}_\beta\\
& = & \left. \frac{i}2 \sum_{j,k}\sum_{\alpha,\beta} \frac{\partial ^2 \log \sum_{l=0}^\infty |\zeta_l|^2}{\partial \zeta_j \partial \bar{\zeta}_k} \frac{\partial \zeta_j}{\partial z_\alpha} \overline{\frac{\partial \zeta_k}{\partial z_\beta}} dz_\alpha \wedge d\bar{z}_\beta \right|_{\zeta_j=h_j(z)}\\
& = & \left.   \frac{i}2  \sum_{j,k}  \frac{\partial ^2 \log \sum_{l=0}^\infty |\zeta_l|^2}{\partial \zeta_j \partial \bar{\zeta}_k} d\zeta_j\wedge d\bar{\zeta}_k\right|_{\zeta_j=h_j(z)}\\
& = & \iota^\ast(\omega_{FS}).
\end{eqnarray*}
This identity and Lemma \ref{lm:imbed} imply the asserted inequality.
\end{proof}

\begin{proposition}
For $\zeta'=(\zeta_0',\zeta_1',\cdots)$ and $\zeta''=(\zeta_0'',\zeta_1'',\cdots)$, we have
$$
d_{FS}([\zeta'],[\zeta'']) \ge \arccos \frac{|\langle \zeta',\zeta''\rangle|}{|\zeta'| \cdot |\zeta''|},
$$
where $\langle\zeta',\zeta''\rangle=\sum\zeta_j'\overline{\zeta_j''}$ and $|\zeta'|=\sqrt{\langle\zeta',\zeta'\rangle}$.
\end{proposition}

\begin{proof}
Assume $|\zeta'|=|\zeta''|=1$. Both sides are invariant under unitary transformations, so we may take $[\zeta']=[1:0:0:\cdots]$. Let $\gamma:[0,1]\rightarrow\mathbb{P}^\infty$ be piecewise $C^1$ with
$\gamma(0)=[\zeta']$ and $\gamma(1)=[\zeta'']$. First assume the entire curve lies in $U_0:=\{[\zeta]:\zeta_0\ne0\}$. In the coordinates $\xi_j=\zeta_j/\zeta_0$,
$$
2 \omega_{FS} =\frac{i\partial\bar{\partial} |\xi|^2}{1+|\xi|^2} -\frac{i\sum_{j=1}^\infty \bar{\xi}_j d\xi_j \wedge \sum_{j=1}^\infty {\xi}_j d\bar{\xi}_j }{(1+|\xi|^2)^2}
 \ge \frac{i\partial\bar{\partial} |\xi|^2}{(1+|\xi|^2)^2}
$$
and $\gamma(0)=0$, $\gamma(1)=(\zeta_1''/\zeta_0'',\zeta_2''/\zeta_0'',\cdots)$. Hence its Fubini--Study length $|\gamma|$ satisfies
\begin{eqnarray*}
|\gamma| & \ge & \int_0^1 \sqrt{\frac{|\gamma'(t)|^2}{(1+|\gamma(t)|^2)^2}} dt\\
& \ge  &  \int_0^1 {\frac{|\gamma(t)|'}{1+|\gamma(t)|^2}} dt\\
& = & \arctan |\gamma(1)|\\
& = & \arccos|\zeta_0''|\ \ \ (\text{since }|\zeta''|=1)\\
& = &  \arccos \frac{|\langle \zeta',\zeta'' \rangle|}{|\zeta'| \cdot |\zeta''|}.
\end{eqnarray*}
If $\gamma$ leaves $U_0$, whether or not its endpoint is in $U_0$, put
$$
t_0=\max\{t:\gamma([0,t))\subset U_0\},  
$$
Then
$$
|\gamma|\ge \int_0^{t_0} |\gamma'(t)|_{FS}\,dt \ge \frac{\pi}2 \ge \arccos \frac{|\langle \zeta',\zeta'' \rangle|}{|\zeta'| \cdot |\zeta''|}.
$$
\end{proof}

\begin{corollary}\label{cor:B-dist}
\begin{eqnarray*}
d_\Omega(z,w) & \ge & \arccos \frac{\left|\sum_{j=0}^\infty h_j(z)\overline{h_j(w)}\right|}{\sqrt{\sum_{j=0}^\infty |h_j(z)|^2}\sqrt{\sum_{j=0}^\infty |h_j(w)|^2}} \\
& = & \arccos \frac{|K_\Omega(z,w)|}{\sqrt{K_\Omega(z)}\sqrt{K_\Omega(w)}}\\
& \ge & \sqrt{1-\frac{|K_\Omega(z,w)|^2}{K_\Omega(z)K_\Omega(w)}}.
\end{eqnarray*}
\end{corollary}

\begin{remark}
The transformation formula shows that $\frac{|K_\Omega(z,w)|^2}{K_\Omega(z)K_\Omega(w)}$ is biholomorphically invariant. It also has a special geometric meaning. In his doctoral dissertation, Calabi \cite{Calabi} introduced a new notion of distance, the \emph{diastasis}. Let $(M,ds^2)$ be a K\"ahler manifold with a local expression
$$
ds^2=\sum_{j,k} g_{j\bar{k}}(z,\bar{z}) dz_j\otimes d\bar{z}_k,
$$ 
where $g_{j\bar k}(z,\bar z)$ is real-analytic. The K\"ahler condition gives a real-analytic function $\Phi(z,\bar z)$ such that
$$
g_{j,\bar{k}}(z,\bar{z})=\frac{\partial^2 \Phi(z,\bar{z})}{\partial z_j\partial\bar{z}_k}.
$$ 
The Calabi diastasis associated with $ds^2$ is defined by
$$
D(z,w):= \Phi(z,\bar{z})+\Phi(w,\bar{w})-\Phi(z,\bar{w})-\Phi(w,\bar{z}),
$$
and is independent of the choice of local coordinates. Calabi introduced it to study isometric holomorphic embeddings of $M$ into $\mathbb{C}^N$ or $\mathbb{P}^N$, with $N\le\infty$. For the Bergman metric of a bounded domain $\Omega$, the diastasis is
$$
D_B(z,w)=-\log \frac{|K_\Omega(z,w)|^2}{K_\Omega(z)K_\Omega(w)}.
$$
Thus Corollary \ref{cor:B-dist} can be rewritten as
$$
d_\Omega(z,w)\ge \sqrt{1-\exp(-D_B(z,w))}.
$$

\end{remark}

\begin{proposition}
Suppose that every discrete sequence $\{a_k\}\subset\Omega$ has a subsequence $\{a_{k_j}\}$ such that
$$
\frac{|K_\Omega(a_{k_1},a_{k_j})|^2}{K_\Omega(a_{k_j})}\rightarrow 0\ \ \ (j\rightarrow \infty),
$$
Then $\Omega$ is Bergman complete.
\end{proposition}

\begin{proof}
If $\Omega$ is not Bergman complete, there is a $d_\Omega$-Cauchy sequence $\{a_k\}$ discrete in $\Omega$. Discard finitely many terms and renumber so that $d_\Omega(a_k,a_l)<1/2$ for all $k,l$.
Apply the hypothesis to this tail and take a subsequence $\{a_{k_j}\}$ such that
$$
\frac{|K_\Omega(a_{k_1},a_{k_j})|^2}{K_\Omega(a_{k_j})}\rightarrow 0\ \ \ (j\rightarrow \infty).
$$
On the other hand, Corollary \ref{cor:B-dist} gives
$$
1/2> d_\Omega(a_{k_1},a_{k_j})\ge \sqrt{1-\frac{|K_\Omega(a_{k_1},a_{k_j})|^2}{K_\Omega(a_{k_1})K_\Omega(a_{k_j})}}\ge 1/2,\ \ \ \forall\,j\gg 1.
$$
which is impossible.
\end{proof}

\begin{theorem}[Kobayashi's criterion \cite{Kobayashi1959,Kobayashi1962}]
Suppose there is a dense subset $\mathcal S\subset A^2(\Omega)$ such that, for every $f\in\mathcal S$ and every discrete sequence $\{a_k\}\subset\Omega$, there is a subsequence $\{a_{k_j}\}$ with
\begin{equation}\label{eq:Kobayashi}
\frac{|f(a_{k_j})|^2}{K_\Omega(a_{k_j})}\rightarrow 0 \ \ \ (j\rightarrow \infty),
\end{equation}
Then $\Omega$ is Bergman complete.
\end{theorem}

\begin{proof}
Suppose otherwise. There is a $d_\Omega$-Cauchy sequence $\{a_k\}$ discrete in $\Omega$. Choose $k_0\gg1$ with $d_\Omega(a_k,a_l)<1/2$ for all $k,l\ge k_0$.
By density of $\mathcal S$, the Gram--Schmidt process gives a complete orthonormal basis $\{\widetilde h_j\}\subset\mathrm{Span\,}\mathcal S$ of $A^2(\Omega)$.
Clearly, \eqref{eq:Kobayashi} holds for each $\widetilde h_j$. Consider
$$
h_0(z)=\frac{K_\Omega(z,a_{k_0})}{\sqrt{K_\Omega(a_{k_0})}}=\sum c_j \widetilde{h}_j(z),
$$
where $\sum|c_j|^2=1$. Choose $j_0$ with $\sum_{j=j_0+1}^\infty|c_j|^2<1/4$. Put $h_{j_0}=\sum_{j=1}^{j_0}c_j\widetilde h_j$. Then
$$
|h_0(z)-h_{j_0}(z)|^2 \le \sum_{j=j_0+1}^\infty |c_j|^2 \sum_{j=j_0+1}^\infty |\widetilde{h}_j(z)|^2 \le \frac{K_\Omega(z)}4.
$$
Therefore
$$
\frac{|h_0(z)|^2}{K_\Omega(z)}\le \frac12+\frac{2|h_{j_0}(z)|^2}{K_\Omega(z)}.
$$
Hence a subsequence $\{a_{k_j}\}$ satisfies $|h_0(a_{k_j})|^2/K_\Omega(a_{k_j})<3/4$ for all sufficiently large $j$. On the other hand, Corollary \ref{cor:B-dist} gives
$$
|h_0(a_{k_j})|^2/K_\Omega(a_{k_j})> 3/4, \ \ \ \forall\,j\gg 1,
$$
a contradiction.
\end{proof}

\begin{corollary}\label{cor:Kobayashi}
Let $\Omega\subset\mathbb{C}^n$ be bounded and satisfy:
\begin{enumerate}
\item $\Omega$ is Bergman exhaustive;
\item for every $a\in\partial\Omega$, the space $A^\infty_a(\Omega)$ is dense in $A^2(\Omega)$, where
$$
A^\infty_a(\Omega)=\left\{f\in A^2(\Omega): \limsup_{z\rightarrow a} |f(z)|<+\infty\right\}.
$$
Then $\Omega$ is Bergman complete.
\end{enumerate}
\end{corollary}

\begin{proof}
Given $f\in A^2(\Omega)$ and a discrete sequence $\{a_k\}\subset\Omega$, choose a subsequence $a_{k_j}\rightarrow a\in\partial\Omega$ and $\{f_\nu\}\subset A^\infty_a(\Omega)$ with
$\|f_\nu-f\|\rightarrow0$. Since
\begin{eqnarray*}
\frac{|f(a_{k_j})|^2}{K_\Omega(a_{k_j})} & \le & \frac{2|f(a_{k_j})-f_\nu(a_{k_j})|^2}{K_\Omega(a_{k_j})} + \frac{2|f_\nu(a_{k_j})|^2}{K_\Omega(a_{k_j})}\\
& \le & 2\|f_\nu-f\|^2 + \frac{2|f_\nu(a_{k_j})|^2}{K_\Omega(a_{k_j})},
\end{eqnarray*}
we obtain
$$
\limsup_{j\rightarrow \infty} |f(a_{k_j})|^2/K_\Omega(a_{k_j})\le 2\|f_\nu-f\|^2 \rightarrow 0, \ \ \ (\nu\rightarrow +\infty).
$$ 
Kobayashi's criterion proves the claim.
\end{proof}

It is natural to study the relationship between  Bergman exhaustiveness and Bergman completeness.

\begin{example}
The Hartogs triangle $\Omega=\{(z_1,z_2):|z_1|<|z_2|<1\}$ is Bergman exhaustive but not Bergman complete.
\end{example}

\begin{proof}
It is easy to verify that
$$
F:\Omega \rightarrow \mathbb{D}\times \mathbb{D}^\ast,\ \ \ (z_1,z_2)\mapsto (z_1/z_2,z_2)
$$
is biholomorphic. Hence
\begin{eqnarray*}
K_\Omega(z) & = & |J_F(z)|^2 K_{ \mathbb{D}\times \mathbb{D}^\ast}(F(z))\\
& = & \frac1{\pi^2 |z_2|^2(1-|z_1/z_2|^2)^2(1-|z_2|^2)^2}\\
&& \rightarrow +\infty \ \ \ (z\rightarrow \partial \Omega).
\end{eqnarray*}
On the other hand, Bergman completeness is biholomorphically invariant, so Proposition \ref{prop:prod} implies that $\Omega$ is not Bergman complete.
\end{proof}

On the other hand,  we have

\begin{theorem}[cf.  {\cite{Chen2000}}]\label{th:exhaus-complete}
Every bounded Bergman exhaustive domain in $\mathbb{C}$ is Bergman complete.
\end{theorem} 

\begin{proof}
Combine Corollary \ref{cor:Kobayashi} with Lemma \ref{lm:point_dense}.
\end{proof}

\begin{remark}
Zwonek \cite{ZwonekNagoya} even constructed a Bergman complete example that does not satisfy Kobayashi's criterion.
\end{remark}

An argument similar to that for Proposition \ref{prop:Zalcman} proves the following:

\begin{proposition}
Let $D$ be a domain of type $(L)$ satisfying
\begin{enumerate}
\item  $x_{k-1}-2\sqrt{r_{k-1}}> x_k+2\sqrt{r_k}$, $\forall\,k\ge 2$,
\item $\sum_k (x_k^4\log 1/r_k)^{-1}<+\infty$,
\end{enumerate}
Then $D$ is not Bergman complete.
\end{proposition}

The localization lemma for the metric also gives the following characterization.

\begin{theorem}[Locality of Bergman completeness]\label{th:BergmanLocalComplete}
Let $\Omega\subset\mathbb{C}^n$ be bounded and pseudoconvex. The following are equivalent:
\begin{enumerate}
\item $\Omega$ is Bergman complete.
\item For every $p\in\partial\Omega$, there is $r_p>0$ such that every component of $\Omega\cap B(p,r_p)$ is Bergman complete.
\item For every Euclidean ball $B\subset\mathbb{C}^n$, every component of $\Omega\cap B$ is Bergman complete.
\end{enumerate}
\end{theorem}

\begin{proof}
First recall two facts. By Hopf--Rinow, a Riemannian metric is complete if and only if every piecewise $C^1$ curve leaving all compact subsets has infinite length. Next, for any bounded domain $G$, choose a ball $B_0$ containing it. Comparison of the Bergman and Carath\'eodory metrics gives
\[
B_G(z;X)\ge C_G(z;X)\ge C_{B_0}(z;X)\ge c|X|,
\]
where $c>0$ is independent of $z,X$. Thus every curve of finite Bergman length has a Euclidean limit at its parameter endpoint; if it leaves all compact subsets of $G$, that limit lies in $\partial G$.

By Lemma \ref{lm:local}, given $p\in\partial\Omega$ and neighborhoods $V\Subset U$ of $p$, there is $A>0$ such that, for each component $D$ of $\Omega\cap U$,
\[
A^{-1}B_D(z;X)\le B_\Omega(z;X)\le A B_D(z;X),
\qquad z\in D\cap V.
\]
For disconnected intersections, the localization lemma is understood componentwise.

To prove $(2)\Rightarrow(1)$, suppose $\Omega$ is incomplete. There is a curve $\gamma:[0,1)\to\Omega$ leaving all compact subsets with finite Bergman length. By the preceding facts, $\gamma(t)\to p\in\partial\Omega$. Take $B(p,r_p)$ as in (2). A tail of $\gamma$ lies in $B(p,r_p/2)$ and in a single component $D$ of $\Omega\cap B(p,r_p)$. Local comparison gives
\[
\int_{t_0}^1 B_D(\gamma(t);\gamma'(t))\,dt
\le A\int_{t_0}^1 B_\Omega(\gamma(t);\gamma'(t))\,dt<\infty.
\]
But $\gamma(t)\to p\in\partial D$, contradicting completeness of $D$.

To prove $(1)\Rightarrow(3)$, fix $B$ and a component $D$ of $\Omega\cap B$. If $D$ is incomplete, there is a curve $\gamma:[0,1)\to D$ of finite Bergman length leaving all compact subsets and tending to some $q\in\partial D$. If $q\in B$, then $q\in\partial\Omega$. Localization near $q$ gives finite $B_\Omega$-length for a tail of $\gamma$, contradicting completeness of $\Omega$. If $q\in\partial B$, then
\[
B_D(z;X)\ge C_D(z;X)\ge C_B(z;X).
\]
The Carath\'eodory metric of the ball is complete, so a curve approaching its sphere has infinite $C_B$-length, hence infinite $B_D$-length, again a contradiction. Thus $D$ is complete.

Finally, $(3)\Rightarrow(2)$ is immediate.
\end{proof}

\begin{remark}
The proof compares the lengths of curve tails directly. It does not turn local comparison of infinitesimal metrics into comparison of global distances, since paths computing global distance may leave the local neighborhood. 
\end{remark}

Next, we use Kobayashi's criterion to establish the Bergman completeness of some domains.

\Needspace{6\baselineskip}
\begin{proposition}\label{prop:KR-complete}
If $\varphi$ is a positive subharmonic function on $\mathbb{D}$, then the Hartogs domain
$$
\Omega=\left\{(z_1,z_2)\in  \mathbb{D}\times \mathbb{C}:|z_2| < e^{-\varphi(z_1)}\right\}
$$
is Bergman complete.
\end{proposition}

\begin{proof}
Suppose otherwise. There is a $d_\Omega$-Cauchy sequence $\{a^{(k)}\}$ with $a^{(k)}\rightarrow a\in\partial\Omega$. Write $a=(a_1,a_2)$. We consider two cases:
1.  $a_1\in \partial \mathbb{D}$.
2.  $a_1\in \mathbb{D}$.

In the first case, note that
$$
B_\Omega(z;X)\ge C_\Omega(z;X)\ge C_{\mathbb{D}}(z_1;X_1),
$$
Hence, as $l\rightarrow\infty$,
$$
d_\Omega(a^{(k)},a^{(l)})\rightarrow \infty,  \ \ \ \forall\,k\in \mathbb{Z}^+, 
$$
contradicting the Cauchy property.

In the second case, choose a small disc $\Delta\Subset\mathbb{D}$ centered at $a_1$ and put $D=\Omega\cap(\Delta\times\mathbb{C})$. By Lemma \ref{lm:local}, the Bergman metrics of $D$ and $\Omega$ are comparable near $a$, so it suffices to prove completeness of $D$. After mapping $\Delta$ linearly onto $\mathbb{D}$, we may assume that $\varphi$ is defined on a neighborhood of $\overline{\mathbb{D}}$.

The power-series approximation in the proof of Proposition \ref{prop:KR} shows that $\Omega$ is a bounded pseudoconvex Runge domain, hence Bergman exhaustive. It remains to prove density of $\mathcal O(\overline\Omega)$ in $A^2(\Omega)$ and apply Kobayashi's criterion.

Given $f\in A^2(\Omega)$, consider its expansion
$$
f(z_1,z_2)=\sum_{j=0}^\infty f_j(z_1)z_2^j,
$$
where
$$
f_j(z_1)=\left.\frac1{j!}\frac{\partial^j f(z_1,z_2)}{\partial z_2^j}\right|_{z_2=0}  \in \mathcal O(\mathbb{D}).
$$
Note that
\begin{eqnarray*}
\int_\Omega |f|^2 & = & \sum_{j=0}^\infty \int_{z_1\in \mathbb{D}} |f_j(z_1)|^2 \int_{|z_2|<e^{-\varphi(z_1)}} |z_2|^{2j}\\
 & = & \sum_{j=0}^\infty \frac{\pi}{j+1} \int_{\mathbb{D}} |f_j|^2 e^{-2(j+1)\varphi}.
\end{eqnarray*}
Thus, if we put
$$
s_m(z_1,z_2)=\sum_{j=0}^m f_j(z_1)z_2^j,
$$
then
\begin{equation}\label{eq:KRC2}
\int_\Omega |f-s_m|^2 = \sum_{j=m+1}^\infty \frac{\pi}{j+1} \int_{\mathbb{D}} |f_j|^2 e^{-2(j+1)\varphi}\rightarrow 0\ \ \ (m\rightarrow \infty).
\end{equation}
If, for each $0\le j\le m$, we can find $\{f_{j,k}\}\subset\mathcal O(\overline{\mathbb{D}})$ such that
$$
\int_\mathbb{D} |f_j-f_{j,k}|^2 e^{-2(j+1)\varphi}\rightarrow 0\ \ \ (k\rightarrow \infty),
$$
then
$$
s_{m,k}(z_1,z_2):=\sum_{j=0}^m f_{j,k}(z_1)z_2^j\in \mathcal O(\overline{\Omega})
$$
and satisfies
\begin{equation}\label{eq:KRC3}
\int_\Omega |s_m-s_{m,k}|^2 =\sum_{j=0}^m \frac{\pi}{j+1} \int_\mathbb{D} |f_j-f_{j,k}|^2 e^{-2(j+1)\varphi}\rightarrow 0\ \ \ (k\rightarrow \infty).
\end{equation}
Equations \eqref{eq:KRC2} and \eqref{eq:KRC3} prove that $\mathcal O(\overline\Omega)$ is dense in $A^2(\Omega)$.

Finally, for fixed $j$, put
$$
f_{j,k}(z_1):=f_j\left(\frac{k z_1}{k+1}\right),\ \ \ k\in \mathbb{Z}^+.
$$
Then $f_{j,k}\in\mathcal O(\overline{\mathbb{D}})$. Since $\varphi>0$ is subharmonic near $\overline{\mathbb{D}}$, upper semicontinuity gives $\varphi\le C$ on $\mathbb{D}$. Thus $f_j\in A^2(\mathbb{D})$. Its power-series expansion readily gives
$$
\int_\mathbb{D} |f_j-f_{j,k}|^2 \rightarrow 0\ \ \ (k\rightarrow \infty).
$$
Hence
$$
\int_\mathbb{D} |f_j-f_{j,k}|^2 e^{-2(j+1)\varphi}\le  \int_\mathbb{D} |f_j-f_{j,k}|^2 \rightarrow 0\ \ \ (k\rightarrow \infty).
$$
\end{proof}

So the real difficulty comes from the singularities of $\varphi$. It is worthwhile to characterize  Bergman completeness of bounded pseudoconvex Hartogs domains in $\mathbb{C}^2$.\footnote{Chen-Forn\ae ss-Wu \cite{ChenFornaessWu} studied this problem in depth. }

\begin{theorem}\label{th:C-boundary}
Every bounded pseudoconvex domain $\Omega\subset\mathbb{C}^n$ with $C^0$ boundary is Bergman complete.
\end{theorem}

\begin{proof}
If $\Omega$ were incomplete, there would be a $d_\Omega$-Cauchy sequence $\{a_k\}$ discrete in $\Omega$. Passing to a subsequence, assume $a_k\rightarrow a\in\partial\Omega$. Theorems \ref{th:approximation}, \ref{th:ChenC^0}, and \ref{th:BergmanLocalComplete} give a contradiction.
\end{proof}

\begin{corollary}[Ohsawa {\cite{Ohsawa1981}}]
Every bounded pseudoconvex domain with $C^1$ boundary is Bergman complete.
\end{corollary}

Pseudoconvex domains with $C^0$ boundary are locally Runge along complex lines. This suggests the question: must every bounded pseudoconvex domain with the Runge along complex lines property be Bergman complete?

\begin{example}\label{ex:lineRungeIncomplete}
Fix $0<\varepsilon<1$ and define
$
\Omega_\varepsilon
 =\{(z,w)\in\mathbb{C}^2:0<|z|<1,\ |w-z^2|<\varepsilon|z|^2\}.
$
Then
\begin{enumerate}
\item $\Omega_\varepsilon$ is a bounded pseudoconvex domain with the Runge along complex lines property, but is not Runge in $\mathbb{C}^2$.
\item $\Omega_\varepsilon$ is Bergman exhaustive, but its Bergman metric is incomplete.
\end{enumerate}
\end{example}

\begin{proof}
Write $\mathbb{D}^\ast=\mathbb{D}\setminus\{0\}$. The map
\[
\Phi:\mathbb{D}^\ast\times\mathbb{D}\longrightarrow\Omega_\varepsilon,
\qquad \Phi(\xi,\zeta)=\bigl(\xi,\xi^2(1+\varepsilon\zeta)\bigr)
\]
is biholomorphic, with inverse
\[
\Phi^{-1}(z,w)=\left(z,\frac{w/z^2-1}{\varepsilon}\right),
\qquad J_\Phi(\xi,\zeta)=\varepsilon \xi^2.
\]
Thus $\Omega_\varepsilon$ is pseudoconvex. Since $|z|<1$ and $|w|<1+\varepsilon$, it is bounded.

First consider complex-line sections. A vertical line $z=c$ cuts out either the empty set or a disc. Every other complex line has the form $w=az+b$, and its section corresponds to the planar open set
\[
D_{a,b}=\{z: 0<|z|<1:|az+b-z^2|<\varepsilon|z|^2\}.
\]
Take a piecewise smooth Jordan curve $\gamma$ in $D_{a,b}$ and let $G$ be its bounded interior. We claim $0\notin G$. Otherwise, since on $\gamma$
\[
|(az+b)-z^2|<\varepsilon|z|^2<|z^2|,
\]
Rouch\'e's theorem implies that $az+b$ and $z^2$ have the same number of zeros in $G$, namely two counting multiplicities. This contradicts the fact that $az+b$ has degree at most one. Since $0\notin\gamma$ as well,
\[
h(z)=\frac{az+b-z^2}{z^2}
\]
is holomorphic near $\overline G$. The maximum modulus principle applied separately to $h$ and $z$ gives $|h|<\varepsilon$ and $|z|<1$ on $G$, so $G\subset D_{a,b}$. The Jordan-curve criterion shows that each section component is simply connected.

Next we show that the domain is not Runge. For $0<r<1$, the compact curve
\[
\Gamma_r=\{(re^{i\theta},r^2e^{2i\theta}):0\le\theta\le2\pi\}
\subset\Omega_\varepsilon
\]
is the boundary of the analytic disc $\zeta\mapsto(r\zeta,r^2\zeta^2)$. Applying the maximum modulus principle to every polynomial shows that the disc lies in the polynomial hull of $\Gamma_r$. In particular, the hull contains $(0,0)\notin\Omega_\varepsilon$, so $\Omega_\varepsilon$ is not Runge.

By Proposition \ref{prop:propertyL}, $K_{\Omega_\varepsilon}\ge c\delta_{\Omega_\varepsilon}^{-2}$, so $\Omega_\varepsilon$ is Bergman exhaustive. Biholomorphic invariance of the metric gives
\[
B_{\Omega_\varepsilon}\bigl((r,r^2);(1,2r)\bigr)
 =B_{\mathbb{D}^\ast\times\mathbb{D}}\bigl((r,0);(1,0)\bigr)
 =\frac{\sqrt2}{1-r^2}.
\]
Thus the curve $p_r=(r,r^2)$, $0<r\le1/2$, tends to the boundary point $(0,0)$ but has finite Bergman length:
\[
\int_0^{1/2}\frac{\sqrt2}{1-r^2}\,dr
 =\frac{\log3}{\sqrt2}<\infty.
\]
This proves incompleteness.
\end{proof}

\begin{problem}
Is every bounded pseudoconvex domain that is (locally) Runge Bergman complete?
\end{problem}

In a later chapter, we will prove that every bounded hyperconvex domain (indeed, any hyperconvex manifold) is Bergman complete; moreover, we will give sharp lower bounds for the Bergman distance. The proofs rely heavily on pluripotential theory.

\section{Curvature Properties}

Extremal problems for point values give the Bergman kernel, those for first derivatives give the metric, and those for second derivatives determine curvature.

\begin{definition}
The holomorphic sectional curvature of the Bergman metric is defined by
$$
{\rm HSC}(z;X):= \frac1{B_\Omega(z;X)^4}\sum_{\alpha,\beta,\gamma,\tau} R_{\alpha\bar{\beta}\gamma\bar{\tau}} X_\alpha \overline{X}_\beta X_\gamma
\overline{X}_\tau,
$$
where $z\in\Omega$, $X=\sum_\alpha X_\alpha\partial/\partial z_\alpha$,
$$
R_{\alpha\bar{\beta}\gamma\bar{\tau}}=-\frac{\partial^2 g_{\alpha\bar{\beta}} }{\partial z_\gamma \partial \bar{z}_\tau} 
+\sum_{\mu,\nu} g^{\bar{\mu}\nu}\frac{\partial g_{\alpha\bar{\mu}}}{\partial z_\gamma} \frac{\partial g_{\nu\bar{\beta}}}{\partial \bar{z}_\tau}
$$
$$
g_{\alpha\bar\beta}=\frac{\partial^2\log K_\Omega(z)}{\partial z_\alpha\partial\bar z_\beta};\ \ \ (g^{\bar\alpha\beta})=(g_{\alpha\bar\beta})^{-1},\ \text{that is, }\sum_\gamma g_{\alpha\bar\gamma}g^{\bar\gamma\beta}=\delta_{\alpha\beta}.
$$
The Ricci curvature of the Bergman metric is defined by
$$
{\rm Ric}:= i \sum_{\alpha,\beta} R_{\alpha\bar{\beta}}\, dz_\alpha\wedge d\bar{z}_\beta,   
$$
where
$$
R_{\alpha\bar{\beta}} := \sum_{\gamma,\tau} g^{\bar{\tau}\gamma} R_{\alpha\bar{\beta}\gamma\bar{\tau}} = -\frac{\partial^2 \log \det (g_{\mu\bar{\nu}})}{\partial z_\alpha \partial \bar{z}_\beta}.
$$
\end{definition}

For $z\in\Omega$, write
$$
A_1=\{f\in A^2(\Omega): f(z)=0\},\ \ \ A_2=\{f\in A_1:\partial f(z)=0\},
$$
$$
\mathcal M_\Omega^{(1)}(z;X)=\sup\{|Xf(z)|^2: f\in A_1, \|f\|=1\},
$$
$$
\mathcal M_\Omega^{(2)}(z;X)=\sup\left\{\left|\sum_{\alpha,\beta} \frac{\partial^2 f(z)}{\partial z_\alpha\partial z_\beta}X_\alpha X_\beta\right|^2: f\in A_2, \|f\|=1\right\}.
$$

\begin{theorem}\label{th:curvature}
\begin{equation}\label{eq:holo-curvature}
{\rm HSC}(z;X) = 2- \frac{K_\Omega(z) \mathcal M_\Omega^{(2)}(z;X)}{[\mathcal M_\Omega^{(1)}(z;X)]^2}.
\end{equation}
\end{theorem}

\begin{corollary}[Bergman, Fuks: $n=2$; Hua Loo-Keng: all $n$]\label{cor:Fuks-Hua}
$$
{\rm HSC}(z;X)<2.
$$
\end{corollary}

Choose a complete orthonormal basis $\{h_j\}_{j=0}^\infty$ of $A^2(\Omega)$ such that $h_0\in A_1^\bot$, $h_1,\cdots,h_n\in A_2^\bot\subset A_1$, and $h_j\in A_2$ for $j\ge n+1$. Then $K_\Omega(z)=|h_0(z)|^2$ and $h_j(z)=0$ for $j\ge1$, so
$$
\omega_\Omega (z) = \frac{i}{2K_\Omega(z)} \sum_{j=1}^n\partial h_j(z)\wedge \overline{\partial h_j(z)}.
$$

Since both sides of \eqref{eq:holo-curvature} are coordinate-independent, a complex affine change gives coordinates $\zeta=(\zeta_1,\cdots,\zeta_n)$ with $\zeta(z)=0$ and $\partial h_j/\partial\zeta_k(z)=h_0(z)\delta_{jk}$. Thus
$$
\omega_\Omega(z)=\frac{i}2 \left.\sum_{j=1}^n d\zeta_j\wedge d\bar{\zeta}_j\right|_z.
$$
For brevity, write $h_\alpha=\partial h/\partial\zeta_\alpha$, $h_{\alpha\beta}=\partial^2h/\partial\zeta_\alpha\partial\zeta_\beta$, and so on.

\begin{lemma}[Kobayashi \cite{Kobayashi1959}]\label{lm:curvature}
$$
R_{\alpha\bar{\beta}\gamma\bar{\tau}}(z)=\delta_{\alpha\beta} \delta_{\gamma\tau} + \delta_{\alpha\tau}\delta_{\beta\gamma}
-\frac1{K_\Omega(z)} \sum_{j=n+1}^\infty h_{j\alpha\gamma}(z) \overline{h_{j\beta\tau}(z)}.
$$
\end{lemma}

\begin{proof}[Lemma \ref{lm:curvature} $\Rightarrow$ Theorem \ref{th:curvature}]
Lemma \ref{lm:curvature} gives
$$
{\rm HSC}(z;X) = 2- \frac1{K_\Omega(z)|X|^4} \sum_{j=n+1}^\infty \left|\sum_{\alpha,\gamma} h_{j\alpha\gamma}(z)X_\alpha X_\gamma\right|^2.
$$
Thus it suffices to verify:
\begin{enumerate}
\item[$(1)$] $\mathcal M_\Omega^{(1)}(z;X)=K_\Omega(z)|X|^2$,
\item[$(2)$] $\mathcal M_\Omega^{(2)}(z;X)= \sum_{j=n+1}^\infty \left|\sum_{\alpha,\gamma} h_{j\alpha\gamma}(z)X_\alpha X_\gamma\right|^2. $
\end{enumerate}

Assertion (1) is the extremal property of the Bergman metric.

For (2), let $f\in A_2$ with $\|f\|=1$. Write $f=\sum_{j=n+1}^\infty c_jh_j$, where
$$
\sum_{j=n+1}^\infty |c_j|^2=1.
$$
Then
\begin{eqnarray*}
\left|\sum_{\alpha,\gamma} f_{\alpha\gamma}(z) X_\alpha X_\gamma\right|^2  
& = & \left|\sum_{j=n+1}^\infty c_j \sum_{\alpha,\gamma} h_{j\alpha\gamma}(z)X_\alpha X_\gamma\right|^2\\
& \le &  \sum_{j=n+1}^\infty \left| \sum_{\alpha,\gamma} h_{j\alpha\gamma}(z)X_\alpha X_\gamma\right|^2
\end{eqnarray*}
Equality is attained by taking $c_j$ proportional to $\overline{\sum_{\alpha,\gamma}h_{j\alpha\gamma}(z)X_\alpha X_\gamma}$ and normalizing.
\end{proof}

\begin{proof}[Proof of Lemma \ref{lm:curvature}]
Put $\phi_j=h_j/h_0$ for $j\ge1$. Near $z$,
\begin{eqnarray*}
\log K_\Omega & = & \log |h_0|^2 + \log\left(1+\sum |\phi_j|^2\right)\\
 & = & \log |h_0|^2 + \sum |\phi_j|^2-\frac12 \left(\sum |\phi_j|^2\right)^2+O(|\zeta|^6).
\end{eqnarray*}
A direct calculation gives
\begin{eqnarray*}
g_{\alpha\bar{\beta}} & = &\sum_j \phi_{j\alpha}\overline{\phi_{j\beta}} -\sum_j  \overline{\phi_j} \phi_{j\alpha}
\sum_j \phi_j\overline{\phi_{j\beta}}\\
&& - \sum_j |\phi_j|^2 \sum_j \phi_{j\alpha} \overline{\phi_{j\beta}} +O(|\zeta|^4),\\
\frac{\partial g_{\alpha\bar{\beta}}}{\partial \bar{\zeta}_\tau} & = & \sum_j \phi_{j\alpha} \overline{\phi_{j\beta\tau}}
-\sum_j \overline{\phi_{j\tau}} \phi_{j\alpha}\sum_j \phi_j \overline{\phi_{j\beta}}\\
&& -\sum_j \overline{\phi_j}\phi_{j\alpha}\sum_j \phi_j \overline{\phi_{j\beta\tau}}
 -\sum_j \phi_j\overline{\phi_{j\tau}} \sum_j \phi_{j\alpha} \overline{\phi_{j\beta}} \\
&& - \sum_j |\phi_j|^2 
\sum_j \phi_{j\alpha}\overline{\phi_{j\beta\tau}}+O(|\zeta|^3).
\end{eqnarray*}
Since, at $z$,
\begin{eqnarray}\label{eq:S3_1}
\phi_{j\alpha} & = & \frac{h_{j\alpha}}{h_0}-  \frac{h_j h_{0\alpha}}{h_0^2}
 =  \delta_{j\alpha},\nonumber\\
 \phi_{j\alpha\gamma} & = & \frac{h_{j\alpha\gamma}}{h_0} -  \frac{h_{j\alpha}h_{0\gamma}}{h_0^2} -  \frac{h_{j\gamma}h_{0\alpha}}{h_0^2}\\
& = &  \frac{h_{j\alpha\gamma}}{h_0} -  \frac{\delta_{j\alpha}h_{0\gamma}}{h_0} -  \frac{\delta_{j\gamma}h_{0\alpha}}{h_0},\nonumber
\end{eqnarray}
we have, at $z$,
\begin{eqnarray*}
\frac{\partial g_{\alpha\bar{\beta}}}{\partial \bar{\zeta}_\tau} & = & \sum_j \phi_{j\alpha} \overline{\phi_{j\beta\tau}}=\overline{\phi_{\alpha\beta\tau}},\\
\frac{\partial^2 g_{\alpha\bar{\beta}}}{\partial \zeta_\gamma\partial\bar{\zeta}_\tau} 
& = & \sum_j \phi_{j\alpha\gamma} \overline{\phi_{j\beta\tau}} -\sum_j \overline{\phi_{j\tau}} \phi_{j\alpha} 
\sum_j  \phi_{j\gamma}\overline{\phi_{j\beta}} \\
&& - \sum_j  \phi_{j\gamma}\overline{\phi_{j\tau}} \sum_j  \phi_{j\alpha}\overline{\phi_{j\beta}},
\end{eqnarray*}
so that
\begin{eqnarray}\label{eq:S3_2}
\sum_{\mu,\nu} g^{\bar{\mu}\nu} \frac{\partial g_{\alpha\bar{\mu}}}{\partial \zeta_\gamma} \frac{\partial g_{\nu\bar{\beta}}}{\partial \bar{\zeta}_\tau}
& = & \sum_\mu \phi_{\mu\alpha\gamma} \overline{\phi_{\mu\beta\tau}}
\end{eqnarray}
\begin{eqnarray}\label{eq:S3_3}
- \frac{\partial^2 g_{\alpha\bar{\beta}}}{\partial \zeta_\gamma\partial\bar{\zeta}_\tau}  & = &
\delta_{\alpha\tau}\delta_{\beta\gamma} + \delta_{\alpha\beta} \delta_{\gamma\tau} -\sum_j \phi_{j\alpha\gamma}\overline{\phi_{j\beta\tau}}\\
& = & \delta_{\alpha\tau}\delta_{\beta\gamma} + \delta_{\alpha\beta} \delta_{\gamma\tau}- |h_0|^{-2} \sum_j h_{j\alpha\gamma}\overline{h_{j\beta\tau}}\nonumber\\
&& +  |h_0|^{-2} \sum_j h_{j\alpha\gamma}\overline{h_{j\beta\tau}} -\sum_j \phi_{j\alpha\gamma}\overline{\phi_{j\beta\tau}}.\nonumber
\end{eqnarray}
On the other hand, \eqref{eq:S3_1} gives
$$
\sum_{j=n+1}^\infty \phi_{j\alpha\gamma}\overline{\phi_{j\beta\tau}}=|h_0|^{-2} \sum_{j=n+1}^\infty  h_{j\alpha\gamma}\overline{h_{j\beta\tau}},
$$
Hence
\begin{equation}\label{eq:S3_4}
 |h_0|^{-2} \sum_{j=1}^\infty h_{j\alpha\gamma}\overline{h_{j\beta\tau}} -\sum_{j=1}^\infty \phi_{j\alpha\gamma}\overline{\phi_{j\beta\tau}}
=  |h_0|^{-2} \sum_{j=1}^n h_{j\alpha\gamma}\overline{h_{j\beta\tau}} -\sum_{j=1}^n \phi_{j\alpha\gamma}\overline{\phi_{j\beta\tau}}.
\end{equation}
Equations (\ref{eq:S3_2})--(\ref{eq:S3_4}) prove the result.
\end{proof}

\begin{corollary}[Kobayashi {\cite{Kobayashi1959}}]
\begin{equation}\label{eq:Ricci_1}
\left(R_{\alpha\bar{\beta}}\right)<  (n+1) \left(g_{\alpha\bar{\beta}} \right);
\end{equation}
In other words,
$$
\log \mathrm{det}\left(g_{\alpha\bar{\beta}}\right)+(n+1)\log K_\Omega 
$$
is strictly plurisubharmonic on $\Omega$.

\end{corollary}

\begin{proof}
Lemma \ref{lm:curvature} gives
$$
R_{\alpha\bar{\beta}}(z) = (n+1)\delta_{\alpha\beta} - \frac{1}{K_\Omega(z)}\sum_\gamma \sum_{j=n+1}^\infty h_{j\alpha\gamma}(z) \overline{h_{j\beta\gamma}(z)},
$$
so that
$$
\sum_{\alpha,\beta} R_{\alpha\bar{\beta}}(z) X_\alpha \overline{X}_\beta - (n+1)|X|^2 = -\frac{1}{K_\Omega(z)}\sum_\gamma \sum_{j=n+1}^\infty \left|  X h_{j\gamma}(z)\right|^2. 
$$
For $f\in A_2$, we have
$$
f=\sum_{j=n+1}^\infty c_j h_j,\ \ \ \sum_{j=n+1}^\infty |c_j|^2 =\|f\|^2.  
$$
Thus
\begin{eqnarray*}
\sum_\gamma |X f_\gamma| &  \le & \sum_\gamma \sum_{j=n+1}^\infty |c_j| |X h_{j\gamma}|\\
& \le & \left[ \sum_\gamma \sum_{j=n+1}^\infty |c_j|^2\right]^{1/2} \left[ \sum_\gamma \sum_{j=n+1}^\infty \left|  X h_{j\gamma}\right|^2 \right]^{1/2},
\end{eqnarray*}
so that
$$
\sum_\gamma \sum_{j=n+1}^\infty \left|  X h_{j\gamma}(z)\right|^2 \ge \frac{|\sum_\gamma |X f_\gamma(z)| |^2}{n\|f\|^2}.
$$
If we choose
$$
f(\zeta)=\frac12\, \overline{X}_1 (\zeta_1-z_1)^2 + \sum_{\alpha=2}^n \overline{X}_\alpha (\zeta_\alpha-z_\alpha)(\zeta_1-z_1),
$$
then $Xf_1(z)=|X|^2$ and $n\|f\|^2\le C|X|^2$, where $C$ depends only on $n$ and the diameter of $\Omega$. Hence
$$
\sum_{\alpha,\beta} R_{\alpha\bar{\beta}}(z) X_\alpha \overline{X}_\beta - (n+1)|X|^2 \le -\frac{|X|^2}{C K_\Omega(z)},
$$
This proves \eqref{eq:Ricci_1}.
\end{proof}

Finally, we prove the celebrated ball theorem. The strategy is to use completeness and constant curvature to obtain a covering from the ball onto the domain, then use the Bergman kernel to rule out multiple sheets.

\begin{theorem}[Lu Qi-Keng's ball theorem]\label{th:lu-ball}
If a bounded domain $\Omega\subset\mathbb{C}^n$ has a complete Bergman metric of constant holomorphic sectional curvature, then $\Omega$ is biholomorphic to the unit ball $\mathbb{B}^n$.
\end{theorem}

The proof requires neither boundary smoothness nor a prior assumption of simple connectivity. Its geometric part invokes a standard result; the rest uses only K\"ahler potentials, polarization of kernels, and basic covering-space theory. See \cite{Lu1966} for the original paper and \cite[Section 7]{Yoo2017} for related arguments involving coverings and Bergman geometry.

\begin{proof}
We use the classical classification theorem in complex differential geometry: a simply connected complete K\"ahler manifold of constant holomorphic sectional curvature is holomorphically isometric to $\mathbb{P}^n$, $\mathbb{C}^n$, or $\mathbb{B}^n$ according as the curvature is positive, zero, or negative, with the standard model metric multiplied by a positive constant.

Lift the Bergman metric of $\Omega$ to its universal cover. The lifted metric remains complete with the same constant holomorphic sectional curvature $c$. If $c>0$, the universal cover is the compact space $\mathbb{P}^n$, so each coordinate function of its holomorphic map to $\Omega\subset\mathbb{C}^n$ is constant. If $c=0$, the universal cover is $\mathbb{C}^n$, and the coordinate functions of the covering map are bounded, hence constant by Liouville's theorem. Both contradict local biholomorphicity of the covering map. Thus $c<0$, and there is a holomorphic covering
\[
\varpi:\mathbb{B}^n\longrightarrow\Omega.
\]
In ball coordinates, a K\"ahler potential for the lifted metric is $-a\log(1-|z|^2)$ for some $a>0$. Since the Bergman metric has potential $\log K_\Omega$,
\begin{equation}\label{eq:lu-potential}
\partial\bar\partial\log K_\Omega(\varpi(z))
=a\,\partial\bar\partial[-\log(1-|z|^2)].
\end{equation}
The constant $a$ is not yet determined; only $a>0$ will be used below.

By \eqref{eq:lu-potential}, the function
\[
\log K_\Omega(\varpi(z))+a\log(1-|z|^2)
\]
is pluriharmonic. Since the ball is simply connected, it is twice the real part of a holomorphic function $F$. Put $h=e^F$. Then $h$ is nowhere zero on the ball, and
\begin{equation}\label{eq:lu-diagonal}
K_\Omega(\varpi(z))
=\frac{|h(z)|^2}{(1-|z|^2)^a}.
\end{equation}

We now use \emph{uniqueness of polarization}: a function holomorphic in $z$ and antiholomorphic in $w$ is uniquely determined by its diagonal values. Indeed, expand locally as
\[
\sum_{\alpha,\beta}c_{\alpha\beta}z^\alpha\overline w^{\,\beta},
\]
Setting $w=z$ gives a real-analytic function determining every coefficient $c_{\alpha\beta}$. The identity theorem extends the conclusion to the connected domain. Thus \eqref{eq:lu-diagonal} yields
\begin{equation}\label{eq:lu-polarized}
K_\Omega(\varpi(z),\varpi(w))
=\frac{h(z)\overline{h(w)}}{(1-\langle z,w\rangle)^a}.
\end{equation}
where $\langle z,w\rangle=\sum_jz_j\overline w_j$. Since $|\langle z,w\rangle|<1$, the factor $(1-\langle z,w\rangle)^{-a}$ is defined single-valuedly by
$\exp[-a\operatorname{Log}(1-\langle z,w\rangle)]$
with the branch equal to $1$ when $\langle z,w\rangle=0$. Hence \eqref{eq:lu-polarized} holds throughout $\mathbb{B}^n\times\mathbb{B}^n$.

Put $p=\varpi(0)$ and suppose $q\in\mathbb{B}^n$ satisfies $\varpi(q)=p$. We show that $q=0$. Setting $(z,w)=(0,0)$ and $(q,0)$ in \eqref{eq:lu-polarized} gives
\[
K_\Omega(p)=|h(0)|^2,
\qquad
K_\Omega(p)=h(q)\overline{h(0)}.
\]
Since $h(0)\ne0$, comparison gives $h(q)=h(0)$. Next set $(z,w)=(q,q)$ to obtain
\[
K_\Omega(p)
=\frac{|h(q)|^2}{(1-|q|^2)^a}
=\frac{K_\Omega(p)}{(1-|q|^2)^a}.
\]
Since $K_\Omega(p)>0$ and $a>0$, necessarily $q=0$. Thus
\[
\varpi^{-1}(p)=\{0\}.
\]
A covering over a connected base has a constant number of sheets. A singleton fiber therefore implies that the entire covering has one sheet. Hence $\varpi$ is biholomorphic.
\end{proof}

\begin{remark}
Completeness and constant curvature give a covering by the ball; the two-variable structure of the Bergman kernel then prevents distinct points from being identified. A general complete K\"ahler metric of constant negative curvature need not have this consequence: an annulus carries a complete Poincar\'e metric of constant negative curvature but is not biholomorphic to the disc.

Completeness cannot be omitted either. Every square-integrable holomorphic function on $\mathbb{B}^n\setminus\{0\}$ extends to the ball with the same norm. Its Bergman kernel and metric are therefore the restrictions of those of the ball, and the holomorphic sectional curvature remains constant. However, the Bergman metric is incomplete.

After the proof, biholomorphic invariance and the formula for the ball's Bergman kernel show that $a=n+1$ in \eqref{eq:lu-potential}. With the curvature convention of this section, the constant is necessarily $-2/(n+1)$.
\end{remark}

An old conjecture of S.-Y. Cheng \cite{Cheng} asserts that for every bounded strongly pseudoconvex domain $\Omega\subset \mathbb C^n$, Bergman metric is K\"ahler-Einstein (i.e., the Ricci curvature is constant) if and only if 
$\Omega$ is biholomorphic to $\mathbb B^n$.  The two-dimensional case was proved through the work of Fu-Wong \cite{FuWong} and Nemirovski-Shafikov \cite{NS2006}. The general-dimensional case was established by Huang-Xiao \cite{HuangXiao}.  An extended version of Cheng's conjecture, known as the Cheng-Yau conjecture, asserts that the strong pseudoconvexity condition can be weakened to weak pseudoconvexity. A recent breakthrough is that Hsiao-Huang-Li \cite{HsiaoHuangLi} proved the case of real-analytic boundary.

\chapter{Parameter Dependence of Weighted Bergman Kernels}

\section{Stability and Continuity}

Let $\Omega\subset\mathbb{C}^n$ be bounded. For $\varphi\in PSH(\Omega)$ with $\varphi\not\equiv-\infty$, define the weighted Bergman space by
$$
A^2_\varphi(\Omega) = \left\{f\in \mathcal O(\Omega): \|f\|_\varphi^2:=\int_
\Omega |f|^2 e^{-\varphi}<\infty\right\}.
$$
The submean-value property gives, for every $f\in A^2_\varphi(\Omega)$ and $z\in\Omega$,
\begin{eqnarray*}
|f(z)|^2 &\le & \frac{4^n n!}{\pi^n\delta_\Omega(z)^{2n}} \int_{B(z,\delta_\Omega(z)/2)} |f|^2\\
 & \le &  \frac{4^n n!}{\pi^n\delta_\Omega(z)^{2n}} \left[ \sup_{B(z,\delta_\Omega(z)/2)} e^\varphi\right] \int_{B(z,\delta_\Omega(z)/2)} |f|^2 e^{-\varphi}\\
& \le &  \frac{4^n n!}{\pi^n\delta_\Omega(z)^{2n}}\left[ \sup_{B(z,\delta_\Omega(z)/2)} e^\varphi\right] \,\|f\|_\varphi^2.
\end{eqnarray*}
Thus the classical theory of Bergman spaces extends in parallel to weighted spaces. Given a complete orthonormal basis $\{h_k\}$ of $A^2_\varphi(\Omega)$, define the weighted Bergman kernel by
$$
K_\varphi(z,w) = K_{\Omega,\varphi}(z,w) =\sum_k h_k(z)\cdot \overline{h_k(w)}.
$$ 
Write $K_\varphi(z)=K_\varphi(z,z)$. Many properties of the ordinary Bergman kernel, with parallel proofs, remain valid in the weighted setting. For example, we have the biholomorphic transformation formula:

\begin{proposition}\label{ch3:prop:W-Trans}
Let $F:\Omega_1\rightarrow\Omega_2$ be biholomorphic and $\varphi_1\in PSH(\Omega_1)$. Put $\varphi_2=\varphi_1\circ F^{-1}$. Then
$$
K_{\Omega_1,\varphi_1}(z,w) = K_{\Omega_2,\varphi_2}(F(z),F(w)) J_F(z)\overline{J_F(w)}.
$$
\end{proposition}

On the open set where $K_\varphi>0$, the form $\partial\bar\partial\log K_\varphi$ is semipositive. It defines a K\"ahler metric when weighted holomorphic functions separate tangent vectors, and is invariant when the domain and weight are transformed together as above.

We next give a very general upper semicontinuity theorem.

\begin{theorem}\label{ch3:th:USC}
Let $\Omega\subset\mathbb{C}^n$ be bounded, $\{\varphi_j\}\subset PSH(\Omega_j)$, and $\varphi\in PSH(\Omega)$. Suppose that, for every compact $S\subset\Omega$, there exists $j_0\in\mathbb{Z}^+$ such that
\begin{enumerate}
\item[$(1)$] $S\subset \Omega_j$, $\forall\,j\ge j_0$; 
\item[$(2)$] $\sup_{j\ge j_0}\sup_S \varphi_j<\infty$; 
\item[$(3)$] $\varphi_j\rightarrow\varphi$ a.e. on $S$.
\end{enumerate}
Then
$$
\limsup_{j\rightarrow \infty} K_{\Omega_j,\varphi_j}(z) \le K_{\Omega,\varphi}(z),\ \ \ \forall\,z\in \Omega.
$$
\end{theorem}

\begin{proof}
Fix $w\in\Omega$. If the limsup on the left is zero, the claim is immediate. Otherwise, choose a subsequence with positive kernel values and put
$$
\lim_{k\rightarrow \infty} K_{\Omega_{j_k},\varphi_{j_k}}(w) = \limsup_{j\rightarrow \infty} K_{\Omega_j,\varphi_j}(w).
$$
Choose $f_k\in\mathcal O(\Omega_{j_k})$ with $\int_{\Omega_{j_k}}|f_k|^2e^{-\varphi_{j_k}}=1$ and $|f_k(w)|^2=K_{\Omega_{j_k},\varphi_{j_k}}(w)$. For every compact $S\subset\Omega$, there are $j_0$ and a neighborhood $U$ of $S$ with $\sup_{j\ge j_0}\sup_U\varphi_j<\infty$. Hence $\|f_k\|_{L^2(U)}$ is uniformly bounded for large $k$, so $\{f_k\}$ is a normal family on $\Omega$. Passing to a subsequence, still denoted by $\{f_k\}$, we obtain local uniform convergence to $f\in\mathcal O(\Omega)$. Fatou's lemma gives
$$
\int_S|f|^2 e^{-\varphi} \le \liminf_{k\rightarrow \infty} \int_S |f_k|^2 e^{-\varphi_{j_k}}\le 1.
$$ 
Since $S$ is arbitrary, $\int_\Omega|f|^2e^{-\varphi}\le1$. Thus
$$
K_{\Omega,\varphi}(w) \ge |f(w)|^2 = \lim_{k\rightarrow \infty} |f_k(w)|^2 =  \limsup_{j\rightarrow \infty} K_{\Omega_j,\varphi_j}(w).
$$
\end{proof}

\begin{corollary}\label{ch3:cor:decreasing}
Let $\{\Omega_j\}$ be an increasing sequence of domains in $\mathbb{C}^n$ with $\Omega=\bigcup_j\Omega_j$. Suppose $\{\varphi_j\}\subset PSH(\Omega_j)$ satisfies $\varphi_j\downarrow\varphi\in PSH(\Omega)$, with $\varphi\not\equiv-\infty$. Then $K_{\Omega_j,\varphi_j}(z)\downarrow K_{\Omega,\varphi}(z)$ for every $z\in\Omega$.
\end{corollary}

\begin{proof}
The extremal property gives $K_{\Omega_j,\varphi_j}(z)\ge K_{\Omega,\varphi}(z)$. By the normal-family argument in Theorem \ref{ch3:th:USC}, it remains to verify monotone decrease. For $f\in A^2_{\varphi_{j+1}}(\Omega_{j+1})$ and $z\in\Omega_j$,
$$
\int_{\Omega_j} |f|^2 e^{-\varphi_j} \le \int_{\Omega_j} |f|^2 e^{-\varphi_{j+1}} \le  \int_{\Omega_{j+1}} |f|^2 e^{-\varphi_{j+1}}.
$$
The extremal property yields
$$
K_{\Omega_j,\varphi_j}(z)\ge \frac{|f(z)|^2}{\int_{\Omega_j} |f|^2 e^{-\varphi_j}} \ge  \frac{|f(z)|^2}{\int_{\Omega_{j+1}} |f|^2 e^{-\varphi_{j+1}}},\ \ \ \forall\, f\in A^2_{\varphi_{j+1}}(\Omega_{j+1}),
$$
so $K_{\Omega_j,\varphi_j}(z)\ge K_{\Omega_{j+1},\varphi_{j+1}}(z)$.
\end{proof}

\begin{example}[Failure of stability on a general bounded domain]
There exist a nonpseudoconvex bounded domain $\Omega\subset\mathbb{C}^2$ and smooth functions
$\varphi_j,\varphi\in PSH(\Omega)$ such that $\varphi_j\le0$ and
$\varphi_j\uparrow\varphi$ locally uniformly, but
\[
K_{\varphi_j}(z)\le K_\varphi(z)-c,
\qquad z\in\Omega,\quad j\ge1,
\]
for a constant $c>0$ independent of $z$ and $j$.
\end{example}

\begin{proof}
Let $\mathbb{B}^2=\{z\in\mathbb{C}^2:|z|<1\}$, and set
\[
\begin{aligned}
S&=\{z\in\mathbb{B}^2:|z|>1/2\},\\
T&=\left\{z\in\mathbb{B}^2:0<|z_1|,\quad
|z_2|<\frac{|z_1|}{\log(e/|z_1|)}\right\},
\qquad \Omega=S\cup T.
\end{aligned}
\]
Thus $\Omega$ consists of a spherical shell together with a thin region approaching
the origin. It is a bounded connected domain and $0\notin\Omega$. By Hartogs'
extension theorem, every $f\in\mathcal O(\Omega)$ has a unique holomorphic extension
$\widetilde f$ to $\mathbb{B}^2$: first extend from the shell $S$, and then use the
identity theorem on $\Omega$.

Define
\[
\varphi_j(z)=\left(4+\frac1j\right)\log|z|,
\qquad \varphi(z)=4\log|z|.
\]
These functions are smooth and plurisubharmonic on $\Omega$, and satisfy
$\varphi_j\le0$ and $\varphi_j\uparrow\varphi$ locally uniformly.

We first examine integrability near the origin. Fix $0<\varepsilon<1/4$ and put
$T_\varepsilon=T\cap\{|z_1|<\varepsilon\}$. On this set, with $r=|z_1|$,
\[
r\le |z|\le {2}\,r,
\]
and the $z_2$-section has area $\pi r^2/[\log(e/r)]^2$. Hence, for $\delta\ge0$,
\[
\int_{T_\varepsilon}|z|^{-4-\delta}
\asymp\int_0^\varepsilon
\frac{dr}{r^{1+\delta}[\log(e/r)]^2}.
\]
The last integral is finite when $\delta=0$ and infinite when $\delta>0$.
Since $|z|\ge\varepsilon$ on $\Omega\setminus T_\varepsilon$, it follows that
\[
M:=\int_\Omega e^{-\varphi}
=\int_\Omega |z|^{-4}<\infty.
\]
In particular, $1\in A^2_\varphi(\Omega)$. On the other hand, every
$f\in A^2_{\varphi_j}(\Omega)$ satisfies $\widetilde f(0)=0$. Indeed, otherwise
$|\widetilde f|$ would have a positive lower bound near the origin, and the
preceding calculation with $\delta=1/j$ would give
$\int_\Omega |f|^2e^{-\varphi_j}=\infty$.

Let
\[
H_0=\{f\in A^2_\varphi(\Omega):\widetilde f(0)=0\}.
\]
We claim that $\langle f,1\rangle_\varphi=M\widetilde f(0)$ for every
$f\in A^2_\varphi(\Omega)$. First, the Cauchy--Schwarz inequality gives
\[
\int_\Omega |f|e^{-\varphi}
\le \left(\int_\Omega |f|^2e^{-\varphi}\right)^{1/2}
     \left(\int_\Omega e^{-\varphi}\right)^{1/2}
=\sqrt{M}\,\|f\|_\varphi<\infty.
\]
For every $\theta\in\mathbb{R}$, the rotation $z\mapsto e^{i\theta}z$ maps
$\Omega$ onto itself, and leaves
$e^{-\varphi(z)}=|z|^{-4}$ unchanged. The change of variables
$w=e^{i\theta}z$ therefore yields
\[
\int_\Omega fe^{-\varphi}
=\int_{z\in \Omega} f(e^{i\theta}z)e^{-\varphi(z)}.
\]
For each fixed $z\in\Omega$, the function
$g_z(\lambda)=\widetilde f(\lambda z)$ is holomorphic on
$\{|\lambda|<1/|z|\}$, a disc of radius greater than $1$. Hence the
one-variable mean-value formula gives
\[
\frac1{2\pi}\int_0^{2\pi}f(e^{i\theta}z)\,d\theta
=\frac1{2\pi}\int_0^{2\pi}g_z(e^{i\theta})\,d\theta
=g_z(0)=\widetilde f(0).
\]
Here the extension $\widetilde f$ allows us to apply the formula on the whole
disc even though $0\notin\Omega$.

To average the preceding change-of-variables identity over $\theta$, note that
the same change of variables with $|f|$, together with Tonelli's theorem, gives
\[
\int_0^{2\pi}\int_\Omega
|f(e^{i\theta}z)|e^{-\varphi(z)}\,dV(z)\,d\theta
=2\pi\int_\Omega |f|e^{-\varphi}<\infty.
\]
Thus Fubini's theorem applies, and we obtain
\[
\begin{aligned}
\langle f,1\rangle_\varphi
&=\frac1{2\pi}\int_0^{2\pi}\int_\Omega
f(e^{i\theta}z)e^{-\varphi(z)}\,dV(z)\,d\theta\\
&=\int_\Omega\left(\frac1{2\pi}\int_0^{2\pi}
f(e^{i\theta}z)\,d\theta\right)e^{-\varphi(z)}\,dV(z)\\
&=\widetilde f(0)\int_\Omega e^{-\varphi(z)}\,dV(z)
=M\widetilde f(0).
\end{aligned}
\]
Consequently, $H_0$ is closed and
\[
A^2_\varphi(\Omega)=\mathbb{C}\cdot1\mathbin{\oplus^\perp}H_0.
\]
If $K_{H_0}$ denotes the reproducing kernel of $H_0$, this orthogonal
decomposition yields
\[
K_{H_0}(z,z)=K_\varphi(z)-\frac1M.
\]
Finally, $A^2_{\varphi_j}(\Omega)\subset H_0$ and
$\|f\|_\varphi\le\|f\|_{\varphi_j}$ because $\varphi_j\le\varphi$.
The extremal characterization of the kernel therefore gives
\[
K_{\varphi_j}(z)\le K_{H_0}(z,z)
=K_\varphi(z)-\frac1M,
\qquad z\in\Omega.
\]
Thus the kernels fail to converge to $K_\varphi$ at every point of $\Omega$.
\end{proof}

The major  stability result is the following

\begin{theorem}[cf.  {\cite{ChenParameter2016}}]\label{ch3:th:singular}
Let $\Omega\subset\mathbb{C}^n$ be bounded and hyperconvex. Let $\{\varphi_j\}\subset PSH(\Omega)$ be uniformly bounded above on compact subsets, and let $\varphi$ be a measurable function finite almost everywhere on $\Omega$.
Suppose there is a closed pluripolar set $E\subset\Omega$ such that
\begin{equation}\label{ch3:eq:S1}
e^{-\varphi_j}\rightarrow e^{-\varphi}\quad\text{in }L^1_{\rm loc}(\Omega\backslash E).
\end{equation}
Then $K_{\varphi_j}=K_{\Omega,\varphi_j}\rightarrow K_\varphi=K_{\Omega,\varphi}$. The assumptions imply that $\varphi$ is locally essentially bounded above, so point-evaluation estimates and the definition of the weighted kernel remain valid.
\end{theorem}

Throughout this section, plurisubharmonic functions are understood {\bf not to be identically $-\infty$}. Here are some common situations satisfying \eqref{ch3:eq:S1}.

\begin{example}
Consider $\{\psi_t:|t|<1\}\subset PSH^-(\Omega)$ such that $e^{\psi_t(z)}$ is continuous in $(z,t)$. Put $E=\psi_0^{-1}(-\infty)$ and $\varphi_j=\psi_{1/j}$. Then $e^{-\varphi_j}$ converges to $e^{-\psi_0}$ uniformly on compact subsets of $\Omega\backslash E$, so \eqref{ch3:eq:S1} holds. A case useful in algebraic geometry is
$$
\psi_t(z)=\log \sum_{j=1}^m c_j(t) |f_j(z)|^{\alpha_j},
$$
where $f_j\in\mathcal O(\Omega)$, $\alpha_j>0$, and the $c_j$ are positive continuous functions.
\end{example}

\begin{example}\label{ch3:ex:common-lower}
Let $\{\varphi_j\}\subset PSH(\Omega)$ be uniformly bounded above on compact subsets and satisfy:
\begin{enumerate}
\item[$(1)$] $\varphi_j\rightarrow\varphi\in PSH(\Omega)$ almost everywhere;
\item[$(2)$] there is $\psi\in PSH(\Omega)$ with $\varphi_j\ge\psi$ for all $j$.
\end{enumerate}
Put
$$
E=\left\{z\in\Omega:e^{-\psi}\ \text{is not integrable in any neighborhood of }z\right\}.
$$
By Bombieri's theorem, $E$ is analytic, hence pluripolar. For every compact $S\subset\Omega\backslash E$, we have $e^{-\varphi_j}\le e^{-\psi}\in L^1(S)$. Dominated convergence gives $e^{-\varphi_j}\rightarrow e^{-\varphi}$ in $L^1(S)$.
\end{example}

\begin{example}
If $\{\varphi_j\}\subset PSH(\Omega)$ satisfies $\varphi_j\downarrow\varphi\in PSH(\Omega)$, then $\varphi\le\varphi_j\le\varphi_1$, so Example \ref{ch3:ex:common-lower} implies (\ref{ch3:eq:S1}).
\end{example}

\begin{example}
If $\{\varphi_j\}\subset PSH(\Omega)$ satisfies $\varphi_j\uparrow\varphi\in PSH(\Omega)$, then $\varphi_1\le\varphi_j\le\varphi$, so Example \ref{ch3:ex:common-lower} implies (\ref{ch3:eq:S1}).
\end{example}

\begin{example}
Let $\Omega\subset\mathbb{C}^n$ be bounded and pseudoconvex, and let $\varphi\in PSH(\Omega)$. Put $\varphi_j=j^{-1}\log K_{j\varphi}$. Demailly's approximation theorem gives, for $r<\delta_\Omega(z)$,
$$
\varphi(z)-C_1 j^{-1} \le \varphi_j(z) \le \sup_{B(z,r)} \varphi +j^{-1} \log \frac{C_2}{r^{2n}}.
$$
Thus $\{\varphi_j\}$ is uniformly bounded above on compact subsets, $\varphi_j\rightarrow\varphi$, and $\varphi_j\ge\varphi-C_1\in PSH(\Omega)$. Example \ref{ch3:ex:common-lower} implies (\ref{ch3:eq:S1}).
\end{example}

In fact,  even $\varphi_j\to\varphi$ in $L^1_{\rm loc}(\Omega)$ also implies (\ref{ch3:eq:S1}) (see  the proof of Theorem \ref{ch3:th:kernel-L1-stability} below).

\begin{corollary}[Berndtsson's openness theorem \cite{BerndtssonOpenness2015}]
Let $\psi\in PSH(\Omega)$ and $0\in\Omega$. Put
$$
c_0(\psi)=\sup\left\{c\ge0:e^{-c\psi}\ \text{is integrable near }0\right\}.
$$
If $c_0(\psi)<\infty$, then $e^{-c_0(\psi)\psi}$ is not integrable in any neighborhood of $0$.
\end{corollary}

\begin{proof}
Shrink $\Omega$ to a ball centered at $0$ and subtract a constant so that $\psi<0$. Put $\varphi_j:=(c_0(\psi)+1/j)\psi$ and $\varphi:=c_0(\psi)\psi$. If $e^{-\varphi}$ were integrable near $0$, the H\"ormander--Bombieri theorem would give $K_\varphi(0)>0$. Example \ref{ch3:ex:common-lower} and Theorem \ref{ch3:th:singular} give $K_{\varphi_j}(0)\to K_\varphi(0)>0$, hence $K_{\varphi_j}(0)>0$ for large $j$. But $e^{-\varphi_j}$ is not integrable in any neighborhood of $0$, so $K_{\varphi_j}(0)=0$, a contradiction.
\end{proof}

The proof of Theorem \ref{ch3:th:singular} is divided into the following lemmas.

\begin{lemma}\label{ch3:lm:S_1}
Let $\Omega\subset\mathbb{C}^n$ be bounded and hyperconvex, let $\rho$ be a negative continuous plurisubharmonic exhaustion, let $\varphi\in PSH(\Omega)$, and let $S\subset\Omega$ be compact. Put $\Omega_\varepsilon=\{\rho<-\varepsilon\}$ for $\varepsilon\ll1$. For every $w\in S$ and $f\in A^2_\varphi(\Omega_\varepsilon)$, there is $g\in A^2_\varphi(\Omega)$ with $g(w)=f(w)$ and
$$
\|g\|_{L^2_\varphi(\Omega)} \le (1+C/|\log\varepsilon|) \|f\|_{L^2_\varphi(\Omega_\varepsilon)},
$$
where $C$ depends only on $S,\Omega,\rho$.
\end{lemma}

\begin{proof}
Replacing $\rho$ by $\max\{\rho,-1/2\}$, assume $-1<\rho<0$. Choose a cutoff $\chi:\mathbb{R}\rightarrow[0,1]$ with $\chi|_{(-\infty,-\log2]}=1$ and $\chi|_{[0,\infty)}=0$. Put
$$
\eta_\varepsilon=\chi(\log(-\log(-\rho))-\log(-\log\varepsilon)).
$$
By the Donnelly--Fefferman estimate, $\bar\partial u=f\,\bar\partial\eta_\varepsilon$ has a solution $u_\varepsilon$ with
\begin{eqnarray*}
\int_\Omega |u_\varepsilon|^2 e^{-\varphi-2n\log |z-w|}  & \le & C_0 \int_\Omega |f|^2 |\bar{\partial}\eta_\varepsilon|^2_{i\partial\bar{\partial}(-\log(-\rho))}e^{-\varphi-2n\log |z-w|}\\
& \le & \frac{C}{|\log\varepsilon|^2} \int_{\Omega_\varepsilon} |f|^2 e^{-\varphi}.
\end{eqnarray*}
For sufficiently small $\varepsilon$, $\eta_\varepsilon=1$ on a fixed neighborhood of $S$. The singular-weight estimate ensures $u_\varepsilon(w)=0$ and $\|u_\varepsilon\|_\varphi\le C|\log\varepsilon|^{-1}\|f\|_{\Omega_\varepsilon,\varphi}$. Thus take $g=\eta_\varepsilon f-u_\varepsilon$.
\end{proof}

\begin{lemma}\label{ch3:lm:S_2}
Let $V\subset\subset\Omega\subset\mathbb{C}^n$ be bounded pseudoconvex domains. If $\varphi_j,\varphi$ satisfy the hypotheses of Theorem \ref{ch3:th:singular}, then for every $f\in A^2_\varphi(\Omega)$ there are $f_j\in A^2_{\varphi_j}(V)$ such that
\begin{enumerate}
\item[$(1)$] $\limsup_{j\rightarrow \infty} \|f_j\|_{L^2_{\varphi_j}(V)} \le \|f\|_{L^2_\varphi(\Omega)}$;
\item[$(2)$]  $\|f_j-f\|_{L^2(V)}\rightarrow 0\ \ \ (j\rightarrow \infty)$.
\end{enumerate}
\end{lemma}

\begin{proof}
Choose $\varrho\in PSH(\Omega)$ with $\varrho|_E=-\infty$ and $\varrho<-2e$ on $\overline V$. Take a nonnegative radial standard mollifier $\chi_\delta$ and put
$$
\varrho_\delta(z) = \varrho\ast \chi_\delta(z)+\delta |z|^2,\ \ \ \delta>0.
$$
For $\delta\ll1$, the function $\varrho_\delta$ is smooth and strictly plurisubharmonic near $\overline V$, with $\varrho_\delta\downarrow\varrho$. Since $E\cap\overline V$ is compact and $\varrho=-\infty$ there, for every $\varepsilon>0$ choose $\delta=\delta(\varepsilon)\downarrow0$ so that $\varrho_\delta<-e^{1/\varepsilon}$ near $E\cap\overline V$. Thus $\overline V\cap\{\varrho_\delta\ge-e^{1/\varepsilon}\}$ is a compact subset of $\Omega\setminus E$.

With $\chi$ as above, put $\eta_\varepsilon=\chi(\log\log(-\varrho_\delta)+\log\varepsilon)$. By the Donnelly--Fefferman estimate,\footnote{In fact, the classical H\"ormander $L^2$ estimate suffices; see \cite{Chen2022}.} the equation $\bar\partial u=f\,\bar\partial\eta_\varepsilon$ has a solution $u_{j,\varepsilon}$ satisfying
$$
\int_V |u_{j,\varepsilon}|^2 e^{-\varphi_j} \le C_0 \int_V |f|^2 |\bar{\partial}\eta_\varepsilon|^2_{i\partial\bar{\partial}(-\log(-\varrho_\delta))} e^{-\varphi_j}
\le C_0\varepsilon^2 \int_{S_\varepsilon} |f|^2 e^{-\varphi_j},
$$
where $S_\varepsilon:=\overline V\cap\{\log(-\varrho_\delta)\le1/\varepsilon\}\subset\Omega\backslash E$ is compact. Since $f\in L^\infty(S_\varepsilon)$ and $e^{-\varphi_j}\rightarrow e^{-\varphi}$ in $L^1(S_\varepsilon)$,
$$
\int_{S_\varepsilon} |f|^2 e^{-\varphi_j}\rightarrow \int_{S_\varepsilon} |f|^2 e^{-\varphi}\ \ \ (j\rightarrow \infty).
$$
Put $f_{j,\varepsilon}=\eta_\varepsilon f-u_{j,\varepsilon}\in\mathcal O(V)$. The preceding integral convergence and the triangle inequality give, for fixed $\varepsilon$, $\limsup_{j\to\infty}\|f_{j,\varepsilon}\|_{L^2_{\varphi_j}(V)}\le(1+\sqrt{C_0}\varepsilon)\|f\|_{L^2_\varphi(\Omega)}$. Thus, for $j\ge j_\varepsilon$ sufficiently large, $\|f_{j,\varepsilon}\|_{L^2_{\varphi_j}(V)}\le[1+(1+\sqrt{C_0})\varepsilon]\|f\|_{L^2_\varphi(\Omega)}$. If $f=0$, take $f_{j,\varepsilon}=0$.

Since the $\varphi_j$ are uniformly bounded above on $V$, for $j\ge j_\varepsilon\gg1$ we have
$$
\|f_{j,\varepsilon}-f\|_{L^2(V)}^2 \le 2 \int_{V\cap \{\log(-\varrho_\delta)\ge 1/(2\varepsilon)\}} |f|^2 +C\varepsilon^2 \int_\Omega |f|^2 e^{-\varphi}.
$$
Since
$$
V\cap\left\{\log(-\varrho_\delta)\ge1/(2\varepsilon)\right\}\subset V\cap\left\{\varrho\le-e^{1/(2\varepsilon)}\right\}\ \ \ (\text{since }\varrho_\delta\ge\varrho),
$$
we obtain
$$
\limsup_{\varepsilon\rightarrow 0} \int_{V\cap \{\log(-\varrho_\delta)\ge 1/(2\varepsilon)\}} |f|^2 \le \int_{V\cap \varrho^{-1}(-\infty)} |f|^2 =0. 
$$
A diagonal choice yields the required full sequence.\footnote{The following argument is due to Yuanpu Xiong. Take $\varepsilon_k=2^{-k}$. We may assume that the corresponding $j_{\varepsilon_k}$ increase to infinity. For $j_{\varepsilon_k}\le j<j_{\varepsilon_{k+1}}$, define $f_j=f_{j,\varepsilon_k}$. This sequence has the required properties.}
\end{proof}

\begin{proof}[Proof of Theorem \ref{ch3:th:singular}]
Let $S\subset\Omega$ be compact and $w\in S$. If $K_\varphi(w)=0$, the liminf inequality is immediate. Assume it is positive and put $f=K_\varphi(\cdot,w)$. Lemma \ref{ch3:lm:S_2} gives $f_{j,\varepsilon}\in A^2_{\varphi_j}(\Omega_\varepsilon)$ such that
$$
\limsup_{j\rightarrow \infty} \|f_{j,\varepsilon}\|_{L^2_{\varphi_j}(\Omega_\varepsilon)}\le \|f\|_{L^2_\varphi(\Omega)} = \sqrt{K_\varphi(w)}, 
$$
For $\varepsilon\ll1$, we also have $f_{j,\varepsilon}(w)\rightarrow f(w)$ as $j\rightarrow\infty$. By Lemma \ref{ch3:lm:S_1}, there exists $g_{j,\varepsilon}\in A^2_{\varphi_j}(\Omega)$ with $g_{j,\varepsilon}(w)=f_{j,\varepsilon}(w)$ and
$$
\|g_{j,\varepsilon}\|_{L^2_{\varphi_j}(\Omega)} \le (1+C/|\log \varepsilon|) \|f_{j,\varepsilon}\|_{L^2_{\varphi_j}(\Omega_\varepsilon)},
$$
where $C$ depends only on $S$, the fixed domain, and the exhaustion, not on $j$. Thus
$$
\liminf_{j\rightarrow \infty} K_{\varphi_j}(w) \ge \liminf_{j\rightarrow \infty} \frac{|g_{j,\varepsilon}(w)|^2}{\|g_{j,\varepsilon}\|_{L^2_{\varphi_j}(\Omega)}^2 }
\ge \frac{K_\varphi(w)}{ (1+C/|\log\varepsilon|)^{2}} \rightarrow K_\varphi(w)\ \ \ (\varepsilon\rightarrow 0).
$$
For the limsup, take a further subsequence of one realizing it. By \eqref{ch3:eq:S1}, arrange $e^{-\varphi_j}\to e^{-\varphi}$ almost everywhere on $\Omega\setminus E$; the pluripolar set $E$ has measure zero. The normal-family and Fatou arguments from Theorem \ref{ch3:th:USC} apply to this measurable limit weight, giving $\limsup_jK_{\varphi_j}(w)\le K_\varphi(w)$. Normal families and uniqueness of polarization then yield uniform convergence on compact subsets for the two-variable kernels.
\end{proof}

Finally, we would like to point out that  Lemma \ref{ch3:lm:S_2} also implies the strong openness theorem below.\footnote{After I shared this short proof with Liyou Zhang, he introduced me to his joint paper with his student Ziyi Zhang \cite{ZhangZhang}. They independently discovered that stability of weighted Bergman kernels implies strong openness. Their method partly overlaps with the one here.}

\begin{theorem}[Guan--Zhou  {\cite{GuanZhou2015}}]\label{th:strong-openness}
Let $\varphi$ be plurisubharmonic near the origin and $f$ holomorphic. If $|f|^2e^{-\varphi}$ is integrable near the origin, then there exists $\varepsilon>0$ such that $|f|^2e^{-(1+\varepsilon)\varphi}$ is also integrable near the origin.
\end{theorem}

We need a few basic facts from local algebra.
Let $\mathcal O$ be the sheaf of holomorphic germs, and let $\mathcal O_0$ be its \emph{stalk} at the origin. We use only the following two facts.
\begin{enumerate}
\item \emph{Finite generation.} Every ideal $I$ of $\mathcal O_0$ is generated by finitely many germs $h_1,\ldots,h_r$. Thus every element of $I$ can be written near the origin as
\[
h=a_1h_1+\cdots+a_rh_r,\qquad a_1,\ldots,a_r\text{ holomorphic}.
\]
\item
\emph{Taylor matching.} For a holomorphic function $g$ near the origin, let $g_0$ be its germ there, and let
\[
T_Ng=\sum_{|\alpha|<N}\frac{\partial^\alpha g(0)}{\alpha!}z^\alpha
\]
be its Taylor polynomial of order less than $N$. If, for every $N$, there is $h_N\in I$ with $T_Ng=T_Nh_N$, then $g_0\in I$. This is Krull's intersection theorem; see \cite{StacksKrull}. Writing $\mathfrak m^N$ for the ideal of germs vanishing to order at least $N$, it says
\[
\bigcap_{N\ge1}(I+\mathfrak m^N)=I.
\]
\end{enumerate}

\begin{lemma}\label{ch3:lm:ideal-closed}
Let $I\subset\mathcal O_0$ be an ideal, and let $g_j,g$ be holomorphic on a common neighborhood $V$ of the origin, with $g_j\to g$ uniformly on compact subsets. If $(g_j)_0\in I$ for all $j$, then $g_0\in I$.
\end{lemma}
\begin{proof}
Fix $N$. Cauchy's formula gives convergence of each coefficient of $T_Ng_j$ to the corresponding coefficient of $T_Ng$. The Taylor polynomials arising from $I$ form the set
\[
T_N(I):=\{T_Nh:h\in I\}.
\]
This is a linear subspace of the finite-dimensional space of polynomials of degree less than $N$, so it is closed. Since $T_Ng_j\in T_N(I)$, passage to the limit gives $T_Ng\in T_N(I)$. Thus some $h_N\in I$ has the same Taylor coefficients as $g$ below order $N$. This holds for every $N$, so Taylor matching gives $g_0\in I$.
\end{proof}

For a plurisubharmonic function $\psi$, define its multiplier ideal at the origin by
\[
\mathcal I(\psi)_0:=\{g_0\in\mathcal O_0:|g|^2e^{-\psi}\text{ is integrable near the origin}\}.
\]
The inequality $|g+h|^2\le2|g|^2+2|h|^2$ and local boundedness of holomorphic functions show that this is indeed an ideal.

\begin{proof}[Proof of Theorem \ref{th:strong-openness}]
Shrinking the neighborhood and subtracting a constant, assume $\varphi<0$ on a ball $U\ni0$ and $\int_U|f|^2e^{-\varphi}<\infty$. Put
\[
\varphi_j=(1+1/j)\varphi.
\]
Since $2\varphi\le\varphi_j<0$ and $\varphi_j\to\varphi$ almost everywhere, Example \ref{ch3:ex:common-lower} shows that these weights satisfy Lemma \ref{ch3:lm:S_2}. Fix a small ball $0\in V\Subset U$. The lemma gives $f_j\in\mathcal O(V)$ such that
\[
\int_V|f_j|^2e^{-\varphi_j}<\infty,
\qquad \|f_j-f\|_{L^2(V)}\longrightarrow0.
\]
The mean-value inequality implies $f_j\to f$ uniformly on compact subsets of $V$.

Put $I_j=\mathcal I(\varphi_j)_0$. Since $\varphi<0$, the integrability condition weakens as $j$ increases, so
\[
I_1\subset I_2\subset\cdots.
\]
These ideals eventually stabilize. Their union $I=\bigcup_jI_j$ is an ideal, since any finite collection of elements belongs to one $I_j$. Choose generators $h_1,\ldots,h_r$ of $I$. Each belongs to some $I_j$, so a common $j_0$ contains all generators. Then $I_{j_0}$ contains all their holomorphic linear combinations, giving
\[
I_j=I_{j_0}=I\qquad(j\ge j_0).
\]
For $j\ge j_0$, all $(f_j)_0$ belong to the same ideal $I_{j_0}$, and $f_j\to f$ uniformly on compact subsets. Closedness of ideals gives $f_0\in I_{j_0}$. By definition, $|f|^2e^{-(1+1/j_0)\varphi}$ is integrable near the origin. Take $\varepsilon=1/j_0$.
\end{proof}

\begin{remark}
Besides Guan and Zhou's original proof, several other proofs are available; interested readers may consult the papers of Berndtsson or Xu Wang.
\end{remark}

We now show that Lemma \ref{ch3:lm:S_2} applies to plurisubharmonic functions converging in $L^1_{\rm loc}$.

\begin{proposition}\label{ch3:prop:DK-weight-convergence}
Let $\varphi_j,\varphi\in PSH(\Omega)$ with $\varphi_j\to\varphi$ in $L^1_{\rm loc}(\Omega)$. For every open $U\subset\subset\Omega$, there is a  closed pluripolar set $E\subset U$ such that
\[
e^{-\varphi_j}\longrightarrow e^{-\varphi}
\quad\text{in }L^1_{\rm loc}(U\setminus E).
\]
\end{proposition}

We first prove two elementary lemmas using one-variable Riesz decomposition and Jensen's inequality.

\begin{lemma}\label{ch3:lm:DK-elementary-exp}
There exists $C_n>0$, depending only on $n$, such that, for every
\[
\varphi\in PSH(B(0,2)),\qquad \varphi\le0,\qquad \varphi(0)\ge-\frac14,
\]
we have
\[
\int_{B(0,1/4)}e^{-\varphi}\,dV\le C_n.
\]
\end{lemma}

\begin{proof}
First let $n=1$. Prove the estimate for smooth convolutions of $\varphi$, then pass to the limit by Fatou's lemma. The approximations still satisfy $\varphi\le0$ and $\varphi(0)\ge-1/4$, so we may assume $\varphi$ is smooth. On the unit disc, use the Riesz decomposition
\[
\varphi(z)=h(z)-\int_{\mathbb{D}}G(z,\zeta)\,d\mu(\zeta),\qquad
G(z,\zeta)=\log\left|\frac{1-\bar\zeta z}{z-\zeta}\right|,
\]
where $\mu=(2\pi)^{-1}\Delta\varphi\,dA\ge0$ and $h$ is the harmonic extension of $\varphi|_{\partial\mathbb{D}}$, so $h\le0$, while $G\ge0$. At the origin,
\[
-h(0)+\int_{\mathbb{D}}\log\frac1{|\zeta|}\,d\mu(\zeta)
=-\varphi(0)\le\frac14.
\]
hence
\[
m:=\mu(B(0,1/2))\le\frac1{4\log2}<1.
\]
For $|z|<1/4$, Harnack's inequality gives $-h(z)\le\frac53[-h(0)]\le\frac5{12}$. For $|\zeta|\ge1/2$, the function $G(\cdot,\zeta)$ is positive harmonic on $B(0,1/2)$, so $G(z,\zeta)\le3G(0,\zeta)$, and this part of the integral is at most $3/4$. For $|\zeta|<1/2$, we have $G(z,\zeta)\le\log(2/|z-\zeta|)$. Thus
\[
-\varphi(z)\le\frac76+\int_{|\zeta|<1/2}\log\frac2{|z-\zeta|}\,d\mu(\zeta).
\]
If $m>0$, put $d\widehat\mu=m^{-1}d\mu|_{B(0,1/2)}$. This is a probability measure, so Jensen's inequality gives
\[
\begin{aligned}
\exp\left(\int_{|\zeta|<1/2}\log\frac2{|z-\zeta|}\,d\mu\right)
&\le\int_{|\zeta|<1/2}\left(\frac2{|z-\zeta|}\right)^m d\widehat\mu\\
&\le\int_{|\zeta|<1/2}\frac2{|z-\zeta|}\,d\widehat\mu.
\end{aligned}
\]
The last step uses $m<1$ and $|z-\zeta|<3/4$. Since
\[
\sup_{|\zeta|<1/2}\int_{|z|<1/4}\frac2{|z-\zeta|}\,dA(z)
\le\int_{|w|<3/4}\frac2{|w|}\,dA(w)=3\pi,
\]
we obtain $\int_{|z|<1/4}e^{-\varphi}\,dA\le3\pi e^{7/6}$. If $m=0$, the same bound follows directly from $-\varphi\le7/6$.

In $\mathbb{C}^n$, apply the one-dimensional result to $\lambda\mapsto\varphi(\lambda\xi)$ for each unit vector $\xi$. Since $\varphi(0)>-\infty$, every such restriction is subharmonic. Let $d\sigma$ be normalized area measure on the unit sphere and $s_{2n-1}$ its area. Polar coordinates and rotation invariance of $d\sigma$ give
\[
\begin{aligned}
\int_{|z|<1/4}e^{-\varphi(z)}\,dV(z)
&=s_{2n-1}\int_{S^{2n-1}}\int_0^{1/4}e^{-\varphi(t\xi)}t^{2n-1}\,dt\,d\sigma(\xi)\\
&\le\frac{s_{2n-1}(1/4)^{2n-2}}{2\pi}
\int_{S^{2n-1}}\int_{|\lambda|<1/4}e^{-\varphi(\lambda\xi)}\,dA(\lambda)\,d\sigma(\xi).
\end{aligned}
\]
The one-dimensional estimate bounds the right-hand side uniformly.
\end{proof}

\begin{lemma}\label{ch3:lm:DK-finite-point}
Let $\varphi_j,\varphi\in PSH(D)$ with $\varphi_j\to\varphi$ in $L^1_{\rm loc}(D)$. If $\varphi(a)>-\infty$, there exist a neighborhood $W$ of $a$ and a positive integer $j_0$ such that
\[
\sup_{j\ge j_0}\int_W e^{-2\varphi_j}<\infty.
\]
\end{lemma}

\begin{proof}
Put $\delta=1/32$. By upper semicontinuity, choose $r>0$ with $B(a,4r)\subset\subset D$ and $\varphi\le\varphi(a)+\delta$ there. For $z\in B(a,3r)$, the submean-value inequality gives
\[
\varphi_j(z)\le \varphi(a)+\delta+
\frac{\|\varphi_j-\varphi\|_{L^1(B(a,4r))}}{|B(0,r)|}.
\]
Thus, for sufficiently large $j$, $\varphi_j\le M:=\varphi(a)+2\delta$ on $B(a,3r)$.
The average of $\varphi$ on $B(a,r/8)$ is at least $\varphi(a)$. By $L^1$ convergence, for large $j$ choose $a_j\in B(a,r/8)$ with $\varphi_j(a_j)\ge\varphi(a)-2\delta$. Put
\[
\psi_j(z)=2\bigl(\varphi_j(a_j+rz)-M\bigr),\qquad |z|<2.
\]
Then $\psi_j\le0$ and $\psi_j(0)\ge-8\delta=-1/4$. Lemma \ref{ch3:lm:DK-elementary-exp} gives
\[
\int_{B(a_j,r/4)}e^{-2\varphi_j}\le C_nr^{2n}e^{-2M}.
\]
Since $B(a,r/8)\subset B(a_j,r/4)$, take the fixed neighborhood $W=B(a,r/8)$.
\end{proof}

\begin{proof}[Proof of Proposition \ref{ch3:prop:DK-weight-convergence}]
The submean-value inequality and $L^1_{\rm loc}$ convergence give a uniform upper bound for $\varphi_j,\varphi$ on $U$ and $\sup_j\int_U|\varphi_j|<\infty$. Subtracting a common constant, assume $\varphi_j,\varphi\le0$ on $U$. Put
\[
\psi=\varphi+\sum_{j=1}^{\infty}2^{-j}\varphi_j.
\]
The partial sums decrease, and
\[
\int_U|\psi|\le\int_U|\varphi|+
\sum_{j\ge1}2^{-j}\int_U|\varphi_j|<\infty.
\]
Thus $\psi\in PSH(U)$ and is finite almost everywhere. In particular, even if $U$ is disconnected, $\psi$ is not identically $-\infty$ on any component.

Let $V$ be the union of all open balls $B\subset\subset U$ satisfying
\[
\sup_j\int_Be^{-2\varphi_j}<\infty,
\]
and put $E=U\setminus V$. Then $E$ is  closed in $U$. We show that
$$
E\subset\{\psi=-\infty\}.
$$
If $\psi(a)>-\infty$, the nonpositive terms in the series imply that $\varphi(a)>-\infty$ and $\varphi_j(a)>-\infty$ for all $j$. Lemma \ref{ch3:lm:DK-finite-point} gives a uniform integral bound near $a$ for sufficiently large $j$. For each of the remaining finitely many $j$, apply the lemma to the constant sequence $\varphi_j,\varphi_j,\ldots$; then $e^{-2\varphi_j}$ is integrable near $a$. Taking a small ball inside the intersection of these neighborhoods gives a uniform bound for all $j$. Thus $a\in V$, so $E\subset\{\psi=-\infty\}$ and $E$ is closed and pluripolar.

For every compact $S\subset\subset U\setminus E$, a finite covering gives
\[
\sup_j\int_S e^{-2\varphi_j}<\infty.
\]
Since $\varphi_j\to\varphi$ in measure and $\varphi$ is finite almost everywhere, $e^{-\varphi_j}\to e^{-\varphi}$ in measure. The $L^2(S)$ bound gives uniform integrability, so Vitali's convergence theorem proves the result.
\end{proof}

Lemma \ref{ch3:lm:S_2} and Proposition \ref{ch3:prop:DK-weight-convergence} imply the following semicontinuity theorem:

\begin{theorem}[Demailly-Koll\'ar {\cite{DemaillyKollar2001}}]\label{ch3:th:DK-semicontinuity}
Let $\Omega\subset\mathbb{C}^n$ be a domain, $K\subset\Omega$ compact, and $\varphi_j,\varphi\in PSH(\Omega)$ with $\varphi_j\to\varphi$ in $L^1_{\rm loc}(\Omega)$. Define the complex singularity exponent
\[
c_K(\varphi)=\sup\left\{c\ge0:
e^{-2c\varphi}\text{ is integrable on a neighborhood of }K\right\}.
\]
Then
\[
c_K(\varphi)\le\liminf_{j\to\infty}c_K(\varphi_j).
\]
Moreover, for every $0<c<c_K(\varphi)$, there is a fixed neighborhood $K\subset W\Subset\Omega$ such that
\[
e^{-2c\varphi_j}\longrightarrow e^{-2c\varphi}
\quad\text{in }L^1(W).
\]
\end{theorem}

\begin{proof}
First take $K=\{a\}$ and write $c_a=c_{\{a\}}$. Fix $0<b<c_a(\varphi)$ and choose a small ball $U\subset\subset\Omega$ centered at $a$ with
\[
A:=\int_Ue^{-2b\varphi}\,dV<\infty.
\]
By Proposition \ref{ch3:prop:DK-weight-convergence} and Lemma \ref{ch3:lm:S_2}, for a ball $a\in V\subset\subset U$ there are $f_j\in\mathcal O(V)$ such that
\[
\|f_j-1\|_{L^2(V)}\longrightarrow0,\qquad
\limsup_{j\to\infty}\int_V|f_j|^2e^{-2b\varphi_j}\,dV\le A.
\]
Take a ball $a\in W\subset\subset V$. The mean-value inequality gives
\[
\sup_{\overline W}|f_j-1|
\le C_{W,V}\|f_j-1\|_{L^2(V)}\longrightarrow0.
\]
Thus $|f_j|\ge1/2$ on $W$ for large $j$, and hence
\[
\sup_{j\ge j_0}\int_W e^{-2b\varphi_j}\,dV
\le4\sup_{j\ge j_0}\int_V|f_j|^2e^{-2b\varphi_j}\,dV<\infty.
\]
This shows that $c_a(\varphi_j)\ge b$ for large $j$. Letting $b\uparrow c_a(\varphi)$ proves semicontinuity; if $c_a(\varphi)=0$, the inequality is immediate.

For $0<c<c_a(\varphi)$, choose $c<b<c_a(\varphi)$. The estimate gives
\[
\sup_{j\ge j_0}\int_W\bigl(e^{-2c\varphi_j}\bigr)^{b/c}\,dV<\infty,
\qquad b/c>1.
\]
Convergence in measure and Vitali's theorem imply $e^{-2c\varphi_j}\to e^{-2c\varphi}$ in $L^1(W)$.

For a general compact $K$, fix $0<b<c_K(\varphi)$, carry out this construction for each $a\in K$, and take a finite cover by balls $W_a$. For large $j$, $e^{-2b\varphi_j}$ has a uniform integral bound on their union, so $c_K(\varphi_j)\ge b$. Let $b\uparrow c_K(\varphi)$ for the first assertion; the case $c_K(\varphi)=0$ is immediate. Finally, fix $0<c<b<c_K(\varphi)$ and use the same finite cover and Vitali's theorem for the second assertion.
\end{proof}

\begin{theorem}\label{ch3:th:kernel-L1-stability}
Let $\Omega\subset\mathbb{C}^n$ be bounded and hyperconvex, and let $\varphi_j,\varphi\in PSH(\Omega)$ with
\[
\varphi_j\longrightarrow\varphi\quad\text{in }L^1_{\rm loc}(\Omega),
\qquad \varphi\not\equiv-\infty.
\]
Write $K_j=K_{\Omega,\varphi_j}$ and $K=K_{\Omega,\varphi}$. Then $K_j$ converges to $K$ uniformly on compact subsets.
\end{theorem}

\begin{proof}
The submean-value inequality gives local uniform upper bounds for $\{\varphi_j\}$. Extend the exceptional-set construction in Proposition \ref{ch3:prop:DK-weight-convergence} to all of $\Omega$. Put
\[
\psi=\varphi+\sum_{j=1}^{\infty}2^{-j}\varphi_j.
\]
On each ball $U\subset\subset\Omega$, choose a common upper bound $M_U$ for $\varphi,\varphi_j$. Then
\[
\psi=(\varphi-M_U)+\sum_{j=1}^{\infty}2^{-j}(\varphi_j-M_U)+2M_U.
\]
The terms on the right are nonpositive, the partial sums decrease, and $L^1(U)$ boundedness implies that their limit is not identically $-\infty$. Thus $\psi\in PSH(\Omega)$ and is finite almost everywhere. The argument also shows that $\psi(a)>-\infty$ implies finiteness of $\varphi(a)$ and every $\varphi_j(a)$.

Let $V$ be the union of all open balls $B\subset\subset\Omega$ satisfying
\[
\sup_j\int_B e^{-2\varphi_j}\,dV<\infty
\]
and put $E=\Omega\setminus V$. If $\psi(a)>-\infty$, Lemma \ref{ch3:lm:DK-finite-point} gives the uniform integral bound near $a$ for large $j$. For each of the finitely many initial terms, apply the lemma to a constant sequence and shrink the neighborhood; this gives a uniform bound for all $j$. Hence $a\in V$, so
\[
E\subset\{\psi=-\infty\}.
\]
Since $V$ is open, $E$ is a closed pluripolar subset of all of $\Omega$.

For every compact $S\subset\Omega\setminus E$, a finite covering gives
\[
\sup_j\int_S e^{-2\varphi_j}\,dV<\infty.
\]
Convergence in measure and Vitali's theorem yield
\[
e^{-\varphi_j}\longrightarrow e^{-\varphi}
\quad\text{in }L^1_{\rm loc}(\Omega\setminus E).
\]
All hypotheses of Theorem \ref{ch3:th:singular} hold, so $K_j(z,w)\to K(z,w)$ uniformly on compact subsets of $\Omega\times\Omega$.
\end{proof}


As an application of Theorem \ref{ch3:th:kernel-L1-stability}, we give a proof of the following classical completeness theorem using stability of weighted Bergman kernels.

\begin{theorem}[B\l ocki--Pflug, Herbort]\label{ch3:th:hyperconvex-completeness-stability}
Every bounded hyperconvex domain $\Omega\subset\mathbb C^n$ is Bergman complete.\footnote{Z. B\l ocki and P. Pflug, \emph{Hyperconvexity and Bergman completeness}, Nagoya Math. J. \textbf{151} (1998), 221--225; G. Herbort, \emph{The Bergman metric on hyperconvex domains}, Math. Z. \textbf{232} (1999), 183--196.}
\end{theorem}

\begin{proof}
Let $\{a_j\}\subset\Omega$ tend to the boundary, and let $g_j(z)=g_\Omega(z,a_j)$ be the negative pluricomplex Green function with pole at $a_j$, normalized by
\[
g_j(z)=\log|z-a_j|+O(1)\qquad(z\to a_j).
\]
By Lemma \ref{ch3:lm:green-poles-boundary} below, $g_j\to0$ in $L^1_{\rm loc}(\Omega)$. Put $\varphi_j=2ng_j$. Theorem \ref{ch3:th:kernel-L1-stability} gives, for every fixed $w\in\Omega$,
\begin{equation}\label{ch3:eq:complete-weighted-kernel-limit}
K_{\varphi_j}(w)\longrightarrow K(w),
\end{equation}
where $K=K_\Omega$. Consider the closed subspace
\[
H_j=\{f\in A^2(\Omega):f(a_j)=0\}.
\]
Since $g_j\le0$, we have $\|f\|_{A^2(\Omega)}\le\|f\|_{\varphi_j}$. Moreover, $e^{-\varphi_j}$ is comparable to $|z-a_j|^{-2n}$ near $a_j$, and this singularity is not integrable. Thus every holomorphic function with finite weighted norm must vanish at $a_j$. Consequently,
\[
A^2_{\varphi_j}(\Omega)\subset H_j,\qquad
\|f\|_{A^2(\Omega)}\le\|f\|_{\varphi_j}.
\]
By the reproducing property, $H_j$ is the orthogonal complement of the unit vector $K(\cdot,a_j)/\sqrt{K(a_j)}$. Removing this one-dimensional component gives the diagonal of its reproducing kernel:
\[
K_{H_j}(w)=K(w)-\frac{|K(w,a_j)|^2}{K(a_j)}.
\]
The extremal characterization gives $K_{H_j}(w)\ge K_{\varphi_j}(w)$. Together with \eqref{ch3:eq:complete-weighted-kernel-limit}, this yields
\begin{equation}\label{ch3:eq:complete-normalized-kernel-limit}
0\le\frac{|K(w,a_j)|^2}{K(a_j)}
\le K(w)-K_{\varphi_j}(w)\longrightarrow0.
\end{equation}
The linear span of the kernel functions $K(\cdot,w)$ is dense in $A^2(\Omega)$. Hence \eqref{ch3:eq:complete-normalized-kernel-limit} and Kobayashi's criterion \eqref{eq:Kobayashi} imply Bergman completeness.
\end{proof}

For completeness, we include the Green-function convergence argument of Carlehed--Cegrell--Wikstr\"om.\footnote{M. Carlehed, U. Cegrell and F. Wikstr\"om, \emph{Jensen measures, hyperconvexity and boundary behaviour of the pluricomplex Green function}, Ann. Polon. Math. \textbf{71} (1999), no. 1, 87--103; see Theorem 4.7 and equation (4.4), p. 100.}

\begin{lemma}\label{ch3:lm:green-poles-boundary}
Let $\Omega$ be bounded and hyperconvex. If $a_j\in\Omega$ tends to the boundary, then $g_j=g_\Omega(\cdot,a_j)\to0$ in $L^1_{\rm loc}(\Omega)$.
\end{lemma}

\begin{proof}
Choose a bounded negative continuous plurisubharmonic exhaustion $\rho$, and use the normalization $dd^c=(i/\pi)\partial\bar\partial$. For a fixed pole, the Green function has boundary values zero and satisfies
\[
(dd^cg_j)^n=\delta_{a_j},
\]
where $\delta_{a_j}$ is the unit point mass.\footnote{J.-P. Demailly, \emph{Mesures de Monge-Amp\`ere et mesures pluriharmoniques}, Math. Z. \textbf{194} (1987), no. 4, 519--564.}
For each $m$, put
\[
v_m=\sup_{j\ge m}g_j,\qquad u_m=v_m^*,
\]
where the star denotes upper semicontinuous regularization. Then $u_m$ is bounded, nonpositive and plurisubharmonic, with boundary values zero. Indeed, every tail contains two distinct poles; the maximum of the corresponding Green functions is bounded below on $\Omega$, and so is $u_m$. Also, $g_m\le u_m\le0$, which gives the boundary values.

Fix $m$. For every $j\ge m$, repeated integration by parts and $u_m\ge g_j$ give
\begin{equation}\label{ch3:eq:green-envelope-mass}
\begin{aligned}
\int_\Omega(-\rho)(dd^cu_m)^n
&=\int_\Omega(-u_m)dd^c\rho\wedge(dd^cu_m)^{n-1}\\
&\le\int_\Omega(-g_j)dd^c\rho\wedge(dd^cu_m)^{n-1}\\
&=\int_\Omega(-\rho)dd^cg_j\wedge(dd^cu_m)^{n-1}\\
&\le\cdots\le\int_\Omega(-\rho)(dd^cg_j)^n=-\rho(a_j).
\end{aligned}
\end{equation}
The integration by parts can first be performed with $g_{j,t}=\max\{g_j,-t\}$. For $t\ge\|u_m\|_\infty$, we still have $u_m\ge g_{j,t}$; then let $t\to\infty$. The Monge--Amp\`ere measures of these truncations converge to $\delta_{a_j}$. This explains the truncation underlying equation (4.4) of Carlehed--Cegrell--Wikstr\"om; the required continuity and comparison principles are those of Bedford--Taylor and Demailly.\footnote{E. Bedford and B. A. Taylor, \emph{A new capacity for plurisubharmonic functions}, Acta Math. \textbf{149} (1982), 1--40; see also Demailly's paper cited above.}

Keeping $m$ fixed, let $j\to\infty$. Since $\rho(a_j)\to0$, \eqref{ch3:eq:green-envelope-mass} implies $(dd^cu_m)^n=0$. The comparison principle and the zero boundary values give $u_m=0$. A supremum of plurisubharmonic functions agrees almost everywhere with its upper semicontinuous regularization. Taking the union of the countably many exceptional sets, we obtain
\begin{equation}\label{ch3:eq:green-envelope-limsup}
\limsup_{j\to\infty}g_j=\lim_{m\to\infty}v_m=0
\qquad\text{almost everywhere}.
\end{equation}
Choose $R>\operatorname{diam}\Omega$. The extremal definition of the Green function gives
\[
\log\frac{|z-a_j|}{R}\le g_j(z)\le0.
\]
Thus, on each compact subset of $\Omega$, the functions $g_j$ are uniformly bounded for sufficiently large $j$. By the $L^1_{\rm loc}$ compactness theorem for plurisubharmonic functions, every subsequence has a further subsequence converging in $L^1_{\rm loc}$; passing to another subsequence gives almost everywhere convergence. Repeating the preceding argument for this subsequence, \eqref{ch3:eq:green-envelope-limsup} shows that its limit must be zero. Therefore the whole sequence converges to zero in $L^1_{\rm loc}(\Omega)$.
\end{proof}

\section{H\"older Continuity}

A natural question is the following.

\begin{problem}\label{ch3:prob:weight-holder}
Let $\Omega$ be a bounded pseudoconvex domain in $\mathbb{C}^n$, and let $\varphi_t,\varphi\in PSH(\Omega)$. If $e^{\varphi_t}$ is H\"older continuous in $t$, is $K_{\varphi_t}(z,w)$ also H\"older continuous in $t$?
\end{problem}

Problem \ref{ch3:prob:weight-holder} has no affirmative answer in general. However, if the interior singularities satisfy a uniform higher integrability condition and the weights do not degenerate near the boundary, the following estimate holds.

\begin{theorem}\label{ch3:th:boundary-controlled-holder}
Let $\Omega\subset\mathbb{C}^n$ be a bounded domain, and let $\varphi_t\in PSH^-(\Omega)$, with $t$ ranging over a real interval. Suppose that there exist $U\subset\subset\Omega$, $p>1$, and constants $C_1,C_2,C_3>0$ such that
\[
 \sup_\Omega|e^{\varphi_t}-e^{\varphi_s}|\le C_1 |t-s|^\alpha,
 \qquad 0<\alpha\le1,
\]
\[
 \sup_t\int_U e^{-p\varphi_t}\,dV\le C_2,
 \qquad e^{\varphi_t}\ge C_3 \quad\hbox{on }\Omega\setminus U.
\]
Set $\theta=\min\{p-1,1\}$ and $K_t=K_{\Omega,\varphi_t}$. Then, for every $S\subset\subset\Omega$, there exists $C_S>0$ such that
\[
 \sup_{z,w\in S}|K_t(z,w)-K_s(z,w)|
 \le C_S|t-s|^{\alpha\theta}.
\]
In particular, when $p\ge2$, the kernel has the same H\"older exponent $\alpha$ as $e^{\varphi_t}$.
\end{theorem}

\begin{proof}
We first show that all the weighted spaces coincide with $A^2(\Omega)$. Since $\varphi_t\le0$, we clearly have $A^2_{\varphi_t}(\Omega)\subset A^2(\Omega)$. Conversely, for $f\in A^2(\Omega)$, the pointwise estimates on compact subsets and the hypotheses give
\[
 \begin{split}
 \int_\Omega |f|^2e^{-\varphi_t}\,dV
 &\le \sup_U|f|^2\int_Ue^{-\varphi_t}\,dV
       +C_3^{-1}\int_{\Omega\setminus U}|f|^2\,dV\\
 &\le(C_U C_2+C_3^{-1})\|f\|_{L^2(\Omega)}^2.
 \end{split}
\]
Here we have used $e^{-\varphi_t}\le e^{-p\varphi_t}$. Thus these norms are uniformly equivalent, and the reproducing property gives
\[
 K_t(z,w)-K_s(z,w)
 =\int_\Omega K_t(z,\zeta)K_s(\zeta,w)
       (e^{-\varphi_s}-e^{-\varphi_t})\,dV(\zeta).
\]

Set $\delta=\|e^{\varphi_t}-e^{\varphi_s}\|_\infty\le1$. If $\delta=0$, the conclusion is immediate. For $a,b>0$ and $0<\theta\le1$, we have the elementary inequality
\[
 |a^{-1}-b^{-1}|
 \le |a-b|^\theta\bigl(a^{-1-\theta}+b^{-1-\theta}\bigr).
\]
Indeed, we may assume that $a\le b$ and use $(b-a)^{1-\theta}a^\theta\le b$. Since $1+\theta\le p$, substituting $a=e^{\varphi_t}$ and $b=e^{\varphi_s}$ and integrating yield
\[
 \int_U|e^{-\varphi_t}-e^{-\varphi_s}|\,dV\le2C_2\delta^\theta.
\]
Moreover, since $\varphi_t\le0$, the Cauchy--Schwarz inequality for the kernel gives
\[
 |K_t(z,\zeta)|^2\le K_t(z)K_t(\zeta)
 \le K_\Omega(z)K_\Omega(\zeta).
\]
Thus, for $z,w\in S$ and $\zeta\in U$, the product of the two kernels in the reproducing formula is uniformly bounded. The integral over $U$ is therefore at most $C_{S,U}C_2\delta^\theta$.

On $\Omega\setminus U$,
\[
 |e^{-\varphi_t}-e^{-\varphi_s}|\le C_3^{-2}\delta.
\]
Applying the Cauchy--Schwarz inequality again and using
\[
 \int_\Omega|K_t(z,\zeta)|^2\,dV(\zeta)
 \le\int_\Omega|K_t(z,\zeta)|^2e^{-\varphi_t(\zeta)}\,dV(\zeta)=K_t(z),
\]
show that this part of the integral is at most $C_SC_3^{-2}\delta$. Combining the two parts and using $\delta\le\delta^\theta$, we obtain
\[
 |K_t(z,w)-K_s(z,w)|\le C_S\delta^\theta
 \le C_S C_1^\theta|t-s|^{\alpha\theta}.
\]
\end{proof}

The conclusion fails if the uniform lower bound for $e^{\varphi_t}$ near the boundary is omitted.

\begin{example}\label{ch3:ex:boundary-weight-slow}
On the unit disc $\mathbb{D}$, set
\[
 F(z)=\exp\left(-\frac{1+z}{1-z}\right),\qquad
 \varphi_t(z)=\log\frac{|F(z)|^2+t}{2},\quad 0\le t<\frac12.
\]
Then $\varphi_t\in PSH^-(\mathbb{D})$, and $e^{\varphi_t}$ is uniformly Lipschitz continuous in $t$. For every $U\subset\subset\mathbb{D}$ and every $p>1$, we have
\[
 \sup_{0\le t<1/2}\int_U e^{-p\varphi_t}\,dA<\infty.
\]
However, writing $K_t=K_{\mathbb{D},\varphi_t}$, there exists $c>0$ such that
\[
 K_t(0)-K_0(0)\ge\frac{c}{(1+\log(1/t))^2},\qquad 0<t<\frac12.
\]
Consequently, $K_t(0)$ satisfies no H\"older estimate with a positive exponent at $t=0$.
\end{example}

\begin{proof}
The function $F$ is holomorphic and nowhere vanishing, and $|F|<1$. Hence $\varphi_t$ is negative and subharmonic, and
\[
 \sup_{\mathbb{D}}|e^{\varphi_t}-e^{\varphi_s}|=\frac{|t-s|}{2}.
\]
Since $\varphi_0$ is harmonic and $\varphi_t\ge\varphi_0$, the integrals over compact subsets above are uniformly bounded.

The identity
\[
 \|Fh\|_{\varphi_0}^2=2\int_{\mathbb{D}}|h|^2\,dA
\]
shows that $A^2_{\varphi_0}(\mathbb{D})=F A^2(\mathbb{D})$, and hence
\[
 K_0(0)=\frac{|F(0)|^2}{2}K_{\mathbb{D}}(0)=\frac{e^{-2}}{2\pi}.
\]
Using $F$ as a test function for $K_t(0)$ and writing
\[
 I_t=\int_{\mathbb{D}}\frac{t}{|F|^2+t}\,dA,
\]
we obtain from the extremal property that
\[
 K_t(0)\ge\frac{e^{-2}}{2(\pi-I_t)},\qquad
 K_t(0)-K_0(0)\ge\frac{e^{-2}}{2\pi^2}I_t.
\]
Set $T=\frac12\log(1/t)$. On the set
\[
 E_T=\left\{z\in\mathbb{D}:\operatorname{Re}\frac{1+z}{1-z}>T\right\}
     =\left\{z:\left|z-\frac{T}{T+1}\right|<\frac1{T+1}\right\}
\]
we have $|F|^2<t$. This set is a small disc internally tangent to the unit circle at $1$, so
\[
 I_t\ge\frac12|E_T|=\frac{\pi}{2(T+1)^2},\qquad
 K_t(0)-K_0(0)\ge\frac{e^{-2}}{4\pi(T+1)^2}.
\]
Since $(1+\log(1/t))^{-2}$ decays more slowly than any $t^\beta$ with $\beta>0$, the conclusion follows.
\end{proof}

\begin{example}\label{ch3:ex:regularization-rates}
On the unit ball $\mathbb{B}\subset\mathbb{C}^n$ with $n\ge2$, set
\[
 \varphi_t(z)= (n-1)\log\frac{|z|^2+t}{2},\quad 0\le t\le \frac12; 
 \qquad K_t=K_{\mathbb{B},\varphi_t}.
\]
The weight is smooth for $t>0$ and has a logarithmic singularity at the origin when $t=0$.
Then
\[
 K_t(0)-K_0(0)\asymp t\log 1/t \qquad(t\downarrow0).
\]
In this example, $e^{\varphi_t(z)}$ is smooth in $(z,t)$, but $K_t$ is not Lipschitz continuous in $t$.
\end{example}

\begin{proof}
Write $\mu_t=2^{n-1}/(|z|^2+t)^{n-1}$ and
$A_t=\int_{\mathbb{B}}\mu_t$. Integration in polar coordinates gives
\[
 A_t=\frac{2^{n-1}\pi^n}{(n-1)!}
             \int_0^1 \frac{s^{n-1}}{(s+t)^{n-1}}\,ds.
\]
Note that $A_0<\infty$ and
$A_t\to A_0>0$. To estimate the difference, split the integral at $s=t$:
\[
 A_0-A_t=\frac{2^{n-1}\pi^n}{(n-1)!}
       \int_0^1s^{n-1}\bigl[s^{1-n}-(s+t)^{1-n}\bigr]\,ds.
\]
On $0<s<t$, the expression in square brackets lies between
$(1-2^{1-n})s^{1-n}$ and $s^{1-n}$;
on $t\le s\le1$, the one-variable mean value theorem shows that it is comparable to $t s^{-n}$.
Thus
\[
 A_0-A_t\asymp t+t\int_t^1s^{-1}\,ds
             \asymp t\log 1/t.
\]

By the rotational invariance of the radial weight, constant functions are orthogonal to all
holomorphic functions vanishing at the origin. Thus, under the constraint $f(0)=1$, the constant function $1$ has the smallest norm, and hence
\[
 K_t(0)=\frac1{A_t},\qquad
 K_t(0)-K_0(0)=\frac{A_0-A_t}{A_0A_t}\asymp t\log 1/t.
\]
\end{proof}

\section{Higher Regularity}
Let $\Omega$ be a bounded domain in $\mathbb{C}^n$, and let $\varphi(z,t)\in C^2(\overline{\Omega}\times[-1,1])$ be real-valued. Set
$
\varphi_t(z):=\varphi(z,t)$, 
$K_t(z,w)=K_{\Omega,\varphi_t}(z,w)$ and write $A^2_t=A^2_{\varphi_t}(\Omega)$.

\begin{theorem}\label{ch3:th:C^2}
The kernel $K_t(z,w)$ is $C^2$ in $t$, and
$$
\frac{\partial^j K_t(\cdot,w)}{\partial t^j} \in A^2_t,\ \ \ \forall\, j\le 2.
$$
\end{theorem}

\begin{proof}
We first verify Lipschitz continuity. Let $S$ be a compact subset of $\Omega$, with $z,w\in S$ and $t_0\in(-1,1)$. Since all the norms $\|\cdot\|_t$ are equivalent, the reproducing property gives
\begin{eqnarray}\label{ch3:eq:high_1}
 K_t(z,w) & = & \int_\Omega K_t(\cdot,w) \overline{K_{t_0}(\cdot,z)} e^{-\varphi_{t_0}}\nonumber\\
 &  = & \int_\Omega K_t(\cdot,w) \overline{K_{t_0}(\cdot,z)} e^{-\varphi_{t}}\nonumber\\
 && + \int_\Omega K_t(\cdot,w) \overline{K_{t_0}(\cdot,z)} (e^{-\varphi_{t_0}}-e^{-\varphi_{t}})\nonumber\\
 & =: & K_{t_0}(z,w)+I,
  \end{eqnarray}
Moreover,
\begin{eqnarray*}
|I|  &  \le &  C|t-t_0| \, \|K_t(\cdot,w)\|_t\, \|K_{t_0}(\cdot,z)\|_{t_0}\\
&  \le  & C|t-t_0|\, K_t(w)^{1/2}\, K_{t_0}(z)^{1/2}
 \le  C|t-t_0|,
\end{eqnarray*}
where the constant $C$ also depends on the corresponding uniform bounds for the weights. Hence $K_t(z,w)$ is Lipschitz continuous in $t$.

Next, we verify $C^1$ regularity. Note that
\begin{eqnarray*}
 I & = & \int_\Omega K_{t_0}(\cdot,w) \overline{K_{t_0}(\cdot,z)} (e^{-\varphi_{t_0}}-e^{-\varphi_t})\\
 & &+ \int_\Omega (K_t(\cdot,w)-K_{t_0}(\cdot,w)) \overline{K_{t_0}(\cdot,z)}(e^{-\varphi_{t_0}}-e^{-\varphi_t})\\
 & =: & II+III.
\end{eqnarray*} 
Since
\begin{eqnarray*}
|III| & \le & C |t-t_0|\, \int_\Omega |K_t(\cdot,w)-K_{t_0}(\cdot,w)|\,|K_{t_0}(\cdot,z)| e^{-\varphi_{t_0}}\\
& \le & C|t-t_0| \, \|K_t(\cdot,w)-K_{t_0}(\cdot,w)\|_{t_0} \, \|K_{t_0}(\cdot,z)\|_{t_0}
 \end{eqnarray*}
and
 \begin{eqnarray}\label{ch3:eq:high_2}
  && \int_\Omega |K_t(\cdot,w)-K_{t_0}(\cdot,w)|^2 e^{-\varphi_{t_0}}\nonumber\\
  & = & \int_\Omega |K_t(\cdot,w)|^2 e^{-\varphi_{t_0}} + \int_\Omega |K_{t_0}(\cdot,w)|^2 e^{-\varphi_{t_0}}\nonumber\\
   && -2{\rm Re\,} \int_\Omega K_t(\cdot,w)\overline{K_{t_0}(\cdot,w)} e^{-\varphi_{t_0}}\nonumber\\
   & = & K_t(w)(1+O(|t-t_0|)) +K_{t_0}(w) - 2 K_t(w) ,
 \end{eqnarray}
we have
 $$
 | III | \le C |t-t_0|^{3/2}.
 $$
The dominated convergence theorem implies that
 \begin{eqnarray*}
 \lim_{t\rightarrow t_0} \frac{K_t(z,w)-K_{t_0}(z,w)}{t-t_0} & = & \lim_{t\rightarrow t_0} \int_\Omega K_{t_0}(\cdot,w)\overline{K_{t_0}(\cdot,z)} \, \frac{e^{-\varphi_{t_0}}-e^{-\varphi_t}}{t-t_0}\\
 & = & \int_\Omega K_{t_0}(\cdot,w) \overline{K_{t_0}(\cdot,z)} \, \left.\frac{\partial \varphi_t}{\partial t}\right|_{t_0}  e^{-\varphi_{t_0}}.
 \end{eqnarray*}
Thus $K_t(z,w)$ is differentiable in $t$, with
 \begin{eqnarray}\label{ch3:eq:high_3}
 \frac{\partial K_t (z,w)}{\partial t}  &  = & \int_\Omega K_t(\cdot,w) \overline{K_t(\cdot,z)} \, \frac{\partial \varphi_t}{\partial t} \, e^{-\varphi_t}\nonumber\\
 & = & P_t\left(K_t(\cdot,w)\frac{\partial \varphi_t}{\partial t}\right)(z) \in A^2_t, 
 \end{eqnarray}
where $P_t$ is the Bergman projection corresponding to $A^2_t$.
 
Finally, we verify $C^2$ regularity. For a function $f_t$ depending on the parameter $t$, define
$$
\Delta f_t:=f_t-f_{t_0},\ \ \ \Delta t:=t-t_0.
$$
By (\ref{ch3:eq:high_3}),
 \begin{eqnarray*}
 \Delta (\partial K_t(z,w)/\partial t)  & = &  \int_\Omega K_t(\cdot,w)\overline{K_t(\cdot,z)}\,\frac{\partial \varphi_t}{\partial t}\, e^{-\varphi_t}\\
 && - \int_\Omega K_{t_0}(\cdot,w)\overline{K_{t_0}(\cdot,z)}\, \left.\frac{\partial \varphi_t}{\partial t}\right|_{t_0}\, e^{-\varphi_{t_0}}.
 \end{eqnarray*}
If we set $g_t:=\frac{\partial\varphi_t}{\partial t}\,e^{-\varphi_t}$, then
 \begin{eqnarray*}
&&  K_t(\zeta,w)\overline{K_t(\zeta,z)} \, g_t(\zeta) \\
 & = & \left(K_{t_0}(\zeta,w)+\Delta K_t(\zeta,w)\right) \overline{(K_{t_0}(\zeta,z)+\Delta K_t(\zeta,z))}\,(g_{t_0}+\Delta g_t)\\
 & = & \left(K_{t_0}(\zeta,w)\overline{K_{t_0}(\zeta,z)}+ K_{t_0}(\zeta,w)\overline{\Delta K_t(\zeta,z)}+\Delta K_t(\zeta,w) 
 \overline{K_{t_0}(\zeta,z)}\right) g_{t_0}\\
 && + K_{t_0}(\zeta,w) \overline{K_{t_0}(\zeta,z)}\, \Delta g_t + R(\zeta,z,w,t),
 \end{eqnarray*}
where $R(\zeta,z,w,t)$ is the remainder. Thus
 \begin{eqnarray*}
 \Delta (\partial K_t(z,w)/\partial t) & = & \int_\Omega K_{t_0}(\cdot,w)\overline{\Delta K_t(\cdot,z)}\, g_{t_0} + \int_\Omega \Delta K_t(\cdot,w) \overline{K_{t_0}(\cdot,z)}\, g_{t_0}\\
 & & + \int_\Omega K_{t_0}(\cdot,w)\overline{K_{t_0}(\cdot,z)}\, \Delta g_t + \int_\Omega R(\cdot,z,w,t).
 \end{eqnarray*}  
By (\ref{ch3:eq:high_1}),
 $$
 \Delta K_t(\zeta,w) = \int_\Omega K_t(\cdot,w) \overline{K_{t_0}(\cdot,\zeta)} (-\Delta e^{-\varphi_t}), 
 $$
hence
 \begin{eqnarray*}
  \frac{\Delta K_t(\zeta,w)}{\Delta t} - \left.\frac{\partial K_t(\zeta,w)}{\partial t}\right|_{t_0}
  & = & \int_\Omega (K_{t_0}(\cdot,w)-K_{t}(\cdot,w)) \overline{K_{t_0}(\cdot,\zeta)} \, \frac{\Delta e^{-\varphi_t}}{\Delta t}\\
  && + \int_\Omega K_{t_0}(\cdot,w) \overline{K_{t_0}(\cdot,\zeta)}\left(\left.\frac{\partial e^{-\varphi_t}}{\partial t}\right|_{t_0}
  -\frac{\Delta e^{-\varphi_t}}{\Delta t}\right).
 \end{eqnarray*}
Since the norm of a Bergman projection is always at most $1$,
 \begin{eqnarray*}
  \left \| \frac{\Delta K_t(\cdot,w)}{\Delta t} - \left.\frac{\partial K_t(\cdot,w)}{\partial t}\right|_{t_0} \right\|_{t_0}
 &  \le & \left\| ( K_t(\cdot,w) -  K_{t_0}(\cdot,w)) \frac{\Delta e^{-\varphi_t}}{\Delta t} \,e^{\varphi_{t_0}} \right\|_{t_0}\\
  && +  \left\|K_{t_0}(\cdot,w)\left(\left.\frac{\partial e^{-\varphi_t}}{\partial t}\right|_{t_0}- \frac{\Delta e^{-\varphi_t}}{\Delta t}
  \right) e^{\varphi_{t_0}}\right\|_{t_0}\\
  & \rightarrow & 0\ \ \ \ \ \ \  (t\rightarrow t_0),
 \end{eqnarray*}
Combining this inequality with (\ref{ch3:eq:high_2}), we obtain
 $$
  \int_\Omega \frac{\Delta K_t(\cdot,w)}{\Delta t} \overline{K_{t_0}(\cdot,z)}\, g_{t_0}
  \rightarrow  \int_\Omega \left.\frac{\partial K_t(\cdot,w)}{\partial t}\right|_{t_0} \overline{K_{t_0}(\cdot,z)} \,g_{t_0}.
 $$
Similarly,
  $$
 \int_\Omega K_{t_0}(\cdot,w) \overline{\frac{\Delta K_t(\cdot,z)}{\Delta t}} \,g_{t_0}
 \rightarrow  \int_\Omega K_{t_0}(\cdot,w) \overline{\left.\frac{\partial K_t(\cdot,z)}{\partial t}\right|_{t_0}} \,g_{t_0}.
 $$
On the other hand, it is easy to show that
 $$
 \int_\Omega K_{t_0}(\cdot,w) \overline{K_{t_0}(\cdot,z)} \, \frac{\Delta g_t}{\Delta t} \rightarrow  \int_\Omega K_{t_0}(\cdot,w) \overline{K_{t_0}(\cdot,z)} \, \left.\frac{\partial g_t}{\partial t}\right|_{t_0}
 $$
and $\int_\Omega R(\cdot,z,w,t)/\Delta t\rightarrow0$ as $t\rightarrow t_0$. Thus $K_t(z,w)$ is $C^2$ in $t$, and
 \begin{eqnarray*}
 \frac{\partial^2 K_t(z,w)}{\partial t^2} & = & \int_\Omega K_t(\cdot,w) \overline{\frac{\partial K_t(\cdot,z)}{\partial t}}\,\frac{\partial \varphi_t}{\partial t} \, e^{-\varphi_t} \\
 & & + \int_\Omega \frac{\partial K_t(\cdot,w)}{\partial t}  \overline{K_t(\cdot,z)}\, \frac{\partial \varphi_t}{\partial t} \, e^{-\varphi_t} \\
 && + \int_\Omega K_t(\cdot,w) \overline{K_t(\cdot,z)} \frac{\partial}{\partial t} \left(\frac{\partial \varphi_t}{\partial t}\,
 e^{-\varphi_t}\right) \\
 & =: & I_1(z,w)+I_2(z,w)+I_3(z,w).
 \end{eqnarray*}
For fixed $w\in\Omega$, we have
 $$
 I_2(z,w)=P_t\left( \frac{\partial K_t(\cdot,w)}{\partial t} \, \frac{\partial \varphi_t}{\partial t} \right) (z)\in A^2_t
 $$
 $$
 I_3(z,w)=P_t \left(K_t(\cdot,w)e^{\varphi_t}\frac{\partial}{\partial t} \left(\frac{\partial \varphi_t}{\partial t}\,
 e^{-\varphi_t}\right)\right)(z)\in A^2_t.
 $$
To prove that $I_1(\cdot,w)\in A^2_t$, first choose an increasing sequence of open sets $\Omega_k\subset\subset\Omega$
such that $\Omega_k\uparrow\Omega$. Let $\chi_k$ be the characteristic function of $\Omega_k$.
Set
  $$
  I_{1,k}(z,w):= \int_\Omega \chi_k\, K_t(\cdot,w) \overline{\frac{\partial K_t(\cdot,z)}{\partial t}}\,\frac{\partial \varphi_t}{\partial t}\,e^{-\varphi_t}.
  $$
Clearly, $I_{1,k}(z,w)\rightarrow I_1(z,w)$ as $k\rightarrow\infty$. On the other hand, since
  \begin{eqnarray*}
  	&& \iint_{\zeta,\xi\in \Omega}\chi_k(\zeta)\, |K_t(\zeta,w)K_t(\xi,z) K_t(\xi,\zeta)|\left| \frac{\partial \varphi_t(\zeta)}{\partial t}\, \frac{\partial \varphi_t(\xi)}{\partial t}\right| e^{-\varphi_t(\zeta)-\varphi_t(\xi)}\\
  	&\le& C_k \left( \iint_{\zeta,\xi\in \Omega} \chi_k(\zeta)|K_t(\xi,z)|^2 e^{-\varphi_t(\xi)}\right)^{1/2}
  	\left( \iint_{\zeta,\xi\in \Omega} \chi_k(\zeta)|K_t(\xi,\zeta)|^2 e^{-\varphi_t(\xi)}\right)^{1/2}\\
  	& < & \infty,
    	\end{eqnarray*}
Fubini's theorem implies that
 \begin{eqnarray*}
 I_{1,k}(z,w) & = & \int_{\zeta\in \Omega} \chi_k(\zeta)\, K_t(\zeta,w) \frac{\partial \varphi_t(\zeta)}{\partial t} e^{-\varphi_t(\zeta)}
 \int_{\xi\in \Omega} \overline{K_t(\xi,z)} K_t(\xi,\zeta) \frac{\partial \varphi_t(\xi)}{\partial t} e^{-\varphi_t(\xi)}\\
 & = &  \int_{\xi\in \Omega} \overline{K_t(\xi,z)}  \frac{\partial \varphi_t(\xi)}{\partial t} e^{-\varphi_t(\xi)}
  \int_{\zeta\in \Omega}\chi_k(\zeta)\, K_t(\zeta,w) \overline{K_t(\zeta,\xi)}\frac{\partial \varphi_t(\zeta)}{\partial t} e^{-\varphi_t(\zeta)}\\
  & = & P_t\left( \frac{\partial \varphi_t}{\partial t}\, P_t\left(\chi_k\,K_t(\cdot,w)\frac{\partial\varphi_t}{\partial t} \right)\right)(z)
  \in A^2_t.
 \end{eqnarray*}
Note that
\begin{eqnarray*}
 && \left\|I_{1,k}(\cdot,w)-P_t\left( \frac{\partial \varphi_t}{\partial t}\, P_t\left(K_t(\cdot,w)\frac{\partial\varphi_t}{\partial t} \right)\right)\right\|_t\\
 & \le & C \left\|(\chi_k-1) K_t(\cdot,w)\frac{\partial\varphi_t}{\partial t}\right\|_t \rightarrow 0\ \ \ (k\rightarrow \infty),
\end{eqnarray*}
Thus
$$
I_1(z,w)=P_t\left( \frac{\partial \varphi_t}{\partial t}\, P_t\left(K_t(\cdot,w)\frac{\partial\varphi_t}{\partial t} \right)\right)(z)
\in A^2_t,
$$
which yields $\frac{\partial^2 K_t(\cdot,w)}{\partial t^2}\in A^2_t$.
\end{proof}

\begin{remark}
The argument above is entirely elementary and also applies to $C^k$ regularity for $2\le k\le\infty$; see Marinescu--Xing \cite{MarinescuXing2026}.
\end{remark}

\begin{remark}
On a fixed bounded domain $\Omega$, take a real-valued function $\psi\in L^\infty(\Omega)$ and set $\varphi_t=t\psi$. If $P_\Omega$ denotes the unweighted Bergman projection, then the operator on $A^2(\Omega)$ defined by
\[
T_tf=P_\Omega(e^{-t\psi}f)
\]
is positive and invertible, and the reproducing property gives
\[
K_{\Omega,\varphi_t}(\cdot,w)=T_t^{-1}K_\Omega(\cdot,w).
\]
Local analytic inversion of operators shows that the kernel is real analytic in $t$. This approach also explains the existence of higher parameter derivatives; for a more general theory of bounded perturbations of weights, see Pasternak--Winiarski \cite{PasternakWiniarski1990}.
\end{remark}

\begin{remark}
Let $\delta_\Omega$ denote the Euclidean distance to the boundary. Marinescu--Xing \cite{MarinescuXing2026} proved that, if $\Omega$ is a bounded pseudoconvex domain with $C^2$ boundary, then the Bergman kernel associated with the measure $\delta_\Omega^t\,dV$, $t>-1$, is real analytic in $t$ and satisfies factorial estimates for its parameter derivatives on compact subsets. In this case, $\varphi_t=-t\log\delta_\Omega$.
\end{remark}

\begin{remark}
For smooth families of strongly pseudoconvex domains, Hadamard-type formulas express the first variation of the kernel as an integral arising from the motion of the boundary, while the second variation also involves curvature. Xu Wang \cite{WangVariation2014} studied such variation formulas and obtained a kernel criterion for triviality in the case of planar holomorphic motions. Parameter dependence can therefore also be used to study the geometry of deformations of domains.
\end{remark}

\section{Plurisubharmonicity}
Let $\Omega\subset\mathbb{C}^n_z\times\mathbb{C}^k_t$ and $\varphi\in PSH(\Omega)$. For a given $t$, set $\Omega_t=\{z:(z,t)\in\Omega\}$ and $A^2_t=A^2_{\varphi_t}(\Omega_t)$. Denote the reproducing kernel of $A^2_t$ by $K_t(z,w)$.

\begin{theorem}[Berndtsson {\cite{Berndtsson2006}}]\label{ch3:th:Berndtsson}
If $\Omega$ is pseudoconvex, then $\log K_t(z)$ is plurisubharmonic on $\Omega$ or identically $-\infty$.
\end{theorem}

\begin{remark}
When $n=1$ and $\varphi=0$, this theorem was first proved by Maitani-Yamaguchi \cite{MaitaniYamaguchi}.
\end{remark}

\begin{lemma}\label{ch3:lm:B_1}
Let $\Omega=D\times\mathbb{D}$, where $D$ is a bounded pseudoconvex domain in $\mathbb{C}^n$, and let $\varphi\in C^\infty(\overline{D}\times\mathbb{D})$ be strictly plurisubharmonic. Then $K_t(z)$ is strictly subharmonic in $t$, and $\log K_t(z)$ is subharmonic in $t$.
\end{lemma}

\begin{proof}
By Theorem \ref{ch3:th:C^2}, $K_t(\zeta,z)$ is $C^2$ in $t$, and
$$
\partial K_t(\cdot,z)/\partial t,\,\partial K_t(\cdot,z)/\partial \bar{t}\in A^2_t.
$$
The proof of Theorem \ref{ch3:th:C^2} also shows that
\begin{equation}\label{ch3:eq:B_0}
\frac{\partial K_t(\zeta,z)}{\partial \bar{t}}=\int_D K_t(\cdot,z)\overline{K_t(\cdot,\zeta)}\,\frac{\partial \varphi_t}{\partial \bar{t}}\,e^{-\varphi_t},
\end{equation}
and
 \begin{eqnarray*}
	\frac{\partial^2 K_t(z)}{\partial t\partial \bar{t}} & = & \int_D K_t(\cdot,z) \overline{\frac{\partial K_t(\cdot,z)}{\partial \bar{t}}}\,\frac{\partial \varphi_t}{\partial \bar{t}} \, e^{-\varphi_t} \\
	& & + \int_D \frac{\partial K_t(\cdot,z)}{\partial t}  \overline{K_t(\cdot,z)}\, \frac{\partial \varphi_t}{\partial \bar{t}} \, e^{-\varphi_t} \\
	&& + \int_D |K_t(\cdot,z)|^2  \left(\frac{\partial^2 \varphi_t}{\partial t\partial \bar{t}}-\left|\frac{\partial \varphi_t}{\partial t}\right|^2\right)
	e^{-\varphi_t}. 
\end{eqnarray*}
The reproducing property gives, for every $f\in A^2_t$,
$$
	f(z)  =  \int_D f\,\overline{K_t(\cdot,z)}\,e^{-\varphi_t},\ \ \ 
\overline{f(z)}  =  \int_D \bar{f}\,{K_t(\cdot,z)}\,e^{-\varphi_t}.
$$
Fixing $f$ and differentiating this identity with respect to $\bar t$, we obtain
\begin{equation}\label{ch3:eq:B_0'}
  \int_D f\,\overline{\partial^\varphi_t K_t(\cdot,z)}\,e^{-\varphi_t}  = 0,\ \ \ 
  \int_D \bar{f}\,{\partial^\varphi_{\bar{t}}K_t(\cdot,z)}\,e^{-\varphi_t}  =  0,
\end{equation}
where $\partial^\varphi_t={\partial}/{\partial t}-{\partial\varphi_t}/{\partial t}\cdot$ and $\partial^\varphi_{\bar t}={\partial}/{\partial\bar t}-{\partial\varphi_t}/{\partial\bar t}\cdot$.
Substituting $f={\partial K_t(\cdot,z)}/{\partial t}$ and $f={\partial K_t(\cdot,z)}/{\partial\bar t}$ into the preceding two identities gives
\begin{equation}\label{ch3:eq:B_1}
\int_D \left|\frac{\partial K_t(\cdot,z)}{\partial t}\right|^2 e^{-\varphi_t}  =   \int_D \frac{\partial K_t(\cdot,z)}{\partial t}  \overline{K_t(\cdot,z)}\, \frac{\partial \varphi_t}{\partial \bar{t}} \, e^{-\varphi_t},
\end{equation}
\begin{equation}\label{ch3:eq:B_2}
\int_D \left|\frac{\partial K_t(\cdot,z)}{\partial \bar{t}}\right|^2 e^{-\varphi_t}  =  \int_D K_t(\cdot,z) \overline{\frac{\partial K_t(\cdot,z)}{\partial \bar{t}}}\,\frac{\partial \varphi_t}{\partial \bar{t}} \, e^{-\varphi_t}.
\end{equation}
Thus
 \begin{eqnarray}\label{ch3:eq:B_3}
	\frac{\partial^2 K_t(z)}{\partial t\partial \bar{t}} & = & \int_D \left|\frac{\partial K_t(\cdot,z)}{\partial t}\right|^2 e^{-\varphi_t}      + \int_D \left|\frac{\partial K_t(\cdot,z)}{\partial \bar{t}}\right|^2 e^{-\varphi_t}\\
		&& + \int_D |K_t(\cdot,z)|^2  \left(\frac{\partial^2 \varphi_t}{\partial t\partial \bar{t}}-\left|\frac{\partial \varphi_t}{\partial t}\right|^2\right)
	e^{-\varphi_t}. \nonumber
\end{eqnarray}
Note that
\begin{eqnarray}\label{ch3:eq:B_3'}
\int_D |\partial^\varphi_t K_t(\cdot,z)|^2 e^{-\varphi_t}
& = & \int_D \left|\frac{\partial K_t(\cdot,z)}{\partial t}\right|^2 e^{-\varphi_t} + \int_D |K_t(\cdot,z)|^2\left|\frac{\partial \varphi_t}{\partial t}\right|^2 e^{-\varphi_t}\nonumber\\
&& -2{\rm Re\,} \int_D \frac{\partial K_t(\cdot,z)}{\partial t}  \overline{K_t(\cdot,z)}\, \frac{\partial \varphi_t}{\partial \bar{t}} \, e^{-\varphi_t}\nonumber\\
& = & \int_D |K_t(\cdot,z)|^2\left|\frac{\partial \varphi_t}{\partial t}\right|^2 e^{-\varphi_t} -  \int_D \left|\frac{\partial K_t(\cdot,z)}{\partial t}\right|^2 e^{-\varphi_t},
\end{eqnarray}
where the second equality follows from \eqref{ch3:eq:B_1}.
Combining \eqref{ch3:eq:B_3'} with (\ref{ch3:eq:B_3}), we obtain the following second variation formula:
\begin{eqnarray}\label{ch3:eq:SecondVariation}
\frac{\partial^2 K_t(z)}{\partial t\partial \bar{t}}  & = &  \int_D \left|\frac{\partial K_t(\cdot,z)}{\partial \bar{t}}\right|^2 e^{-\varphi_t}
		 + \int_D |K_t(\cdot,z)|^2 \,\frac{\partial^2 \varphi_t}{\partial t\partial \bar{t}}\,e^{-\varphi_t}\\
		&& -\int_D |\partial^\varphi_t K_t(\cdot,z)|^2 e^{-\varphi_t}. \nonumber
	\end{eqnarray}
Set $u_t=\partial^\varphi_t K_t(\cdot,z)$. Note that $u_t$ is the solution of
$$
\bar{\partial} u = \bar{\partial}(\partial^\varphi_t K_t(\cdot,z)) = - K_t(\cdot,z) \bar{\partial} (\partial \varphi_t/\partial t)
$$
with minimal $L^2_{\varphi_t}(D)$ norm. Hence H\"ormander's estimate gives
$$
\int_D |u_t|^2 e^{-\varphi_t} \le \int_D |\bar{\partial}(\partial \varphi_t/\partial t)|^2_{i\partial\bar{\partial}\varphi_t} |K_t(\cdot,z)|^2 e^{-\varphi_t}.
$$
Combining this inequality with (\ref{ch3:eq:SecondVariation}), we obtain
\begin{eqnarray}\label{ch3:eq:B_4}
\frac{\partial^2 K_t(z)}{\partial t\partial \bar{t}}  & \ge &   \int_D \left|\frac{\partial K_t(\cdot,z)}{\partial \bar{t}}\right|^2 e^{-\varphi_t}\\
&&	 + \int_D |K_t(\cdot,z)|^2 \left[\frac{\partial^2 \varphi_t}{\partial t\partial \bar{t}} - |\bar{\partial}(\partial \varphi_t/\partial t)|^2_{i\partial\bar{\partial}\varphi_t} \right] e^{-\varphi_t}. \nonumber
\end{eqnarray}
Note that
$$
 |\bar{\partial}(\partial \varphi_t/\partial t)|^2_{i\partial\bar{\partial}\varphi_t} =\sum_{j,k} \left(\frac{\partial^2 \varphi_t}{\partial z_j\partial\bar{z}_k}\right)^{-1}
\frac{\partial^2 \varphi_t}{\partial t\partial\bar{z}_j} \overline{\frac{\partial^2 \varphi_t}{\partial t\partial\bar{z}_k}}.
$$
The block determinant formula gives, when $B$ is invertible,
$$
\left|
\begin{array}{cc} 
A & D\\
C & B
\end{array}
\right|
=\left|A- D B^{-1} C\right| |B|,
$$
hence
\begin{equation}\label{ch3:eq:B_5}
\frac{\partial^2 \varphi_t}{\partial t\partial \bar{t}} - |\bar{\partial}(\partial \varphi_t/\partial t)|^2_{i\partial\bar{\partial}\varphi_t}= \frac{|H(\varphi)|}{|H(\varphi_t)|}>0,
\end{equation}
where $H(\varphi)$ and $H(\varphi_t)$ are the complex Hessian matrices of $\varphi$ and $\varphi_t$, respectively.

Next, substituting $f=K_t(\cdot,z)$ into the first identity in (\ref{ch3:eq:B_0'}) and combining it with (\ref{ch3:eq:B_0}), we obtain
\begin{eqnarray*}
\left|\frac{\partial K_t(z)}{\partial \bar{t}}\right|^2 & = & \left| \int_D \frac{\partial K_t(\cdot,z)}{\partial \bar{t}}\,\overline{K_t(\cdot,z)} e^{-\varphi_t} \right|^2\\
& \le &  \int_D \left|\frac{\partial K_t(\cdot,z)}{\partial \bar{t}}\right|^2 e^{-\varphi_t}\cdot \int_D |K_t(\cdot,z)|^2 e^{-\varphi_t}\\
& \le & \frac{\partial^2 K_t(z)}{\partial t\partial\bar t}\cdot K_t(z)\quad\text{by (\ref{ch3:eq:B_4}) and (\ref{ch3:eq:B_5})},
\end{eqnarray*}
that is,
$$
\frac{\partial^2 \log K_t(z)}{\partial t\partial \bar{t}}\ge 0.
$$
\end{proof}

\begin{lemma}\label{ch3:lm:B_2}
Let $\Omega_1$ be a bounded pseudoconvex domain in $\mathbb{C}^n$, and let $\Omega_0\subset\subset\Omega_1$ be strongly hyperconvex: there exists a continuous plurisubharmonic function $\rho$ on a neighborhood of $\overline{\Omega}_0$ such that $\Omega_0=\{\rho<0\}$ and $\overline{\Omega}_0=\{\rho\le0\}$, with both sublevel sets taken in that neighborhood. Suppose that $\{\varphi_j\}\subset PSH(\Omega_1)$ satisfies the following conditions:
     \begin{enumerate}
  \item[$(1)$] $\varphi_j=\varphi$ on $\overline{\Omega}_0$;
  \item[$(2)$] $\varphi_j\uparrow+\infty$ a.e.\ on $\Omega_1\backslash\overline{\Omega}_0$.
    \end{enumerate}   
Then $K_{\Omega_1,\varphi_j}(z)$ converges to $K_{\Omega_0,\varphi}(z)$ uniformly on compact subsets of $\Omega_0$.
        \end{lemma}

\begin{proof}
Given $\varepsilon\ll1$, set $\Omega_\varepsilon=\{\rho<\varepsilon\}$. Clearly, $\Omega_\varepsilon$ is pseudoconvex. We may assume that $\rho>-1/e$. Choose a cutoff function $\chi:\mathbb{R}\rightarrow[0,1]$ such that $\chi|_{(-\infty,-\log2]}=1$ and $\chi|_{[0,\infty)}=0$. Let $f\in A^2_\varphi(\Omega_0)$. The Donnelly--Fefferman estimate implies that the equation
$$
\bar{\partial} u = \bar{\partial} \left[f\chi(\log(-\log(-\rho+\varepsilon))-\log(-\log \varepsilon))\right]=:v_\varepsilon
$$
has a solution $u_\varepsilon$ on
$\Omega_\varepsilon$ satisfying
$$
\int_{\Omega_\varepsilon} |u_\varepsilon|^2 e^{-\varphi_1} 
 \le  C_0 \int_{\Omega_\varepsilon} |v_\varepsilon|^2_{i\partial\bar{\partial}(-\log(-\rho+\varepsilon))} e^{-\varphi_1}
 \le  \frac{C_0}{|\log \varepsilon|^2} \int_{\Omega_0} |f|^2 e^{-\varphi}.
$$
Thus
$$
f_\varepsilon=f\chi(\log(-\log(-\rho+\varepsilon))-\log(-\log \varepsilon))-u_\varepsilon
$$
is holomorphic on $\Omega_\varepsilon$ and satisfies
\begin{eqnarray*}
\int_{\Omega_0} |f_\varepsilon-f|^2 e^{-\varphi} & \le & 2 \int_{-\rho\le \sqrt{\varepsilon}-\varepsilon} |f|^2 e^{-\varphi}+ 2 \int_{\Omega_0} |u_\varepsilon|^2 e^{-\varphi} \\
& \le & 2 \int_{-\rho\le \sqrt{\varepsilon}-\varepsilon} |f|^2 e^{-\varphi}+  \frac{2C_0}{|\log \varepsilon|^2} \int_{\Omega_0} |f|^2 e^{-\varphi}\\
& \rightarrow & 0 \ \ \ (\varepsilon\rightarrow 0),
\end{eqnarray*}
and
$$
\|f_\varepsilon\|_{\Omega_\varepsilon,\varphi_1}\le \|f\|_{\Omega_0,\varphi}+ \|u_\varepsilon\|_{\Omega_\varepsilon,\varphi_1}.
$$
Next, choose $\kappa_\varepsilon\in C^\infty_0(\Omega_\varepsilon)$ such that $\kappa_\varepsilon=1$ on a neighborhood of $\overline{\Omega}_0$.
By H\"ormander's $L^2$ estimate, the equation
$$
\bar{\partial} u =\bar{\partial}(\kappa_\varepsilon f_\varepsilon)
$$
has a solution $u_{\varepsilon,j}$ on $\Omega_1$ such that
$$
\int_{\Omega_1} |u_{\varepsilon,j}|^2 e^{-\varphi_j} \le C \int_{\Omega_1} |f_\varepsilon|^2 |\bar{\partial}\kappa_\varepsilon|^2 e^{-\varphi_j}
\rightarrow 0\ \ \ (j\rightarrow \infty),  
$$
Here the constant $C$ depends only on the diameter of $\Omega_1$.
Thus the function $f_{\varepsilon,j}=\kappa_\varepsilon f_\varepsilon-u_{\varepsilon,j}$ is holomorphic on $\Omega_1$ and satisfies
$$
\int_{\Omega_0} |f_{\varepsilon,j}-f_\varepsilon|^2 e^{-\varphi} \le  \int_{\Omega_1} |u_{\varepsilon,j}|^2 e^{-\varphi_j},
$$
and
$$
 \|f_{\varepsilon,j}\|_{\Omega_1,\varphi_j}  \le \|f_\varepsilon\|_{\Omega_\varepsilon,\varphi_1} +  \|u_{\varepsilon,j}\|_{\Omega_1,\varphi_j}.
$$
Consequently, for every fixed $w\in\Omega_0$,
$$
\liminf_{j\rightarrow\infty} K_{\Omega_1,\varphi_j}(w)\ge \liminf_{j\rightarrow\infty} \frac{|f_{\varepsilon,j}(w)|^2}{\|f_{\varepsilon,j}\|^2_{\Omega_1,\varphi_j}} \ge \frac{|f_\varepsilon(w)|^2}{\|f_\varepsilon\|^2_{\Omega_\varepsilon,\varphi_1}}.
$$
Letting $\varepsilon\rightarrow0$, we obtain
$$
\liminf_{j\rightarrow\infty} K_{\Omega_1,\varphi_j}(w)\ge |f(w)|^2/\|f\|^2_{\Omega_0,\varphi},\ \ \ \forall\, f\in A^2_{\varphi}(\Omega_0),
$$
and hence
$$
\liminf_{j\rightarrow\infty} K_{\Omega_1,\varphi_j}(w)\ge K_{\Omega_0,\varphi}(w).
$$
On the other hand, the extremal property implies that $K_{\Omega_1,\varphi_j}(w)\le K_{\Omega_0,\varphi}(w)$.
Thus
$$
\lim_{j\rightarrow\infty} K_{\Omega_1,\varphi_j}(w) = K_{\Omega_0,\varphi}(w).
$$
\end{proof}

\begin{proof}[Proof of Theorem \ref{ch3:th:Berndtsson}]
We divide the proof into three steps.

Step 1: Plurisubharmonicity in $t$ in a special case. First suppose that $\Omega=\{(z,t):\rho(z,t)<0\}$ is bounded and that $\rho$ is smooth and strictly plurisubharmonic on a neighborhood of its closure. Also assume that $\varphi$ is strictly plurisubharmonic on a neighborhood of $\overline{\Omega}$. We first consider $k=1$. Fix $t$, say $t=0$, and choose $\varepsilon\ll1$ and a neighborhood $U$ of $0\in\mathbb{C}$ such that
$$
\overline{\Omega}_t\subset \widetilde{\Omega}:=\{z:\rho(z,0)<\varepsilon\},\ \ \ \forall\,t\in U
$$
and $\varphi$ is strictly plurisubharmonic on a neighborhood of the closure of $\widetilde{\Omega}$. Choose a sequence of increasing convex functions $\lambda_j$ satisfying $\lambda_j|_{(-\infty,0]}=0$ and $\lambda_j\uparrow+\infty$ on $(0,+\infty)$. Set $\varphi_j=\varphi+\lambda_j\circ\rho$. Lemma \ref{ch3:lm:B_1} shows that $\log K_{\widetilde{\Omega},\varphi_j(\cdot,t)}(z)$ is subharmonic in $t$. On the other hand, Lemma \ref{ch3:lm:B_2} gives
$$
 \log K_{\widetilde{\Omega},\varphi_j(\cdot,t)}(z) \uparrow \log K_t(z),
$$
Theorem \ref{ch3:th:USC} implies that $\log K_t(z)$ is upper semicontinuous in $t$. Hence $\log K_t(z)$ is subharmonic in $t$.

For general $k\ge1$, the function $\log K_t(z)$ is upper semicontinuous in $t$ and subharmonic on every complex line in $\mathbb{C}^k_t$.
It is therefore plurisubharmonic in $t$.

Step 2: Plurisubharmonicity in $t$ in the general case. Choose a strictly plurisubharmonic exhaustion function $\rho$ on $\Omega$ and a sequence of positive numbers
$N_j\uparrow+\infty$ such that $\Omega_j:=\{\rho<N_j\}$ has smooth boundary. Choose smooth strictly plurisubharmonic functions $\varphi_j$ on $\Omega_{j+1}$ such that $\varphi_j\downarrow\varphi$. We have already proved that $\log K_{\Omega_{j,t},\varphi_j(\cdot,t)}(z)$ is plurisubharmonic in $t$. On the other hand, Corollary \ref{ch3:cor:decreasing} gives
$$
\log K_{\Omega_{j,t},\varphi_j(\cdot,t)}(z)\downarrow \log K_t(z),
$$ 
If the limiting weight on a fiber is identically $-\infty$, its kernel is defined to be zero; convergence of the kernels then follows directly from the pointwise estimates. Thus $\log K_t(z)$ is plurisubharmonic in $t$.

Step 3: Plurisubharmonicity in $(z,t)$. We use a trick due to Oka. Clearly, it suffices to show that, for fixed $t_0,z_0$ and $(a,b)\in\mathbb{C}^n\times\mathbb{C}^k$, the function $\log K_{t_0+\zeta b}(z_0+\zeta a)$ is subharmonic in $\zeta$. Since the complex Jacobian of a translation is identically $1$, Proposition \ref{ch3:prop:W-Trans} gives
$$
K_{t_0+\zeta b}(z_0+\zeta a) = K_{\Omega_{t_0+\zeta b}-\zeta a,\varphi_{t_0+\zeta b}(\cdot+\zeta a)}(z_0),
$$
The family $\{\Omega_{t_0+\zeta b}-\zeta a:|\zeta|<\varepsilon\}$ is still pseudoconvex, and $\varphi(z+\zeta a,t_0+\zeta b)$ is plurisubharmonic in $(z,\zeta)$.
The conclusion therefore follows from Step 2.
\end{proof}

We now give some applications of Berndtsson's theorem.

\begin{theorem}[Complex Pr\'ekopa theorem]\label{ch3:th:Prekopa}
Let $\Omega\subset\mathbb{C}_z^n\times\mathbb{C}_t^k$ be a pseudoconvex domain, and let $\varphi\in PSH(\Omega)$. Suppose that every $\Omega_t$ is a complete circular domain
and that every $\varphi_t$ is invariant under
$z\mapsto e^{i\theta}z$ for all $\theta\in\mathbb{R}$. Then
$$
\widetilde{\varphi}(t):= -\log \int_{\Omega_t} e^{-\varphi_t}
$$
is plurisubharmonic.
\end{theorem}

\begin{proof}
Since $\Omega_t$ is a complete circular domain, every $f\in\mathcal O(\Omega_t)$
has a homogeneous expansion converging uniformly on compact subsets:
\[
f(z)=\sum_{m=0}^{\infty}\sum_{|\alpha|=m}c_\alpha z^\alpha
=:\sum_{m=0}^{\infty}P_m(z).
\]
For $f\in A^2_t$, the homogeneous components satisfy
\[
P_m(z)=\frac1{2\pi}\int_0^{2\pi}f(e^{i\theta}z)e^{-im\theta}\,d\theta.
\]
Rotational invariance of the weight, Parseval's identity, and Tonelli's theorem
give
\[
\|f\|_t^2=\sum_{m=0}^{\infty}\|P_m\|_t^2.
\]
The components of different degrees are orthogonal. Thus $P_m\in A^2_t$,
and the expansion also converges in $A^2_t$.

Write $M_t=\int_{\Omega_t}e^{-\varphi_t}$ and first suppose that $M_t<\infty$.
For every square-integrable homogeneous polynomial $P_m$ of degree $m\ge1$,
rotational invariance gives
\[
\int_{\Omega_t}P_m e^{-\varphi_t}
=e^{im\theta}\int_{\Omega_t}P_m e^{-\varphi_t},
\qquad \theta\in\mathbb{R}.
\]
Hence $\int_{\Omega_t}P_m e^{-\varphi_t}=0$, so $P_m$ is orthogonal to $1$.
Applying Gram--Schmidt to bases of the spaces of square-integrable homogeneous
polynomials, starting with the constant function, gives a complete orthonormal
basis $\{h_j\}$ of $A^2_t$ with
\[
h_0=M_t^{-1/2},\qquad h_j(0)=0\quad(j\ge1).
\]
Consequently, $K_t(0)=\sum_j|h_j(0)|^2=M_t^{-1}$.
If $M_t=\infty$, the fact that the constant component $P_0=f(0)$ belongs to
$A^2_t$ forces $f(0)=0$ for every $f\in A^2_t$. In this case also,
$K_t(0)=0=M_t^{-1}$. Therefore
\[
\widetilde\varphi(t)=\log K_t(0),
\]
and the conclusion follows from Theorem \ref{ch3:th:Berndtsson}, allowing the
case where the function is identically $-\infty$.
\end{proof}

Taking $\varphi=0$ in Theorem \ref{ch3:th:Prekopa}, we obtain

\begin{corollary}
If $\Omega$ is pseudoconvex and every $\Omega_t$ is a complete circular domain, then $-\log|\Omega_t|$ is plurisubharmonic.
\end{corollary}

\begin{corollary}[Kiselman's minimum principle]
Under the hypotheses of Theorem \ref{ch3:th:Prekopa},
$$
\widetilde{\varphi}_\infty(t):=\inf_{\Omega_t} \varphi_t
$$
is plurisubharmonic.
\end{corollary}

\begin{proof}
First suppose that all fibers lie in a fixed ball $B(0,R)$ and that
$\varphi\ge-N$ for some $N>0$. By Theorem \ref{ch3:th:Prekopa}, for every $p>0$,
\[
\widetilde\varphi_p(t):=-\frac1p\log\int_{\Omega_t}e^{-p\varphi_t}
\]
is plurisubharmonic. The usual limit of $L^p$ norms on a finite measure space
gives
\[
\widetilde\varphi_p(t)\longrightarrow
\operatorname*{ess\,inf}_{\Omega_t}\varphi_t
=\inf_{\Omega_t}\varphi_t
=\widetilde\varphi_\infty(t)
\qquad(p\longrightarrow\infty).
\]
The essential infimum equals the infimum because $\varphi_t$ is upper
semicontinuous. Moreover, upper semicontinuity of $\varphi$ and openness of
$\Omega$ imply that $\widetilde\varphi_\infty$ is upper semicontinuous.

Locally in the parameter, there is a fixed ball $B\subset\Omega_t$ on which
$\varphi_t\le C$, since $(0,t)\in\Omega$. For $p\ge1$ this gives
\[
-N-\frac{\log|B(0,R)|}{p}
\le\widetilde\varphi_p(t)
\le C-\frac{\log|B|}{p}.
\]
Thus the functions $\widetilde\varphi_p$ are locally uniformly bounded above
and below. Dominated convergence allows us to pass to the limit in their
submean inequalities on complex lines. Together with upper semicontinuity,
this proves that $\widetilde\varphi_\infty$ is plurisubharmonic.

For the general case, apply the preceding argument to
\[
\Omega^{(j)}=\Omega\cap\bigl(B(0,j)\times\mathbb{C}^k\bigr),
\qquad \varphi^{(j)}=\max\{\varphi,-j\}.
\]
These domains and weights retain pseudoconvexity, complete circularity of the
fibers, and rotational invariance. The resulting plurisubharmonic functions
\[
t\longmapsto\inf_{z\in\Omega_t\cap B(0,j)}\max\{\varphi(z,t),-j\}
\]
decrease to $\widetilde\varphi_\infty(t)$. Hence the latter is plurisubharmonic
or identically $-\infty$.
\end{proof}

\begin{corollary}\label{ch3:cor:convex}
Let $\Omega$ be a pseudoconvex domain, and let $\rho$ be a plurisubharmonic function on $\Omega$. Set $\Omega_t=\{\rho<t\}$ for $t\in\mathbb{R}$ and $K_t=K_{\Omega_t}$. Then $\log K_t(z)$ is convex in $t$. On any parameter interval where $K_t(z)>0$, the function $\log K_t(z)$ is finite and convex. The case where the kernel vanishes identically is treated separately as a trivial case.
\end{corollary}

\begin{proof}
It is easy to verify that
$$
\widehat{\Omega}:=\{(z,\zeta)\in \Omega \times \mathbb{C}^*: \rho(z)-\log |\zeta|<0\}
$$
is a pseudoconvex domain in $\mathbb{C}^{n+1}$. Theorem \ref{ch3:th:Berndtsson} therefore implies that $\log\widehat K_\zeta(z):=\log K_{\widehat\Omega_\zeta}(z)$ is subharmonic in $\zeta$. Hence
$$
\log \widehat{K}_{r} (z) = \frac1{2\pi} \int_0^{2\pi} \log \widehat{K}_{re^{i\theta}}(z) d\theta
$$ 
is convex in $\log r$; that is, $\log K_t(z)$ is convex in $t$.
\end{proof}

This convexity also yields H\"older continuity in the parameter for unweighted Bergman kernels on varying domains.

\begin{definition}
Let $\{\Omega_t\}_{t\in(-1,1)}$ be a family of bounded domains in $\mathbb{C}^n$ whose total space $\mathcal D=\{(z,t):z\in\Omega_t,\ -1<t<1\}$ is open. Suppose that $\rho(z,t)$ is continuous on $\mathcal D$ and that $\rho_t=\rho(\cdot,t)$ is a negative exhaustion function on $\Omega_t$. Given $0<\alpha\le1$, if for every $\gamma>0$ there exist $b_\gamma\gg1$ and $c_\gamma\ll1$ such that, whenever $t,s\in(-1,1)$ satisfy $|t-s|\le c_\gamma$,
$$
\left\{z\in \Omega_t: \rho_t<-b_\gamma |t-s|^\alpha\right\}\subset \left\{z\in \Omega_s: \rho_s<-\gamma |t-s|^\alpha\right\},
$$
then the family $\{\Omega_t\}$ is said to be $\rho_t$-H\"older continuous of order $\alpha$.
\end{definition}

\begin{example}
Suppose that $\rho(z,t)$ is H\"older continuous of order $\alpha$ on $\mathbb{C}^n\times(-1,1)$. Set $\rho_t=\rho(\cdot,t)$ and assume that $\Omega_t=\{\rho_t<0\}$ is a bounded domain. Then $\Omega_t$ is $\rho_t$-H\"older continuous of order $\alpha$.
\end{example}

It is proved in  \cite{ChenParameter2016}  that, if $\Omega_t$ is $\rho_t$-H\"older continuous of order $\alpha$ and $\rho_t\in PSH(\Omega_t)$, then $K_{\Omega_t}(z,w)$ is H\"older continuous in $t$ of every order $0<\beta<\alpha$.
This result can be improved as follows.

\begin{proposition}\label{ch3:prop:Holder-endpoint}
If $\Omega_t$ is $\rho_t$-H\"older continuous of order $\alpha$ and $\rho_t\in PSH(\Omega_t)$, then $K_t(z,w)$ is H\"older continuous in $t$ of order $\alpha$.
In particular, $\alpha=1$ gives local Lipschitz continuity.
\end{proposition}

\begin{proof}
First consider a fixed bounded pseudoconvex domain $D$ with a negative continuous plurisubharmonic exhaustion function $\rho$. Write $D_\varepsilon=\{\rho<-\varepsilon\}$. Fix $a>0$ and $z\in D_a$. By Corollary \ref{ch3:cor:convex}, the function
\[
t \longmapsto \log K_{\{\rho<t\}}(z)
\]
is convex, so $t\mapsto K_{\{\rho<t\}}(z)$ is also convex. Note that $\{\rho<0\}=D$ and, by Ramadanov's theorem, $K_{\{\rho<t\}}(z)\downarrow K_D(z)$ as $t\uparrow0$. Thus the convexity inequality extends to the endpoint $t=0$. For $0<\varepsilon<a$, writing $-\varepsilon$ as a convex combination of $-a$ and $0$ gives
\[
K_{D_\varepsilon}(z)
\le\frac{\varepsilon}{a}K_{D_a}(z)
+\left(1-\frac{\varepsilon}{a}\right)K_D(z).
\]
Combining this with the monotonicity of the kernel yields the linear approximation estimate
\begin{equation}\label{ch3:eq:sublevel-linear}
0\le K_{D_\varepsilon}(z)-K_D(z)
\le\frac{\varepsilon}{a}\bigl(K_{D_a}(z)-K_D(z)\bigr).
\end{equation}

We now return to the family $\Omega_t$. Fix $t_0\in(-1,1)$ and a compact set $S\subset\subset\Omega_{t_0}$, and write $D_t(\varepsilon)=\{\rho_t<-\varepsilon\}$. By the openness of the total space and the joint continuity of $\rho(z,t)$, we can choose fixed $a,r,\eta>0$ such that
\[
B(z,r)\subset D_t(a)\qquad(z\in S,\ |t-t_0|<\eta).
\]
Thus $K_{D_t(a)}(z)\le K_{B(z,r)}(z)$. By \eqref{ch3:eq:sublevel-linear},
\begin{equation}\label{ch3:eq:uniform-sublevel-linear}
0\le K_{D_t(\varepsilon)}(z)-K_t(z)\le (K_{B(z,r)}(z)/a)\varepsilon
\qquad(z\in S,\ 0<\varepsilon<a).
\end{equation}

For $t\ne s$, set $\delta=|t-s|^\alpha$. Taking $\gamma=1$ in the definition of $\rho_t$-H\"older continuity gives a constant $b>1$. After shrinking $\eta$, we may assume that $b\delta<a$ and
\[
D_t(b\delta)\subset D_s(\delta)\subset\Omega_s,
\qquad D_s(b\delta)\subset\Omega_t.
\]
Set $G=\Omega_t\cap\Omega_s$. For $q=t,s$, we have
\[
D_q(b\delta)\subset G\subset\Omega_q.
\]
Hence, by monotonicity and \eqref{ch3:eq:uniform-sublevel-linear},
\begin{equation}\label{ch3:eq:intersection-linear}
0\le K_G(z)-K_q(z)\le ( K_{B(z,r)}(z) b/a)\delta\qquad(z\in S,\ q=t,s).
\end{equation}

Finally, Corollary \ref{cor:Wang} extends the diagonal estimate to the off-diagonal case.
\end{proof}

A particularly challenging question is the following.

\begin{problem}
For a smoothly varying family of bounded pseudoconvex domains $\Omega_t$, is the Bergman kernel $K_{\Omega_t}(z)$ jointly smooth in $(z,t)$?
\end{problem}

When all $\Omega_t$ are strongly pseudoconvex, Hamilton's theory (see Hamilton \cite{Hamilton1979}) gives an affirmative answer. The case of general weakly pseudoconvex families requires further study.

\needspace{6\baselineskip}
Using Corollary \ref{ch3:cor:convex}, Lempert gave an elegant proof of the Suita conjecture\footnote{First proved by B\l ocki \cite{Blocki2013}.} based on convexity. We supply the notation and limiting argument below. Let $\Omega$ be a domain in $\mathbb{C}$ admitting a Green function, denoted by $g_\Omega$. For $w\in\Omega$, define 
$$
c_\Omega(w)=\exp\bigl(\lim_{z\to w}[g_\Omega(z,w)-\log|z-w|]\bigr).
$$
 Set
$$
\Omega_t=\{z\in \Omega: g_\Omega(z,w)<t\},\ \ \ t\le 0.
$$
Note that $\Omega_t\uparrow\Omega$ as $t\uparrow0$, so Ramadanov's theorem gives $K_{\Omega_t}(w)\downarrow K_\Omega(w)$. Corollary \ref{ch3:cor:convex} shows that $\log K_t(z):=\log K_{\Omega_t}(z)$ is convex in $t$. Since $K_{\Omega_t}(w)$ depends only on the connected component $U_t$ of $\Omega_t$ containing $w$, and
$$
U_t \sim \left\{ |z-w|<e^t/c_\Omega(w)\right\}\ \ \ (t\rightarrow -\infty),
$$
we have
$$
K_t(w) \sim \frac{c_\Omega(w)^2}{\pi e^{2t}}\ \ \ (t\rightarrow -\infty).
$$
Thus $k(t):=\log K_t(w)+2t$ is uniformly bounded above as $t\rightarrow-\infty$. Since $k$ is convex, it is increasing. Therefore
$$
\log K_\Omega(w) = k(0) \ge \lim_{t\rightarrow -\infty} k(t) = \log \frac{c_\Omega(w)^2}{\pi},
$$
that is, $K_\Omega(w)\ge c_\Omega(w)^2/\pi$.

\begin{remark}
The square of the dual norm of the evaluation functional is the Bergman kernel. Replacing the evaluation functional by a suitable continuous linear functional $\ell$, consider
\[
\|\ell\|_{t,*}^{\,2}=\sup_{f\ne0}\frac{|\ell(f)|^2}{\|f\|_t^2}.
\]
Berndtsson--Lempert \cite{BerndtssonLempert2016} constructed suitable parameter families and used the logarithmic convexity of these dual norms, together with degeneration limits, to prove the Ohsawa--Takegoshi extension theorem with the optimal constant,  which was first obtained by B\l ocki \cite{Blocki2013}. Conversely, Guan-Zhou \cite{GuanZhouExtension} found that the Ohsawa--Takegoshi extension theorem with the optimal constant also implies Berndtsson's theorem. 
\end{remark}

\section{Convexity}\label{ch3:sec:real-convexity}

In the preceding section, Berndtsson's theorem was used to prove the logarithmic convexity of Bergman kernels of sublevel sets with respect to a real parameter. This section concerns another type of convexity: when both the domain and the weight are real convex, the Bergman kernel is also logarithmically convex as a function of the spatial variable.

Here we identify $\mathbb{C}^n$ with $\mathbb{R}^{2n}$. By a convex function we always mean a real convex function, that is,
\[
u((1-t)z_0+tz_1)\le (1-t)u(z_0)+tu(z_1),\qquad 0\le t\le1.
\]
Strict convexity means that the inequality is strict whenever $z_0\ne z_1$ and $0<t<1$. We allow a convex function to be identically $-\infty$; if a convex function takes values in $[-\infty,+\infty)$ and equals $-\infty$ at one point, it is identically $-\infty$ on the convex domain. Apart from this case, the weights $\varphi$ below are finite convex functions and are therefore continuous and plurisubharmonic.

\begin{theorem}[Xiong \cite{Xiong2026}]\label{ch3:th:Xiong-log}
Let $\Omega\subset\mathbb{C}^n$ be a convex domain, and let $\varphi$ be a convex function on $\Omega$. Then $\log K_{\Omega,\varphi}$ is convex, with the identically $-\infty$ case allowed.
\end{theorem}

Before proving the theorem, we review its historical background. In 1985, Caffarelli--Friedman \cite{CF1985} used partial differential equations to prove that $\log\lambda_\Omega$ is convex whenever $\Omega\subsetneq\mathbb{C}$ is a convex domain. Here $\lambda_\Omega$ denotes the density of the hyperbolic metric $\lambda_\Omega(z)|dz|$, normalized by $\lambda_{\mathbb{D}}(z)=(1-|z|^2)^{-1}$. By the Riemann mapping theorem and the transformation formula for the kernel, $K_\Omega=\lambda_\Omega^2/\pi$, so $\log K_\Omega=2\log\lambda_\Omega-\log\pi$ is also convex. In 1990, Gustafsson \cite{Gustafsson1990} gave a proof using only complex analysis. Interestingly, this result was already contained in the 1982 paper of Minda--Wright \cite{MindaWright1982}, with a particularly simple proof: given any two points $z_0,z_1\in\Omega$, set $z_t=(1-t)z_0+tz_1$ for $t\in[0,1]$. Choose conformal maps $F_t:\mathbb{D}\to\Omega$ with $F_t(0)=z_t$ and $F_t'(0)>0$. Set $g_t=(1-t)F_0+tF_1$. Then $g_t\in\mathcal O(\mathbb{D})$ and $g_t(0)=z_t$, and the convexity of $\Omega$ gives $g_t(\mathbb{D})\subset\Omega$. The decreasing property of the hyperbolic metric under holomorphic maps yields
$$
\lambda_\Omega (g_t(0))|g_t'(0)| \le \lambda_{\mathbb{D}}(0)=1, 
$$
or, equivalently,
\begin{equation}
\frac1{\lambda_\Omega(z_t)}\ge (1-t) F_0'(0) + tF_1'(0) = \frac{1-t}{\lambda_\Omega(z_0)}+\frac{t}{\lambda_\Omega(z_1)},
\end{equation}
Thus $1/\lambda_\Omega$ is concave. Since $-\log(-t)$ is increasing and convex on $(-\infty,0)$,
$$
\log \lambda_\Omega = -\log(-(-1/\lambda_\Omega))
$$
is convex.

In higher dimensions, there is no tool as powerful as the Riemann mapping theorem, so a new method is needed. The following lemma uses real linear transformations to reduce convexity to subharmonicity, thereby connecting the problem to Berndtsson's theorem.

\begin{lemma}\label{ch3:lm:real-convex-test}
Let $U\subset\mathbb{C}$ be a convex domain and $u\in C^2(U)$. For each $|\lambda|<1$, set
\[
t_\lambda(s)=s+\lambda^2\bar s,\qquad V_\lambda=t_\lambda^{-1}(U).
\]
If $u\circ t_\lambda$ is subharmonic on $V_\lambda$, then $u$ is convex.
\end{lemma}
\begin{proof}
A direct computation shows that the determinant of the real linear map $t:=t_\lambda$ is $1-|\lambda|^4>0$, so it is invertible. The chain rule gives
$$
\frac{\partial u(t)}{\partial \bar{s}}  =  u_t \frac{\partial t}{\partial \bar{s}} + u_{\bar{t}} \frac{\partial \bar{t}}{\partial \bar{s}} = \lambda^2 u_t + u_{\bar{t}},
$$
\begin{eqnarray*}
0\le \frac{\partial^2 u( t)}{\partial s\partial\bar s} & = & \lambda^2 u_{tt} \frac{\partial t}{\partial s} + \lambda^2 u_{t\bar{t}} \frac{\partial \bar{t}}{\partial s} + u_{t\bar{t}}\frac{\partial t}{\partial s} + u_{\bar{t}\, \bar{t}}\frac{\partial \bar{t}}{\partial s} \\
& = & (1+|\lambda|^4)u_{t\bar t}
 +2\operatorname{Re}(\lambda^2u_{tt}).
\end{eqnarray*}
Writing $\lambda=r\eta$ with $|\eta|=1$ and letting $r\uparrow1$, we obtain
\[
u_{t\bar t}+\operatorname{Re}(\eta^2u_{tt})\ge0.
\]
On the other hand, the second derivative in the real direction $\eta$ is
\[
\left.\frac{d^2}{dx^2}u(t+x\eta)\right|_{x=0}
=2\bigl(u_{t\bar t}+\operatorname{Re}(\eta^2u_{tt})\bigr).
\]
Thus the real Hessian is positive semidefinite, and $u$ is convex.
\end{proof}

\begin{remark}
This criterion also applies to nonsmooth functions: if every $u\circ t_\lambda$ is subharmonic, then $u$ is convex or identically $-\infty$. In the nontrivial case, taking $\lambda=0$ first shows that $u$ is subharmonic. Convolve with a nonnegative smooth radial kernel on relatively compact subdomains. Every translate of $u\circ t_\lambda$ remains subharmonic, so convolution preserves the condition above. Apply the lemma to the smooth approximations and then pass to the pointwise limit. The application to Bergman kernels below requires only the smooth case.
\end{remark}

\begin{proof}[Proof of Theorem \ref{ch3:th:Xiong-log}]
Clearly, it suffices to consider the case where $\varphi$ is finite everywhere.
First assume that $\Omega$ is bounded and that $\varphi$ is defined on a neighborhood of $\overline{\Omega}$. Then constant functions belong to $A^2_\varphi(\Omega)$, so $K_{\Omega,\varphi}>0$ and its logarithm is real analytic.

Fix $a\in\Omega$ and $v\in\mathbb{C}^n\setminus\{0\}$, and write
\[
U=\{t\in\mathbb{C}:a+tv\in\Omega\},\qquad
u(t)=\log K_{\Omega,\varphi}(a+tv).
\]
The set $U$ is a convex domain. For $|\lambda|<1$, consider
\[
\mathcal D_\lambda
 =\{(w,s):s\in V_\lambda,\ w+a+t_\lambda(s)v\in\Omega\},
 \]
 \[
\psi_\lambda(w,s) =\varphi(w+a+t_\lambda(s)v).
\]
This domain is convex, and the weight is real convex, so they are pseudoconvex and plurisubharmonic, respectively. Although the change of variables here is not holomorphic, it preserves real convexity; this is precisely what gives the pseudoconvexity and plurisubharmonicity needed to apply Berndtsson's theorem. The fibers and their weights are
\[
(\mathcal D_\lambda)_s=\Omega-a-t_\lambda(s)v,\qquad
(\psi_\lambda)_s(w)=\varphi(w+a+t_\lambda(s)v).
\]
By the transformation formula for the kernel under translations and Theorem \ref{ch3:th:Berndtsson},
\[
s\longmapsto\log K_{(\mathcal D_\lambda)_s,(\psi_\lambda)_s}(0)
=u(t_\lambda(s))
\]
is subharmonic. Lemma \ref{ch3:lm:real-convex-test} shows that $u$ is convex. Every real line segment lies in a complex line, so $\log K_{\Omega,\varphi}$ is convex on $\Omega$.

In the general case, choose an increasing sequence of relatively compact bounded convex domains $\Omega_j\Subset\Omega$ exhausting $\Omega$. The weighted Ramadanov theorem gives
\[
K_{\Omega_j,\varphi}(z)\downarrow K_{\Omega,\varphi}(z).
\]
The logarithmic convexity inequality is preserved in this limit. If the limiting kernel vanishes at one point, convexity of its logarithm forces it to vanish everywhere; otherwise, we obtain a finite convex function.
\end{proof}

The same method also gives joint convexity in the parameter and spatial variables.
\begin{corollary}\label{ch3:cor:Xiong-family}
Let $\mathcal D\subset\mathbb{C}_z^n\times\mathbb{C}_t^k$ be a convex domain, and let $\Phi$ be a finite convex function on it. Writing $\mathcal D_t=\{z:(z,t)\in\mathcal D\}$, the function
\[
(z,t)\longmapsto\log K_{\mathcal D_t,\Phi(\cdot,t)}(z)
\]
is convex on $\mathcal D$, with the identically $-\infty$ case allowed.
\end{corollary}
\begin{proof}
First assume that the total space is bounded. Restrict to any complex affine line $(z,t)=(z_0+\tau a,t_0+\tau b)$ and then set $\tau=s+\lambda^2\bar s$. Translating the fiber variable by $z_0+\tau a$ gives a total space and weight that are still convex. Applying Berndtsson's theorem at the origin and the criterion above yields real convexity along this line. For a general total space, use an increasing exhaustion by bounded convex domains and apply the weighted Ramadanov theorem on each fiber.
\end{proof}

We also have the following strict convexity theorem.

\begin{theorem}[Xiong \cite{Xiong2026}]\footnote{In the same paper, it is shown that for a domain $\Omega\subset \mathbb C^n$,  $\log K_\Omega$ is convex if and only if $\Omega$ does not contain any real line.}\label{ch3:th:Xiong-strict}
If $\Omega\subset\mathbb{C}^n$ is a bounded convex domain and $\varphi$ is a finite convex function, then $\log K_{\Omega,\varphi}$ is strictly convex. In particular, $\log K_\Omega$ is strictly convex.
\end{theorem}

Before proving the theorem, we establish a criterion for strict convexity.

\begin{lemma}\label{ch3:lm:Xiong-strict}
Let $\Omega$ be a convex domain, with $K_{\Omega,\varphi}>0$ and $\log K_{\Omega,\varphi}$ convex. Suppose that every real line $L$ intersecting $\Omega$ has a finite endpoint $p\in\partial\Omega$ of $L\cap\Omega$ at which
\[
K_{\Omega,\varphi}(z)\longrightarrow+\infty
\qquad(L\cap\Omega\ni z\to p),
\]
Then $\log K_{\Omega,\varphi}$ is strictly convex.
\end{lemma}

\begin{proof}
Write $u=\log K_{\Omega,\varphi}$. Since the kernel is strictly positive and real analytic, $u$ is also finite and real analytic. Suppose that $u$ is not strictly convex; we shall derive a contradiction.

First, we show that equality at one point in the convexity inequality implies affinity on the entire segment. Suppose that there exist $z_0\ne z_1$ and $t_0\in(0,1)$ such that
\[
u((1-t_0)z_0+t_0z_1)=(1-t_0)u(z_0)+t_0u(z_1).
\]
Set $v=z_1-z_0$ and define on $[0,1]$
\[
f(t)=u(z_0+tv)-\bigl((1-t)u(z_0)+tu(z_1)\bigr).
\]
Then $f$ is convex, $f(0)=f(t_0)=f(1)=0$, and convexity gives $f(t)\le0$. If $0<t<t_0$, write
\[
t_0=\frac{1-t_0}{1-t}\,t+\frac{t_0-t}{1-t}\,1.
\]
By convexity,
\[
0=f(t_0)\le\frac{1-t_0}{1-t}f(t)+\frac{t_0-t}{1-t}f(1)
=\frac{1-t_0}{1-t}f(t),
\]
so $f(t)\ge0$. If $t_0<t<1$, the identity
$t_0=(1-t_0/t)\,0+(t_0/t)\,t$ similarly gives
$0=f(t_0)\le(t_0/t)f(t)$, so again $f(t)\ge0$.
Thus $f\equiv0$; that is, $u(z_0+tv)=At+B$ for $0\le t\le1$, where
$A=u(z_1)-u(z_0)$, $B=u(z_0)$.

Next, we extend affinity to the entire line section. Set
\[
I=\{t\in\mathbb{R}:z_0+tv\in\Omega\}.
\]
By the openness and convexity of $\Omega$, the set $I$ is an open interval containing $[0,1]$. The function
$F(t)=u(z_0+tv)-At-B$ is real analytic on $I$ and vanishes identically on $(0,1)$. The identity principle for real analytic functions gives $F\equiv0$ on all of $I$. Hence
\[
K_{\Omega,\varphi}(z_0+tv)=e^{At+B},\qquad t\in I.
\]

By hypothesis, $I$ has at least one finite endpoint $b$, and the kernel tends to $+\infty$ as that endpoint is approached. However, the preceding identity gives
\[
\lim_{I\ni t\to b}K_{\Omega,\varphi}(z_0+tv)=e^{Ab+B}<+\infty,
\]
a contradiction. Thus $u$ is strictly convex.
\end{proof}

\begin{proof}[Proof of Theorem \ref{ch3:th:Xiong-strict}]
A bounded convex domain is pseudoconvex, and every nonempty intersection with a complex line is convex. It is therefore a domain with the Runge along complex lines  property, so Proposition \ref{prop:propertyL} shows that $K_\Omega(z)$ is an exhaustion function on $\Omega$.
A finite convex weight on a bounded convex domain is bounded below, say $\varphi\ge-C$. By the extremal property,
\[
K_{\Omega,\varphi}(z)\ge e^{-C}K_\Omega(z).
\]
Thus the weighted kernel is also an exhaustion function. The conclusion follows from Lemma \ref{ch3:lm:Xiong-strict}.
\end{proof}

In what follows, we use the convention $0^{-1/(2n)}=+\infty$.

\begin{theorem}[Brunn--Minkowski-type inequality]\label{ch3:th:Xiong-BM}
Let $\Omega_0,\Omega_1\subset\mathbb{C}^n$ be convex domains, and let $z_j\in\Omega_j$. Then
\begin{equation}\label{ch3:eq:Xiong-BM}
K_{\Omega_0+\Omega_1}(z_0+z_1)^{-1/(2n)}
\ge K_{\Omega_0}(z_0)^{-1/(2n)}
 +K_{\Omega_1}(z_1)^{-1/(2n)},
\end{equation}
where $\Omega_0+\Omega_1=\{w_0+w_1:w_j\in\Omega_j\}$.
\end{theorem}
\begin{proof}
First assume that both domains are bounded. After translation, we may assume that $z_0=z_1=0$. For $0<t<1$, set
\[
\Omega_t=(1-t)\Omega_0+t\Omega_1.
\]
The total space
\[
\mathcal D=\{(z,\tau):0<\operatorname{Re}\tau<1,\
 |\operatorname{Im}\tau|<1,\ z\in\Omega_{\operatorname{Re}\tau}\}
\]
is a bounded convex domain. By Corollary \ref{ch3:cor:Xiong-family}, the function $h(t)=\log K_{\Omega_t}(0)$ is convex. We must check the endpoints: since $0\in\Omega_0\cap\Omega_1$,
\[
(1-t)\Omega_0\subset\Omega_t,\qquad t\Omega_1\subset\Omega_t.
\]
Monotonicity and the scaling formula $K_{cD}(0)=c^{-2n}K_D(0)$ give
\[
\limsup_{t\downarrow0}h(t)\le\log K_{\Omega_0}(0),\qquad
\limsup_{t\uparrow1}h(t)\le\log K_{\Omega_1}(0).
\]
Applying convexity on closed intervals in the interior and then approaching the endpoints, we obtain
\begin{equation}\label{ch3:eq:Xiong-max}
K_{\Omega_t}(0)\le\max\{K_{\Omega_0}(0),K_{\Omega_1}(0)\}.
\end{equation}

Set $c_j=K_{\Omega_j}(0)^{-1/(2n)}>0$. Scaling gives $K_{c_j^{-1}\Omega_j}(0)=1$. Applying \eqref{ch3:eq:Xiong-max} to these two normalized domains and taking $t=c_1/(c_0+c_1)$ yields
\[
1\ge K_{(\Omega_0+\Omega_1)/(c_0+c_1)}(0)
=(c_0+c_1)^{2n}K_{\Omega_0+\Omega_1}(0),
\]
This is the desired inequality.

In the general case, exhaust each $\Omega_j$ by an increasing sequence of bounded convex domains containing $z_j$. Their Minkowski sums increase to exhaust $\Omega_0+\Omega_1$. Passing to the limit by Ramadanov's theorem gives \eqref{ch3:eq:Xiong-BM}, including the case where the right-hand side is $+\infty$.
\end{proof}

\begin{corollary}[Convexity of negative powers]\label{ch3:cor:Xiong-power}
If $\Omega\subset\mathbb{C}^n$ is a convex domain and $K_\Omega>0$, then
\[
z\longmapsto -K_\Omega(z)^{-1/(2n)}
\]
is convex; equivalently, $K_\Omega^{-1/(2n)}$ is concave. If $K_\Omega\equiv0$, the negative power is interpreted in the degenerate sense as identically $-\infty$.
\end{corollary}
\begin{proof}
For $0<t<1$, apply the theorem to $(1-t)\Omega$ and $t\Omega$, evaluated at $(1-t)z_0$ and $tz_1$. By convexity,
$(1-t)\Omega+t\Omega=\Omega$. The scaling formula then gives
\[
K_\Omega((1-t)z_0+tz_1)^{-1/(2n)}
\ge (1-t)K_\Omega(z_0)^{-1/(2n)}
+tK_\Omega(z_1)^{-1/(2n)}.
\]
\end{proof}

\begin{example}[Optimality of the exponent]
For the unit polydisc $\mathbb{D}^n$, set $z=(t,\ldots,t)$ with $-1<t<1$. Then
\[
K_{\mathbb{D}^n}(z)=\frac{1}{\pi^n(1-t^2)^{2n}}.
\]
If $\alpha>1/(2n)$, set $p=2n\alpha>1$. Then
\[
\frac{d^2}{dt^2}\bigl[-\pi^{n\alpha}(1-t^2)^p\bigr]
=2p\pi^{n\alpha}(1-t^2)^{p-2}\bigl[1-(2p-1)t^2\bigr]
\]
is negative for $t$ sufficiently close to $1$. Hence convexity of $-K_\Omega^{-\alpha}$ cannot hold for all convex domains, and the exponent $1/(2n)$ is optimal uniformly over such domains. Stronger conclusions may hold for particular domains: for the unit ball $\mathbb{B}^n$, the formula
$K_{\mathbb{B}^n}(z)=\frac{n!}{\pi^n}(1-|z|^2)^{-n-1}$ shows that $-K_{\mathbb{B}^n}^{-1/(n+1)}$ is convex.
\end{example}

\begin{example}[Limitations in the weighted case]
Convexity of negative powers does not extend directly to arbitrary convex weights. Take $\Omega=\mathbb{C}$ and $\varphi(z)=(\operatorname{Re}z)^2=x^2$. With the Euclidean area normalization used in these lectures,
\[
K_{\mathbb{C},\varphi}(z)=\frac{1}{2\pi}e^{x^2}.
\]
Indeed, $x^2=|z|^2/2+\operatorname{Re}(z^2)/2$, so multiplication by $e^{-z^2/4}$ maps this weighted space isometrically onto the Fock space with weight $e^{-|z|^2/2}$, whose kernel is $(2\pi)^{-1}e^{z\bar w/2}$. The transformation formula for the kernel immediately gives the identity above. For every $\alpha>0$,
\[
\frac{d^2}{dx^2}\bigl[-(2\pi)^\alpha e^{-\alpha x^2}\bigr]
=2\alpha(2\pi)^\alpha e^{-\alpha x^2}(1-2\alpha x^2)
\]
takes negative values, so $-K_{\mathbb{C},\varphi}^{-\alpha}$ is not convex. However, $\log K_{\mathbb{C},\varphi}=x^2-\log(2\pi)$ is still convex, in agreement with Theorem \ref{ch3:th:Xiong-log}.
\end{example}

\backmatter

\end{document}